\documentclass[a4paper,12pt]{amsart}

\usepackage{kantlipsum} 
\calclayout

\usepackage[authoryear]{natbib}    
\usepackage{enumitem}
\usepackage[utf8]{inputenc}
\usepackage[T1]{fontenc}

\usepackage{amssymb,amsthm,amsmath}
\usepackage{bm}
\usepackage{enumitem}
\usepackage{graphicx}
\usepackage{xr}

\newtheorem{cor}{Corollary}
\newtheorem{defn}{Definition}
\newtheorem{lemma}{Lemma}
\newtheorem{prop}{Proposition}
\newtheorem{thm}{Theorem}

\newcommand{\beginsupplement}{%
        \setcounter{lemma}{0}
        \renewcommand{\thelemma}{S\arabic{lemma}}%
        \setcounter{prop}{0}
        \renewcommand{\theprop}{S\arabic{prop}}%
        \setcounter{cor}{0}
        \renewcommand{\thecor}{S\arabic{cor}}%
        \setcounter{table}{0}
        \renewcommand{\thetable}{S\arabic{table}}%
        \setcounter{figure}{0}
        \renewcommand{\thefigure}{S\arabic{figure}}%
        \setcounter{section}{0}
        \renewcommand{\thesection}{S\arabic{section}}%
        \setcounter{equation}{0}
        \renewcommand{\theequation}{S.\arabic{equation}}%
     }

\newcommand{\bx}{{\bf x}}

\newcommand{\by}{{\bf y}}

\newcommand{\bbeta}{\bm{\beta}}
\newcommand{\bgamma}{\bm{\gamma}}
\newcommand{\bmu}{\bm{\mu}}

\newcommand{\btheta}{\bm{\theta}}

\begin{document}

\title{Concentration of posterior model probabilities and normalized $L_0$ criteria}

\author{David Rossell \\
Universitat Pompeu Fabra, Department of Business and Economics, Barcelona (Spain)}

\begin{abstract}
We study frequentist properties of Bayesian and $L_0$ model selection, with a focus on (potentially non-linear) high-dimensional regression. We  propose a construction to study  how posterior probabilities and normalized $L_0$ criteria concentrate on the (Kullback-Leibler) optimal model and other subsets of the model space. When such concentration occurs, one also bounds the frequentist probabilities of selecting the correct model, type I and type II errors. These results hold generally, and help validate the use of posterior probabilities and $L_0$ criteria  to control frequentist error probabilities associated to model selection and hypothesis tests.  Regarding regression, we help understand the effect of the sparsity imposed by the prior or the $L_0$ penalty, and of problem characteristics such as the sample size, signal-to-noise, dimension and true sparsity. 
 A particular finding is that one may use less sparse formulations than would be asymptotically optimal, but still attain consistency and often also significantly better finite-sample performance. 
We also prove new results related to misspecifying the mean or covariance structures, and give tighter rates for certain non-local priors than currently available. 

\vspace{3mm}
\end{abstract}

\maketitle

Selecting a probability model and quantifying the associated uncertainty are two fundamental tasks in Statistics.
In Bayesian model selection (BMS), given models and priors one obtains posterior model probabilities that guide model choice and measure the (Bayesian) certainty on that choice.
It is interesting to understand how such posterior probabilities are related to the frequentist probability of selecting the optimal model (defined below).
$L_0$ penalties are also powerful selection criteria, but it is less clear how to portray uncertainty.
Suppose one selects the model optimizing the Bayesian information criterion (BIC, \cite{schwarz:1978}), how is one to measure the certainty about that choice? Given the connection between the BIC and Bayes factors, it is tempting to define a pseudo-posterior probability via a normalized $L_0$ criterion (defined below).
Again the question is how do these pseudo-probabilities relate to frequentist selection probabilities.

Our goals are two-fold. First, we present a general framework to study the $L_1$ convergence of posterior model probabilities and normalized $L_0$ criteria, and show that the obtained rates bound the frequentist probabilities of choosing the wrong model and making type I-II errors.
There exists previous work studying $L_1$ convergence (see below), our specific construction however
is novel (to our knowledge) and reduces the problem to integrating certain Bayes factor tail probabilities.
The result on bounding frequentist error probabilities is also new (although elementary),
and validates using posterior probabilities and normalized $L_0$ criteria to quantify model choice uncertainty from a frequentist standpoint.
Our second goal is to apply our framework to Gaussian regression to synthesize and extend current theoretical results.
We show that posterior model probabilities in high dimensions depend on the same three elements that drive their behavior in finite dimensions. These are the sparsity of the prior on the models, the dispersion of the prior on the parameters, and whether the latter is a local or a non-local prior (\cite{johnson:2010}).
We impose fairly mild conditions on these prior elements so that, in contrast to current results proving consistency for a specific prior sparsity regime, we portray the impact of the prior sparsity.
As novel aspects, we consider model misspecification within (possibly non-linear) regression and we obtain tighter rates for the non-local product MOM prior (pMOM, \cite{johnson:2012}) than currently available.
A practical implication of our results is that, by using less sparse priors than those leading to optimal asymptotic rates, one can still get consistency and sometimes attain significantly better finite $n$ performance. Some of our examples may be striking in that regard.


The introduction is organized as follows. First, we lay out the minimal notation needed to define the problem. We subsequently review existing results for fixed- and high-dimensional settings, and finally we outline the paper.

Let $\by=(\by_1,\ldots,\by_n)$ be an observed outcome and $n$ the sample size.
One wishes to consider a set of $K$ candidate models $M_1,\ldots,M_K$ for $\by$. 
Each model is defined by a density $p(\by \mid \btheta_k,\phi,M_k)$ for $k=1,\ldots,K$,
where $\btheta_k \in \Theta_k$ is a parameter of interest and $\phi \in \Phi$ a (potential) nuisance parameter.
The model dimension is given by $p_k=\mbox{dim}(\Theta_k)$ and $d=\mbox{dim}(\Phi)$.
Densities are in the Radon-Nikodym sense, 
in particular allowing discrete and continuous $\by$.
Without loss of generality let $\Theta_k \subseteq \Theta \subseteq \mathbb{R}^p$ for $k=1,\ldots,K$, i.e. models are nested
within a larger model of dimension $p+d$.
Although not denoted explicitly in high-dimensional problems 
 both the number of parameters  $p$ and  models  $K$ may grow with $n$, and this is precisely our main focus.
 In BMS each model is equipped with  a prior density $p(\btheta_k,\phi \mid M_k)$, and one obtains
posterior probabilities
\begin{align}
p(M_k \mid \by)= \left( 1 + \sum_{l \neq k}^{}\frac{p(\by \mid M_l)}{p(\by \mid M_k)}
  \frac{p(M_l)}{p(M_k)} \right)^{-1}=
\left( 1 + \sum_{l \neq k}^{} B_{lk} \frac{p(M_l)}{p(M_k)} \right)^{-1}
\label{eq:pp}
\end{align}
where $p(\by \mid M_k)= \int_{}^{}\!\, p(\by \mid \btheta_k,\phi,M_k) dP(\btheta_k,\phi\mid M_k)$ is the integrated likelihood under model $M_k$,
$p(M_k)$ its prior probability and $B_{lk}= p(\by \mid M_l)/p(\by \mid M_k)$ the Bayes factor between $(M_l,M_k)$.
We focus our discussion on BMS, but one can obtain analogous expressions for normalized $L_0$ criteria, see Section \ref{sec:l0penalties}.

 To fix ideas, consider a Gaussian regression where 
$\by \in \mathbb{R}^n$ and the data analyst assumes the model
$p(\by \mid \btheta,\phi)= N(\by; X\btheta, \phi I)$ where $X$ is an $n \times p$ matrix,
$\btheta \in \mathbb{R}^p$ the regression coefficients, and $\phi>0$ the error variance.
Models $M_k$ are defined by selecting subsets of columns in $X$,
 i.e. $p(\by \mid \btheta_k,\phi)= N(\by; X_k\btheta_k, \phi I)$ where 
$X_k$ is the $n \times p_k$ matrix containing the columns selected by $M_k$, and $\btheta_k \in \mathbb{R}^{p_k}$.
Note that $X$ may contain non-linear effects and interactions  such as wavelets, splines, or tensor-products. 

Here for example one might set a conjugate Normal-inverse Gamma prior
$p(\btheta_k, \phi \mid M_k)= N(\btheta_k; {\bf 0}, \tau \phi I) \mbox{IG}(\phi; a_\phi/2,b_\phi/2)$,
where $(\tau,a_\phi,b_\phi)$ are prior parameters.

 We consider the following question. 
Suppose that $\by$ arises from some data-generating density $f^*$, which may be outside the considered models (model misspecification).
Let $M_t$ be the optimal model in that it is the smallest model minimizing Kullback-Leibler (KL) divergence to $f^*$ (see Section \ref{sec:approach}).
 For example,  in Gaussian regression
$M_t$ is the smallest model minimizing mean squared prediction error under $f^*$.
 
If the models are well-specified, then $M_t$ is simply the smallest model containing $f^*$, and is often referred to as {\it true model}.
Ideally one wants to assign large $p(M_t \mid \by)$, so that one not only selects the optimal model but is also confident about that choice.

Our goal is to study if $p(M_t \mid \by)$ converges to 1 as $n$ grows, and at what rate.

 This problem has been well-studied in 
finite dimensions where $(p,K)$ do not grow with $n$.

Consider a model $M_k$ that includes the optimal $M_t$ ($\Theta_t \subset \Theta_k$), i.e. $M_k$ contains spurious parameters.
We refer to such $M_k$ as a {\it spurious model}.
For fairly general models and priors, the Bayes factor $B_{kt}$ converges in probability to 0 at a polynomial rate in $n$ and in the prior dispersion ($\tau$ in our regression example).
See Theorem 1 in \cite{dawid:1999}, 
Propositions 3, 4 and 7 in \cite{rossell:2019t} for misspecified Gaussian, binary and survival regression,
and the proof of our Theorem \ref{THM:NECESSARY_CONDITIONS_PC} for models with concave log-likelihood.
This polynomial rate holds when $p(\btheta_k \mid \phi,M_k)$ is a local prior, for non-local priors the rates are faster \citep{johnson:2010}.
In contrast, if $M_k$ is a non-spurious model ($\Theta_t \not\in \Theta_k$, i.e. missing parameters from $M_t$), then $B_{kt}$ vanishes exponentially in $n$ (more precisely, in a non-centrality parameter that is proportional to $n$, under suitable assumptions).
In summary, to help discard spurious models one may either set large $\tau$ (i.e. a diffuse prior on parameters),
set sparse model priors $p(M_k)$ that penalize model size, and/or set a non-local prior.
 
A caveat with diffuse and sparse priors is that they penalize complexity purely a priori, which can lead to a drop in statistical power.


Extensions of such precise rates to high dimensions are of fundamental interest yet hard to come by.
Most results focus on a prior satisfying relatively rigid sparsity conditions and either only study consistency (with no rates)
or focus attention on asymptotic optimality.
We show that in high-dimensional regression the finite-dimensional rates discussed above still hold, up to lower-order terms, for the stronger form of $L_1$ convergence.
 In particular the main prior features driving posterior consistency remain the same (the use of diffuse, sparse and non-local). 

We review selected high-dimensional BMS literature.
\cite{johnson:2012} proved that $p(M_t \mid \by)$ converges to 1 in linear regression with $p \ll n$ under NLPs and uniform $p(M_k)$.
\cite{narisetty:2014} showed that if $p \ll e^n$ then certain diffuse priors $p(\btheta_k \mid M_k)$ also attain consistency.
In fact, the RIC of \cite{foster:1994} is a related early advocate for diffuse priors,  and can also be shown to attain consistency for $p \ll e^n$ (Section \ref{sec:l0penalties}). 
\cite{shin:2018} extended \cite{johnson:2012} 
to $p \ll e^n$ under certain diffuse NLPs.
\cite{yangyun:2017} also used diffuse priors  (defined implicitly via a prior anti-concentration condition) , in a more general framework  that allows for non-parametric models. 
These results proved consistency but no specific rates were given.
\cite{castillo:2015} showed that, by using so-called {\it Complexity priors} $p(M_k)$ and Laplace priors on parameters,
one can consistently select the data-generating model in regression.
\cite{chae:2016} proved  that the same prior structure attains consistency in regression with non-parametric symmetric errors. 
\cite{gao:2015} extended these results to general structured linear models under misspecified sub-Gaussian errors,
and \cite{rockova:2017} to regression trees.
\cite{yangyun:2016} studied a regression setting
 where one uses diffuse priors on parameters and a type of Complexity prior on models. 
These contributions provide significant insights, the focus however is showing that complex models are asymptotically discarded a posteriori
under a given, sufficiently sparse, prior setting. 

In summary, the diffuse priors as in \cite{narisetty:2014} and Complexity priors as in \cite{castillo:2015} underlie much of the state-of-the-art literature. These priors excel at discarding spurious models, and do not require one to restrict the maximum model complexity.
Our results suggest that they should be used with care, however, and that there can be advantages to setting less sparse priors by placing mild restrictions on the model complexity.
As an illustration, we preview

Figures \ref{fig:margpp_ortho}-\ref{fig:pp_corpred} where three prior formulations were used.
Although the Complexity prior attains better asymptotic rates,
the combined pMOM and Beta-Binomial priors attained better finite $n$ power/sparsity tradeoffs.
Although in this example $M_t$ has small dimension $p_t=5,10,20$ (sparse truth), the losses in power due to setting sparse priors are substantial.

It is therefore of interest to study consistency allowing for less sparse priors.

Note that we focus on BMS where one entertains a collection of models.
An alternative is to use shrinkage priors, i.e. set a single model and a continuous
prior on $\btheta$ concentrating most mass on subsets of $\Theta$. 
While interesting shrinkage priors are fundamentally different
as any zero Lebesgue measure subset of $\Theta$ has zero posterior probability,
hence in our view they are more suitable for parameter estimation than for structural learning.
For results on shrinkage priors see
\cite{bhattacharya:2012} or
\cite{song:2017}, for example.

Also, we adopt a fully Bayesian framework where no priors are data-dependent.
Extending our framework to empirical Bayes approaches where prior features are learned from data is interesting
but requires a delicate treatment beyond our scope to avoid certain posterior degeneracy issues, we refer the reader to \cite{petrone:2014}.

The paper is organized as follows.
Section \ref{sec:approach}  sets notation,  presents our general framework,  and shows that the expectation of posterior probabilities such as $E_{f^*}(p(M_t \mid \by))$ bound relevant frequentist error probabilities.
Section \ref{sec:priors_and_sparsity} discusses the priors that we focus attention on, and important technical conditions related to the model complexity and sparsity embedded in the prior.
 It also outlines necessary conditions that fairly general priors need to satisfy, if one wishes to attain consistency. 
Section \ref{sec:modelrates} characterizes the posterior probability of individual models in Gaussian regression, under essentially any prior on the models and the priors on parameters from Section \ref{sec:priors_and_sparsity}.
Specifically, we consider Zellner priors (with known and unknown error variance) and more general Normal priors, for which Bayes factors have tractable expressions and hence simplify our exposition. We also include the pMOM prior where such an expression is unavailable, and misspecified (possibly non-linear) models where tail probabilities are harder to bound.

We show that failing to include true non-linearities or omitting relevant variables
causes an exponential drop in power, whereas misspecifying the error covariance (truly correlated and/or heteroskedastic errors) need not do so but may inflate false positives.
Section \ref{sec:l0penalties} extends Section \ref{sec:modelrates} to normalized $L_0$ penalties, including the BIC, EBIC and RIC.
Section \ref{sec:varsel} obtains global rates for $p(M_t \mid \by)$ and other interesting model subsets.
 The results show that it is often possible to discard spurious parameters, even when not using particularly sparse priors in problems of fairly large dimension, e.g. by combining a Beta-Binomial prior on models with non-local priors on parameters.  
Section \ref{sec:varsel_examples} offers examples, and Section \ref{sec:discussion} concludes.
A significant number of auxiliary lemmas, technical results and all proofs are in the supplementary material.

\section{Approach}
\label{sec:approach}

We first formalize the notion of optimal model $M_t$ and introduce notation used throughout the paper in Section \ref{ssec:notation}.
Then Section \ref{ssec:l1convergence} presents a framework to study $L_1$ convergence of $p(M_t \mid \by)$ in fully general settings, discusses its tightness,
and shows that the associated rates bound relevant frequentist error probabilities. 

\subsection{Definitions and notation}
\label{ssec:notation}

We define the optimal model $M_t$ to be that of smallest dimension $p_t$ among all models minimizing Kullback-Leibler (KL) divergence to $f^*$.
For simplicity we assume $M_t$ to be unique, but our results hold more generally by defining $M_t$ to be the union of all smallest KL-optimal models.
\begin{defn}
Let $(\btheta^*,\phi^*)= \arg\min_{\theta \in \Theta, \phi \in \Phi} \mbox{KL}(f^*,p(\by \mid \btheta,\phi))$.
Define $t= \arg\min_{k \in \mathcal{M}^*} p_k$, where
\vspace{-2mm}
\begin{align}
\mathcal{M}^*= \left\{k: \exists (\btheta_k^*,\phi_k^*) \in \Theta_k \times \Phi:
\mbox{KL}(f^*,p(\by \mid \btheta_k^*,\phi_k^*))= \mbox{KL}(f^*,p(\by \mid \btheta^*,\phi^*)) \right\}
\nonumber
\end{align}
is the set of all models minimizing KL-divergence to $f^*$.
\label{def:Mt}
\end{defn}


We denote by $(\btheta^*,\phi^*)$ the global optimal parameter value minimizing KL-divergence to $f^*$, and by $(\btheta_k^*,\phi_k^*)$ that under a model $M_k$. If $f^*$ lies in the assumed model family (well-specified case), then $(\btheta^*,\phi^*)$ is the true parameter value.

We shall study the posterior probability assigned to models other than $M_t$.
To that end, it is convenient to 
denote the set of $l$-dimensional models that contain $M_t$ plus some spurious parameters by $S_l=\left\{k: \Theta_t \subset \Theta_k, p_k=l \right\}$.
We refer to $S_l$ as {\it spurious models} of dimension $l$.
Similarly, let $S_l^c=\left\{k: \Theta_t \not\subset \Theta_k, p_k=l \right\}$ be the size $l$ {\it non-spurious models},

and let $S= \bigcup_{l=p_t+1}^{\bar{p}} S_l$ and $S^c= \bigcup_{l=0}^{\bar{p}} S_l^c$ 
 the complete set of spurious and non-spurious models. 
Denote by $|S|$ the cardinality of $S$.

In our study it is often convenient to express certain conditions and results in terms of their asymptotic order as $n$ grows. To this end, 

$a_n \ll b_n$  denotes $\lim_{n \rightarrow \infty} a_n/b_n=0$ for two deterministic sequences $a_n,b_n>0$,
and similarly $a_n \preceq b_n$ denotes $\lim_{n \rightarrow \infty} a_n/b_n \leq c$ for some constant $c>0$.
Finally, $a_n \asymp b_n$ denotes that both $a_n \preceq b_n$ and $a_n \succeq b_n$.

As we discuss later, although $p$ could potentially grow exponentially with $n$, for certain prior/$L_0$ penalty settings
to achieve consistency it may be necessary to impose restrictions on the model complexity.
We assume that the analyst specifies a maximum model size that we denote by $\bar{p}=\max_k p_k$, and describe rates as a function of $\bar{p}$
For instance, in regression one may have $p \gg n$ but restrict attention to models selecting at most $\bar{p}=\mbox{min}\{n,p\}$ out of the $p$ variables, as choosing a model with $p_k>n$ parameters may not be desirable.
The number of models is then $K= \sum_{j=0}^{\bar{p}} {p \choose j}$, which is still $\gg n$. 


\subsection{$L_1$ convergence}
\label{ssec:l1convergence}

From \eqref{eq:pp}, posterior consistency requires $\sum_{k \neq t}^{} B_{kt} p(M_k)/p(M_t) \stackrel{P}{\longrightarrow} 0$.
The difficulty in high dimensions is that the number of models $K-1$ grows with $n$,
hence the sum can only vanish if each term $B_{kt} p(M_k)/p(M_t)$ converges to 0 quickly enough.
This intuition is clear, but obtaining probabilistic bounds for this stochastic sum is non-trivial,
since the $B_{kt}$'s may exhibit complex dependencies. 
To avoid dealing with such high-dimensional stochastic sums, it is simpler to study deterministic expectations.
Specifically, we study when $p(M_t \mid \by) \stackrel{L_1}{\longrightarrow} 1$, which by definition of $L_1$ convergence is equivalent to
\begin{align}
\mathop {\lim }\limits_{n \to \infty } \sum_{k \neq t}^{}
  E_{f^*} \left( p(M_k \mid \by) \right)= 0,
\label{eq:mean_ppwrong}
\end{align}
where $E_{f^*}(\cdot)$ is the expectation under $f^*$.
Some remarks are in order.
First, $L_1$ convergence in \eqref{eq:mean_ppwrong} implies convergence in probability.
Let $b_n >0$ be a sequence such that $\mathop {\lim }\limits_{n \to \infty } b_n=0$,
if $E_{f^*}(1-p(M_t\mid \by)) \preceq b_n$ then $1-p(M_t \mid \by)= O_p(b_n)$.
Naturally \eqref{eq:mean_ppwrong} may require more stringent conditions than convergence in probability,
but in regression we obtain essentially tight rates and the gains in clarity are substantial.
Second, one can evaluate the sum on the left-hand side of \eqref{eq:mean_ppwrong} for fixed $n$, $p$ and $K$,
i.e. the expression can be used in non-asymptotic regimes.

An advantage of studying $L_1$ convergence is that one automatically obtains a form of frequentist validity, in the sense of bounding relevant error probabilities, which helps justify the use of posterior probabilities (or normalized $L_0$ criteria) to quantify model choice uncertainty. 

Proposition \ref{PROP:CORRECTSEL_FROM_PP} shows that $E_{f^*}( p(M_t \mid \by))$ bounds the (frequentist) model selection probability
$P_{f^*}(\hat{k} \neq t)$, where $\hat{k}$ is the highest posterior probability model.
Proposition \ref{PROP:CORRECTSEL_FROM_PP} also holds when $\hat{k}$ is the median probability model of \cite{barbieri:2004},
 that is when $\hat{k}$ selects parameters with marginal posterior inclusion probability 
$P(\theta_j \neq 0 \mid \by) > 0.5$  (see the proof).

\begin{prop}
Let $\hat{k}= \arg\max_k p(M_k \mid \by)$ be the posterior mode, then
$$
P_{f^*}(\hat{k} \neq t) \leq 
2 E_{f^*}(1 - p(M_t \mid \by) )=
2 \sum_{k \neq t}^{} E_{f^*} (p(M_k \mid \by)).
$$
\label{PROP:CORRECTSEL_FROM_PP}
\end{prop}

Corollaries \ref{COR:POW_TYPEI_FROM_PP}-\ref{COR:POW_TYPEI_FROM_MARGPP} relate type I-II error probabilities to expected posterior model probabilities.

 Corollary \ref{COR:POW_TYPEI_FROM_PP}  is based on the trivial observation that family-wise type I-II error rates are both $\leq P_{f^*}(\hat{k} \neq t)$, and hence also bounded by Proposition \ref{PROP:CORRECTSEL_FROM_PP}.
Alternatively, suppose that $\hat{k}$ is obtained by selecting parameters with $P(\theta_j \neq 0 \mid \by)>t$, for some threshold $t$.

For instance, one may control the Bayesian False Discovery rate below some level $\alpha$ by setting a certain $t \leq 1- \alpha$ \citep{mueller:2004}.

Then, by Corollary \ref{COR:POW_TYPEI_FROM_MARGPP} the type I-II errors for individual coefficients are bounded by   $E_{f^*}(P(\theta_j \neq 0 \mid \by))$, times a factor that depends on $t$. 


\begin{cor}
Let $S(\hat{k})$ the set of non-zero parameters in model $\hat{k}= \arg\max_k p(M_k \mid \by)$.
\begin{itemize}[leftmargin=*]
\item The family-wise type I error is
$
P_{f^*} \left( \bigcup_{j: \theta_j^*= 0} \{ j \in S(\hat{k}) \}  \right) \leq P_{f^*}( \hat{k} \neq t).
$

\item The family-wise type II error is
$
P_{f^*} \left( \bigcup_{j: \theta_j^*= 1} \{ j \not\in S(\hat{k}) \}  \right) \leq P_{f^*}( \hat{k} \neq t).
$
\end{itemize}
  \label{COR:POW_TYPEI_FROM_PP}
\end{cor}

\begin{cor}
  Let $S(\hat{k})= \{j: P(\theta_j \neq 0 \mid \by) > t \}$ for a given threshold $t$.

\begin{itemize}[leftmargin=*]
\item False positives. Assume that $\theta_j^*=0$. Then
  $
P_{f^*}(j \in S(\hat{k})) \leq 
\frac{1}{t} E_{f^*}( P(\theta_j \neq 0 \mid \by)).
$

\item Power. Assume that $\theta_j^* \neq 0$. Then
  $
P_{f^*}( j \not\in S(\hat{k}) ) \leq 
\frac{1}{1-t} E_{f^*}( P(\theta_j=0 \mid \by) ).
  $
\end{itemize}
  \label{COR:POW_TYPEI_FROM_MARGPP}
\end{cor}

To summarize, if one can bound sums of expectations $E_{f^*}(p(M_k \mid \by))$ across models one can then prove that the posterior probability of $M_t$ converges to 1 via \eqref{eq:mean_ppwrong}, as well as bound the frequentist probability of selecting $M_t$, and of the selected model including type I-II errors.
The question is therefore how to bound the right-hand side in Proposition \ref{PROP:CORRECTSEL_FROM_PP}, which we discuss next.

Our strategy is to use that
\begin{align}
1 - p(M_t \mid \by)= \sum_{k \neq t} p(M_k \mid \by) \leq \sum_{k \neq t} \left( 1 + B_{tk} \frac{p(M_t)}{p(M_k)} \right)^{-1}.
\nonumber
\end{align}
Per Lemma \ref{LEM:TAILPROB_TO_MEAN} below, the $L_1$ convergence of the right-hand side can be proven by integrating tail probabilities that, conveniently, only involve pairwise Bayes factors $B_{kt}$. 
A natural question is whether said right-hand side provides a sufficiently tight bound. Lemma \ref{LEM:TIGHTNESS_PAIRWISEBOUND} shows that, indeed, whenever $p(M_t \mid \by)$ converges to 1 the right-hand side is asymptotically equivalent to $1- p(M_t \mid \by)$.


\begin{lemma}




$$
E_{f^*}(p(M_k \mid \by)) \leq
E_{f^*} \left( \left(1 + B_{kt} \frac{p(M_k)}{p(M_t)} \right)^{-1} \right) =
\int_0^1 P_{f^*}\left(B_{kt}> \frac{p(M_k)}{p(M_t) (1/u-1)} \right) du.
$$

\vspace{-4mm}
\label{LEM:TAILPROB_TO_MEAN}
\end{lemma}


\begin{lemma}
Suppose that $p(M_t \mid \by) \stackrel{L_1}{\longrightarrow} 1$. Then
\begin{align}
\frac{1 - p(M_t \mid \by)}{\sum_{k \neq t} (1 + B_{tk} p(M_t)/p(M_k))^{-1}} \stackrel{L_1}{\longrightarrow} 1.
\nonumber
\end{align}
\label{LEM:TIGHTNESS_PAIRWISEBOUND}
\end{lemma}



Our strategy is based on two steps.
First, we use Lemma \ref{LEM:TAILPROB_TO_MEAN} to bound the posterior probability assigned to an individual model $E_{f^*} (p(M_k \mid \by) )$.
This is achieved by bounding tail probabilities for $B_{kt}$, for all $n \geq n_{k0}$ and some fixed  $n_{k0}$.
Sections \ref{sec:modelrates}-\ref{sec:l0penalties} use such bounds for (possibly non-linear) Gaussian regression for Bayesian and normalized $L_0$ methods, respectively.
The key is that Bayes factors can be bounded by quadratic forms involving least-squares estimators (or Bayesian analogues), for which we derived tail inequalities (Section \ref{ssec:tailbounds}).

To facilitate applying our framework to other models, Section \ref{ssec:tailbounds} also gives finite-$n$ bounds for $E_{f^*}(p(M_k \mid \by))$
in more general cases where suitably re-scaled $\log(B_{tk})$ have exponential or polynomial tails.

The second step is to bound the right-hand side in Proposition \ref{PROP:CORRECTSEL_FROM_PP} for all $n \geq n_0$ and fixed $n_0= \max_k n_{k0}$ by adding the model-specific bounds. 
Note that one can similarly bound the posterior probability of other interesting model subsets, e.g. adding spurious parameters to $M_t$.
Section \ref{sec:varsel} performs this task for Gaussian regression (and implicitly for other settings where rates for $E_{f^*}(p(M_k \mid \by))$ take a similar form).
As a technical remark, one must ensure that such fixed $n_0$ exists. This need not hold in general, since the number of models $k \neq t$ grows with $n$, but in our regression examples such $n_0$ indeed exists.
We refer the reader to Section \ref{ssec:bound_seriessum} (A4) for further discussion.

\section{Conditions for consistency}
\label{sec:priors_and_sparsity}

We outline priors and conditions that are related to the extent to which they encourage sparsity.
The conditions feature non-centrality parameters that measure the signal strength when comparing $M_t$ versus another model $M_m$. 
We generically denote these by $\lambda_{tm}$, and define their precise meaning in each setting below.
Section \ref{ssec:necessary_cond} states necessary conditions for a wide class of models and priors for $p(M_t \mid \by)$ to converge to 1.
Section \ref{ssec:priors} lists the priors used in our Gaussian regression examples.
Section \ref{ssec:conditions_complexity} sets conditions on these priors
(see Section \ref{sec:l0penalties} for analogous conditions on $L_0$ criteria),
and discusses connections to related literature. 

The main difference to earlier work is that we do not restrict attention to situations where $p(M_k)$ is a sparse prior or one sets diffuse parameter priors.
By restricting the maximum model complexity, our study includes the use of less sparse priors, to provide a wider depiction of when one can hope $p(M_t \mid \by)$ to converge to 1.


\subsection{Necessary conditions}
\label{ssec:necessary_cond}

We list two necessary conditions for $p(M_t \mid \by)$ to converge to 1.
For simplicity we state them for models with concave log-likelihood, 
but in much more general models one obtains similar Bayes factor rates, and hence necessary conditions (Theorem 1 in \cite{dawid:1999}).
Without loss of generality, suppose that the prior on $\btheta_k$ is defined by taking a scale transformation of a random variable following some distribution $\tilde{p}()$, that is
$$
p(\btheta_k \mid \phi,M_k)= \tau^{-p_k/2} \tilde{p}(\btheta_k/\tau^{1/2} \mid \phi, M_k).
$$
where $\tau > 0$ is the scale parameter (see Section \ref{ssec:priors} for examples).

As discussed earlier, large $\tau$ leads to diffuse priors that favor sparsity.
Specifically, the Bayes factor $B_{tm}$ to compare $M_t$ versus some other $M_m$ asymptotically includes a term $(n\tau)^{(p_m-p_t)/2}$.
Theorem \ref{THM:NECESSARY_CONDITIONS_PC} states that, if the combined sparsity induced by this term and 
model prior probabilities is too strong (relative to the signal strength), then $p(M_t \mid \by)$ cannot converge to 1. 
Signal strength is measured by a parameter
\begin{align}
 \lambda_{tm}= E_{f^*}\left[ \log \frac{p(\by \mid \btheta_t^*, \phi_t^*)}{p(\by \mid \btheta_m^*, \phi_m^*)} \right]
= \mbox{KL}(f^*, p(\by \mid \btheta_m^*, \phi_m^*)) - \mbox{KL}(f^*, p(\by \mid \btheta_t^*, \phi_t^*))
\label{eq:ncp_logconcave}
\end{align}
which, in Gaussian regression, is given by differences in mean-squared prediction errors (Section \ref{ssec:conditions_complexity}).
Theorem \ref{THM:NECESSARY_CONDITIONS_PC} assumes near-minimal conditions used by \cite{hjort:2011} to prove asymptotic normality, allowing for misspecification
(see Section \ref{proof:necessary_conditions_pc}). 

\begin{thm}
Consider models $M_t$ and $M_m$ such that, for $k \in \{t,m\}$, $p_k$ is fixed,
$\log p(\by \mid \btheta_k,\phi)$ is continuous and strictly concave 
and $\tilde{p}(\btheta_k/\tau^{1/2}, \phi \mid M_k) \in (0,\infty)$ for all $(\btheta_k,\phi)$.
Assume the conditions of Theorem 4.1 in \cite{hjort:2011}.

\begin{enumerate}[leftmargin=*,label=(\roman*)]
\item Let $m \in S$ be a spurious model. 
If $(\tau n)^{(p_m-p_t)/2} p(M_t)/p(M_m) \preceq 1$, then $p(M_t \mid \by)$ does not converge in probability to 1.

\item Let $m \in S^c$ be a non-spurious model and let $\lambda_{tm}$ in \eqref{eq:ncp_logconcave}.
If
\begin{align}
\frac{\lambda_{tm}}{2} +  \frac{p_m - p_t}{2} \log(\tau n) + \log \left( \frac{p(M_t)}{p(M_m)} \right) \preceq 1
\nonumber
\end{align}
then $p(M_t \mid \by)$ does not converge in probability to 1.
\end{enumerate}
\label{THM:NECESSARY_CONDITIONS_PC}
\end{thm}

Theorem \ref{THM:NECESSARY_CONDITIONS_PC} allows $p$ to grow with $n$, but to ease exposition it only lists necessary conditions for fixed-dimensional models.
From Part (i), for any spurious model it is necessary that $(\tau n)^{(p_m-p_t)/2} p(M_t)/p(M_m) \gg 1$, which
is identical to Condition C1 in Section \ref{ssec:conditions_complexity} for high-dimensional regression.
Part (ii) is also analogous to Condition C2.

\subsection{Prior distributions for regression}
\label{ssec:priors}

The framework from Section \ref{sec:approach} applies to any prior but the required algebra varies, for illustration we focus on several popular priors.
First we consider Zellner's prior
\begin{align}
p(\btheta_k \mid M_k,\phi)= N(\btheta; {\bf 0}, \tau n \phi (X_k'X_k)^{-1}),
\label{eq:zellnerprior}
\end{align}
where $\tau > 0$ is a known prior dispersion.
For simplicity $X_k'X_k$ is assumed invertible for $p_k \leq \bar{p}$.
We then extend results to Normal priors
\begin{align}
p(\btheta_k \mid M_k,\phi)= N(\btheta_k; {\bf 0}, \tau n \phi V_k)
\label{eq:normalprior}
\end{align}
with general covariance $V_k$ and to the pMOM prior \citep{johnson:2012}
\begin{align}
p(\btheta_k \mid \phi, M_k)=
\prod_{j \in M_k}^{} \theta_j^2 \bx_j'\bx_j/(\tau n \phi) N(\theta_j; {\bf 0}, \tau n \phi/ \bx_j'\bx_j) 
\label{eq:pmomprior}
\end{align}
where $\bx_j$ is the $j^{th}$ column in $X_k$. 

The idea is that constant $\tau$ leads to roughly constant prior variance,

e.g. for the pMOM prior if $X_k$ has zero column means and unit column variances then $n/\bx_j'\bx_j=1$.
Such constant $\tau$ may be desirable from a foundational Bayesian point of view, where the prior does not to depend on $n$.

In fact, the default choice $\tau=1$ leads to the unit information prior, which in turn leads to the BIC \citep{schwarz:1978}.
An alternative is to set $\tau$ growing with $n$, which leads to diffuse priors.
For example, one may set $\tau=p^2/n$ \citep{foster:1994},
$\tau=\mbox{max}\{1,p^2/n\}$ \citep{fernandez:2001}
and $\tau \gg p^2$ \citep{narisetty:2014}. 

As discussed, diffuse priors are used by many high-dimensional methods to induce sparsity.

Regarding the error variance $\phi$, whenever we treat it as unknown, we set $p(\phi \mid M_k)= \mbox{IG}(\phi; a_\phi/2,l_\phi/2)$ for fixed $a_\phi,l_\phi>0$.

For the prior on the models, in Section \ref{sec:modelrates} we allow for a general prior.
For concreteness, when discussing prior sparsity conditions below and when providing global rates in Section \ref{sec:varsel}, we focus on three popular choices.
These assume that all models with the same dimension $p_k$ receive equal prior probability, that is
\begin{align}
p(M_k)= P(p_k = l)/{p \choose l},
\label{eq:prior_model}
\end{align}
where $P(p_k=l)$ is the prior on the model size, and ${p \choose l}$ the number of models selecting $l$ parameters out of $p$.

First, we consider the uniform prior where $P(p_k=l)= {p \choose l}$, so that $p(M_k)=1/K$ for all $k=1,\ldots,K$.
Second, we consider the Beta-Binomial(1,1) prior where $P(p_k=l)= 1/\bar{p}$ \citep{scott:2010},
and finally and a so-called Complexity prior where $P(p_k=l) \propto 1/p^{cl}$  for $c>0$ \citep{castillo:2015}.
Note that the Beta-Binomial corresponds to $c=0$.

\subsection{Conditions on model complexity and prior sparsity}
\label{ssec:conditions_complexity}

We state two sets of conditions. First, B1-B2 constrain the sizes of the optimal and largest allowed models.
\begin{enumerate}[leftmargin=*,label=(B\arabic*)]
\item The maximum model size satisfies $\bar{p} \ll \min \{n,p, n\tau\}$.

\item  The optimal model size satisfies $p_t \ll \min\{\bar{p}, n\}$. 
\end{enumerate}

If one assigns non-vanishing $\tau$, as in all default choices above, B1 simplifies to $\bar{p} \ll \min\{n,p\}$.

%
%
One can allow for larger $\bar{p}=p$,
e.g. under Zellner's prior $B_{mt}=1$ for $p_m \geq n$ and one can immediately bound $E_{f^*}(p(M_m \mid \by))$,
 but then one must impose stricter prior sparsity conditions that our C1-C2 below. 
Setting $\bar{p} \ll n$ seems natural, however, as $p_m \geq n$ results in data interpolation.

See \cite{martin:2017} (Section 2.1) and references therein for further arguments for setting $\bar{p} \leq n$.

The second set of Conditions C1-C2 restrict the sparsity induced by the model prior and the prior dispersion $\tau$ in \eqref{eq:zellnerprior}-\eqref{eq:pmomprior}, and are related to the non-centrality parameter measuring the signal strength.

Specifically, let $H_m=X_m(X_m'X_m)^{-1}X_m'$ be the projection matrix onto the column space of $X_m$. For any non-spurious model $m \in S^c$, denote by
\begin{align}
\lambda_{tm}=(X_t \btheta_t^*)'(I-H_m)X_t \btheta_t^*/\phi^*
\label{eq:ncp_regression}
\end{align}
the non-centrality parameter measuring the difference in mean squared prediction error between $M_t$ and $M_m$ under the KL-optimal $(\btheta_t^*,\phi^*)$.

Equivalently, $\lambda_{tm}$ is the difference between the $L_2$ norm of the optimal predictor $X_t\btheta_t^*$ relative to its projection onto $X_m$.

This non-centrality parameter can be lower-bounded by
$$
\lambda_{tm} \geq n v_{tm} (\btheta_t^*)' \btheta_t^*/\phi^*
$$
where $v_{tm}$ is the smallest non-zero eigenvalue of $X_t'(I-H_m)X_t/n$.

Conditions C1-C2 suffice for $p(M_m \mid \by) \stackrel{L_1}{\longrightarrow} 0$ in high-dimensional regression.
As we shall see in Section \ref{sec:varsel}, uniform versions of C1-C2 also guarantee that $p(M_t \mid \by) \stackrel{L_1}{\longrightarrow} 1$. 
C1-C2 are stated for a generic $p(M_k)$, see Section \ref{supplsec:regcond} for concrete expressions for the uniform, Beta-Binomial and Complexity priors in \eqref{eq:prior_model}.


\begin{enumerate}[leftmargin=*,label=(C\arabic*)]
\item Let $m \in S$ be a spurious model. As $n\rightarrow\infty$, $(\tau n)^{(p_m-p_t)/2} p(M_t)/p(M_m) \gg 1$.

\item Let $m \in S^c$ be a non-spurious model. As $n\rightarrow\infty$,

\begin{align}
\frac{\lambda_{tm}}{2\log(\lambda_{tm})} + {\frac{p_m-p_t}{2}} \log(\tau n) + \log \left( \frac{p(M_t)}{p(M_m)} \right)  - \log(p_m) \gg 1.
\nonumber
\end{align}

\end{enumerate}

C1-C2 ensure that $p(M_k)$ and $\tau$ do not favor $M_m$ over $M_t$ too strongly a priori

and, per Theorem \ref{THM:NECESSARY_CONDITIONS_PC}, are near-necessary.
See also Section \ref{tightness:zellner_phiknown} for an extension of Theorem \ref{THM:NECESSARY_CONDITIONS_PC} to high-dimensional models for Zellner's prior.
 For the pMOM prior one can relax slightly C1 to $(\tau n)^{3(p_m-p_t)/2} p(M_t)/p(M_m) \gg 1$ under certain conditions, see Section \ref{ssec:mom}. 

We compare our conditions to those in \cite{narisetty:2014}, \cite{castillo:2015}, \cite{yangyun:2016} and \cite{yangyun:2017}.
We offer a summary, and discuss further details in Section \ref{supplsec:regcond}.
A main difference is on the prior setup. These authors restricted attention to diffuse priors (large $\tau$) and/or Complexity priors akin to that in \eqref{eq:prior_model}.
Specifically \cite{narisetty:2014} and \cite{yangyun:2016} required $\tau n \gg p^2$, whereas \cite{yangyun:2017} set a prior anti-concentration condition that also leads to $\tau$ growing with $n$.
\cite{castillo:2015} and \cite{yangyun:2016} required $p(M_k)$ to be a Complexity prior, and \cite{narisetty:2014} also used $p(M_k)$ that converges to a Complexity prior as $p$ grows. 

Our C1-C2 in principle allow more general $\tau$ and $p(M_k)$, such as fixed $\tau$ and $p(M_k)$ that do not penalize model size exponentially, e.g. the Beta-Binomial. 
For such choices the asymptotic rates for $p(M_t \mid \by)$ are then usually slower,
further one may need to restrict the maximum model complexity $\bar{p}$ (see Section \ref{sec:varsel}).
Nevertheless, there can be significant improvements for finite $n$, as illustrated in Section \ref{sec:varsel_examples}.

Regarding conditions on the data-generating truth,
\cite{narisetty:2014} require that $p_t$ is fixed, \cite{castillo:2015} that $p_t \leq \sqrt{n/ \log p}$, \cite{yangyun:2016} that $p_t \leq n/\log p$ and \cite{yangyun:2017} that $p_t \log(p/p_t) \leq n$.
These are related to our B1-B2, which require $p_t \ll n$ and $\bar{p} \ll n$, though these authors did not restrict $\bar{p}$.
Rather, they set $\bar{p}=p$ and priors that strongly penalize complexity.

Finally, these authors also set assumptions which, under restricted eigenvalue conditions, are related to beta-min conditions.
Specifically, \cite{castillo:2015} essentially required that
$\min_j |\theta_j^*|^2/\phi^* > p_t (\log p)/n,$
\cite{yangyun:2016} that $\min_j |\theta_j^*|^2/\phi^* > (c + p_t) (\log p)/n$, where $c$ is the Complexity prior's parameter,
and \cite{narisetty:2014} and \cite{yangyun:2017} that $\min_j |\theta_j^*|^2/\phi^* > (\log p)/n$.
Under such eigenvalue conditions, if $p(M_k)$ is the Complexity prior then for our C2 to hold it suffices that
\begin{align}
\min_j |\theta_j^*|^2/\phi^* \gg [\log(\tau n) + (1+c) \log p]/n,
\label{eq:suffcond_c2}
\end{align}
which is similar to these conditions above. Recall that $c=0$ corresponds to the Beta-Binomial prior, illustrating that using less sparse $p(M_k)$ lowers the required signal strength.
These conditions are mild, e.g. \cite{wainwright:2009information} showed that 
$\min_j |\theta_j^*|^2/\phi^* > [\log(p/p_t)]/n$ is a necessary condition for any method to consistently select $M_t$.

\section{Model-specific rates for regression}
\label{sec:modelrates}

In this section we bound $E_{f^*}(p(M_m\mid\by))$ for a single model $M_m$

for Gaussian regression and the priors in Section \ref{ssec:priors}.

Per Lemma \ref{LEM:TAILPROB_TO_MEAN} the proof strategy is to bound tail probabilities for pairwise Bayes factors $B_{mt}$.
Sections \ref{ssec:zellner_phiknown_spur}-\ref{ssec:mom} consider the case where the model is well-specified,
that is they assume a data-generating $f^*(\by)= N(\by; X_t \btheta_t^*, \phi^* I)$.
Then $B_{tm}$ is bounded by chi-square and F distribution tails 
for Normal priors, and by a slightly more involved term for pMOM priors.
Section \ref{ssec:misspec_mean} considers 
 the situation where the mean structure has been misspecified, e.g. $E_{f^*}(y)$ features 
variables or non-linear terms that were omitted from $X$.

Finally, Section \ref{ssec:misspec_covar} considers a misspecified covariance case, i.e. $f^*$ has heteroskedastic and/or correlated errors.

 The rates in Sections \ref{ssec:zellner_phiknown_spur}-\ref{ssec:normalprior} 
are similar to standard finite-dimensional rates (Theorem 1 in \cite{dawid:1999}, proof of our Theorem \ref{THM:NECESSARY_CONDITIONS_PC}).
Roughly speaking, non-spurious models are discarded at an exponential rate in $n$ (more precisely, in the non-centrality parameter $\lambda_{tm}$ in \eqref{eq:ncp_regression}, proportional to $n$ under restricted eigenvalue conditions).
More critically, spurious models are discarded at a rate that is essentially
\begin{align}
 E_{f^*}(p(M_m \mid \by)) \preceq \frac{p(M_m)}{p(M_t) (\tau n)^{(p_m-p_t)/2}},
\nonumber
\end{align}
up to lower-order terms. 
This result portrays the effect of the prior dispersion $\tau$ and model prior probabilities to encourage sparsity in more general regimes that in current high-dimensional literature (see Section \ref{sec:priors_and_sparsity}). As shown in Section \ref{sec:varsel}, the implication is that one can often allow for fixed $\tau$ and/or $p(M_k)$ that are not particularly sparse (e.g. the Beta-Binomial prior in \eqref{eq:prior_model}), and still attain consistency.

The pMOM rates to discard spurious models are faster, as is standard for non-local priors, but we provide tighter rates than currently available (see Section \ref{ssec:mom}).

\subsection{{\bf Zellner's prior with known variance}}
\label{ssec:zellner_phiknown_spur}

Under Zellner's prior and known error variance $\phi^*$, simple algebra gives
\begin{align}
B_{tm}= \exp \left\{ -\frac{\tau n}{2\phi^*(1+\tau n)} W_{mt} \right\}
(1+\tau n)^{\frac{p_m-p_t}{2}}
\label{eq:bf_zellner_known}
\end{align}
and hence in Lemma \ref{LEM:TAILPROB_TO_MEAN}
\begin{align}
P_{f^*}\left(B_{mt}> \frac{p(M_t)}{p(M_m)(1/u-1)}  \right)= P_{f^*} \left( \frac{W_{mt}}{\phi^*} >
\frac{1+\tau n}{\tau n} 2 \log \left[ \frac{(1+\tau n)^{\frac{p_m-p_t}{2}} p(M_t)}{p(M_m)(1/u-1)} \right]
 \right),
\label{eq:tailprob_zellner_known}
\end{align}
where $W_{mt}=\hat{\btheta}_m'X_m'X_m\hat{\btheta}_m - \hat{\btheta}_t'X_t'X_t\hat{\btheta}_t$
is the difference between residual sums of squares
under $M_t$ and $M_m$ and $\hat{\btheta}_m= (X_m'X_m)^{-1} X_m' \by$ the least-squares estimate.

Proposition \ref{PROP:ZELLNER_PHIKNOWN} gives a simple asymptotic expression for the $L_1$ rate at which $p(M_m \mid \by)$ vanishes. 
Spurious models are discarded at a rate that depends on $p(M_m)$ and $\tau n$.
Non-spurious models are discarded near-exponentially in the non-centrality parameter, times a factor driven by $p(M_m)$ and $\tau n$. 
The result portrays the effect of favoring sparse models either via $p(M_k)$ or by setting large $\tau$, namely a faster rate in Part (i) at the cost of a slower rate in Part (ii) for any model of size $p_m < p_t$.

\begin{prop}
Assume that $f^*(\by)= N(\by; X_t \btheta_t^*; \phi^* I)$ and consider $m \neq t$.
\begin{enumerate}[leftmargin=*,label=(\roman*)]
\item Let $m \in S$ be a spurious model, and $g= (\tau n)^{(p_m-p_t)/2} p(M_t)/p(M_m)$.
Assume Conditions B1 and C1. Then, for all fixed $\alpha <1$,

$$E_{f^*}(p(M_m\mid\by)) \preceq 
\frac{[\log(g)]^{(p_m-p_t)/2}}{g}  
\ll \left[ \frac{p(M_m)}{p(M_t) (\tau n)^{(p_m-p_t)/2} } \right]^\alpha
$$

\item Let $m \in S^c$ be a non-spurious model. Assume Conditions B2 and C2. Then
\begin{align}
 E_{f^*}(p(M_m \mid \by)) \ll
e^{-\lambda_{tm}^\gamma / 2} 
\left[ \frac{p(M_m) }{p(M_t) (\tau n)^{(p_m-p_t)/2} } \right]^\gamma,
\nonumber
\end{align}
for all fixed $\gamma <1$, where $\lambda_{tm}$ is as in \eqref{eq:ncp_regression}.
\end{enumerate}
\label{PROP:ZELLNER_PHIKNOWN}
\end{prop}

We remark that $\alpha,\gamma$ are taken arbitrarily close to 1, and are introduced to provide simpler expressions, see the proof for slightly tighter bounds where $\alpha=\gamma=1$, after adding lower-order terms.

Proposition \ref{PROP:ZELLNER_PHIKNOWN} gives upper-bounds. To see that they are reasonably tight, Section \ref{tightness:zellner_phiknown} 
shows that for spurious models $E_{f^*}(p(M_m \mid \by)) \succeq [p(M_m)/p(M_t)] (\tau n)^{-(p_m-p_t)/2}$, which equals the rate in Part (i), up to a log term.
A similar argument is made for Part (ii). 

\subsection{{\bf Zellner's prior with unknown variance}}
\label{ssec:zellner_unknown_spur}

Proposition \ref{PROP:ZELLNER_PHIUNKNOWN} extends Proposition \ref{PROP:ZELLNER_PHIKNOWN} to the case where $\phi^*$ is unknown, and one sets a prior $\phi \sim \mbox{IG}(a_\phi/2, l_\phi/2)$.
The rates are essentially equivalent, up to lower-order terms.

Let $s_k= \by'\by - \by'X_k(X_k'X_k)^{-1}X_k'\by$ be the residual sum of squares under $M_k$, then
\begin{align}
B_{tm}= \left( \frac{\tilde{s}_m}{\tilde{s}_t}\right)^{\frac{a_\phi+n}{2}} (1+ \tau n)^{\frac{p_m-p_t}{2}}=
\left( 1 + \frac{p_m-p_t}{n-p_m} \tilde{F}_{mt} \right)^{-\frac{a_\phi+n}{2}} (1+ \tau n)^{\frac{p_m-p_t}{2}}
\label{eq:bf_zellner}
\end{align}
where $\tilde{s}_m= l_\phi + \by'\by - \frac{\tau n}{\tau n+1} \by'X_m(X_m'X_m)^{-1}X_m'\by$  is a Bayesian analogue of $s_m$  and
\begin{align}
\tilde{F}_{mt}= \frac{(\tilde{s}_t-\tilde{s}_m)/(p_m-p_t)}{\tilde{s}_m/(n-p_m)} \leq \frac{(s_t-s_m)/(p_m-p_t)}{s_m/(n-p_m)}= F_{mt}.
\label{eq:ftest_stat}
\end{align}
$F_{mt}$ is the F-statistic to test $M_t$ versus $M_m$,
 $\tilde{F}_{mt}$ is its Bayesian analogue, 
and the inequality in \eqref{eq:ftest_stat} follows from trivial algebra.

\begin{prop}
Assume $f^*(\by)= N(\by; X_t \btheta_t^*; \phi^* I)$ and Conditions B1-B2, C1-C2.
\begin{enumerate}[leftmargin=*,label=(\roman*)]
\item Let $m \in S$  be a spurious model  and $g= (\tau n)^{(p_m-p_t)/2} p(M_t)/p(M_m)$.
If $\log(g) \ll n-p_m$, then
\begin{align}
E_{f^*}(p(M_m\mid\by))
\preceq \left(\frac{1}{g} \right)^{1 - 4 \sqrt{\frac{\log(g)}{n-p_m}}}
\ll \frac{[p(M_m)/p(M_t)]^\alpha}{(\tau n)^{\alpha \frac{p_m-p_t}{2}}}
\nonumber
\end{align}
for any fixed $\alpha <1$.
If $\log(g) \gg n-p_m$ then
$$E_{f^*}(p(M_m\mid\by)) \ll \exp\left\{- \frac{(n-p_t-5)}{2} \log \left( \frac{\log^\gamma (g)}{n-p_t-6} \right) \right\}
\ll e^{-\kappa n}$$
for any fixed $\gamma <1$, $\kappa >1$.

\item Let $m \in S^c$  be a non-spurious model  and $\lambda_{tm}$ as in \eqref{eq:ncp_regression}. Then
\begin{align}
 E_{f^*}(p(M_m \mid \by)) \ll
\max \left\{  
e^{-\lambda_{tm}^\gamma / 2}
 \left[ \frac{p(M_m) }{p(M_t) (\tau n)^{(p_m-p_t)/2} } \right]^\gamma 
,  e^{- \kappa n} \right\},
\nonumber
\end{align}
for any fixed $\gamma <1$, $\kappa >0$.
\end{enumerate}
\label{PROP:ZELLNER_PHIUNKNOWN}
\end{prop}

\subsection{{\bf Normal prior with general covariance}}
\label{ssec:normalprior}

 We extend Proposition \ref{PROP:ZELLNER_PHIUNKNOWN} to more general normal priors
$p(\btheta_k \mid \phi, M_k)= N(\btheta_k; {\bf 0}, \tau n \phi V_k)$.
The rates are essentially equivalent, subject to mild eigenvalue conditions.

Let $\rho_{k1} \geq \ldots \geq \rho_{kp_k} > 0$ be the $p_k$ non-zero eigenvalues of $V_k X_k'X_k$,
$F_{mt}$ the F-test statistic in \eqref{eq:ftest_stat},
and $\tilde{F}_{mt}$ be as in \eqref{eq:ftest_stat}
after replacing $\tilde{s}_k= l_\phi + \by'\by - \by'X_k(X_k'X_k+ (\tau n)^{-1}V_k^{-1})^{-1} X_k'\by$.
Then simple algebra gives the following expression for Bayes factors
\begin{align}
B_{tm}=
\left( 1 + \frac{p_m-p_t}{n-p_m} \tilde{F}_{mt} \right)^{-\frac{a_\phi+n}{2}}
\frac{\prod_{j=1}^{p_m} (\tau n \rho_{mj} + 1)^{\frac{1}{2}}}{\prod_{j=1}^{p_t} (\tau n \rho_{tj} + 1)^{\frac{1}{2}}}.
\label{eq:bf_normalprior}
\end{align}

 Proposition \ref{PROP:NORMALPRIOR} assumes two further technical conditions D1-D2, beyond those in Proposition \ref{PROP:ZELLNER_PHIUNKNOWN}. 
Both can be relaxed, but they simplify exposition.
D1 allows  interpreting $\tau$ as driving the prior variance in a similar fashion than for Zellner's prior. 
D2 ensures that the Bayesian-flavoured F-statistic $\tilde{F}_{mt}$ is close to the classical $F_{mt}$,
and is a mild requirement since typically $\tau n \succeq n \succeq \lambda_{t0}$.

\begin{enumerate}[leftmargin=*,label=(D\arabic*)]
\item For some constant $c_{mt}>0$,
$\prod_{j=1}^{p_m} (\tau n \rho_{mj}+1)^{\frac{1}{2}}/ \prod_{j=1}^{p_t} (\tau n \rho_{tj}+1)^{\frac{1}{2}} \asymp (c_{mt}\tau n)^{(p_m-p_t)/2}$.

\item As $n \rightarrow \infty$, $\lambda_{t0} \ll \tau n \rho_{tp_t}$, where $\lambda_{t0}$ is as in \eqref{eq:ncp_regression}.
\end{enumerate}

\begin{prop}
Assume that $f^*(\by)= N(\by; X_t \btheta_t^*; \phi^* I)$.
Consider $m \neq t$ and that Conditions B1, B2, C1, C2, D1 and D2 hold.
\begin{enumerate}[leftmargin=*,label=(\roman*)]
\item Let $m \in S$ be a spurious model. Then, for any fixed $\alpha \in (0,1)$ and $\kappa>0$,
$$E_{f^*}(p(M_m \mid \by)) \ll \max \left\{ [p(M_m)/p(M_t)]^\alpha (\tau n)^{\alpha (p_t-p_m)/2}, e^{-\kappa n}, e^{-\tau n\rho_{tp_t}/2} \right\}$$

\item Let $m \in S^c$ be a non-spurious model and $\lambda_{tm}$ as in \eqref{eq:ncp_regression}. Then, for any fixed $\gamma < 1$ and $\kappa >0$,
\begin{align}
  E_{f^*}(p(M_m \mid \by)) \ll
\max \left\{  
e^{-\lambda_{tm}^\gamma / 2}
 \left[ \frac{p(M_m) }{p(M_t) (\tau n)^{(p_m-p_t)/2} } \right]^\gamma , e^{- \kappa n}, e^{-\tau n \rho_{tp_t}/2}
\right\}.
\nonumber
\end{align}
\end{enumerate}
\label{PROP:NORMALPRIOR}
\end{prop}

\subsection{{\bf pMOM prior}}
\label{ssec:mom}

Proposition \ref{PROP:MOMPRIOR_SIMPLE} below states that, under suitable conditions, the pMOM prior attains a rate to discard spurious models featuring a term that is essentially $(\tau n)^{3(p_m-p_t)/2}$, and hence faster than the $(\tau n)^{(p_m-p_t)/2}$ shown for Normal priors.

To ease the algebra we assume that $X_k$ has zero column means and unit variances.
By Proposition 1 in \cite{rossell:2017} the Bayes factor under the pMOM prior in \eqref{eq:pmomprior} is
\begin{align}
B_{tm}= D_{tm} \left( 1 + \frac{p_m-p_t}{n-p_m} \tilde{F}_{mt} \right)^{-\frac{a_\phi+n}{2}}
\frac{\prod_{j=1}^{p_m} (\tau n \rho_{mj} + 1)^{\frac{1}{2}}}{\prod_{j=1}^{p_t} (\tau n \rho_{tj} + 1)^{\frac{1}{2}}},
\label{eq:bf_momprior}
\end{align}
where
$$
D_{tm}
=\frac{\int\int N(\btheta_t; \tilde{\btheta}_t, \phi \tilde{V}_t) \mbox{IG} \left(\phi; \frac{a_\phi+n}{2}, \frac{\tilde{s}_t}{2} \right)
  \prod_{j \in M_t}^{}   d(\theta_{tj}/\sqrt{\phi}) d\btheta_t d\phi}
{\int\int N(\btheta_m; \tilde{\btheta}_m, \phi \tilde{V}_m) \mbox{IG} \left(\phi; \frac{a_\phi+n}{2}, \frac{\tilde{s}_m}{2} \right)
  \prod_{j \in M_m}^{}  d(\theta_{mj}/\sqrt{\phi}) d\btheta_m d\phi},
$$
$d(z)=z^2/\tau$,
$\tilde{V}_k^{-1}= X_k'X_k + V_k^{-1}/(\tau n)$,
$\tilde{\btheta}_k= \tilde{V}_k X_k' \by$
and
$\tilde{F}_{mt}$, $\tilde{s}_k$ and $\rho_{kj}$ are as in \eqref{eq:bf_normalprior} for the particular case $V_k= \mbox{diag}(X_k'X_k)^{-1}$.

The Bayes factor in \eqref{eq:bf_momprior} is hence equal to that in \eqref{eq:bf_normalprior} times a penalty term $D_{tm}$ that helps penalize spurious models $m \in S$.
Intuitively, this is because the posterior distribution of $d(\theta_{mj}/\sqrt{\phi})= \theta_{mj}^2/(\phi \tau)$ concentrates at 0 for  truly spurious $\theta_{mj}^*=0$, at a rate that is at most $\sigma/\tau$, where $\sigma$ is the largest (posterior) variance in $\tilde{V}_m$.
In order to state a simple rate, Proposition \ref{PROP:MOMPRIOR_SIMPLE} assumes technical conditions E1-E5 discussed in Section \ref{proof:momprior_simple}. These can be relaxed, at the cost of a more involved expression for the bound on $E_{f^*}(p(M_m \mid \by))$.

\begin{prop}
Assume that $f^*(\by)= N(\by; X_t \btheta_t^*; \phi^* I)$.
Let $m \in S$ be a spurious model and assume that Conditions B1, C1, D1 and E1-E5 hold. Then
$$
E_{f^*}(p(M_m \mid \by)) \ll 
\max \left\{ \left(\frac{p(M_m)}{p(M_t)}\right)^\alpha \left(\frac{\tau^{3} n}{\sigma^2} \frac{\rho_{m p_m}}{\rho_{m 1}}  \right)^{-\alpha \frac{p_m-p_t}{2}}, e^{-\kappa n}, e^{-\tau n\rho_{tp_t}/2} \right\}
$$
for any fixed $\kappa >0$ and $\alpha < 1$, where $\sigma$ is the largest diagonal element in $\tilde{V}_m$.
\label{PROP:MOMPRIOR_SIMPLE}
\end{prop}

Relative to Sections \ref{ssec:zellner_phiknown_spur}-\ref{ssec:normalprior}, Proposition \ref{PROP:MOMPRIOR_SIMPLE} features an acceleration factor $\tau/\sigma$ for each truly spurious variable in $M_m$ and a term $\rho_{m p_m}^{1/2}/\rho_{m 1}^{1/2}$ involving eigenvalues. If the latter is bounded and $\sigma \asymp 1/n$ (e.g. under restricted eigenvalue conditions), the acceleration is of order $(\tau n)^{p_m-p_t}$.

Proposition \ref{PROP:MOMPRIOR_SIMPLE} is tighter than results in \cite{johnson:2012},
e.g. under uniform $p(M_k)$ we prove consistency when $p \ll (\tau n)^{\alpha/2}$ for any $\alpha<3$ (Section \ref{sec:varsel})
whereas \cite{johnson:2012} required $p \ll n$.


\subsection{{\bf Misspecified mean structure}}
\label{ssec:misspec_mean}

So far 
 we assumed that the data analyst poses a model $p(\by \mid \btheta, \phi)= N(\by; X\btheta, \phi)$ and that 
the data-generating $f^*(\by)=N(\by; X\btheta^*, \phi^* I)$ lies in the considered family.
 Although $X$ may contain non-linear basis expansions, e.g. splines or tensor products, there are practically-relevant situations where either the mean or the error structure are misspecified.
Proposition \ref{PROP:BF_NONPARAM} considers the mean misspecification case. Specifically, it considers that

$f^*(\by)= N(\by; W\bbeta^*, \xi^*I)$
for some $n \times q$ matrix $W$, $\bbeta^* \in \mathbb{R}^q$ and $\xi^* \geq 0$.
This includes situations where one did not record truly relevant variables
($X$ misses columns from $W$) or the mean of $\by$ depends on non-linearly 
 in ways that are not captured by $X$ (e.g. $X$ assumes an additive structure, whereas $W$ contains non-linear interactions). 
For simplicity we state Proposition \ref{PROP:BF_NONPARAM} for Zellner's prior but extensions to other priors follow similar lines.
Proposition \ref{PROP:BF_HETEROSKEDASTICITY_KNOWNPHI}
considers $f^*(\by)=N(\by; X_t\btheta_t^*, \phi^* \Sigma^*)$
for general $\Sigma^*$, allowing for heteroskedastic and potentially correlated errors.

 The proof strategy is as follows. The framework in Section \ref{sec:approach} applies to any $f^*$, as long as one can bound Bayes factor tail probabilities in Lemma \ref{LEM:TAILPROB_TO_MEAN}.
In Propositions \ref{PROP:BF_NONPARAM}-\ref{PROP:BF_HETEROSKEDASTICITY_KNOWNPHI} 
it is possible to bound said tails, using that $f^*$ has Gaussian errors and eigenvalues of $\Sigma^*$.
Further extensions are possible, e.g. Propositions S1-S2 in \cite{rossell:2020} deploy Lemma \ref{LEM:TAILPROB_TO_MEAN} to the case where $f^*$ has sub-Gaussian errors, e.g. when $\by$ is a binary outcome.

 Proposition \ref{PROP:BF_NONPARAM} says that the rate to discard spurious models is similar to the well-specified case (slightly sped-up by a factor $e^{\phi_t^*/\xi^*} \geq 1$). The rate for non-spurious models $m \in S^c$ vanishes exponentially in a non-centrality parameter $\lambda_{tm}$, but is exponentially slower than
a certain $\lambda_m^*$ obtained when using the correct mean structure. Specifically, 
denote by $M^*$ the true model class $N(\by; W\bbeta, \xi I)$ indexed by $(\bbeta,\xi)$.
Let $H_m=X_m(X_m'X_m)^{-1}X_m'$ be the projection matrix associated to a model $M_m$, 
and define the non-centrality parameter
\begin{align}
\lambda_{tm}= (W\bbeta^*)' H_t (I-H_m) H_t W\bbeta^*/\xi^*.
\label{eq:ncp_misspec_mean}
\end{align}
 Note that $\lambda_{tm}$ extends the non-centrality parameter in \eqref{eq:ncp_regression} to the misspecified case, by projecting the true mean $W\bbeta^*$ onto the column space of $X_t$. Similarly, 

let $\lambda_m^*= (W\bbeta^*)' (I-H_m) W\bbeta^*/\xi^*$. 
Denote the KL-optimal parameters under $M_m$ by $\btheta_m^*= (X_m'X_m)^{-1} X_m'W\bbeta^*$ (assuming full-rank $X_m$)
and  the optimal error variance by 
$$\phi_m^*=\xi^* + \frac{1}{n} (W\bbeta^*)' (I - H_m) W\bbeta^*.$$

If one were to compare $M^*$ versus $M_m$, then by Proposition \ref{PROP:ZELLNER_PHIUNKNOWN} one would select $M^*$ at an exponential rate in $\lambda_m^*$.
However, under misspecification the best one can hope for is to select $M_t$. 
When comparing $M_t$ and $M_m$, the Bayes factor for $M_m$ vanishes at an exponential rate in $\lambda_{tm} \leq \lambda_m^*$, with equality if and only if $W \bbeta^*= X\btheta_t^*$ (the mean is well-specified).

%

\begin{prop}
Let $p(\btheta_k \mid \phi_k, M_k)= N(\btheta_k; 0, \phi_k \tau n (X_k'X_k)^{-1}$ be Zellner's prior and $\alpha$, $\kappa$ be any constants satisfying $\alpha \in (0,1)$, $\kappa >0$.
Assume that $f^*(\by)=N(\by; W\bbeta^*, \xi^* I)$.
Further assume B1, B2, C1, C2 for $\lambda_{tm}$ be as in \eqref{eq:ncp_misspec_mean}, and $\phi_t^*/\xi^* \ll \log(\lambda_{tm})$.

\begin{enumerate}[leftmargin=*,label=(\roman*)]
\item Let $m \in S$. If $\log((\tau n)^{\frac{p_m-p_t}{2}} e^{\phi_t^*/\xi^*} p(M_t)/p(M_m)) \ll n-p_m$ then

$$
E_{f^*}(p(M_m\mid\by)) \ll \left[ (p(M_t)/p(M_m)) (\tau n)^{\frac{p_m-p_t}{2}} e^{\phi_t^*/\xi^*} \right]^{-\alpha}.
$$
If $\log([p(M_m)/p(M_t)] (\tau n)^{\frac{p_t-p_m}{2}}e^{\phi_t^*/\xi^*}) \gg n-p_m$ then $E_{f^*}(p(M_m\mid\by)) \ll e^{-\kappa n}$.

\item Let $m \in S^c$. If $\lambda_{tm} +\log((\tau n)^{\frac{p_m-p_t}{2}} p(M_t)/p(M_m)) \ll n-p_q$ then
\begin{align}
E_{f^*}(p(M_m \mid \by)) \ll
 \max \left\{  
e^{- \lambda_{tm}^\gamma / 2}   \left[ \frac{p(M_m) }{p(M_t) (1+\tau n)^{(p_m-p_t)/2} } \right]^\gamma, 
e^{- \kappa n}
\right\},
\nonumber
\end{align}
for any fixed $\gamma <1$, $\kappa >0$.

Further, $\lambda_{tm} \leq \lambda_m^*$, with equality if and only if $W\bbeta^*=X\btheta_t^*$.
\end{enumerate}
\label{PROP:BF_NONPARAM}
\end{prop}

As a technical remark, Proposition \ref{PROP:BF_NONPARAM} uses the minimal assumption that $\phi_t^*/\xi^* \ll \log(\lambda_{tm})$. Since the latter grows with $n$, this assumption holds in standard cases where $\xi^*$ and $\phi_t^*$ are constant and when $\phi_t^*$ decreases with $n$ 
 (e.g. $X_t$ is a non-parametric basis with growing dimension and hence lower error variance as $\phi_t^*$ as $n$ grows), but also allows for pathological cases where $\phi_t^*/\xi^*$ grows slowly with $n$. 

\subsection{{\bf Misspecified covariance structure}}
\label{ssec:misspec_covar}

We now consider the misspecified covariance case, i.e. $f^*(\by)=N(\by; X_t\btheta_t^*, \phi^* \Sigma^*)$ for positive-definite $\Sigma^*$.
Without loss of generality we constrain $\mbox{tr}(\Sigma^*)=n$, so that $\phi^*=\sum_{i=1}^{n} \mbox{Var}_{f^*}(y_i)/n$ can be interpreted as the average variance.

For simplicity we assume that $\phi^*$ is known, in analogy to Section \ref{ssec:zellner_phiknown_spur}. Extensions to unknown $\phi^*$ are possible, along the lines of the proof of Proposition \ref{PROP:ZELLNER_PHIUNKNOWN}.

We obtain rates that resemble the well-specified case, but there are potentially important differences related to 
certain eigenvalues  and an adjusted non-centrality parameter $\tilde{\lambda}_{tm}$.
Specifically, for any model $M_k$ with design matrix $X_k$ denote by $\tilde{X}_t=(I-H_k) X_t$
and by $(\underline{\omega}_{tk},\bar{\omega}_{tk})$ the smallest and largest eigenvalues of 
$\tilde{X}_t' \Sigma^* \tilde{X}_t (\tilde{X}_t' \tilde{X}_t)^{-1}$.
Consider the non-centrality parameter
\begin{align}
\tilde{\lambda}_{tm}= (\btheta_t^*)' \tilde{X}_t' \tilde{X}_t (\tilde{X}_t' \Sigma^* \tilde{X}_t)^{-1} \tilde{X}_t' \tilde{X}_t \btheta_t^*/\phi^*,
\label{eq:ncp_misspec_covar}
\end{align}
where $\tilde{X}_t= (I - H_m)X_t$. 
To gain intuition, in the well-specified case where $\Sigma^*=I$ then $\tilde{\lambda}_{tm}$ simplifies to $\lambda_{tm}$ in \eqref{eq:ncp_regression}, and also $\underline{\omega}_{tm}=\bar{\omega}_{tm}=1$.
More generally,
$
\underline{\omega}_{tm} \tilde{\lambda}_{tm} \leq \lambda_{tm} \leq \bar{\omega}_{tm} \tilde{\lambda}_{tm}.
$

Proposition \ref{PROP:BF_HETEROSKEDASTICITY_KNOWNPHI} says that spurious models are discarded at the same rate as in the well-specified case, raised to a power $1/\bar{\omega}_{tm}$. Hence, when $\bar{\omega}_{tm}$ is large, misspecifying $\Sigma^*$ can lead to a significantly slower rate.
The intuition is that $\bar{\omega}_{tm}$ measures the discrepancy between the model-based least-squares covariance $(\tilde{X}_t'\tilde{X}_t)^{-1}$ and its actual sampling covariance $(\tilde{X}_t'\tilde{X}_t)^{-1} \tilde{X}_t' \Sigma^* \tilde{X}_t  (\tilde{X}_t'\tilde{X}_t)^{-1}$.
In contrast, non-spurious models are discarded exponentially in $\tilde{\lambda}_{tm}$ so, provided $\underline{\omega}_{tm}$ is bounded, the rate remains exponential in $\lambda_{tm}$.

In contrast to Proposition \ref{PROP:BF_NONPARAM} where misspecifying the mean was guaranteed to decrease power, this need not happen when misspecifying $\Sigma^*$.

Proposition \ref{PROP:BF_HETEROSKEDASTICITY_KNOWNPHI} requires adjusting Condition C2 into C2' below.

\begin{enumerate}[leftmargin=*,label=(C\arabic*')]
\setcounter{enumi}{1}
\item Let $m \in S^c$, $\tilde{\lambda}_{tm}$ as in \eqref{eq:ncp_misspec_covar} and $M_q=M_t \cup M_m$ be the model with design matrix $X_q$
combining $X_t$ and $X_m$.
As $n\rightarrow\infty$, $[\underline{\omega}_{mq}/\bar{\omega}_{tq}] \log(\tilde{\lambda}_{tm}) \gg 1$ and
$$
\frac{\tilde{\lambda}_{tm}}{2 \log(\tilde{\lambda}_{tm})} 
+ \frac{1}{\bar{\omega}_{tq}} \left[ \frac{p_m-p_t}{2} \log(\tau n) + \log \left( \frac{p(M_t)}{p(M_m)} \right)   \right] 
 - \log p_m \gg 1.
$$
\end{enumerate}

The interpretation of C2' is similar to C2, albeit with the incorporation of eigenvalues.
The presence of the eigenvalues can be relaxed somewhat, at the expense of obtaining slower rates in Proposition \ref{PROP:BF_HETEROSKEDASTICITY_KNOWNPHI} (see the proof).

We avoid a detailed study, but note that
$n^{-1} \tilde{X}_t' \Sigma^* \tilde{X}_t$ and $n^{-1} \tilde{X}_t' \tilde{X}_t$
are sample covariance matrices. 
Under suitable assumptions (e.g. the rows of $X$ are independent draws from a Normal distribution) one can show that $\underline{\omega}_{tm}/\bar{\omega}_{tm}$ are bounded by constants with high probability,
see \cite{wainwright:2019} (Chapter 6).

\begin{prop}

Assume that $f^*(\by)=N(\by; X_t\btheta_t^*, \phi^* \Sigma^*)$ where $\Sigma^*$ is positive-definite, $\mbox{tr}(\Sigma^*)=n$ and $\phi^*$ is known.
Let $p(\btheta_k \mid \phi_k, M_k)=N(\btheta_k; 0, \phi^* \tau n (X_k'X_k)^{-1})$ be Zellner's prior.

\begin{enumerate}[leftmargin=*,label=(\roman*)]
\item Let $m \in S$ be a spurious model. Assume Conditions B1 and C1. Then
\begin{align}
E_{f^*}(p(M_m\mid\by))
\ll
\left[ \frac{[p(M_m)/p(M_t)]}{(1+\tau n)^{(p_m-p_t)/2} } \right]^{\alpha \min\{1, 1 / \bar{\omega}_{tm} \}}
\nonumber
\end{align}
for any fixed $\alpha \in (0,1)$.

\item Let $m \in S^c$ be a non-spurious model. Assume Conditions B2 and C2'. Then
\begin{align}
 E_{f^*}(p(M_k \mid \by)) \ll
\max \left\{  e^{ -\gamma \tilde{\lambda}_{tm}/2},
\left[ \frac{[p(M_m)/p(M_t)]}{(\tau n)^{(p_m-p_t)/2} } \right]^{\gamma \min\{1, 1 / \bar{\omega}_{tq} \}}
e^{-\frac{\gamma \min\{1, \bar{\omega}_{tq} \} \tilde{\lambda}_{tm}}{2}  }
\right\}.
\nonumber
\end{align}
for any fixed $\gamma < 1$,
where $\bar{\omega}_{tq}$ is the largest eigenvalue of 
$\tilde{X}_q' \Sigma^* \tilde{X}_q (\tilde{X}_q' \tilde{X}_q)^{-1}$.
\end{enumerate}

\label{PROP:BF_HETEROSKEDASTICITY_KNOWNPHI}
\end{prop}

\section{Normalized $L_0$ penalties}
\label{sec:l0penalties}

 An $L_0$ criterion proceeds by selecting the model
\begin{align}
 \hat{k}= \arg\max_k \log p(\by \mid \hat{\btheta}_k,\hat{\phi}_k) - \eta_k,
\nonumber
\end{align}

where $(\hat{\btheta}_k,\hat{\phi}_k)= \arg\max_{\btheta_k \in \Theta_k,\phi \in \Phi} p(\by \mid \btheta,\phi)$
is the maximum likelihood estimator under model $M_k$, and $\eta_k$ is a penalty that may depend on the model size $p_k$, $n$ and $p$.
For example the BIC corresponds to $\eta_k= 0.5 p_k \log(n)$, the RIC to $\eta_k= p_k \log(p)$  
and the EBIC to $\eta_k= 0.5 p_k \log(n) + \xi \log {p \choose p_k}$ for some $\xi \in (0,1)$.

We provide results analogous to Section \ref{sec:modelrates} for normalized $L_0$ methods. By normalized we refer to equivalently defining
\begin{align}
\hat{k}= \arg\max_k \frac{h(\by,k)}{\sum_{l=1}^{K} h(\by,l)},
\nonumber
\end{align}

where $h(\by,k)= p(\by \mid \hat{\btheta}_k,\hat{\phi}_k) e^{-\eta_k}$.

We refer to $\tilde{h}(\by,k)= h(\by, k)/\sum_{l=1}^{K} h(\by,l)$ as a normalized $L_0$ criterion. The idea is that, given the connection between BMS and $L_0$ penalties (see below), one could view $\tilde{h}(\by,k)$ as a pseudo-posterior probability for $M_k$ that quantifies the certainty in $\hat{k}$.
Let $M_t$ be the optimal model defined in Section \ref{sec:approach}.
Akin to \eqref{eq:mean_ppwrong}, our goal is to show that $\tilde{h}(\by,t)$ converges to 1 in the $L_1$ sense by studying
\begin{align}
\sum_{k \neq t}^{} E_{f^*} \left( \tilde{h}(\by,k) \right) \leq
\sum_{k \neq t}^{} E_{f^*} \left( \left[ 1 + h(\by,t)/h(\by,k) \right]^{-1} \right).
\label{eq:mean_ppwrong_l0penalty}
\end{align}
Note that $h(\by,t)/h(\by,k)$ is analogous to $B_{tk} p(M_t)/p(M_k)$, a product of Bayes factors and prior model probabilities.
From Proposition \ref{PROP:CORRECTSEL_FROM_PP} and Corollaries \ref{COR:POW_TYPEI_FROM_PP}-\ref{COR:POW_TYPEI_FROM_MARGPP},
\eqref{eq:mean_ppwrong_l0penalty} bounds the frequentist probability of selecting $M_t$, type I error and power.

This section is organized as follows. First, we discuss the connection between Zellner's prior and normalized $L_0$ criteria. We subsequently show our main result, Proposition \ref{PROP:L0_PENALTY}, 
which bounds $E_{f^*} ( \tilde{h}(\by,k) )$ for an individual model, analogously to Section \ref{sec:modelrates} where we bounded $E_{f^*}( p(M_k \mid \by) )$. 

Section \ref{sec:varsel} combines these bounds across models to obtain global variable selection rates.


 To see the connection between Zellner's prior and normalized $L_0$ criteria, 
in Gaussian regression straightforward algebra shows that
\begin{align}
\frac{h(\by,t)}{h(\by,k)}=
\left( 1 +  \frac{p_k-p_t}{n-p_k} F_{kt} \right)^{-\frac{n}{2}} e^{\eta_k-\eta_t} ,
\label{eq:bf_l0penalty}
\end{align}
where $F_{kt}$ is the F-test statistic in \eqref{eq:ftest_stat}.
The resemblance of \eqref{eq:bf_l0penalty} to \eqref{eq:bf_zellner} allows extending
Proposition \ref{PROP:ZELLNER_PHIUNKNOWN} to $L_0$ penalties.
Briefly, proceeding as in the proof of Proposition \ref{PROP:ZELLNER_PHIUNKNOWN} gives
\begin{align}
E_{f^*} \left( \left[ 1 + \frac{h(\by,t)}{h(\by,k)} \right]^{-1} \right) <
\int_0^1 P_{f^*} \left( (p_k-p_t) F_{mt} > 2 \frac{n-p_k}{n} \log \left[\frac{e^{\eta_k-\eta_t}}{1/u-1} \right] \right) du.
\label{eq:meanbound_l0}
\end{align}
Expression \eqref{eq:meanbound_l0} is identical to
replacing $g=(\tau n)^{(p_k-p_t)/2}p(M_t)/p(M_k)$ by $e^{\eta_k-\eta_t}$ in the proof of Proposition \ref{PROP:ZELLNER_PHIUNKNOWN}.
Proposition \ref{PROP:L0_PENALTY}(i) follows immediately,
and Proposition \ref{PROP:L0_PENALTY}(ii) is also obtained directly by proceeding as in \eqref{eq:ftail_bound_nonspur}.

We state two technical conditions C1''-C2'' required by Proposition \ref{PROP:L0_PENALTY}, which are trivial modifications of Conditions C1-C2 from Section \ref{ssec:conditions_complexity}.

\begin{enumerate}[leftmargin=*,label=(C\arabic*'')]
\item Let $m \in S$. As $n\rightarrow\infty$, $\eta_m-\eta_t \gg 1$.

\item Let $m \in S^c$. As $n\rightarrow\infty$,

\begin{align}
\frac{\lambda_{tm}}{2\log(\lambda_{tm})} + \eta_m - \eta_t - \log(p_m) \gg 1.
\nonumber
\end{align}

\end{enumerate}

Condition C1'' holds for the BIC, RIC and EBIC, and for any penalty $\eta_k$ that increases with model size $p_k$ and diverges to infinity as $n \rightarrow \infty$.
Condition C2'' is also mild.

For example, for the BIC it suffices that $\lambda_{tm}/[\log(\lambda_{tm}) p_t \log(n)] \gg 1$,
for the RIC that $\lambda_{tm}/[2 \log(\lambda_{tm}) p_t \log(p)  ] \gg 1$
and for the EBIC that $\lambda_{tm}/[\log(\lambda_{tm}) p_t \log(n^{1/2} p^\xi)] \gg 1$.
See Section \ref{ssec:conditions_complexity} for discussion why these conditions are near-minimal.

\begin{prop}
Assume that $f^*(\by)= N(\by; X_t \btheta_t^*; \phi^* I)$.
Consider $m \neq t$ and that Conditions B1, B2, C1'' and C2'' hold.
\begin{enumerate}[leftmargin=*,label=(\roman*)]
\item Let $m \in S$.
If $\eta_m-\eta_t \ll n-p_m$ then 
\begin{align}
E_{f^*} \left( \tilde{h}(\by,m) \right) \leq e^{-(\eta_m-\eta_t) \left(1 - 4 \sqrt{\frac{\eta_m-\eta_t}{n-p_m}}\right)}.
\nonumber
\end{align}
for all $n \geq n_0$, where $n_0$ is fixed and does not depend on $m$.
If $\eta_m-\eta_t \gg n-p_m$ then $E_{f^*} \left( \tilde{h}(\by,m) \right) < e^{-\kappa n}$ for any fixed $\kappa>0$ and $n \geq n_0$.

\item Let $m \in S^c$. Then, for any fixed $\gamma <1$, $\kappa >0$,
$$E_{f^*} \left( \tilde{h}(\by,m)  \right) < 
\max \left\{
e^{- \lambda_{tm}^\gamma / 2} e^{-\gamma (\eta_m - \eta_t)}   ,  e^{-\kappa n} \right\}
$$
for all $n \geq n_0$ where $n_0$  is fixed and does not depend on $m$.
\end{enumerate}
\label{PROP:L0_PENALTY}
\end{prop}

For example, for the BIC $\eta_m= 0.5 p_m \log(n)$,
then $E_{f^*}(\tilde{h}(\by,m))$ vanishes essentially at a rate $n^{-(p_m-p_t)/2}$ for spurious models $m \in S$,
and a faster $e^{-\lambda_{tm}/2} n^{-(p_m-p_t)/2}$ for non-spurious models.
This is no surprise, as discussed the BIC is essentially identical to using a uniform model prior and setting Zellner's prior dispersion to $\tau=1$, hence one obtains the same rates.
Similarly, the RIC is essentially identical to $\tau n=p^2$ and uniform $p(M_m)$,
and the EBIC (for the particular choice $\xi=1$) to $\tau=1$ and Beta-Binomial $p(M_m)$.
Note also that the model misspecification results for Zellner's prior in Sections \ref{ssec:misspec_mean}-\ref{ssec:misspec_covar} extend directly to normalized $L_0$ penalties.

\section{Global rates for regression}
\label{sec:varsel}

We now use the model-specific bounds from Sections \ref{sec:modelrates}-\ref{sec:l0penalties} to obtain global bounds. 
We saw that, under suitable conditions,
$E_{f^*}(1 - p(M_t \mid \by))=E_{f^*} (P(S \mid \by)) + E_{f^*}(P(S^c \mid \by))$
\begin{align}
\leq
 \sum_{l= p_t+1}^{\bar{p}} \sum_{k \in S_l} \left[ \frac{p(M_k)}{p(M_t) (\tau n)^{\frac{(p_k-p_t)}{2}}} \right]^\alpha
  +  \sum_{l= 0}^{\bar{p}} \sum_{k \in S_l^c} e^{-\frac{\lambda_{kt}^\alpha}{2}} \left[ \frac{p(M_k) }{p(M_t) (\tau n)^{\frac{(p_k-p_t)}{2}}} \right]^\alpha,
\label{eq:pp_growingp}
\end{align}
for sufficiently large $n$ and some fixed $\alpha$, where $\tau$ is the prior dispersion.

In well-specified Gaussian regression, as well as with a misspecified mean structure, we showed that one can take essentially take $\alpha=1$ (up to lower-order terms). For the pMOM prior in the first term of \eqref{eq:pp_growingp} one may take a larger $\alpha <3$.
Similarly, for normalized $L_0$ criteria,

\begin{align}
E_{f^*}(1 - \tilde{h}(t, \by) ) \leq
 \sum_{l= p_t+1}^{\bar{p}} \sum_{k \in S_l}
e^{-(\eta_k-\eta_t) \left(1 - 4 \sqrt{\frac{\eta_k-\eta_t}{n-p_k}}\right)}
+  \sum_{l= 0}^{\bar{p}} \sum_{k \in S_l^c}  e^{- \frac{\lambda_{tk}^\alpha}{2} -\alpha (\eta_k - \eta_t)},
\label{eq:l0_growingp}
\end{align}
where $\eta_k$ is the $L_0$ penalty, e.g. for the BIC $\eta_k= 0.5 p_k \log(n)$.

The bounds in Sections \ref{sec:modelrates}-\ref{sec:l0penalties} also feature terms such as $e^{- \kappa n}$ that vanish exponentially with $n$. For simplicity we omitted these, since they are typically of a smaller order, but they can easily be plugged into \eqref{eq:pp_growingp}-\eqref{eq:l0_growingp}.

This section derives simpler asymptotic expressions for \eqref{eq:pp_growingp}-\eqref{eq:l0_growingp} for the uniform, Beta-Binomial and Complexity priors in \eqref{eq:prior_model}, and for the BIC, RIC and EBIC.


We study separately spurious and non-spurious models, i.e. $E_{f^*}(p(S \mid \by))$ and $E_{f^*}(p(S^c \mid \by))$,
and we also discuss the effect of setting priors or $L_0$ penalties that are not particularly sparse.
Such priors attain worse asymptotic rates to discard spurious models, but they can also significantly improve finite $n$ performance.

The reason for the mismatch between asymptotic and finite $n$ performance is that $E_{f^*}(p(S^c \mid \by))$ is typically negligible for large $n$,
as it vanishes exponentially under eigenvalue conditions.
However, $E_{f^*}(p(S^c \mid \by))$ can be large for finite $n$, particularly when optimal model is not sparse.
See Section \ref{sec:varsel_examples} for examples.

\subsection{Uniform prior, spurious models}
\label{ssec:varsel_spur_unif}

The uniform prior sets $p(M_k)/p(M_t)=1$.
From the first term in \eqref{eq:pp_growingp},
using that there are $|S_l| = {p - p_t \choose l-p_t}$ spurious models of size $l$ 
and the geometric series, one obtains
\begin{align}
E_{f^*}(P(S \mid \by)) \leq
  \frac{\frac{p-p_t}{(\tau n)^{\alpha/2}} - \left(\frac{p-p_t}{(\tau n)^{\alpha/2}}\right)^{\bar{p}-p_t+1} }{1-(p-p_t)/(\tau n)^{\alpha/2}}
  \asymp
  \frac{p-p_t}{(\tau n)^{\alpha/2}},
\label{eq:Sbound_unif}
\end{align}
for sufficiently large $n$.
The asymptotic expression in the right-hand side of \eqref{eq:Sbound_unif} holds if $p - p_t \ll (\tau n)^{\alpha/2}$, i.e. when $E_{f^*}(P(S \mid \by))$ converges to 0.
Rates for the BIC and RIC are obtained by plugging $\tau=1$ and $\tau n=p^2$ into \eqref{eq:Sbound_unif}.

%

Expression \eqref{eq:Sbound_unif} describes the effect of the prior dispersion $\tau$ on sparsity.
For example, for $\tau=1$ then $P(S \mid \by)$ vanishes as long as $p-p_t \ll n^{1/2}$, under Zellner and Normal priors. Under the pMOM one can handle $p \ll n^{3/2}$.
Another default is  $\tau =\mbox{max}\{1,p^{2+a}/n\}$ for some small $a > 0$ \citep{fernandez:2001},
 which effectively sets a diffuse prior ($\tau$ grows with $n$, whenever $p \gg \sqrt{n}$).
Under such a diffuse prior,  $p - p_t \ll (\tau n)^{\alpha/2}$ and $P(S\mid\by)$ vanishes under Zellner's, Normal and pMOM priors, regardless of the magnitude of $p$.

\subsection{Beta-Binomial prior, spurious models}
\label{ssec:varsel_spur_betabin}

The Beta-Binomial prior sets $p(M_m)/p(M_t)= {p \choose p_t} / {p \choose p_m}$.
Using simple algebra and the binomial coefficient's ordinary generating function,
\begin{align}
E_{f^*}(P(S \mid \by)) 
  <\left[ 1 - \frac{(p-p_t)^{1-\alpha}}{(\tau n)^{\frac{\alpha}{2}}} \right]^{-p_t-1} -1
 \asymp \frac{(p_t+1) (p-p_t)^{1-\alpha}}{(\tau n)^{\alpha/2}},
\nonumber
\end{align}
where the right-hand side holds if $(\tau n)^{\alpha/2} \gg (p_t+1) (p-p_t)^{1-\alpha}$, by l'Hopital's rule.
If $\alpha$ is arbitrarily close to 1,
$P(S \mid \by)$ vanishes as long as $p_t^{a + \epsilon} (p-p_t) \ll (\tau n)^{a/2}$ for arbitrarily large but fixed $a>0$ and any small $\epsilon>0$.
For instance, under $\tau=1$ one can handle $p-p_t \ll n^{a/2}$ variables, i.e. $p$ can grow polynomially with $n$ (provided $p_t \ll n$ grows sub-linearly in $n$, as in Condition B1.

Therefore, despite not necessarily leading to the optimal asymptotic rate, the Beta-Binomial prior can handle problems of fairly large dimension and still discard all spurious models.

We remark that one can obtain slightly tighter rates for $L_0$ penalties and for specific priors.
For Zellner's prior and known $\phi^*$ Lemma \ref{LEM:CHISQTAIL_INTBOUND} gives
$E_{f^*}(p(M_m\mid\by)) \preceq [\log(g)]^{(p_m-p_t)/2+1} / g$, where $g=(\tau n)^{(p_m-p_t)/2} p(M_t)/p(M_m)$,
then Lemma \ref{LEM:SUMSPUR_ZELLNERKNOWN} shows that
\begin{align}
E_{f^*}(P(S\mid\by))
  \preceq \frac{(p_t+1) (\bar{p}-p_t)^{a/2} \log^{3/2}((\tau n)^{1/2} (p-p_t))}{(\tau n)^{1/2}}
  \nonumber
\end{align}
for any fixed $a>1$, i.e. the dependence on $p$ is now logarithmic.
Similarly, for unknown $\phi^*$ and Zellner's prior
Lemma \ref{LEM:SUMSPUR_ZELLNERUNKNOWN} gives that
\begin{align}
  E_{f^*}( P(S \mid \by)) \preceq
  \frac{(p_t+1)}{(\tau n)^{1/2}} e^{2 [\log^{3/2} ((\tau n)^{1/2} (p-p_t))] \sqrt{(p-p_t)/(n-\bar{p})  }}.
\nonumber
\end{align}
Rates for the EBIC are obtained by plugging $\tau=1$ into this last expression.

\subsection{Complexity prior, spurious models}
\label{ssec:varsel_complexity}

Here
$p(M_m)/p(M_t) \asymp p^{c(p_t-p_m)} {p \choose p_t} / {p \choose p_m}$, where $c$ is the Complexity prior's parameter in \eqref{eq:prior_model}. Simple algebra shows that
\begin{align}
E_{f^*}(P(S \mid \by)) \preceq
  \sum_{l=p_t+1}^{\bar{p}} {l \choose p_t} \left( \frac{(p-p_t)^{1-\alpha}}{(\tau n)^{\frac{\alpha}{2}}p^c } \right)^{l-p_t} 
\preceq \frac{(p_t+1) (p-p_t)^{1-\alpha}}{(\tau n)^{\alpha/2} p^c}.
\nonumber
\end{align}
Since $\alpha$ is arbitrarily close to 1, $P(S \mid \by)$ vanishes
under the minimal requirement that $p_t^{1+\epsilon} \ll p^{c} (\tau n)^{1/2}$ for some (small) fixed $\epsilon >0$.
That is, the complexity prior can handle almost any $p$ and still excel at discarding spurious models, even with moderately small $c>0$.

However, as illustrated next, $c$ also plays a role in slowing down the rate at which one discards small non-spurious models, which can reduce the statistical power to detect non-zero coefficients.

\subsection{Non-spurious models}
\label{ssec:varsel_nonspur}

Our main result is Proposition \ref{PROP:VARSEL_NONSPUR}, which gives rates for the total posterior probability assigned to models of size $p_m <p_t$ (smaller than the optimal $M_t$) and to those of size $p_m \geq p_t$.
The rates depends on two parameters $(\underline{\lambda},\bar{\lambda})$ that bound uniformly the non-centrality parameters $\lambda_{tm}$.
In Gaussian regression, subject to restricted eigenvalue conditions, $(\underline{\lambda},\bar{\lambda})$ are roughly proportional to $n$ and a beta-min parameter.
We first define $(\underline{\lambda},\bar{\lambda})$, and then present Proposition \ref{PROP:VARSEL_NONSPUR}.

First, define $\underline{\lambda}= \min_{p_m < p_t} \lambda_{tm}^\alpha / (p_t-p_m)$, where $\alpha$ is as in \eqref{eq:pp_growingp}.
In Gaussian regression Lemma \ref{LEM:LAMBDABOUND_SMALLMODELS} shows that
one may set $\underline{\lambda}= [n \underline{v} \min_j (\theta_j^*)^2/\phi^*]^\alpha$, where $\underline{v}$ is the smallest eigenvalue $v_{tm}$ across models of size $p_m < p_t$.
Regarding $\bar{\lambda}$, let $S_{l,j}^c \subseteq S_l^c$ be the set of non-spurious models $M_m$ of size $p_m=l$ that contain $j$ truly active parameters (non-zero elements in $\btheta^*$).
Let $\bar{\lambda} = \min_{j \geq p_t, m \in S_{l,j}^c} \lambda_{tm}^\alpha/(p_t-j)$ be an analogous quantity to $\underline{\lambda}$, when taking the minimum over $m \in S_{l,j}^c$.
Lemma \ref{LEM:LAMBDABOUND_SMALLMODELS} shows that one may set $\bar{\lambda}= [n \bar{v} \min_j (\theta_j^*)^2/\phi^*]^\alpha$, where $\bar{v}$ is the smallest $v_{tm}$ across models of size $p_m \in [p_t, \bar{p}]$.


\begin{lemma}
Let $\lambda_{tm}=(\btheta_t^*)' X_t' (I-H_m) X_t \btheta_t^*/\phi^*$ in \eqref{eq:ncp_regression}, where $H_m=X_m(X_m'X_m)^{-1}X_m'$,
and $v_{tm}$ be the smallest non-zero eigenvalue of $X_t' (I-H_m) X_t/n$.

\begin{enumerate}[leftmargin=*,label=(\roman*)]
\item Let $m \in S^c$ be a model of size $p_m<p_t \leq n$.
If $X_t$ has rank $p_t$, then
$\lambda_{tm} \geq n v_{tm} (p_t-p_m) \min_{j} (\theta_{tj}^*)^2/\phi^*$,
where 

\item Let $m \in S_{p_m,j}^c$ be a model of size $p_m \geq p_t$ with design matrix $X_m$ containing $j$ columns from $X_t$.
Let $X_q$ be the union of the columns in $X_t$ and the $p_m - j$ remaining columns in $X_m$.
If $X_q$ has full rank, then
$\lambda_{tm} \geq n v_{tm} (p_t - j) \min_j (\theta_{j}^*)^2/\phi^*$.

\end{enumerate}

\label{LEM:LAMBDABOUND_SMALLMODELS}
\end{lemma}

Proposition \ref{PROP:VARSEL_NONSPUR} Part (i) requires Condition (F1) below,
which ensures that Condition C2 in Section \ref{ssec:conditions_complexity} holds uniformly across models smaller than $M_t$.
For example, for the uniform prior $p(M_k)$, (F1) is the mild requirement that $\underline{\lambda}/2 + \log p - 0.5 \log(n \tau) \gg 1$.
For the Beta-Binomial (F1) requires $\underline{\lambda}/2 + (1 - \alpha) \log p - 0.5 \log(n \tau) \gg 1$.
(F1) is more stringent for the Complexity and for large prior dispersion $\tau$.
These priors give higher support to small models, and hence require a stronger signal.
Interestingly, (F1) is not needed for Part (ii). There, by setting sufficiently sparse priors (large $\tau$ or $c$) one may discard models of size $> p_t$. In particular, one could potentially set the maximum model complexity to $\bar{p}>n$ and still attain convergence in Part (ii).

\begin{enumerate}[leftmargin=*,label=(F\arabic*)]

\item 
Assume that
$\lim_{n \rightarrow \infty} \underline{\lambda}/2 - (\alpha (1+c) - 1) \log p - 0.5 \log(n\tau) = \infty$
holds for $c= -1$ when $p(M_k)$ is the uniform prior, $c=0$ when it is the Beta-Binomial and $c>0$ when it is the Complexity(c) prior in \eqref{eq:prior_model}.

\end{enumerate}

\begin{prop}
Let $p(M_k)$ be either the uniform or the Complexity(c) prior in \eqref{eq:prior_model}, where $c=0$ corresponds to the Beta-Binomial prior.
Assume that for all non-spurious $m \in S^c$ it holds that
\begin{align}
E_{f^*}( p(M_m \mid \by)) \leq  e^{- \frac{\lambda_{tm}^\alpha}{2}}  (n\tau)^{-\alpha(p_m - p_t)/2}  \left(\frac{p(M_m)}{p(M_t)}\right)^\alpha 
\nonumber
\end{align}
for some $\alpha <1$ and all $n \geq n_0$, where $n_0$ is fixed.

\begin{enumerate}[leftmargin=*,label=(\roman*)]
\item Assume that (F1) holds. Then, for the Complexity prior 
\begin{align}
\lim_{n \rightarrow \infty} E_{f^*}\left( \sum_{p_m=0}^{p_t-1} P(S_l^c \mid \by) \right) \leq
e^{-\frac{\underline{\lambda}}{2} + [p_t -1 +  \alpha (1+c)  ] \log p   + \frac{\alpha}{2} \log(n\tau)   }
\nonumber
\end{align}
for all $n \geq n_0$.  The result for the uniform prior is obtained by setting $c=-1$. 

\item Suppose that $\lim_{n \rightarrow \infty} \bar{\lambda}/2 + \log p_t - \log(p - p_t)= \infty$. Then
\begin{align}
\lim_{n \rightarrow \infty} E_{f^*}\left( \sum_{p_m=p_t}^{\bar{p}} P(S_l^c \mid \by) \right) \leq
 e^{- \frac{\bar{\lambda}}{2} + p_t \log(pe)} + 
 \frac{ e^{-\bar{\lambda}/2 + p_t \log p_t + \log p}}{[ (n\tau)^{\alpha/2}  p^{\alpha(c+1)-1}]^{\bar{p} - p_t}}.
\nonumber
\end{align}
for all $n \geq n_0$.
 
If $\lim_{n \rightarrow \infty} \bar{\lambda}/2 + \log p_t - \log(p - p_t)= -\infty$, then
\begin{align}
\lim_{n \rightarrow \infty} E_{f^*}\left( \sum_{p_m=p_t}^{\bar{p}} P(S_l^c \mid \by) \right)
 \leq e^{-p_t\bar{\lambda}/2}
 \left( \frac{1}{n \tau} \right)^{\frac{\alpha (\bar{p} - p_t)}{2}} \left( \frac{1}{p} \right)^{\alpha (c+1) (\bar{p}-p_t) -1}.
\nonumber
\end{align}
The results for the uniform prior are obtained by setting $c=-1$ above. 
\end{enumerate}
\label{PROP:VARSEL_NONSPUR}
\end{prop}

\section{Empirical examples}
\label{sec:varsel_examples}

We illustrate the effect of the prior formulation and signal strength on linear regression rates with two simple studies.
Section \ref{ssec:lm_ortho} shows simulated data under orthogonal $X'X$
and Section \ref{ssec:lm_corpred} a setting where all pairwise correlations are 0.5,
in both cases covariates are normally distributed with zero mean and unit variance.
We considered three prior formulations:
Zellner's prior $(\tau=1)$ coupled with either a Complexity($c=1$) or Beta-Binomial(1,1) priors on the model space,
and the pMOM prior (default $\tau=0.348$ from \cite{johnson:2010}) coupled with a Beta-Binomial(1,1).
For the error variance we set $p(\phi \mid M_k) \sim \mbox{IG}(0.005,0.005)$.
In Section \ref{ssec:lm_ortho} we used the methodology in \cite{papaspiliopoulos:2017}
to obtain exact posterior probabilities,
and in Section \ref{ssec:lm_corpred} the Gibbs sampling algorithm from \cite{johnson:2012}
(functions \texttt{postModeOrtho} and \texttt{modelSelection} in R package \texttt{mombf}, respectively)
with 10,000 iterations (i.e. $10^4\times p$ variable updates) after a 1,000 burnin.

\subsection{Orthogonal design}
\label{ssec:lm_ortho}

\begin{figure}
\begin{center}
\begin{tabular}{cc}
$p_t=5$, $p=100$, $n=110$ & $p_t=20$, $p=100$, $n=110$ \\
\includegraphics[width=0.48\textwidth]{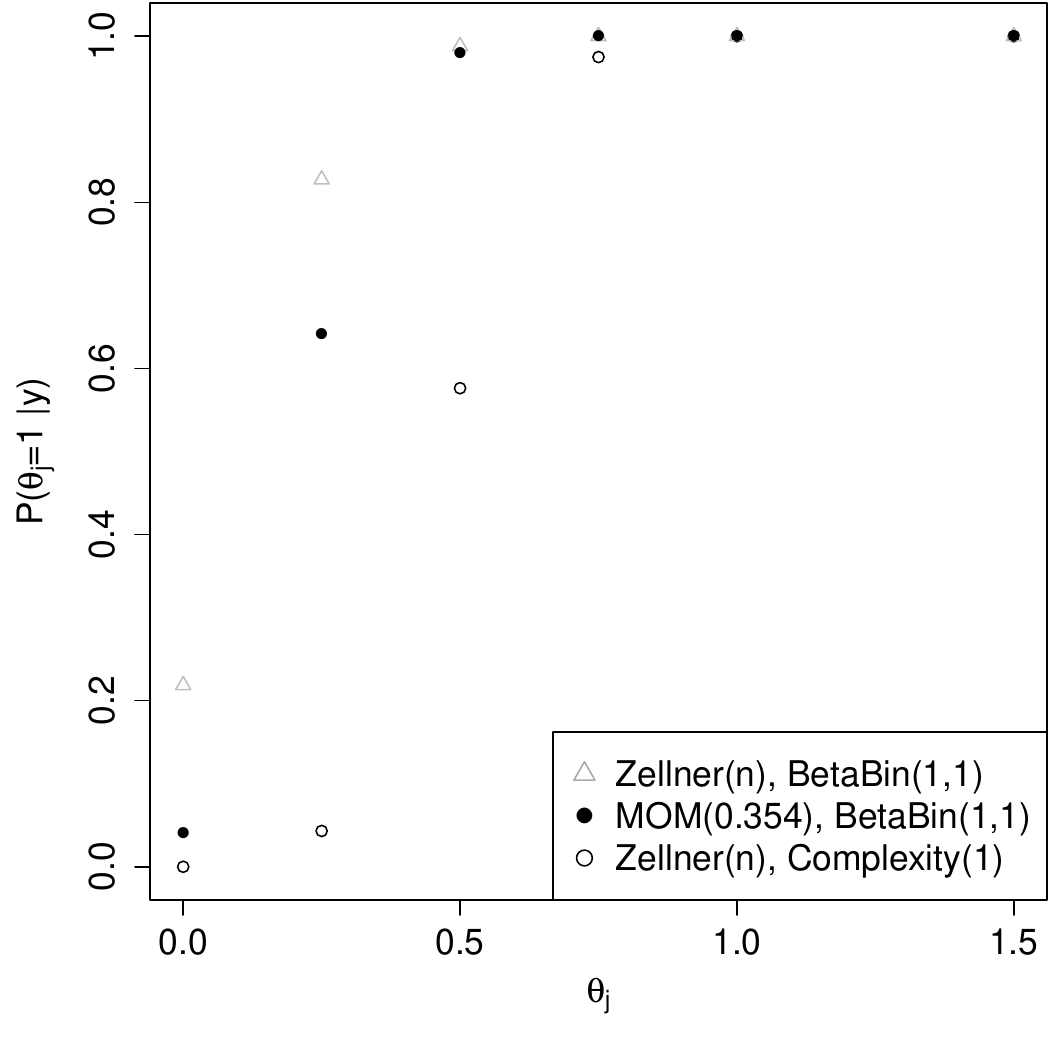} &
\includegraphics[width=0.48\textwidth]{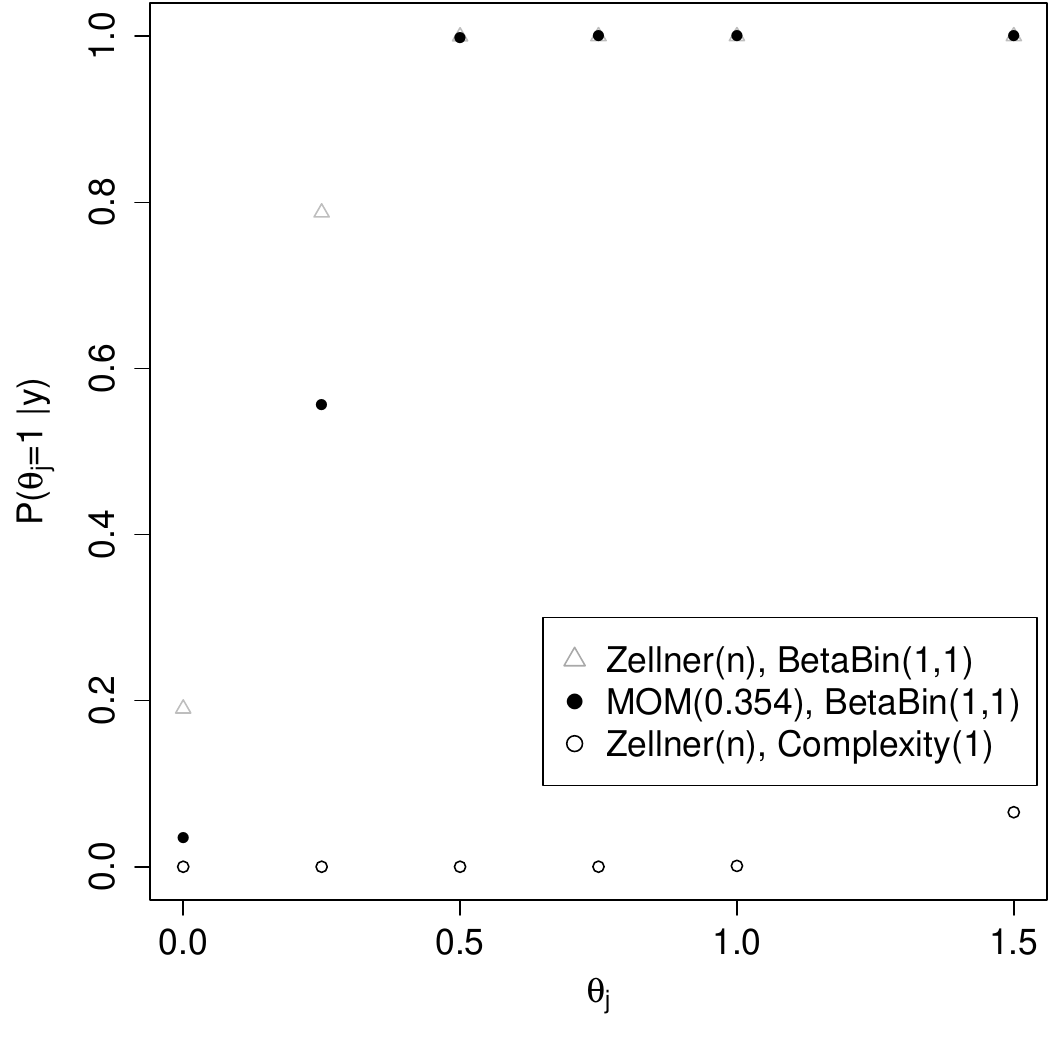} \\
$p_t=5$, $p=500$, $n=510$ & $p_t=20$, $p=500$, $n=510$ \\
\includegraphics[width=0.48\textwidth]{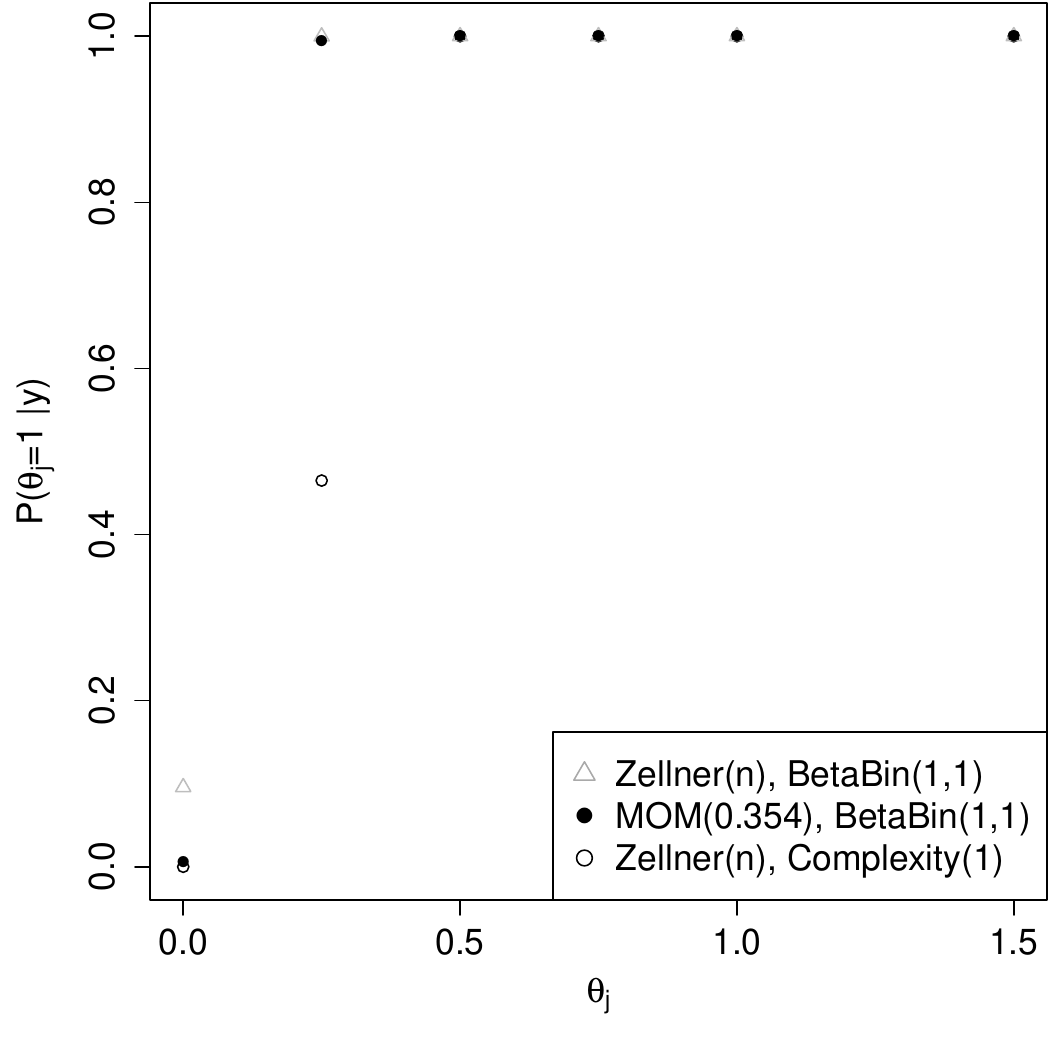} &
\includegraphics[width=0.48\textwidth]{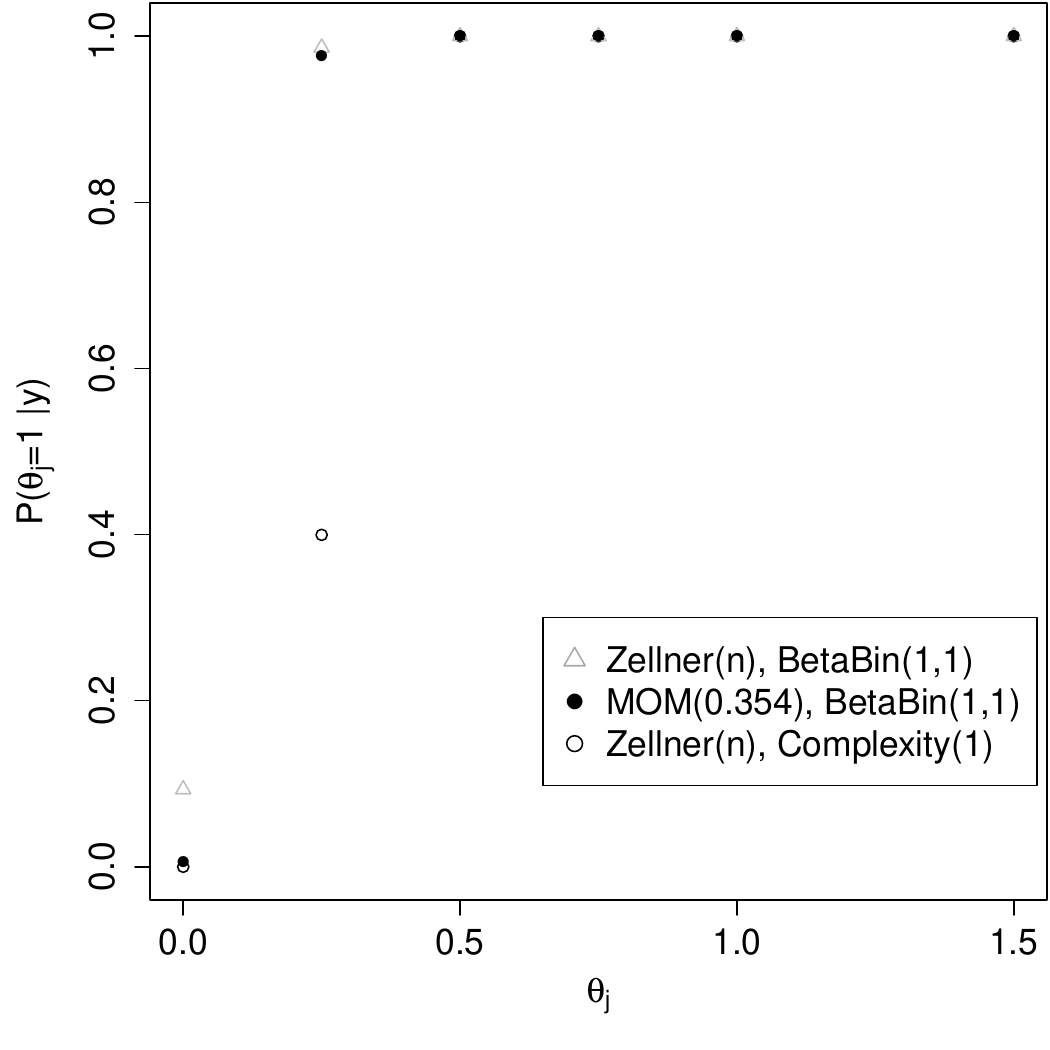}
\end{tabular}
\end{center}
\caption{Average marginal inclusion probabilities under orthogonal $X'X$ and $\phi^*=1$
for three prior formulations:
Zellner-Complexity(1), Zellner-Beta-Binomial(1,1), pMOM-Beta-Binomial(1,1).
For both Zellner's and pMOM priors $\tau$ was set to obtain unit prior variance ($\tau=1$, $\tau=0.348$)}
\label{fig:margpp_ortho}
\end{figure}

We considered four scenarios and simulated 100 independent datasets under each.
In Scenario 1 we set $p=100$, $n=105$ and $p_t=5$ truly active variables with coefficients
$\theta_j^*=0.25,0.5,0.75,1,1.5$ for $j=1,\ldots,p_t$.
In Scenario 2 again $p=100$, $n=105$ but coefficients were less sparse,
we set $p_t=20$ by repeating four times each coefficient in Scenario 1,
i.e. $\theta_j^*=0.25,0.25,0.25,0.25,\ldots,1.5,1.5,1.5,1.5$ for $j=1,\ldots,p_t$.
Scenarios 3-4 were identical to Scenarios 1-2 (respectively) setting $p=500$ and $n=510$.
The true error variance was $\phi^*=1$ under all scenarios.

Figure \ref{fig:margpp_ortho} shows marginal inclusion probabilities $P(\theta_j \neq 0 \mid \by)$.
The Zellner-Complexity prior gave the smallest inclusion probabilities to truly inactive variables ($\theta_j^*=0$),
but incurred a significant loss in power to detect truly active variables.
In agreement with our theory this drop was particularly severe for $p_t=20$,
e.g. when $n=110$ inclusion probabilities were close to 0 even for fairly large coefficients.
Also as predicted by the theory the power increased for $(n,p)=(510,500)$ under all priors,
but under the Zellner-Complexity prior it remained low for $\theta_j^*=0.25$.
The MOM-Beta-Binomial prior showed a good balance between power and sparsity,
although for $n=100$ it had slightly lower power to detect $\theta_j^*=0.25$ relative to the Zellner-Beta-Binomial.

\subsection{Correlated predictors}
\label{ssec:lm_corpred}

\begin{figure}
\begin{center}
\begin{tabular}{cc}
Scenario1: $p_t=10$, $p=n$, $\theta_j^*=0.5$ & Scenario 2: $p_t=10$, $p=n$, $\theta_j^*=0.25$ \\
\includegraphics[width=0.45\textwidth,height=0.4\textwidth]{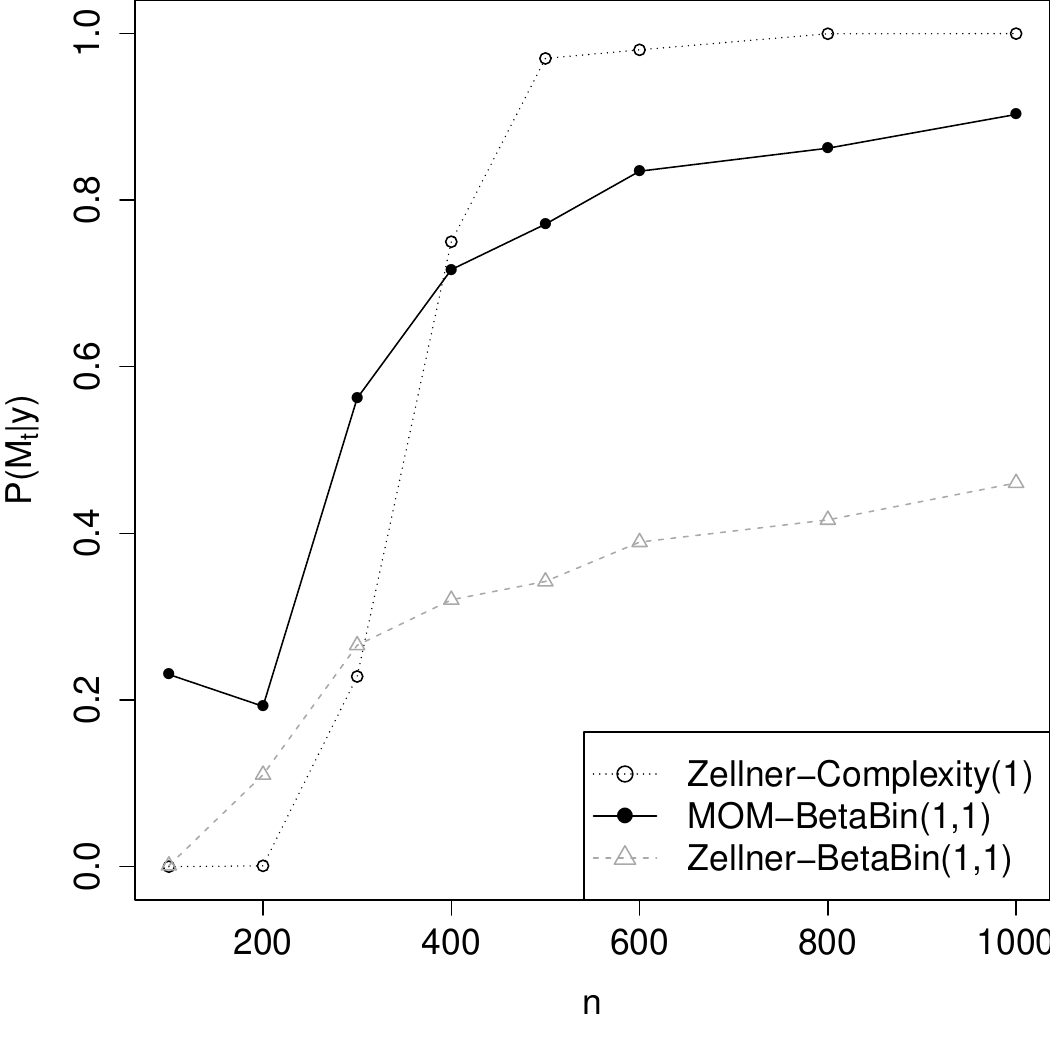} &
\includegraphics[width=0.45\textwidth,height=0.4\textwidth]{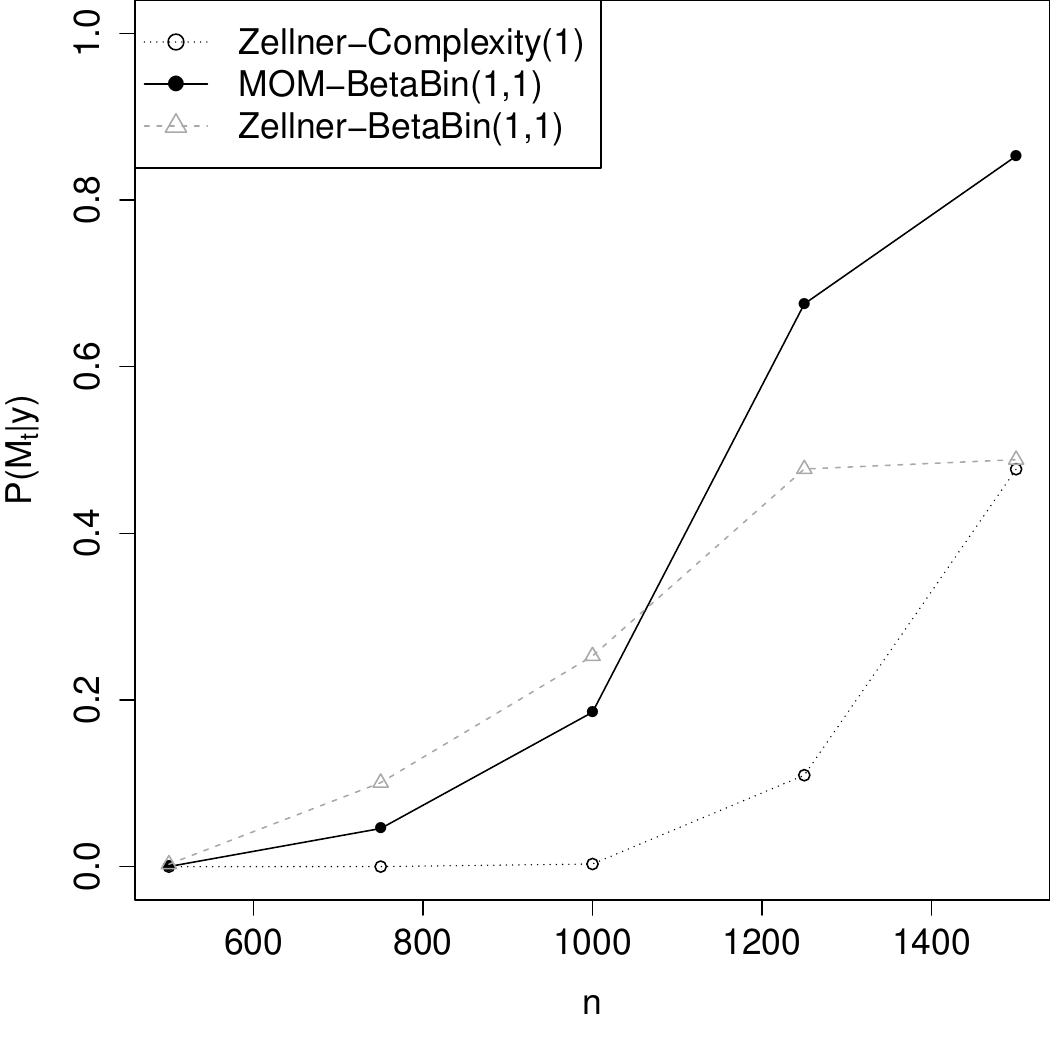} \\
\includegraphics[width=0.45\textwidth,height=0.4\textwidth]{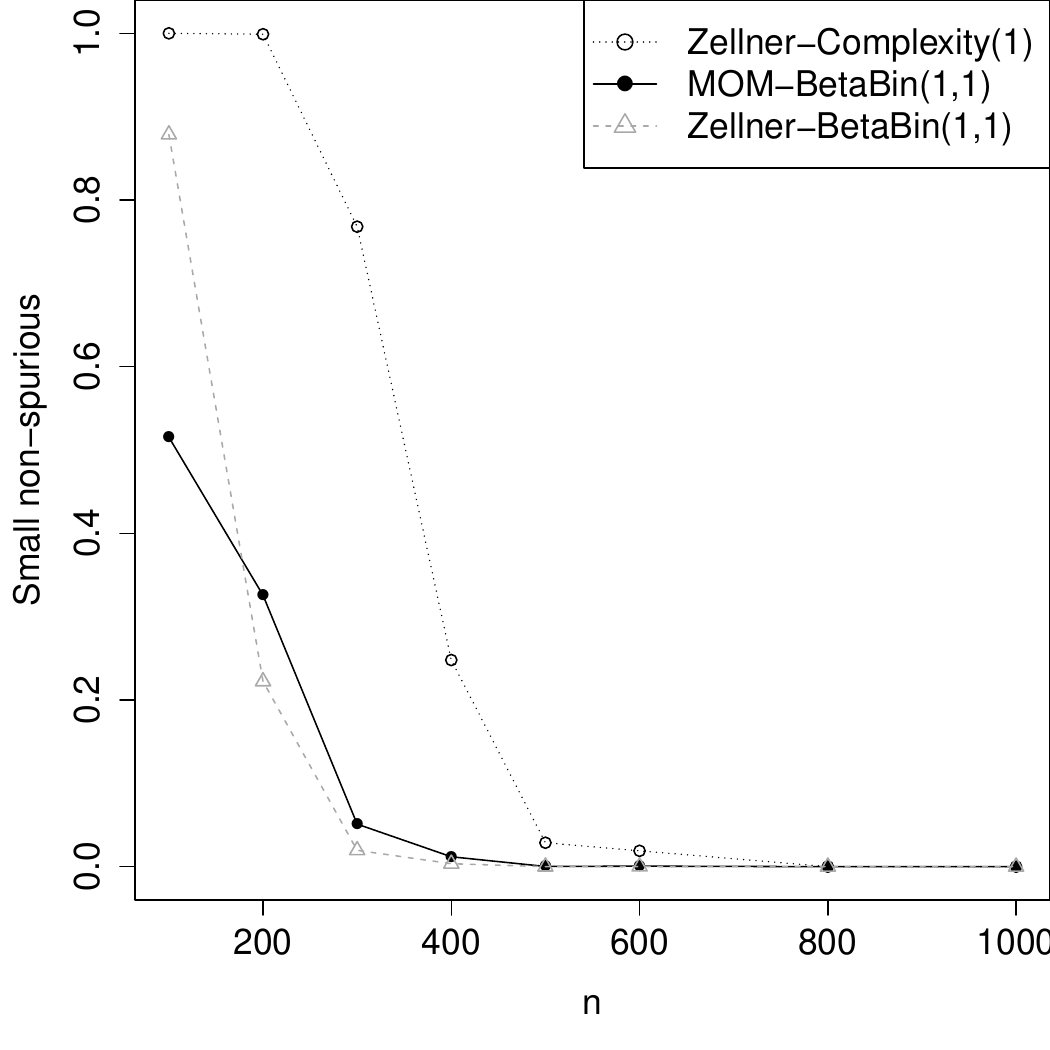} &
\includegraphics[width=0.45\textwidth,height=0.4\textwidth]{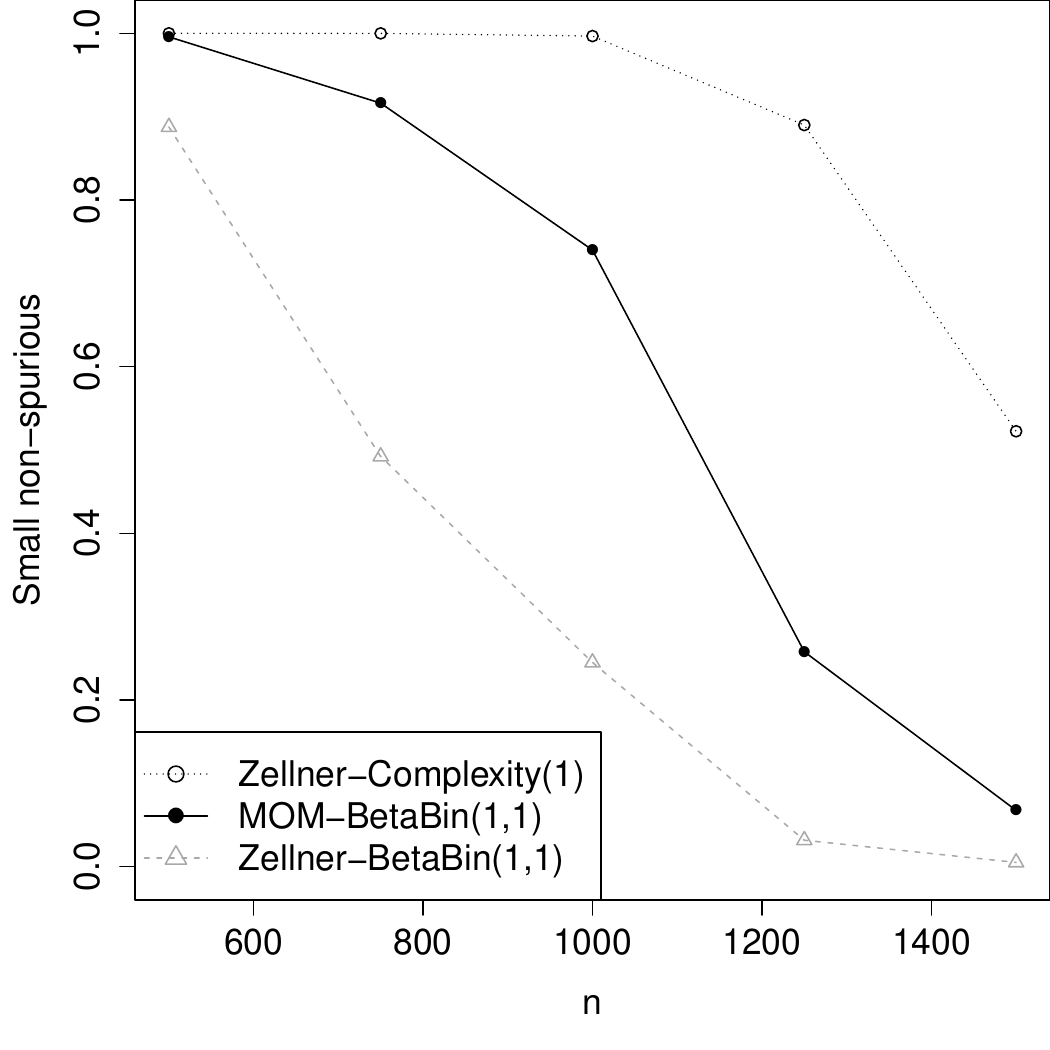} \\
\includegraphics[width=0.45\textwidth,height=0.4\textwidth]{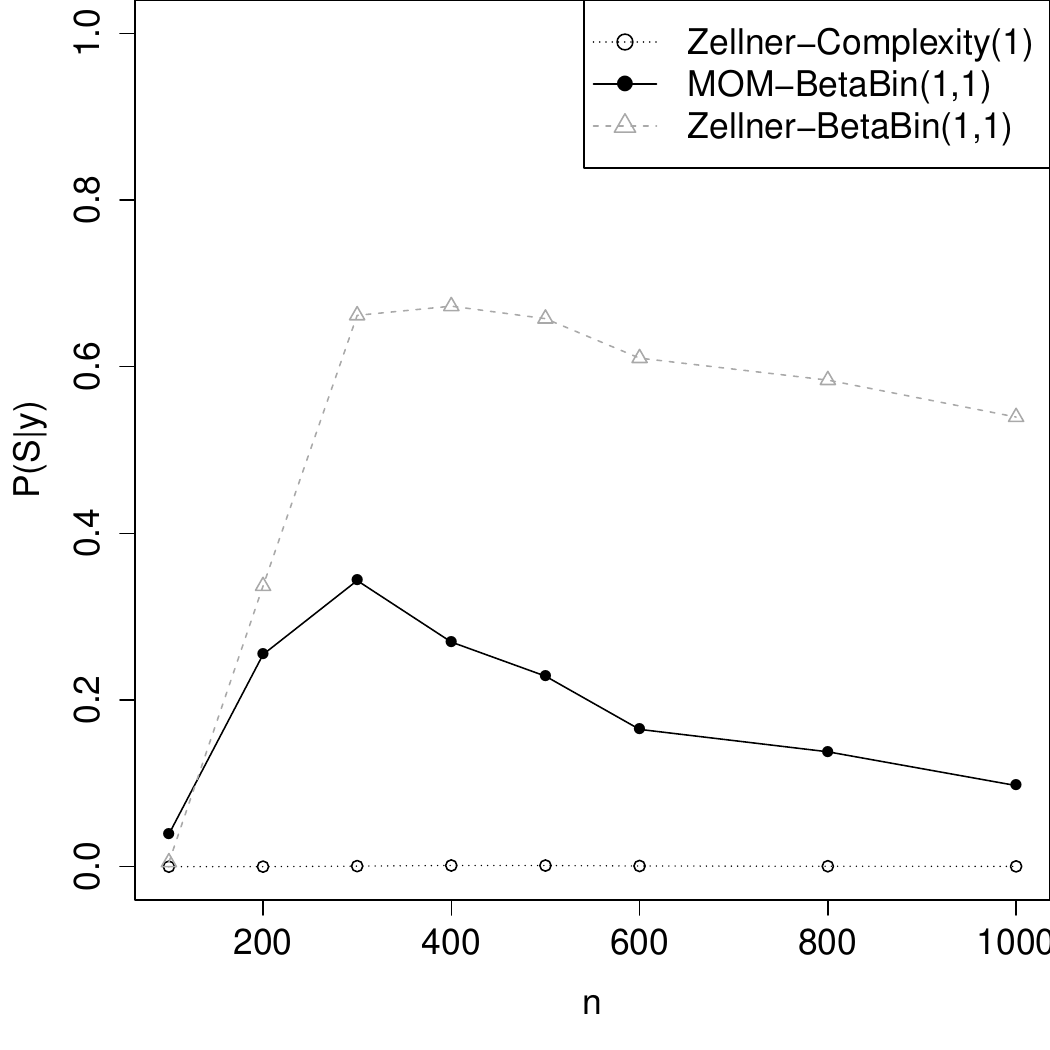} &
\includegraphics[width=0.45\textwidth,height=0.4\textwidth]{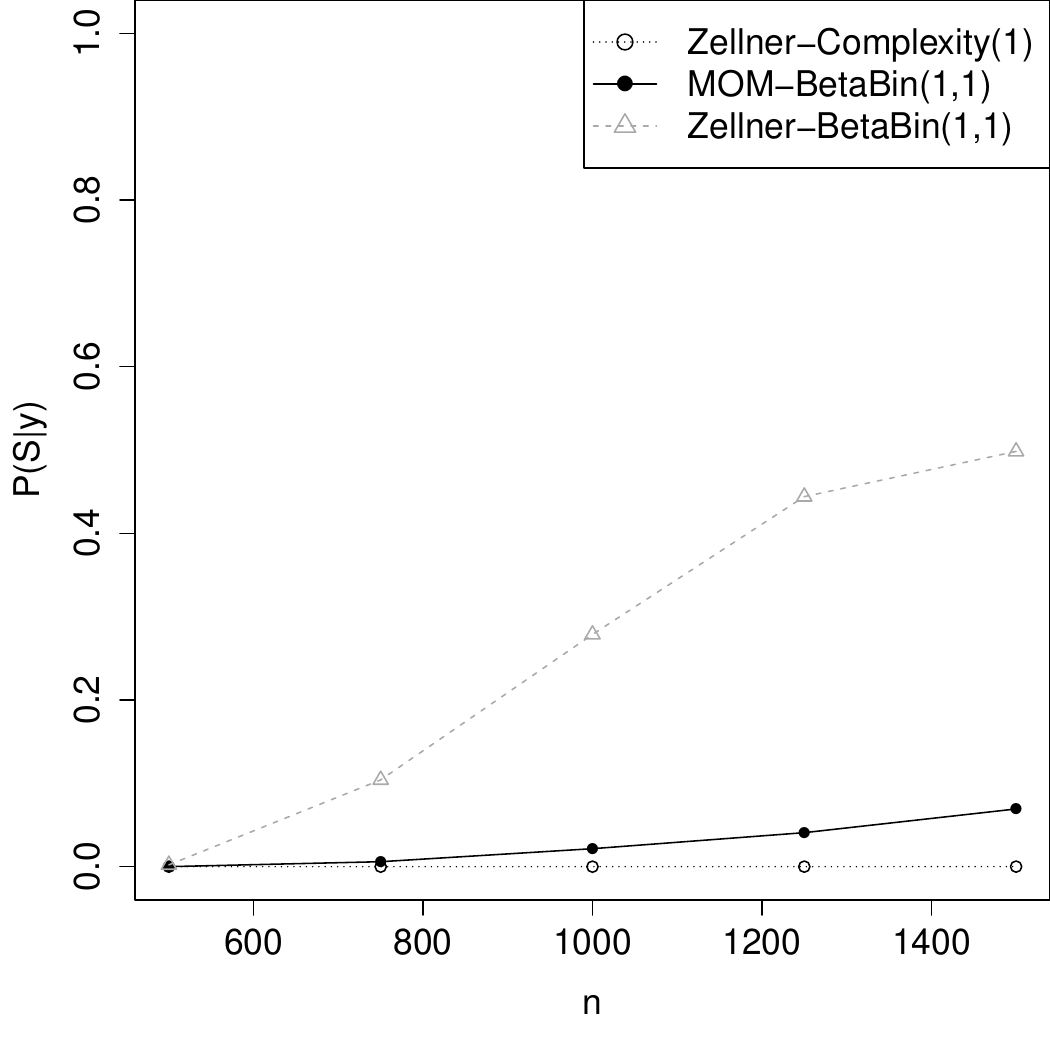}
\end{tabular}
\end{center}
\caption{Linear regression simulation with pairwise correlations$=0.5$.
Average $p(M_t \mid \by)$, $P(S \mid \by)$ and $\sum_{l<p_t} P(S_l^c \mid \by)$
under Zellner-Complexity(1), Zellner-Beta-Binomial(1,1), pMOM-Beta-Binomial(1,1) priors}
\label{fig:pp_corpred}
\end{figure}

We considered normally-distributed covariates with all pairwise correlations equal to 0.5.
We set $p=n$, $p_t=10$ and considered two scenarios.
In Scenario 1 $\theta_j^*=0.5$ for all active variables $j=1,\ldots,p_t$,
whereas Scenario 2 considered weaker signals $\theta_j^*=0.25$ again for $j=1,\ldots,p_t$.
Figure \ref{fig:pp_corpred} shows that whichever prior achieved largest $p(M_t\mid\by)$
depended on $n$ and the signal strength.
For large enough $n$ all three priors discarded small non-spurious models, i.e. $\sum_{l<p_t} P(S_l^c \mid \by)$ vanished, but the required $n$ can be fairly large.
Overall, the MOM-Beta-Binomial prior achieved a reasonable compromise between discarding spurious $m \in S$
and detecting truly active variables.

\section{Discussion}
\label{sec:discussion}

We outlined a strategy to study the $L_1$ convergence of posterior model probabilities and normalized $L_0$ criteria
and showed that, when such convergence occurs,
it is possible to bound frequentist probabilities of correct model selection, and type I-II errors.
The strategy is in principle generic in that it applies to any probability model and prior or $L_0$ penalty,
but requires non-negligible work in bounding the tails of Bayes factors and likelihood-ratio test statistics.
To address this issue, we developed a significant amount of supplementary material to deploy the framework to Gaussian models.
The obtained rates for regression unify current literature and  clarify the consequences of setting sparse priors and $L_0$ penalties. They also clarify 
how convergence depends on the prior dispersion, model prior probabilities, and whether the prior is local or non-local,
as well as on problem characteristics such as $n$, $p$, true sparsity $p_t$ and the signal strength.
Model misspecification also plays a role,
in particular misspecifying the mean structure in (potentially non-linear) regression causes an exponential drop in power,
whereas choosing the wrong error correlation structure can hamper the type I error control.

We gave simple asymptotic expressions for several popular priors and $L_0$ criteria.
We did not study thick-tailed parameter priors,
but such variations typically affect model selection rates only up to lower-order terms.
For instance for a wide class of local priors it is known that for spurious models
$B_{mt}= O_p((\tau n)^{(-p_m-p_t)/2})$ \citep{dawid:1999},
which implies that $L_1$ convergence rates cannot be any faster.
Since our obtained $L_1$ rates are $(\tau n)^{-\alpha (p_m-p_t)/2}$ (or tighter) for any fixed $\alpha<1$, one cannot attain significantly faster rates with other prior families.

Throughout we avoided a detailed study of eigenvalues,  and referred to restricted eigenvalue conditions common in the literature.
This was to highlight the main principles (the role of non-centrality parameters) and keep the results as general as possible.
For a study on eigenvalues see \cite{narisetty:2014} (Remarks 4-5 and Lemma 6.1), for example.

An interesting observation is that, depending on how large $p$ is relative to $n$ one can consider less sparse priors,
this opens a venue to detect smaller signals and may have implications for parameter estimation.
 Also, by restricting the model complexity, one can often use less sparse formulations within the set of allowed models. 
This is particularly relevant when the truth is non-sparse, effect sizes are small or the model's mean structure is strongly misspecified.
Per our examples in this situation it can be helpful to consider strategies that exercise moderation at enforcing sparsity (e.g. the Beta-Binomial prior or the EBIC), 
or that do so in a data-adaptive manner (e.g. using non-local priors on parameters or empirical Bayes).
Such strategies are an interesting venue for future research.

\section*{Acknowledgments}
The author thanks Gabor Lugosi and James O. Berger for helpful discussions.
DR was partially funded by the Europa Excelencia grant EUR2020-112096, the NIH grant R01 CA158113-01, the Ram\'on y Cajal Fellowship RYC-2015-18544,
Plan Estatal PGC2018-101643-B-I00 and
Ayudas Fundaci\'on BBVA a equipos de investigaci\'on cient\'ifica en Big Data 2017.

\beginsupplement

\section*{Supplementary material}

Section \ref{suplsec:auxiliary} provides a number of auxiliary results required for our derivations.
These include bounds on central and non-central chi-square and F distributions,
obtaining non-centrality parameters and bounding Bayesian F-test statistics for nested linear models,
and bounding certain high-dimensional deterministic sums
(e.g. as arising in establishing posterior consistency for variable selection under Zellner's prior).
Section \ref{ssec:tailbounds} provides bounds for the integral of Bayes factor tail probabilities featuring in Lemma \ref{LEM:TAILPROB_TO_MEAN}.
The remaining sections provide proofs for all our main and auxiliary results.

\setcounter{tocdepth}{1}
\tableofcontents

\section{Auxiliary results}
\label{suplsec:auxiliary}

\subsection*{A1. Chi-square and F-distribution bounds}

For convenience Lemma \ref{LEM:INEQ_CHISQ} states well-known chi-square tail bounds.
Lemmas \ref{LEM:INEQ_NCCHISQ_LEFT}-\ref{LEM:INEQ_NCCHISQ_RIGHT} provide Chernoff
and useful related bounds for left and right non-central chi-square tails.
In Lemmas \ref{LEM:INEQ_NCF_LEFT}-\ref{LEM:INEQ_NCF_RIGHT} we derived convenient
moment-generating function-based bounds for the ratio of a non-central divided by a central
chi-square variables, in particular including the F-distribution when the two variables are independent.
Finally, Lemma \ref{LEM:INEQ_F} gives moment bounds used in our theorems to characterize extreme events.

\begin{lemma}
{\bf Chernoff bounds for chi-square tails}

Let $W \sim \chi_\nu^2$. For any $w>\nu$
$$P(W>w) \leq \left( \frac{ew}{\nu} \right)^{\frac{\nu}{2}} e^{-w/2}.$$
Further, for any $w<\nu$
$$
P(W<w) \leq \left( \frac{ew}{\nu} \right)^{\frac{\nu}{2}} e^{-w/2}.
$$
\label{LEM:INEQ_CHISQ}
\end{lemma}

\begin{lemma}
  {\bf Chernoff bounds for non-central chi-square left tails}

Let $W \sim \chi_\nu^2(\lambda)$ be a chi-square with non-centrality parameter $\lambda$ and $w<\lambda$.
Then
$$
P(W<w) \leq \frac{\exp\{\frac{\lambda s}{1-2s} -sw\}}{(1-2s)^{\nu/2}}
$$
for any $s<0$, and the right hand side is minimized for $s= \frac{1}{2} - \frac{\nu}{4w} - \frac{1}{2} \sqrt{\frac{\nu^2}{4w^2} + \frac{\lambda}{w}}$.
In particular,
setting $s=\frac{1}{2} - \frac{1}{2}\sqrt{\lambda/w}$ gives
$$
P(W<w) \leq \frac{\exp \{ -\frac{1}{2} (\sqrt{\lambda} - \sqrt{w})^2\} }{(\lambda / w)^{\nu/4}}.
$$
\label{LEM:INEQ_NCCHISQ_LEFT}
\end{lemma}

\begin{lemma}
  {\bf Chernoff bounds for non-central chi-square right tails}

Let $W \sim \chi_\nu^2(\lambda)$. Let $w>\lambda + \nu$, then
$$
P(W>w) \leq \frac{\exp\{\frac{\lambda s}{1-2s} -sw\}}{(1-2s)^{\nu/2}}
$$
for any $s \in (0,\frac{1}{2})$, and the right hand side is minimized for $s= \frac{1}{2} - \frac{\nu}{4w} - \frac{1}{2} \sqrt{\frac{\nu^2}{4w^2} + \frac{\lambda}{w}}$.

In particular, setting $s=\frac{1}{2} - \frac{1}{2}\sqrt{\lambda/w}$ gives
$$
P(W>w) \leq
e^{ -\frac{w}{2} \left(1 - \sqrt{\frac{\lambda}{w}} \right)^2} \left(\frac{w}{\lambda}\right)^{\frac{\nu}{4}}.
$$

Alternatively, one may set $s=\frac{1}{2} - \frac{\nu}{2w}$ to obtain
$$
P(W>w) \leq
\left( \frac{ew}{\nu} \right)^{\frac{\nu}{2}} e^{-\frac{\lambda}{2}} e^{-\frac{w}{2}(1-\frac{\lambda}{\nu})}
$$
\label{LEM:INEQ_NCCHISQ_RIGHT}
\end{lemma}


\begin{lemma}
{\bf Moment-generating function-based bounds for right F tails}

Let $W= U_1\nu_2/(U_2\nu_1)$ where $U_1 \sim \chi_{\nu_1}^2(\lambda)$, $U_2 \sim \chi_{\nu_2}^2$, $\lambda \geq 0$, $\nu_1\geq 1$ and $\nu_2\geq 1$.
In particular, if $U_1$ and $U_2$ are independent then $W \sim \mathcal{F}(\lambda,\nu_1,\nu_2)$.
Let $w>\lambda + \nu_1$.

(i) Consider the case $\lambda=0$. Then for any $s \in (\nu_1/w,1)$,
$$
P(\nu_1W > w) \leq \left( \frac{ews}{\nu_1} \right)^{\frac{\nu_1}{2}} e^{-ws/2} + (es)^{\frac{\nu_2}{2}} e^{-\frac{s\nu_2}{2}}.
$$

(ii) Consider the case $\lambda > 0$. Then for any $s \in ((\lambda+\nu_1)/w, 1)$
$$
P(\nu_1W > w) \leq \frac{e^{\frac{\lambda t}{1-2t} - tws}}{(1-2t)^{\frac{\nu_1}{2}}} + (es)^{\frac{\nu_2}{2}} e^{-\frac{s\nu_2}{2}},
$$
where $t \in (0,1/2)$, and the right hand side is minimized by setting
$t= \frac{1}{2} - \frac{\nu_1}{4ws} - \frac{1}{2} \sqrt{\frac{\nu_1^2}{4w^2s^2} + \frac{\lambda}{ws}}$.

In particular, we may set $t=\frac{1}{2} - \frac{1}{2} \sqrt{\frac{\lambda}{ws}}$
to obtain
$$
P(\nu_1W>w) \leq e^{-\frac{ws}{2} \left( 1- \sqrt{\frac{\lambda}{ws}} \right)^2 } \left( \frac{ws}{\lambda} \right)^{\frac{\nu_1}{4}}
+ (es)^{\nu_2/2} e^{- \frac{s\nu_2}{2}}.
$$
Alternatively, we may also set $t= \frac{1}{2}  - \frac{\nu_1}{2ws}$ to obtain
$$
P(\nu_1W>w) \leq \left( \frac{ews}{\nu_1} \right)^{\frac{\nu_1}{2}} e^{-\lambda/2}
e^{-\frac{ws}{2}(1- \lambda/\nu_1)} + (es)^{\frac{\nu_2}{2}} e^{-s\nu_2/2}.
$$
\label{LEM:INEQ_NCF_RIGHT}
\end{lemma}

\begin{cor}
Let $W= U_1\nu_2/(U_2\nu_1)$ where $U_1 \sim \chi_{\nu_1}^2(\lambda)$, $U_2 \sim \chi_{\nu_2}^2$, $\nu_1\geq 1$ and $\nu_2> \nu_1/(2-\sqrt{3})$.
Consider $w \in ((\nu_1+\lambda)/(2-\sqrt{3}),\nu_2)$.

(i) If $\lambda=0$, then
$$
P(\nu_1 W > w) \leq \left( \frac{ew}{\nu_1} \right)^{\frac{\nu_1}{2}}
e^{-\frac{w}{2} \left(1 - \sqrt{2w/\nu_2}\right)} + e^{-\frac{w}{2}}.
$$

(ii) If $\lambda>0$, then
$$
P(\nu_1 W > w) \leq \left( \frac{w}{\lambda} \right)^{\frac{\nu_1}{4}}
e^{-\frac{w}{2} \left(1 - \sqrt{2w/\nu_2}\right) \left(1-\sqrt{\lambda/\left[w(1 - \sqrt{2w/\nu_2})\right]}  \right)^2}
+ e^{-\frac{w}{2}}.
$$

\label{COR:INEQ_F_RIGHT}
\end{cor}

\begin{lemma}
{\bf Moment-generating function-based bounds for left F tails}

Let $W= U_1\nu_2/(U_2\nu_1)$ where $U_1 \sim \chi_{\nu_1}^2(\lambda)$, $U_2 \sim \chi_{\nu_2}^2$, $\lambda \geq 0$, $\nu_1\geq 1$ and $\nu_2\geq 1$.
Then for any $s \geq 1$ and $t<0$
$$
P(\nu_1W < w) \leq \frac{\exp\{\frac{\lambda t}{1-2t} -tws\}}{(1-2t)^{\nu_1/2}}
+ e^{-\frac{\nu_2}{2} (s - 1 -\log(s))},
$$
and the right hand side is minimized for
$t= \frac{1}{2} - \frac{\nu_1}{4ws} - \frac{1}{2} \sqrt{\frac{\nu_1^2}{4ws^2} + \frac{\lambda}{ws}}$.
Further, if $s< \lambda/w$ we may set $t=\frac{1}{2} - \frac{1}{2} \sqrt{\frac{\lambda}{sw}}$ to obtain
$$
P(\nu_1W < w) \leq
\frac{e^{-\frac{\lambda}{2}\left(1 - \sqrt{ws/\lambda}\right)^2}}{\left(\lambda/(ws)\right)^{\frac{\nu_1}{4}}}
+ e^{-\frac{\nu_2}{2} (s - 1 - \log(s))}.
$$
\label{LEM:INEQ_NCF_LEFT}
\end{lemma}

\begin{lemma}
{\bf Moment bounds for F distribution}

Let $W \sim \mathcal{F}(\nu_1,\nu_2)$ be an F-distributed random variable with degrees of freedom $\nu_1 \geq 1, \nu_2 > 4$.
Then for any $w> \nu_2/(\nu_2-2)$, $s \in [1,\nu_2/2-2]$
\begin{align}
P(W>w) < a \left( \frac{\nu_2 (s+\nu_1/2-1)}{\nu_1w(\nu_2/2-s-1)}  \right)^s
\frac{(s+\nu_1/2-1)^{\frac{\nu_1-1}{2}}}{(\nu_1/2-1)^{\frac{\nu_1-1}{2} \mbox{I}(\nu_1>2)} }
\left( 1 - \frac{s}{\nu_2/2-1} \right)^{\frac{\nu_2-1}{2}}.
\label{eq:ftail_momentbound}
\end{align}
where $a= e^{\frac{5}{2}}/(\pi\sqrt{2})$ if $\nu_1=1$,
$a= e^2/\sqrt{2\pi}$ if $\nu_1=2$ and $a= e^2 / (2\pi)$ if $\nu_1>2$.

Consider $s=\min\{ (w-1)\nu_1/2+1, \nu_2/2-2 \}$.
If $w \leq (\nu_1+\nu_2-6)/\nu_1$ then $P(W>w)<$
$$
a
\left( 1 + \frac{\nu_1(w-1)+4}{\nu_2-\nu_1(w-1)-4} \right)^{\frac{(w-1)\nu_1}{2}+1}
\frac{(\nu_1 w/2)^{\frac{\nu_1-1}{2}}}{(\nu_1/2-1)^{\frac{\nu_1-1}{2} \mbox{I}(\nu_1>2)} }
\left( 1 - \frac{(w-1)\nu_1/2-1}{\nu_2/2-1} \right)^{\frac{\nu_2-1}{2}}
$$
where
$$
\left( 1 + \frac{\nu_1(w-1)+4}{\nu_2-\nu_1(w-1)-4} \right)^{\frac{(w-1)\nu_1}{2}}
\leq e^{\frac{\nu_1(w-1) (\nu_1(w-1)+2)}{2(\nu_2 - \nu_1(w-1)-2)}}
$$
and
$$
\left( 1 - \frac{(w-1)\nu_1/2-1}{\nu_2/2-1} \right)^{\frac{\nu_2-1}{2}} < e^{-\frac{w \nu_1}{2}} e^{\frac{\nu_1-1}{2}+\frac{3}{2}}
.
$$
If $w > (\nu_1+\nu_2-6)/\nu_1$ then
$$
P(W>w) < ae \left( \frac{\nu_1+\nu_2-6}{\nu_1w} \right)^{\frac{\nu_2}{2}-2}
\frac{(\nu_1/2+\nu_2/2-3)^{\frac{\nu_1-1}{2}}}{(\nu_1/2-1)^{\frac{\nu_1-1}{2}\mbox{I}(\nu_1>2)} }
\frac{1}{(\nu_2/2-1)^{\frac{3}{2}}}.
$$
\label{LEM:INEQ_F}
\end{lemma}

\subsection*{A2. Non-centrality parameter for nested models}

For convenience Lemma \ref{LEM:NESTEDLM_CHITEST} states
a known result on the difference of sum of squares between nested linear models
(proven in the supplementary material).
Lemma \ref{LEM:NESTEDLM_NONLINEAR} is an extension to misspecified linear models
and Lemma \ref{LEM:NESTEDLM_CHITEST_HETERO} is an extension to heteroskedastic errors.
Lemmas \ref{LEM:NCCHISQ_BAYES_UNIV}-\ref{LEM:BAYES_UNIVFTEST_INEQ}
characterize the distribution and tails of Normal shrinkage estimators and related F-test statistics.

\begin{lemma}
Consider two nested Normal linear regression models $M_m \subset M_q$
with respective full-rank design matrices $X_m$ and $X_q=(X_m,X_s)$,
where $X_s$ are the columns in $X_q$ not contained in $X_m$.
Let $H_m=X_m(X_m'X_m)^{-1}X_m'$ and $H_q=X_q(X_q'X_q)^{-1}X_q'$ be the projection matrices for $M_m$ and $M_q$.
Assume that truly $f^*(\by)= N(\by; X_q \btheta_q^*, \phi^* I)$,
where potentially some or all of the entries in $\btheta_q^*$ can be zero
and $\btheta_q^*=((\btheta_m^*)',(\btheta_s^*))'$.
Let $\hat{\btheta}_q$ and $\hat{\btheta}_m$ be the least squares estimates, then
$$
\frac{1}{\phi^*}(\hat{\btheta}_q' X_q'X_q \hat{\btheta}_q - \hat{\btheta}_m' X_m'X_m \hat{\btheta}_m) \sim \chi^2_{p_q-p_m}(\lambda_{qm}),
$$
where $\lambda_{qm}=(X_s\btheta_s^*)'(I-H_m)X_s \btheta_s^*/\phi^*= (X_t\btheta_t^*)'H_q(I-H_m)H_qX_t\btheta_t^*/\phi^*$.
\label{LEM:NESTEDLM_CHITEST}
\end{lemma}

\begin{lemma}
Consider two linear regression models $M_m \subset M_q$ as in Lemma \ref{LEM:NESTEDLM_CHITEST}
and let $\hat{\btheta}_m$ and $\hat{\btheta}_q$ be the least squares estimates.
Assume that truly $f^*(\by)= N(\by; W \bbeta^*, \xi^* I)$
and let $\btheta_q^*=(X_q'X_q)^{-1}X_q'W\bbeta^*$ be the KL-optimal regression coefficient.
Then
$$
\frac{1}{\xi^*}(\hat{\btheta}_q' X_q'X_q \hat{\btheta}_q - \hat{\btheta}_m' X_m'X_m \hat{\btheta}_m) \sim \chi^2_{p_q-p_m}(\lambda_{qm}),
$$
where $\lambda_{qm}=(X_s\btheta_s^*)'(I-H_m)X_s \btheta_s^*/\xi^*=
(W\bbeta^*)'H_q (I-H_m)H_q W\bbeta^*/\xi^*$.
\label{LEM:NESTEDLM_NONLINEAR}
\end{lemma}

\begin{lemma}
Let $X_q=(X_m,X_s)$ as in Lemma \ref{LEM:NESTEDLM_CHITEST}, where $p_q \leq n$.
Assume that truly $\by \sim N(X_q \btheta_q^*, \phi^* \Sigma^*)$
where potentially some or all of the entries in $\btheta_q^*$ can be zero
and $\btheta_q^*=((\btheta_m^*)',(\btheta_s^*))'$.
Assume that $\sum_{i=1}^{n} \Sigma_{ii}^*=n$,
so that $n^{-1} \sum_{i=1}^{n} \mbox{Var}(y_i)=\phi^*$.
Let $\hat{\btheta}_q$ and $\hat{\btheta}_m$ be the least squares estimates,
$\tilde{X}_s= X_s - X_m(X_m'X_m)^{-1}X_m'X_s$.
Let $\underline{\varrho}(A)$ and $\bar{\varrho}(A)$ be the smallest and largest eigenvalue of $A$.

\begin{enumerate}[leftmargin=*,label=(\roman*)]
\item
Let $\underline{\omega}_{mq}=\underline{\varrho}(\tilde{X}_s' \Sigma \tilde{X}_s (\tilde{X}_s'\tilde{X}_s)^{-1})$,
$\bar{\omega}_{mq}= \bar{\varrho}(\tilde{X}_s' \Sigma \tilde{X}_s (\tilde{X}_s'\tilde{X}_s)^{-1})$.
Then
$$
\underline{\omega}_{mq} Z_1
\leq \frac{\hat{\btheta}_q' X_q'X_q \hat{\btheta}_q - \hat{\btheta}_m' X_m'X_m \hat{\btheta}_m}{\phi^*} \leq
\bar{\omega}_{mq} Z_1
$$
where $Z_1 \sim \chi_{p_q-p_m}^2(\tilde{\lambda}_{qm} )$,
$\tilde{\lambda}_{qm}=(\btheta_s^*)' W^{-1} \btheta_s^*$ and
$W=(\tilde{X}_s'\tilde{X}_s)^{-1} \tilde{X}_s' \Sigma \tilde{X}_s (\tilde{X}_s'\tilde{X}_s)^{-1}$.
Further,
$$
\frac{1}{\bar{\omega}_{mq}} \lambda_{qm} \leq \tilde{\lambda}_{qm} \leq \frac{1}{\underline{\omega}_{mq}} \lambda_{qm},
$$
where $\lambda_{qm}=(\btheta_s^*)'X_s'(I-X_m(X_m'X_m)^{-1}X_m)X_s \btheta_s^*/\phi^*$.

\item

Let $L_q$ be the set of $n \times (n-p_q)$ matrices $T$
such that the matrix $(X_q,T)$ is full-rank.
Define $\tilde{T}= (I - X_q(X_q'X_q)^{-1}X_q')T$,
$\underline{\omega}_q= \max_{L_q} \underline{\varrho}(\tilde{T}_l' \Sigma^* \tilde{T}_l (\tilde{T}_l'\tilde{T}_l)^{-1})$
and $\bar{\omega}_q= \min_{L_q} \bar{\varrho}(\tilde{T}_l' \Sigma^* \tilde{T}_l (\tilde{T}_l'\tilde{T}_l)^{-1})$.
Then
$$
\frac{\underline{\omega}_{mq} Z_1}{\bar{\omega}_q Z_2}
\leq \frac{\hat{\btheta}_q' X_q'X_q \hat{\btheta}_q - \hat{\btheta}_m' X_m'X_m \hat{\btheta}_m}{\by'\by - \hat{\btheta}_q' X_q'X_q \hat{\btheta}_q} \leq
\frac{\bar{\omega}_{mq} Z_1}{\underline{\omega}_q Z_2},
$$
where $Z_1$ is as in Part (i) and $Z_2 \sim \chi_{n-p_q}^2$.
\end{enumerate}

\label{LEM:NESTEDLM_CHITEST_HETERO}
\end{lemma}

\begin{lemma}
Assume that truly $\by \sim N(X_q \btheta_q^*, \phi^* I)$.
Let $\hat{\btheta}_q= (X_q'X_q+ V_q^{-1}/(\tau n)) X_q'\by$ where $V_q$ is positive-definite,
$\bmu= (X_q'X_q+ V_q^{-1}/(\tau n))^{-1} X_q'X_q \btheta_q^*$ and
$\Sigma= (X_q'X_q+ V_q^{-1}/(\tau n))^{-1} X_q'X_q (X_q'X_q+ V_q^{-1}/(\tau n))^{-1}$.
Denote by $\hat{\theta}_{qi}$, $\theta_{qi}^*$ and $\mu_i$ the $i^{th}$ entry in $\hat{\btheta}_q$, $\btheta_q^*$ and $\bm{\mu}$ respectively,
by $\sigma_{ii}$ the element $(i,i)$ in $\Sigma$,
${\bf x}_{qi}$ the $i^{th}$ column in $X_q$,
and by $X_m$ the matrix obtained by removing the $i^{th}$ column in $X_q$.
Then
$$
\frac{\hat{\theta}_{qi}^2}{\phi^* \sigma_{ii}} \sim \chi_1^2 \left(\mu_i^2/(\phi^* \sigma_{ii}) \right).
$$
Let $\tilde{\sigma}_{ii}= ({\bf x}_{qi}' (I - X_m (X_m'X_m)^{-1}X_m') {\bf x}_{qi})^{-1}$ and
$\lambda_{qi}= (\theta_{qi}^*)^2/(\tilde{\sigma}_{ii} \phi^*)$, then
\begin{align}
\frac{(\tau n)^2\rho_{qp_q}^2}{(\tau n)^2\rho^2_{qp_q} +1}
&\leq \frac{\sigma_{ii}}{\tilde{\sigma}_{ii}} \leq
\frac{(\tau n)^2\rho_{q1}^2}{(\tau n)^2\rho_{q1}^2 +1}
\nonumber \\
(1 - 2 \delta_{qi}^*) &\leq \frac{\mu_i^2}{(\theta_{qi}^*)^2} \leq
\left( 1 + 2 \delta_{qi}^* + (\delta_{qi}^*)^2 \right)
\nonumber \\
\lambda_{qi} (1 - 2 \delta_{qi}^*) \left(1 + \frac{1}{(\tau n)^2\rho_{q1}^2} \right)
&\leq \frac{\mu_i^2}{\phi^* \sigma_{ii}} \leq
\lambda_{qi}
\left( 1 + 2 \delta_{qi}^* + (\delta_{qi}^*)^2 \right) \left(1+ \frac{1}{(\tau n)^2\rho_{qp_q}^2} \right)
\nonumber
\end{align}
where $\delta_{qi}^*= \sqrt{(\btheta_q^*)' \btheta_q^*}/(((\tau n)\rho_{q1} +1)|\theta_{qi}^*|)$
and $\rho_{q1} \geq \ldots \geq \rho_{qp_q} >0$ are the eigenvalues of $V_qX_q'X_q$.
\label{LEM:NCCHISQ_BAYES_UNIV}
\end{lemma}

\begin{lemma}
Assume that truly $\by \sim N(X_q \btheta_q^*, \phi^* I)$
and let $\hat{\btheta}_q$, $\lambda_{qi}$, $V_q$
and $\rho_{q1} \geq \ldots \geq \rho_{qp_q} >0$ be as in Lemma \ref{LEM:NCCHISQ_BAYES_UNIV}.
Let $\tilde{s}_q= l_\phi + \by'\by - \by'X_q(X_q'X_q+ (\tau n)^{-1}V_q^{-1})^{-1} X_q'\by$.
Assume that, as $n \rightarrow \infty$,
$(\tau n) \rho_{qp_q} \gg 1$ and
$(\tau n) \rho_{q1} \gg \sqrt{(\btheta_q^*)' \btheta_q^*} / |\theta_{qi}^*|$.

\noindent
(i) Let $\theta_{qi}^*\neq 0$ and $h_n>0$ be a sequence satisfying $h_n \gg \phi^*/(\theta_{qi}^*)^2$.
Then for any fixed $\gamma \in (0,1)$, as $n \rightarrow \infty$,
\begin{align}
P \left( \frac{\hat{\theta}_{qi}^2}{\tilde{s}_q/(n-p_q)} < \frac{1}{h_n} \right) &\ll
e^{-\gamma \lambda_{qi}/2} + e^{-(n-p_q)^\gamma}
\nonumber \\
P \left( \prod_{i=1}^{p_q} \frac{\hat{\theta}_{qi}^2}{\tilde{s}_q/(n-p_q)} < \frac{1}{h_n^{p_q}} \right) &\ll
p_q( e^{-\gamma \min_i \lambda_{qi}/2} + e^{-(n-p_q)^\gamma}).
\nonumber
\end{align}

\noindent
(ii) Let $\theta_{qi}^* \neq 0$ and $\tilde{h}_n \gg (\theta_{qi}^*)^2/\phi^*$.
Then for any fixed $\gamma \in (0,1)$, as $n \rightarrow \infty$,
$$
P \left( \frac{\hat{\theta}_{qi}^2}{\tilde{s}_q/(n-p_q)} > \tilde{h}_n \right)
  \ll  e^{-\gamma \tilde{h}_n/(2\sigma_{ii})} + e^{-(n-p_q)^\gamma/2}
\ll  e^{-\gamma \lambda_{qi}/2} + e^{-(n-p_q)^\gamma/2}.
$$

\noindent
(iii) Let $\theta_{qi}^*=0$ and $\tilde{h}_n \gg \sigma_{ii}$.
Then for any fixed $\gamma \in (0,1)$, as $n \rightarrow \infty$,
$$
P \left( \frac{\hat{\theta}_{qi}^2}{\tilde{s}_q/(n-p_q)} > \tilde{h}_n \right)
  \ll  e^{-\gamma \tilde{h}_n/(2\sigma_{ii})} + e^{-(n-p_q)^\gamma/2}
\ll  e^{-\gamma \lambda_{qi}/2} + e^{-(n-p_q)^\gamma/2}.
$$

\label{LEM:BAYES_UNIVFTEST_INEQ}
\end{lemma}

\subsection*{A3. Bound for a Bayesian F-test statistic under a general Normal prior}

\begin{lemma}
For any model $k$ let $s_k= \by'\by - \by'X_k(X_k'X_k)^{-1} X_k'\by$
and $\tilde{s}_k= l_\phi + \by'\by - \by'X_k(X_k'X_k+ (\tau n)^{-1}V_k^{-1})^{-1} X_k'\by$
where $(\tau n)>0$, $l_\phi>0$, and $V_k$ is a symmetric positive-definite matrix.
In particular $s_0=\by'\by$ denotes the sum of squared residuals under the model with no covariates.
Let $m$ be a given model
and for any other model $k$ define $F_{mk}=\frac{(s_k-s_m)/(p_m-p_k)}{s_m/(n-p_m)}$ and
$\tilde{F}_{mk}= \frac{(\tilde{s}_k-\tilde{s}_m)/(p_m-p_k)}{\tilde{s}_m/(n-p_m)}$.
Then, for any model $t$, it holds that
$$
1 \leq \frac{\tilde{s}_t}{s_t} \leq 1 + \frac{(s_0-s_t)/s_t}{1+(\tau n)\rho_{t,p_t}}
$$
and that
$$(p_m-p_t) \tilde{F}_{mt} \leq
\frac{(\tau n) \rho_{t,p_t}}{1+(\tau n)\rho_{t,p_t}} (p_m-p_t) F_{mt} +
\frac{1}{1+(\tau n)\rho_{t,p_t}} p_m F_{m0},
$$
where $\rho_{tp_t}$ denotes the smallest non-zero eigenvalue of $V_t X_t'X_t$.
\label{LEM:FSTAT_NORMALPRIOR_BOUND}
\end{lemma}

\subsection*{A4. Asymptotic bounds on sums of deterministic sequences}
\label{ssec:bound_seriessum}

A convenient strategy to bound the right-hand side in Proposition \ref{PROP:CORRECTSEL_FROM_PP} is to first obtain model-specific bounds, often expressed in as asymptotic bounds of the type $E_{f^*}( p(M_k \mid \by)) \ll a_n^{(k)}$ for some sequence $a_n^{(k)}$.
An important technical remark is that even if such a bound exists for each $k \neq t$,
this does not necessarily imply that
$\sum_{k \neq t} E_{f^*}(p(M_k \mid \by)) \ll \sum_{k \neq t} a_n^{(k)}$.
This issue can be addressed by finding finite-$n$ bounds such that
$E_{f^*}( p(M_k \mid \by)) \leq a_n^{(k)}$ for all $n \geq n_0$ where $n_0$ does not depend on $k$,
then clearly $\sum_{k \neq t} E_{f^*}(p(M_k \mid \by)) \leq \sum_{k \neq t} a_n^{(k)}$ for $n \geq n_0$.
A refinement is to show that there is a single $n_0$ that holds for all models within a subset,
and then showing that these $n_0$ are strictly ordered across subsets, so one may take the largest $n_0$ as a uniform bound across all models.
This is the strategy used for our regression examples in Sections \ref{sec:modelrates}-\ref{sec:l0penalties}.
There is a common $n_{k0}$ for all spurious models $M_k \in S_l$ of size $l$, and these $n_{k0}$ are strictly increasing in $l$, so that they are all bounded by the largest model size $l$.
Similarly, for non-spurious models $k \in S_l^c$ a uniform fixed $n_0$ is obtained by taking the minimum over certain non-centrality parameters.

In some general settings beyond our regression examples it can be hard to find such finite-$n$ or uniform bounds.
The rest of this section offers some discussion on how to bound the right-hand side in Proposition \ref{PROP:CORRECTSEL_FROM_PP} in these situations.

We provide sufficient conditions
to bound the posterior probability of model subsets defined by their size in Lemma \ref{LEM:SUMO_SIZE},
or by more general groupings in Lemma \ref{LEM:SUMO}.
We first outline the idea. 
Let $a_n^{(l)}$ and $\tilde{a}_n^{(l)}$ be sequences such that
$E_{f^*}(p(M_k \mid \by)) \ll a_n^{(l)}$ for all $l \in S_l$  and
$E_{f^*}(p(M_k \mid \by)) \ll \tilde{a}_n^{(l)}$ for all $l \in S_l^c$.
Then under suitable conditions (Lemmas \ref{LEM:SUMO_SIZE}-\ref{LEM:SUMO})
one can obtain global bounds by adding up model-specific bounds, i.e.
\begin{align}
E_{f^*}(P(S\mid\by)) \preceq \sum_{l=p_t+1}^{\bar{p}} a_n^{(l)} |S_l|
\nonumber \\
E_{f^*}(P(S^c \mid \by)) \preceq \sum_{l=0}^{p_t} \tilde{a}_n^{(l)} |S_l^c| + \sum_{l=p_t+1}^{\bar{p}} \tilde{a}_n^{(l)} |S_l^c|.
\label{eq:sumo}
\end{align}


The rates $(a_n^{(l)},\tilde{a}_n^{(l)})$ are associated to Bayes factors and $p(M_k)$
as in \eqref{eq:pp_growingp}, see Sections \ref{sec:modelrates}-\ref{sec:l0penalties}.
Expression \eqref{eq:sumo} splits the sum to emphasize the role of $p_t$, i.e. the sparsity of $f^*$.

Specifically, Lemma \ref{LEM:SUMO_SIZE} provides sufficient conditions to bound the total posterior probability
assigned to $S_l$ and $S_l^c$, and to unions thereof such as $S$ and $S^c$.
In turn, sufficient conditions for Lemma \ref{LEM:SUMO_SIZE}(i)
are that $\sigma_{l,n}/b_n$ is a (non-strictly) decreasing series in $n$ for all $n \geq n_0$ with common $n_0$ across all $l$,
or alternatively that
$\lim_{n \rightarrow \infty} \sum_{l=p_t + 1}^{\bar{p}} \sigma_{l,n+1}/b_{n+1} \leq \lim_{n \rightarrow \infty} \sum_{l=p_t+1}^{\bar{p}} \sigma_{l,n}/b_n < \infty$.

Lemma \ref{LEM:SUMO} provides sufficient conditions to bound the total posterior probability
across different model subsets $A_{l,n}$ that are indexed by some model characteristic $l$,
say its total number of variables or the number of active variables,
by adding asymptotic bounds for each specific $A_{l,n}$.
Corollary \ref{COR:SUMO} is a specialization to the case where all models have a common asymptotic bound,
for instance in our variable selection examples all spurious models $k \in S$ of equal size $p_k$ share such a bound.
Lemmas \ref{LEM:SUMSPUR_ZELLNERKNOWN}-\ref{LEM:SUMSPUR_ZELLNERUNKNOWN}
are auxiliary results to bound the total posterior probability assigned to spurious models
under Zellner's prior when the residual variance is assumed either known or unknown, respectively.

\begin{lemma}
Let $a_n^{(l)}$ and $\tilde{a}_n^{(l)}$ be sequences such that, as $n \rightarrow \infty$,
$E_{f^*}(p(M_k \mid \by)) \ll a_n^{(l)}$ for all $l \in S_l$  and
$E_{f^*}(p(M_k \mid \by)) \ll \tilde{a}_n^{(l)}$ for all $l \in S_l^c$.
Denote by $\sigma_{l,n}= \sum_{k \in S_l} E_{f^*}( p(M_k \mid \by) )$,
$\tilde{\sigma}_{l,n}= \sum_{k \in S_l^c} E_{f^*}( p(M_k \mid \by) )$
the mean posterior probability assigned to size $l$ spurious and non-spurious models respectively.

\begin{enumerate}[leftmargin=*,label=(\roman*)]
\item Let $b_n= \sum_{l=p_t+1}^{\bar{p}} a_n^{(l)} |S_l|/(\bar{p}-p_t)$.
Suppose that the following two conditions hold
\begin{align}
\sigma_{\bar{p},n} \ll b_n
\nonumber \\
\lim_{n \rightarrow \infty} \sum_{l=p_t+1}^{\bar{p}} \frac{\sigma_{l,n+1}}{b_{n+1}} - \frac{\sigma_{l,n}}{b_n}=0.
\nonumber
\end{align}
Then
$
E_{f^*}(P(S\mid\by)) \preceq \sum_{l=p_t+1}^{\bar{p}} a_n^{(l)} |S_l|.
$
\item Let $\tilde{b}_n= \sum_{l=0}^{\bar{p}} \tilde{a}_n^{(l)}|S_l^c|/(\bar{p}+1)$.
Suppose that the following two conditions hold
\begin{align}
\tilde{\sigma}_{\bar{p},n} \ll \tilde{b}_n
\nonumber \\
\lim_{n \rightarrow \infty} \sum_{l=0}^{\bar{p}} \frac{\tilde{\sigma}_{l,n+1}}{\tilde{b}_{n+1}} - \frac{\tilde{\sigma}_{l,n}}{\tilde{b}_n}=0.
\nonumber
\end{align}
Then
$
E_{f^*}(P(S^c \mid \by)) \preceq \sum_{l=0}^{p_t} \tilde{a}_n^{(l)} |S_l^c| + \sum_{l=p_t+1}^{\bar{p}} \tilde{a}_n^{(l)} |S_l^c|
$

\end{enumerate}
\label{LEM:SUMO_SIZE}
\end{lemma}

\begin{lemma}
Let $A_{l,n} \subseteq \{1,\ldots,K\}$ be subsets of the model space indexed by $l=l_n^{(0)},\ldots,l_n^{(1)}$,
where their size $|A_{l,n}|$ may grow with $n$.
Let $g_n^{(l)}>0$ be a set of decreasing series also indexed by $l=l_n^{(0)},\ldots,l_n^{(1)}$.
Denote by $\sigma_{l,n}= \sum_{k \in A_{l,n}}^{} E_{f^*}(p(M_k \mid \by))$
and by $\bar{g}_n= \sum_{l=l_n^{(0)}}^{l_n^{(1)}} g_n^{(l)}/(l_n^{(1)}-l_n^{(0)}+1)$.
Assume that the following two conditions hold
\begin{enumerate}[leftmargin=*,label=(\roman*)]
\item  $\sigma_{l_n^{(1)},n} \ll \bar{g}_n$
\item $\lim_{n\rightarrow \infty} \sum_{l=l_n^{(0)}}^{l_n^{(1)}} \sigma_{l,n+1} / \bar{g}_{n+1} - \sigma_{l,n} / \bar{g}_n \leq 0$
\end{enumerate}

Then $\sum_{l=l_n^{(0)}}^{l_n^{(1)}} \sum_{k \in A_{l,n}}^{} E_{f^*}(p(M_k \mid \by)) \preceq \sum_{l=l_0^{(n)}}^{l_1^{(n)}} g_n^{(l)}$.
\label{LEM:SUMO}
\end{lemma}

\begin{cor}
Let $A_n \subseteq \{1,\ldots,K\}$ be a subset of the model space, where $|A_n|$ may grow with $n$.
Assume that $\mu_{k,n}=E_{f^*}(p(M_k \mid \by)) \ll b_n$ as $n \rightarrow \infty$ for all $k \in A_n$, where $b_n>0$.
Assume that the following condition holds
$$\lim_{n\rightarrow \infty} \sum_{k \in A_n} \mu_{k,n+1} / b_{n+1} - \mu_{k,n} / b_n \leq 0.$$
Then $\sum_{k \in A_n} \mu_{k,n} \preceq |A_n| b_n$.
\label{COR:SUMO}
\end{cor}

\begin{lemma}
  Let $(\tau n)$, $p$, $\bar{p} \leq p$, $p_t \leq \bar{p}$ be such that
  as $n \rightarrow \infty$, $(\tau n)^{1/2} \gg (p_t+1) (\bar{p}-p_t)^{a/2} \log^{3/2}((\tau n)^{1/2} (p-p_t))$,
  $p_t \succeq 1$, $p \succeq 1$, $\bar{p} \succeq 1$. Then
  for any fixed $a>1$ it holds that
  \begin{align}
    \sum_{l=p_t+1}^{\bar{p}} {l \choose p_t} \frac{\left[ (l-p_t) \log((\tau n)^{1/2} (p-p_t)) \right]^{\frac{l-p_t}{2}+1}}{(\tau n)^{\frac{l-p_t}{2}}}
    \preceq \frac{(p_t+1) (\bar{p}-p_t)^{a/2} \log^{3/2}((\tau n)^{1/2} (p-p_t))}{(\tau n)^{1/2}}
    \nonumber
   \end{align}
  \label{LEM:SUMSPUR_ZELLNERKNOWN}
\end{lemma}

\begin{lemma}
  Let $P(S \mid \by)$ be the posterior probability assigned to spurious models
  with $p_k \leq \bar{p}$ variables under Zellner's prior
  $p(\btheta_k,\phi \mid M_k)= N(\btheta_k; {\bf 0}, (\tau n) \phi (X_k'X_k)^{-1}) \mbox{IG}(\phi; a_\phi,l_\phi)$,
  and let $p$ be the total number of variables.
  Let $p(M_k)= \left(  (\bar{p}+1) {p \choose p_k} \right)^{-1}$ be Beta-Binomial(1,1) prior probabilities on the models.
  Assume that truly $\by \sim N(X_t\btheta_t; \phi I)$ for some $t$ with $p_t \leq \bar{p}$.
  Assume also that as $n \rightarrow \infty$ it holds that $(\bar{p}-p_t) \log((\tau n) (p-p_t)) \ll n-\bar{p}$
  and that $[(\tau n)^{1/2} (p-p_t)]^{2c_n} \ll (\tau n)^{1/2}$,
  where $c_n= \sqrt{ (\bar{p}-p_t) [\log((\tau n)^{1/2} (p-p_t))]/(n-\bar{p})}$.
  Then
  $$
  E_{f^*}( P(S \mid \by)) \preceq
  \frac{(p_t+1)}{(\tau n)^{1/2}} e^{2 [\log^{3/2} ((\tau n)^{1/2} (p-p_t))] \sqrt{(p-p_t)/(n-\bar{p})  }}.
  $$
\label{LEM:SUMSPUR_ZELLNERUNKNOWN}
\end{lemma}

\section{Tail integral bounds}
\label{ssec:tailbounds}

\subsection*{Generic bounds for exponential and polynomial tails}

A generic strategy to apply Lemma \ref{LEM:TAILPROB_TO_MEAN} is to upper-bound $B_{kt}p(M_k)/p(M_t)$ or a suitable transformation,
e.g. in Sections \ref{sec:modelrates}-~\ref{sec:l0penalties} we bound $d \log(B_{kt} g p(M_k)/p(M_t)) \leq W$ for some $d,g>0$,
where $W$ is a random variable for which one can characterize the tails.
Then trivially $P(B_{kt} > p(M_t)/p(M_k)/(1/u-1)) \leq P(W > d \log(g/(1/u-1)))$
and $E_{f^*}(p(M_k\mid\by))$ can be bounded by integrating the tails of $W$.
One can use this strategy on a case-by-case basis,
i.e. for a given model/prior, but to facilitate applying our framework
Lemmas \ref{LEM:INTBOUND_EXPTAILS}-\ref{LEM:INTBOUND_POLYTAILS} give non-asymptotic bounds
for the common cases where $W$ has exponential or polynomial tails, respectively.
In Sections \ref{sec:modelrates}-\ref{sec:l0penalties}, as $n \rightarrow \infty$ we let $g \rightarrow \infty$
and set the other parameters such that
the integrals in Lemmas \ref{LEM:INTBOUND_EXPTAILS}-\ref{LEM:INTBOUND_POLYTAILS} and hence $E_{f^*}(p(M_k\mid\by))$
are bounded by $b/g^\alpha$ times a term of smaller order,
where $b$ is a constant and $\alpha$ should be thought of as a constant close to 1.
Lemmas \ref{LEM:CHISQTAIL_INTBOUND}-\ref{LEM:FTAIL_INTBOUND} are adaptations to chi-square and F distributions useful for linear regression,
for simplicity they state asymptotic bounds as $g \rightarrow \infty$, but the proofs also provide non-asymptotic bounds.

\begin{lemma}

  Let $\underline{u},\bar{u} \in (0,1)$ such that $\underline{u}<\bar{u}$, $d>0$ and $g \geq 1/\underline{u}-1$.
  Let $W>0$ be a random variable satisfying $P(W>w) \leq b w^c e^{-lw}$ for $w\in (\underline{u},\bar{u})$ and some $b>0$, $c\geq 0$, $l>0$.

If $ld=1$, then
  $$
  \int_{\underline{u}}^{\bar{u}} P\left(W > d \log \left( \frac{g}{1/u-1} \right)   \right) du <
\frac{b}{g^{ld}} \left[ d \log \left( \frac{g}{1/\bar{u} -1} \right)  \right]^c
\log(1/\underline{u}).
  $$

If $ld < 1$, then
$$
  \int_{\underline{u}}^{\bar{u}} P\left(W > d \log \left( \frac{g}{1/u-1} \right)   \right) du <
\frac{b}{g^{ld}} \left[ d \log \left( \frac{g}{1/\bar{u} -1} \right)  \right]^c
\frac{1}{1-ld} \left( \frac{\bar{u}}{1 - \bar{u}} \right)^{1-ld}.
$$

If $ld > 1$, then
$$
  \int_{\underline{u}}^{\bar{u}} P\left(W > d \log \left( \frac{g}{1/u-1} \right)   \right) du <
\frac{b}{g^{ld}} \left[ d \log \left( \frac{g}{1/\bar{u} -1} \right)  \right]^c
\frac{( 1/\underline{u}-1 )^{ld-1}}{ld-1} .
$$
\label{LEM:INTBOUND_EXPTAILS}
\end{lemma}

\begin{lemma}
  Let $\underline{u},\bar{u} \in (0,1)$ such that $\underline{u}<\bar{u}$ and let $d>0$, $g \geq 1/\underline{u}-1$.
  Let $W>0$ be a random variable satisfying $P(W>w) \leq b/w^c$ for all $w \in (\underline{u},\bar{u})$ and some $b>0$, $c>1$.

  If $c-1 > \log \left( \frac{g}{1/\underline{u}-1} \right)$ then
  $$
  \int_{\underline{u}}^{\bar{u}} P\left(W > d \log \left( \frac{g}{1/u-1} \right)   \right) du <
\frac{b}{d^c (c-1)}
\left[ \frac{1}{\log \left(\frac{g}{1/\underline{u}-1}\right)} \right]^{c-1}.
  $$
  and otherwise
  $$
  \int_{\underline{u}}^{\bar{u}} P\left(W > d \log \left( \frac{g}{1/u-1} \right)   \right) du <
\frac{b}{d^c \left[\log \left( \frac{g}{1/\underline{u}-1} \right) \right]^c}.
  $$
\label{LEM:INTBOUND_POLYTAILS}
\end{lemma}

\subsection*{Tail integral bounds for chi-square and F distributions}

Lemmas \ref{LEM:CHISQTAIL_INTBOUND}-\ref{LEM:FTAIL_INTBOUND} adapt
Lemmas \ref{LEM:INTBOUND_EXPTAILS}-\ref{LEM:INTBOUND_POLYTAILS} to chi-square and F-distributed random variables $W$.
For simplicity Lemmas \ref{LEM:CHISQTAIL_INTBOUND}-\ref{LEM:FTAIL_INTBOUND}
state asymptotic bounds as $g \rightarrow \infty$
and one should think of $\alpha$ and $d$ as being arbitrarily close to 1 and 2 (respectively),
but the proofs also provide non-asymptotic bounds.

\begin{lemma}
Let $g>0$ and $W \sim \chi^2_\nu$ where $\nu$ may depend on $g$ and,
as $g \rightarrow \infty$, $\log g \gg \nu$.

\noindent (i) Let $d \geq 2$ and $\alpha \in (0,1)$ be fixed constants. Then, as $g \rightarrow \infty$,
$$
\int_0^1 P \left( W > d \log \left( \frac{g}{1/u-1} \right) \right) du
\preceq \frac{1}{g} \left(\frac{4e}{\nu}\right)^{\nu/2} \left[ \log \left( \frac{g}{e^{\nu/4}} \right) \right]^{\frac{\nu}{2}+1}  \ll \frac{1}{g^\alpha}.
$$

\noindent (ii) Let $d \in (1,2)$ and $\alpha \in (0,1)$ be fixed constants. Then, 
$$
\int_0^1 P \left( W > d \log \left( \frac{g}{1/u-1} \right) \right) du
\leq \left(\frac{1}{g}\right)^{\frac{d}{2}(2-\frac{d}{2})}
\left[ 2 + \frac{d e}{\nu (1-\frac{d}{2})} \log \left( g^{\frac{d}{2}(2-\frac{d}{2}) + 1} \right) \right].
$$
Further, let $g \geq 1$. Since $d \in (1,2)$, then $\frac{d}{2}(2-\frac{d}{2}) \in (d-1,d-1/4)$ and the right-hand side above is
$$
\leq \frac{2}{g^{d-1}}\left[1+ \frac{d e}{\nu (2-d)} \log( g^{d+3/4} ) \right],
$$
which is $ \ll  \frac{1}{g^{\alpha (d-1)}}$ as $g \rightarrow \infty$.
\label{LEM:CHISQTAIL_INTBOUND}
\end{lemma}

%

\begin{lemma}
  Let $F \sim \mathcal{F}_{\nu_1,\nu_2}$, $d>1$ be a fixed constant and $g>0$ be a function of $(\nu_1,\nu_2)$.

(i) Assume that $\nu_1 \ll \log(g) \ll \nu_2$ as $\nu_2 \rightarrow \infty$.
Then, for any fixed $\alpha \in (0,d-1)$,
$$
\int_{0}^{1} P \left( \nu_1 F > d \log \left( \frac{g}{1/u-1}  \right) \right) du
\preceq e^{\frac{\nu_1}{d(2-\sqrt{3})}} \frac{1}{g} + \frac{e^{-\nu_1}}{g^{d-1 - d\sqrt{ \frac{4}{\nu_2} \log \left(\frac{g}{e^{\nu_1/2}} \right)  }}}
\ll \frac{1}{g^\alpha}.
$$

(ii) Assume that $\nu_1 \ll \nu_2 \ll \log^\gamma(g)$ as $\nu_2 \rightarrow \infty$ for all fixed $\gamma <1$. Then,
for any fixed $\alpha <1$,
$$
\int_{0}^{1} P \left( \nu_1 F > d \log \left( \frac{g}{1/u-1}  \right) \right) du
\ll \exp\left\{-\frac{\nu_1+\nu_2-5}{2} \log \left( \frac{\log^\alpha(g)}{\nu_1+\nu_2-6} \right) \right\}.
$$
\label{LEM:FTAIL_INTBOUND}
\end{lemma}

\section{Discussion of regularity conditions}
\label{supplsec:regcond}

In \ref{suppsec:regcond_modelpriors} we specialize the regularity conditions C1-C2 on prior sparsity to the cases where one sets the model prior $p(M_k)$ to be either the uniform, Beta-Binomial or Complexity priors.

\subsection{Conditions C1-C2 for the uniform, Beta-Binomial and Complexity prior}
\label{suppsec:regcond_modelpriors}

Lemma \ref{LEM:CONSISTENCYCOND_MODELPRIORS} gives sufficient conditions for C1-C2 to hold when $p(M_k)$ is the uniform, Beta-Binomial and Complexity priors defined in \eqref{eq:prior_model}. The conditions are given separately for spurious models $m \in S$ and non-spurious models, and involve the sample size $n$, the problem dimension $p$, the size of the optimal model $p_t$, and the size of the signal $\btheta_t^*$ as measured by the non-centrality parameter in \eqref{eq:ncp_regression}.

\begin{lemma}
Let $M_m$ be a model $p_m \in [p_t,\bar{p}]$ variables,
 and $\tau$ be the prior dispersion parameter in Conditions (C1)-(C2).
  \begin{enumerate}[label=(\roman*)]
  \item Uniform prior.   If $\tau n \gg 1$ then (C1)-(C2) hold.

  \item Beta-Binomial(1,1) prior. If $\tau n \gg [\bar{p}/(p-p_t)]^2$ then (C1)-(C2) hold.

  \item Complexity prior. If $\tau n \gg 1$ then (C1)-(C2) hold.
  \end{enumerate}

Let $M_m \in S^c$ be a model with $p_m < p_t$ variables.
\begin{enumerate}[label=(\roman*)]
\item Uniform prior. (C2) holds if and only if $\frac{\lambda_{tm}}{2 \log(\lambda_{tm})} \gg (p_t-p_m) \log(\tau n)$.
\item Beta-Binomial(1,1) prior. If 
$$\frac{\lambda_{tm}}{2 \log(\lambda_{tm})} \gg   (p_t-p_m) \log\left(\frac{\sqrt{n\tau} (p-p_t)}{ep_m}\right)$$ then (C2) holds.
\item Complexity$(c)$ prior. If 
$$\frac{\lambda_{tm}}{2 \log(\lambda_{tm})} \gg   (p_t-p_m) \log\left(\frac{\sqrt{n\tau} p^c (p-p_t)}{ep_m}\right)$$ then (C2) holds.

\end{enumerate}
\label{LEM:CONSISTENCYCOND_MODELPRIORS}
\end{lemma}

The interpretation of Lemma \ref{LEM:CONSISTENCYCOND_MODELPRIORS} is as follows.
Consider first a model of size at least as that of $M_t$, i.e. $p_m \geq p_t$, then one may take any $\tau n \gg 1$. 

This is a truly minimal requirement that even allows $\tau$ to decrease in $n$, although as discussed in high dimensions the custom is that $\tau \succeq 1$ is either fixed or grows with $n$.
Consider now models of smaller size than $p(M_t)$, that is models favored over $M_t$ when setting sparse $p(M_k)$ and $\tau$.
Then there is a limit to the prior sparsity, dictated by the signal strength (the non-centrality parameter $\lambda_{tm}$).

For example, under the Complexity prior (C2) holds when $\lambda_{tm}/\log(\lambda_{tm}) \gg p_t \log(\sqrt{\tau} n p^{1+c})$, 

and note that the Beta-Binomial prior corresponds to $c=0$.
The larger $c$, the larger $\lambda_{tm}$ needs to be if one wishes to attain pairwise consistency.

\subsection{Comparison to conditions in existing literature}
\label{supplseq:regcond_comparison}

We discuss connections between the conditions assumed by \cite{narisetty:2014}, \cite{castillo:2015}, \cite{yangyun:2016} and \cite{yangyun:2017}
and our model complexity conditions (B1)-(B2) and prior sparsity conditions (C1)-(C2).
We start by explaining where the score of our results differs from this literature.

\subsubsection*{{\bf Scope}}

The main result of \cite{castillo:2015} on model selection consistency (Corollary 1) 
proves that $E_{f^*}(p(M_t \mid \by))$ converges to 1 but, rather than giving the rate of said convergence, gives a minimax analysis that focuses on a worst-case $\btheta_t^*$.
Also, the results are restricted to the case where the prior on the models $p(M_k)$ is a Complexity prior, which is critically needed to prove that supersets of $M_t$ receive vanishing posterior probability $P(S \mid \by)$. 
In contrast, we give rates as a function of $\btheta_t^*$, consider some settings where the model may be misspecified,
and allow for more general $p(M_k)$ (made possible by restricting the largest model size $\bar{p}$).
A similar comment regarding the restrictiveness of the prior sparsity structure applies to \cite{narisetty:2014}, \cite{yangyun:2016} and \cite{yangyun:2017}.
\cite{narisetty:2014} prove that when the linear regression error variance $\phi$ is known the posterior probability of $M_t$ converges to 1, and give separate rates for spurious models $m \in S$, large and small non-spurious models $m \in S^c$ (Theorem 4.1 and Lemma 4.2). The results are for a spike-and-slab prior on the coefficients where the prior variance $\tau$ must grow with $n$ at a fairly fast rate, to attain consistency (see below). The result for the unknown variance $\phi$ case requires a restriction on the maximum model size $\bar{p} \preceq n/\log p $, analogous to our Condition B1.
\cite{yangyun:2016} consider a restrictive prior setting where $p(M_k)$ is a complexity prior and $\tau$ grows sufficiently fast with $n$ (see below). Their analysis is restricted to Zellner's prior, but this is less critical.
\cite{yangyun:2017} allow for very general families of likelihoods, including for example non-iid Gaussian regression, non-parametric regression and density estimation. The priors can also be fairly general in terms of the chosen distributional family, but are subject to a key requirement (their prior anti-concentration Assumption B2) that leads to diffuse parameter priors ($\tau$ growing with $n$, in our notation, see below). Also in terms of scope, the main result (Theorem 4) shows that $p(M_t \mid \by) \stackrel{L_1}{\longrightarrow} 1$, but does not describe the associated convergence rates.

\subsubsection*{{\bf \cite{castillo:2015}}}

The main technical are so-called compatibility, smallest sparse singular value (SSV) and mutual coherence conditions. 
We focus on the latter two, as they lead to simpler interpretation.
The SSV condition basically says that the smallest non-zero singular value of sub-matrices $X_k$ for models $M_k$ of size $p_k$ equal to (a multiple of) $M_t$ is bounded away from zero, as $n$ grows.
Mutual coherence is a stronger condition on the largest absolute pairwise correlation between columns in $X$. For example, if the rows in $X$ are iid random variables, then their framework can recover models of dimension $p_t \leq \sqrt{n/\log n}$ (\cite{castillo:2015}, Section 2).
Provided these conditions hold and $\btheta^*$ is sufficiently large (essentially, a beta-min condition $\min_{j \in M_t} |\theta_j^*|^2 >  p_t (\log p)/n$),
the authors prove $E_{f^*} p(M_t \mid \by) \stackrel{L_1}{\longrightarrow} 1$.

There are connections to our assumptions (B2) and (C2).
Our Condition (B2) requires a milder $p_t \ll n$.
Regarding (C2), recall that 
\begin{align}
  \lambda_{tm} \geq n v_{tm} (\btheta^*)' \btheta^*/\phi^* \geq n v_{tm} p_t \min_j |\theta_j^*|^2/\phi^*,
\label{eq:boundncp_suppl}
\end{align}
where $v_{tm}$ is the smallest non-zero eigenvalue of $X_t'(I-H_m)X_t/n$, and $H_m=X_m(X_m'X_m)^{-1}X_m'$ is the projection matrix onto the column space of $X_m$. 
Under the SSV condition $v_{tm}$ is bounded away from zero for models of size up to $p_t$.
A sufficient condition for our (C2) to hold under the Complexity(c) prior is that
$\lambda_{tm}/\log(\lambda_{tm}) \gg p_t \log(\sqrt{\tau n}) + (1+c)p_t \log p$ (Section \ref{suppsec:regcond_modelpriors}).
Hence, it essentially suffices that
\begin{align}
\min_{j \in M_t} (\theta_j^*)^2/\phi^* \gg [\log(\sqrt{\tau n}) + (1+c) \log p]/n,
\label{eq:suffcond_c2}
\end{align}
and recall that $c=0$ corresponds to the Beta-Binomial prior. Hence our condition on the signal strength is slightly milder when  $p_t \gg \log(\tau n)$, and slightly stronger when $p_t \ll \log(\tau n)$.


\subsubsection{{\bf \cite{yangyun:2016}}} 

The authors focus on a restrictive setting where $p(M_k)$ is a Complexity prior and the prior dispersion $\tau$ must be fairly large.
Specifically, by combining their Assumption C and their sparse projection condition in Assumption B, they require $\tau n \asymp p^{2 \alpha}$, for some fixed $\alpha \geq 1/2$,  and $c + \alpha \succeq p_t$, where $c$ is the parameter in the Complexity prior. That is, $\tau$ and $c$ must grow with $p$ and $p_t$, respectively.
Our Conditions (C1)-(C2) allow $\tau$ and $c$ to be constant, for example, i.e. a significantly less sparse prior setting.
In particular, we allow for $c=0$, which corresponds to the Beta-Binomial prior.

\cite{yangyun:2016} make a further Assumption D that basically requires that $p_t \leq n/\log p$, which is similar to our Condition (B2) that $p_t \ll n$.
Further technical conditions include a lower restricted eigenvalue condition (their Assumption B) that $X_k'X_k/n$ has smallest eigenvalue bounded away from 0 (for models up to a certain size),
and that the signal strength satisfies
$\min_j  (\theta_j^*)^2/\phi^* \geq (\alpha + c + p_t) (\log p) / n$.
This is stronger than \eqref{eq:suffcond_c2} (given the authors' assumption that $\tau n \asymp p^{2 \alpha}$), which suffices to guarantee our Condition (C2) under the Complexity prior (the setting in \cite{yangyun:2016}),



\subsubsection*{{\bf \cite{narisetty:2014}}}

The main assumptions made by the authors regarding the sparsity of the prior are as follows.
First, they assume that $n \tau v \asymp \max\{n , p^{2+\delta} \}$ for some $\delta>0$ (their Condition 4.2),
where $v$ is the smallest non-zero eigenvalue of $X_k'X_k/n$ among models of size $p_k \leq \min\{ n/(2\log p), p \}$.
Said $v$ is assumed $\gg 1/p^\kappa$, for some small power $\kappa>0$ (Condition 4.5).
This limits attention to diffuse priors where $\tau$ grows fairly quickly with $n$.
A second assumption is that the marginal prior inclusion probabilities satisfy $P(\theta_j \neq 0) \asymp 1/p$, limiting attention to sparse $p(M_k)$.
For example, as $p \rightarrow \infty$ then the prior distribution on the model size is approximately $P(p_k=l)=1/p^l$, i.e. roughly equivalent to a Complexity prior with $c=1$.
The authors also assume that $p_t$, the dimension of $M_t$, is fixed (Condition 4.3).
We relax these assumptions by allowing for smaller $\tau$, less sparse $p(M_k)$, and $p_t$ that may grow with $n$.
In terms of signal strength, Condition 4.4 in \cite{narisetty:2014} assumes that 
$\lambda_{kt} > 5 p_t (1+\delta) \log (\max\{\sqrt{n},p\})$
for all models $M_k$ of size $p_k$ at most a multiple of $p_t$.
This is very similar to our Condition C2. 
For example, using \eqref{eq:boundncp_suppl} and the restricted eigenvalue condition that led to our sufficient condition for C2 given in \eqref{eq:suffcond_c2} (for the particular case of the Complexity prior),
their signal strength requirement is basically satisfied when
$\min_{j \in M_t} (\theta_j^*)^2/\phi^* \succeq n^{-1} \log (\max\{\sqrt{n},p\})$,
which is similar to \eqref{eq:suffcond_c2}.

\subsection{{\bf \cite{yangyun:2017}}}

The main assumption made by the authors related to prior sparsity (anti-concentration Assumption B2) is that
that the prior $p(\btheta_k \mid M_k)$ under spurious models $k \in S$ assign sufficiently small mass to a neighborhood of $\btheta_t^*$.
A direct comparison to our prior setting is not straightforward. This is due to the neighborhood being defined in terms of general distance measures,
and the specific size of the neighborhood not being explicitly stated (rather it is assumed a neighborhood of suitable size exists).
It is possible to draw connections, however.
First, note that the prior probability assigned to neighborhoods of $\btheta_t^*$ is controlled by the prior dispersion $\tau$ (in our notation), hence the anti-concentration assumption refers to the range of allowable $\tau$.
Then, for Gaussian priors the anti-concentration condition can be interpreted as $\tau$ behaving in a way such that, as one approaches $\btheta_t^*$, the prior probability of the neighborhood decreases at a slower-than-exponential rate in the distance to $\btheta_t^*$. If $\tau$ were fixed then the prior probability would decay at the exponential rate, implying that $\tau$ must grow with $n$.
Also, the authors state that the size of the neighborhood is expected to increase with the model size $p_k$, implying that $\tau$ would need to grow faster for larger models.
This effectively leads to setting increasingly diffuse priors, i.e. $\tau \gg 1$ in our notation.

In their Gaussian process variable selection application, \cite{yangyun:2017} also use a sparse model space prior $p(M_k)$ with marginal inclusion probabilities $P(\theta_j \neq 0)$ of order $1/p$, i.e. $p(M_k)$ is a Complexity prior with parameter $c=1$. Such a sparse $p(M_k)$ is not required by their main model selection consistency result (Theorem 4), however.

Regarding assumptions on the signal strength, their Assumption B3 requires that the Bayes factor between any non-spurious $M_k$ of size $p_k < p_t$ decreases exponentially in $n$. Our Condition C2 is analogous, and also guarantees that said Bayes factor is exponential in $n$.
Specifically in Gaussian regression their Assumption B3 requires the beta-min conditions in \cite{wainwright:2009information}.
If $\min_j |\theta_j^*|^2/\phi^* \gg 1/(n-p_t) $ a sufficient beta-min condition is that $p_t \log(p/p_t) < n$
\begin{align}
\min_j |\theta_j^*|^2/\phi^* > n^{-1} \log(p - p_t),
\nonumber
\end{align}
which is analogous to \eqref{eq:suffcond_c2} sufficing for our Condition (C2) to hold.



\section{Proof of Proposition \ref{PROP:CORRECTSEL_FROM_PP}}

The result follows from simple observations and Markov's inequality.
We first prove the result when $\hat{k}= \arg\max_k p(M_k \mid \by)$ is the posterior mode.
If $p(M_t \mid \by) > 1/2$ then $\hat{k}=t$, therefore
\begin{align}
P_{f^*}(\hat{k}=t) \geq P_{f^*} \left( p(M_t \mid \by) \geq 1/2 \right).
\label{eq:probphat_pp}
\end{align}
This implies
\begin{align}
P_{f^*}(\hat{k} \neq t) \leq
P_{f^*}\left( p(M_t \mid \by) < 1/2 \right)=
P_{f^*}\left( \sum_{k \neq t} p(M_k \mid \by) \geq 1/2 \right)
\leq 2 E_{f^*} \left( \sum_{k \neq t} p(M_k \mid \by) \right)
\nonumber
\end{align}
where the right-hand side follows from Markov's inequality.
Therefore
$$
P_{f^*}(\hat{k}=t) \geq 1 - 2 E_{f^*} \left( \sum_{k \neq t} p(M_k \mid \by) \right),
$$
as we wished to prove.

For the case where $\hat{k}$ is the median probability model it suffices to prove that
if $p(M_t \mid \by) > 1/2$ then $\hat{k}=t$,
since then \eqref{eq:probphat_pp} holds and all subsequent arguments remain valid.
Let $\theta_j$ be the $j^{th}$ element in $\btheta$,
$\Theta_t= \Theta_1^* \times \ldots \times \Theta_p^*$ be the parameter space for the KL-optimal $M_t$,
and recall that the median probability model decides that $\theta_j \in \Theta_j^*$ if and only if $P(\theta_j \in \Theta_j^* \mid \by) \geq 1/2$.
Since $p(M_t \mid \by)= P \left( \bigcap_{j=1}^p \theta_j \in \Theta_j^* \mid \by \right)$, we have that
$$
1-p(M_t \mid \by)= P \left( \bigcup_{j=1}^p \theta_j \not\in \Theta_j^* \mid \by \right) \geq P(\theta_j \not\in \Theta_j^* \mid \by)
$$
for all $j=1,\ldots,p$.
Therefore if $P(\theta_j \not\in \Theta_j^* \mid \by)>1/2$ for any $j=1,\ldots,p$ that implies that $p(M_t \mid \by) < 1/2$.
Equivalently,
$p(M_t \mid \by)\geq 1/2$ implies that
$P(\theta_j \in \Theta_j^* \mid \by)\geq 1/2$ for all $j=1,\ldots,p$, as we wished to prove.

\section{Proof of Corollary \ref{COR:POW_TYPEI_FROM_PP}}

The event that one makes any type I or any type II error implies that $\hat{k} \neq t$, hence the probability of making a type I or II error is $\leq P_{f^*}(\hat{k} \neq t)$.

%

\section{Proof of Corollary \ref{COR:POW_TYPEI_FROM_MARGPP}}

Consider first the case where truly $\gamma_j^*=0$.
\begin{align}
 P_{f^*}(\hat{\gamma}_j=1) = P_{f^*} \left( P(\gamma_j=1 \mid \by) > t \right)=
P_{f^*} \left( \sum_{\gamma_j=1} p(\bgamma \mid \by) > t \right)
\leq \frac{1}{t} \sum_{\gamma_j=1} E_{f^*} \left( p(\bgamma \mid \by)\right),
\nonumber
\end{align}
by Markov's inequality.

Consider now the case where truly $\gamma_j^*=1$. Then
\begin{align}
  P_{f^*}(\hat{\gamma}_j=0) = P_{f^*} \left( P(\gamma_j=1 \mid \by) \leq t \right)=
  P_{f^*} \left( \sum_{\gamma_j=0} p(\bgamma \mid \by) > 1-t \right)
  \leq \frac{1}{1-t} \sum_{\gamma_j=0} E_{f^*}( p(\bgamma \mid \by) ),
\nonumber
\end{align}
again by Markov's inequality, as we wished to prove.

\section{Proof of Lemma \ref{LEM:TAILPROB_TO_MEAN}}

%
%


By definition 
$$p(M_k \mid \by)=\left(1+ \sum_{l \neq k}^{} B_{lk}p(M_l)/p(M_k) \right)^{-1} \leq \left(1+ B_{tk}p(M_t)/p(M_k) \right)^{-1},
$$
the right-hand side following from observing that $B_{lk} \geq 0$ and $p(M_l)/p(M_k) \geq 0$ for all $l$.
Denote by $U_k=(1+ B_{tk}p(M_t)/p(M_k))^{-1}$, so that $E_{f^*}(p(M_k \mid \by)) \leq E_{f^*}(U_k)$. 
Since $U_k \in [0,1]$ is a positive random variable its expectation is given by the integral of its survival (or right-tail probability) function, that is
\begin{align}
E_{f^*}(U_k)= \int_0^1 P_{f^*}(U_k>u) du \leq
\int_0^1 P_{f^*} \left(B_{kt} > \frac{p(M_t)}{(1/u-1) p(M_k)} \right) du,
\label{eq:boundu_cdf}
\end{align}
since $P_{f^*}(U_k>u)= P_{f^*}\left(B_{tk}< (1/u-1) p(M_k)/p(M_t)\right)$, as we wished to prove.



\section{Proof of Lemma \ref{LEM:TIGHTNESS_PAIRWISEBOUND}}
\label{proof:tightness_pairwisebound}

The goal is to prove that
\begin{align}
\frac{1 - p(M_t \mid \by)}{\sum_{k \neq t} (1 + B_{tk} p(M_k)/p(M_t))^{-1}} \stackrel{L_1}{\longrightarrow} 1.
\nonumber
\end{align}
The proof strategy is to show that the term above is bounded between $p(M_t \mid \by)$ and 1 so that, if $p(M_t \mid \by)$ converges to 1, then so does the term above.

First, note that $1- p(M_t \mid \by)= \sum_{k \neq t} p(M_k \mid \by)$ and that
$$
p(M_k \mid \by) \leq \left(1 + B_{tk} \frac{p(M_t)}{p(M_k)} \right)^{-1} 
$$
and hence
\begin{align}
 \frac{1 - p(M_t \mid \by)}{\sum_{k \neq t} (1 + B_{tk} p(M_k)/p(M_t))^{-1}} =
 \frac{\sum_{k\neq t} p(M_k \mid \by)}{\sum_{k \neq t} (1 + B_{tk} p(M_k)/p(M_t))^{-1}} 
\leq 1.
\nonumber
\end{align}

Further, using that $p(M_k \mid \by)= p(\by \mid M_k) p(M_k)/p(\by)$ and upper-bounding the denominator gives
\begin{align}
\frac{1 - p(M_t \mid \by)}{\sum_{k \neq t} (1 + B_{tk} p(M_k)/p(M_t))^{-1}} \geq
\frac{\frac{1}{p(\by)} \sum_{k \neq t} p(\by \mid M_k) p(M_k)}
{ \sum_{k \neq t} \frac{p(\by \mid M_k) p(M_k)}{p(\by \mid M_t) p(M_t)} }= p(M_t \mid \by).
\nonumber
\end{align}
The result follows from noting that $p(M_t \mid \by) \stackrel{L_1}{\longrightarrow} 1$ by assumption.

\section{Proof of Theorem \ref{THM:NECESSARY_CONDITIONS_PC}}
\label{proof:necessary_conditions_pc}

We first outline and discuss the technical conditions required by Theorem 4.1 in \cite{hjort:2011}, which establishes the asymptotic Normality of the maximum likelihood estimator (allowing for model misspecification) for models where the log-likelihood is concave.
These conditions are near-minimal, and basically require the existence of certain hessian matrices.
We also discuss the technical conditions of Proposition 8 in \cite{rossell:2019t} establishing the asymptotic equivalence between the marginal likelihood and its Laplace approximation. These follow from the conditions of Theorem 4.1 in \cite{hjort:2011}, plus our assumption that $\tilde{p}(\btheta_k/\tau^{1/2}, \phi \mid M_k)$ is bounded away from 0 and $\infty$.
Finally, we elaborate on some details required to complete the proof outline.

\subsection{Technical conditions}
\label{ssec:conditions_hjort}

Theorem 4.1 in \cite{hjort:2011} assumes that the log-likelihood $\log p(\by \mid \btheta_k,\phi)$ is concave in $\eta_k=(\btheta_k,\phi)$, and that the expected log-likelihood under the data-generating $f^*(\by)$ has a unique global maximum
\begin{align}
\bm{\eta}_k^*= E_{f^*}(\log p(\by \mid \btheta_k,\phi))
\nonumber
\end{align}
that lies outside the boundary of the parameter space $\Theta_k \times \Phi$.
The existence of such a maximum is guaranteed if $\Theta_k \times \Phi$ is a compact set, see \cite{pollard:1991}.

The theorem also assumes a condition on the remainder term $R()$ associated to a local expansion around $\bm{\eta}_k^*$.
Specifically, suppose that $\log p(\by, \bm{\eta}_k^* + t ) - \log p(\by \mid \bm{\eta}_k^*)= D(\by)' t + R(\by,t)$
where $t \in \mathbb{R}^{p_k}$ and $E_{f^*}(D(y))$ has zero mean and finite covariance matrix $W_k$, and that the remainder term satisfies
\begin{align}
 E_{f^*} \left[ \log p(\by, \bm{\eta}_k^* + t ) - \log p(\by \mid \bm{\eta}_k^*) \right]= E_{f^*} [ R(\by,t) ]= -\frac{1}{2} t' H_k t + o(t' t)
\nonumber
\end{align}
as well as $\mbox{Var}_{f^*} R(y_i, t)= o(t' t)$, where $H_k$ is symmetric and positive-definite.
In ordinary cases where $\log p(\by \mid \btheta_k,\phi)$ is smooth then $D(\by)$ is its gradient at $\btheta_k^*$, and one can find a remainder $R(\by,t)$ of a quadratic Taylor expansion satisfying the conditions above for $H_k$ being the expected hessian and $W_k$ the covariance of the log-likelihood gradient under $f^*$. That is, for
\begin{align}
H_k= - E_{f^*} \nabla_{\eta_k}^2 \log p(y_i \mid \btheta_k^*)
\nonumber 
W_k= \mbox{Cov}_{f^*} \nabla_{\eta_k}^2 \log p(y_i \mid \btheta_k^*)
\nonumber
\end{align}
Intuitively, in regression models, $H_k$ involves the expectation of $\bx_{ki} \bx_{ki}'$ times the model-predicted variance for $y_i$, where $\bx_i \in \mathbb{R}^{p_k}$ is the vector of covariates for individual $i$ under model $M_k$, whereas $H$ involves $\bx_{ki} \bx_{ki}'$ times the quadratic error when predicting $y_i$ from the model-based mean at $\bm{\eta}_k^*$.
See \cite{hjort:2011} (Sections 5-6) for specific expressions for $H_k$ and $W_k$ for logistic and Cox regression, and \cite{rossell:2019t} for accelerated failure time models and probit regression.

Proposition 8 in \cite{rossell:2019t} requires three conditions, numbered D1-D3 in that paper. Conditions D1-D2 follow immediately from Theorem 4.1 in \cite{hjort:2011}. Condition D3 is satisfied by the fact that $p(\by \mid \btheta_k,\phi,M_k)$ and 
$p(\btheta_k \mid \phi,M_k)= \tau^{p_k/2} \tilde{p}(\btheta_k/\tau^{1/2} \mid \phi, M_k)$ 
are continuous functions in $(\btheta_k,\phi)$, and that $\tilde{p}(\btheta_k/\tau^{1/2}, \phi, M_k)$ is bounded away from 0 and $\infty$ by assumption.

\subsection{Proof outline}
\label{ssec:outline_necessary_conditions_pc}

For any model $M_k$, denote the whole parameter vector by $\eta_k= (\btheta_k,\phi)$, let $\bm{\eta}_k^*= E_{f^*}(\log p(\by \mid \btheta_k,\phi))$ be its optimal value minimizing KL-divergence to $f^*$, and $\hat{\bm{\eta}}_k$ be its maximum likelihood estimator (which is unique, from the assumption of strict concavity).
The proof strategy is to note that
\begin{align}
 p(M_t \mid \by) = \left( 1 + \sum_{k \neq t} B_{kt} \frac{p(M_k)}{p(M_t)} \right)^{-1} \leq \left( 1 + B_{mt} \frac{p(M_m)}{p(M_t)} \right)^{-1},
\nonumber
\end{align}
where $B_{mt}$ is the Bayes factor between $M_t$ and $M_m$. 
Hence, for $p(M_t \mid \by)$ to converge in probability to 1, it is necessary that $B_{mt} p(M_m)/p(M_t)$ converges in probability to 0.
To show that the latter does not hold, we use that, by Proposition 8 in \cite{rossell:2019t} we have that $B_{mt}/B_{mt}^* \stackrel{P}{\longrightarrow} 1$, where
$B_{mt}^*$ can be interpreted as an asymptotic Laplace approximation to $B_{mt}$. Specifically,
\begin{align}
 B_{mt}^*= e^{L_{mt}/2}
\tau^{(p_t-p_m)/2} \frac{\tilde{p}(\bm{\eta}_m^* \mid M_m)}{\tilde{p}(\bm{\eta}_t^* \mid M_t)} 
\left( \frac{2\pi}{n} \right)^{\frac{p_m-p_t}{2}} \frac{|H_t|^{1/2}}{|H_m|^{1/2}}
\nonumber
\end{align}
where $L_{mt}=2[\log p(\by \mid \hat{\bm{\eta}}_m, M_m) - \log p(\by \mid \hat{\bm{\eta}}_t, M_t)]$ is the likelihood-ratio statistic to test $M_t$ versus $M_m$.
Note that $|H_t|/|H_m|$ is a non-zero finite constant that does not depend on $n$, since both $H_t$ and $H_m$ are positive-definite by assumption,
and that $\tilde{p}(\bm{\eta}_m^* \mid M_m) / \tilde{p}(\bm{\eta}_t^* \mid M_t)$ is also bounded above and below by a finite constant by assumption.
Hence, it suffices to show that
\begin{align}
e^{L_{mt}/2} \left( \frac{2\pi}{\tau n} \right)^{\frac{p_m-p_t}{2}}
\label{eq:asymp_bf_logconcave}
\end{align}
converges to 0 in probability.

Consider first the case when $m \in S$ is a spurious model.
In Section \ref{ssec:further_necessary_conditions_pc} we show that, then $L_{mt}= O_p(1)$ under the data-generating $f^*$. Hence
\begin{align}
B_{mt}^* \frac{p(M_m)}{p(M_t)}=  \left( \frac{2\pi}{\tau n} \right)^{\frac{p_m-p_t}{2}}  \frac{p(M_m)}{p(M_t)} O_p(1).
\label{eq:asymp_bf_logconcave}
\end{align}
If $\lim_{n \rightarrow \infty} (\tau n)^{-(p_m-p_t)/2} p(M_m)/p(M_t) \neq 0$, 
that is $(\tau n)^{-(p_m-p_t)/2} p(M_m)/p(M_t) \succeq 1$ in our notation, then 
$B_{mt}^* p(M_m)/p(M_t)$ does not converge in probability to 0. This implies that $p(M_t \mid \by)$ does not converge in probability to 1, as we wished to prove..

Consider now the case when $m \in S^c$ is a non-spurious model.
Then the law of large numbers and the fact that $\log p(\by \mid \bm{\eta}_k,M_k)/n$ converges to $E_{f^*}(p(\by \mid \bm{\eta}_k,M_k))$ uniformly in $\bm{\eta}_k$ (see the Convexity Lemma in \cite{pollard:1991}) imply that
\begin{align}
&\frac{L_{mt}}{n}= \frac{\lambda_{tm}}{n} (1 + o_p(1))
\nonumber
\end{align}
where
\begin{align}
\lambda_{mt}&= E_{f^*}[\log p(\by \mid \bm{\eta}_m^*,M_m) - \log p(\by \mid \bm{\eta}_t^*,M_t)]
\nonumber \\ 
&=\mbox{KL}(f^*, \log p(\by \mid \bm{\eta}_t^*,M_t)) - \mbox{KL}(f^*, \log p(\by \mid \bm{\eta}_m^*,M_m)) < 0.
\nonumber
\end{align}

Note that this difference between KL-divergences $\lambda_{mt}$ is strictly negative (since by definition $M_t$ is the smallest model minimizing KL-divergence to $f^*$, and also by definition any non-spurious $m \in S^c$ has a larger KL-divergence to $f^*$). Hence we obtain that

\begin{align}
\frac{1}{n} \log \left( B_{mt}^*  \frac{p(M_m)}{p(M_t)} \right) =
\frac{\lambda_{mt}}{2 n} (1+o_p(1))  - \frac{p_m - p_t}{2n} \log(\tau n) + \frac{1}{n} \log \left( \frac{p(M_m)}{p(M_t)} \right) + O(1/n).
\nonumber
\end{align}

Therefore, since $\lambda_{tm}=-\lambda_{mt}$ by definition, if
\begin{align}
\lim_{n \rightarrow \infty} -\frac{\lambda_{tm}}{2} - \frac{p_m - p_t}{2} \log(\tau n) + \log \left( \frac{p(M_m)}{p(M_t)} \right) \neq -\infty,
\label{eq:condA2_equiv}
\end{align}
then 
\begin{align}
 \log \left( B_{mt}^*  \frac{p(M_m)}{p(M_t)} \right)
\nonumber
\end{align}
does not converge in probability to $-\infty$.
That is, $B_{mt}^*  \frac{p(M_m)}{p(M_t)}$ does not converge in probability to 0, which implies that $p(M_t \mid \by)$ does not converge in probability to 1, as we wished to prove.
Note that \eqref{eq:condA2_equiv} can be equivalently stated as
\begin{align}
 \frac{\lambda_{tm}}{2} \preceq -\frac{p_m - p_t}{2} \log(\tau n) + \log \left( \frac{p(M_m)}{p(M_t)} \right),
\nonumber
\end{align}
completing the proof.

\subsection{Further arguments}
\label{ssec:further_necessary_conditions_pc}

To complete the proof we show that, when $m \in S$ is a spurious model, then $L_{mt}= O_p(1)$.
Under the conditions of Theorem 4.1 in \cite{hjort:2011}, it follows that
$\sqrt{n}( \hat{\eta}_k - \eta_k^*) \stackrel{D}{\longrightarrow} N(0, V_k^{-1})$,
where $V_k= H_k^{-1} W_k H_k^{-1}$, for both models $k \in \{t,m\}$.
Further, it also implies that the log-likelihood admits a quadratic expansion around $\bm{\eta}_k^*$ such that
\begin{align}
 \log p(\by \mid \hat{\bm{\eta}}_k, M_k) - \log p(\by \mid \bm{\eta}_k^*, M_k)= 
\frac{n}{2} (\hat{\bm{\eta}}_k - \bm{\eta}_k^*)' \Sigma_k^{-1} (\hat{\bm{\eta}}_k - \bm{\eta}_k^*) [ 1 + o_p(1) ]
\nonumber
\end{align}
where $\Sigma_k= -H_k^{-1}$.

Since $M_t$ is nested in $M_m$, without loss of generality let $\bm{\eta}_m'=(\bm{\eta}_{m1}',\bm{\eta}_{m2}')$ where
$\bm{\eta}_{m2}$ are the parameters shared with $M_t$, and $\bm{\eta}_{m1}$ those only featuring in $M_m$.
By definition, $\bm{\eta}_m^*= (0, \bm{\eta}_t^*)$ and 
$\log p(\by \mid \bm{\eta}_m^*, M_m)= \log p(\by \mid \bm{\eta}_t^*, M_t)$, hence

\begin{align}
L_{mt}=2[\log p(\by \mid \hat{\bm{\eta}}_m, M_m) - \log p(\by \mid \hat{\bm{\eta}}_t, M_t) \pm \log p(\by \mid \hat{\bm{\eta}}_t^*, M_t)]=
\nonumber \\
n [(\hat{\bm{\eta}}_m - \bm{\eta}_m^*)' \Sigma_m^{-1} (\hat{\bm{\eta}}_m - \bm{\eta}_m^*)
- (\hat{\bm{\eta}}_t - \bm{\eta}_t^*)' \Sigma_t^{-1} (\hat{\bm{\eta}}_t - \bm{\eta}_t^*)] [ 1 + o_p(1) ]
\label{eq:quadformdif_logconcave_nested}.
\end{align}
The goal is hence to show that the difference between quadratic forms in the right-hand side of \eqref{eq:quadformdif_logconcave_nested} is $O_p(1)$.
Decompose $\Sigma_m$ into the blocks defined by $(\bm{\eta}_{m1}',\bm{\eta}_{m2}')$, that is
\begin{align}
 \Sigma_m = \begin{pmatrix} \Sigma_{11} & \Sigma_{12} \\ \Sigma_{21} & \Sigma_{22} \end{pmatrix}.
\nonumber
\end{align}
To ease notation let $A= (\Sigma_{11} - \Sigma_{12} \Sigma_{22}^{-1} \Sigma_{21})^{-1}$. Using the block-wise inversion formula,
\begin{align}
 \Sigma_m^{-1}= \begin{pmatrix} 
A & -\Sigma_{22}^{-1} \Sigma_{21} A \\
- A' \Sigma_{12} \Sigma_{22}^{-1} & \Sigma_{22}^{-1} + \Sigma_{22}^{-1} \Sigma_{21} A \Sigma_{12} \Sigma_{22}^{-1}
\end{pmatrix}
\nonumber
\end{align}
it is possible to re-arrange terms as follows
\begin{align}
 [(\hat{\bm{\eta}}_m - \hat{\bm{\eta}}_m^*)' \Sigma_m^{-1} (\hat{\bm{\eta}}_m - \hat{\bm{\eta}}_m^*)=
(\hat{\bm{\eta}}_{m2} - \hat{\bm{\eta}}_t^*)' \Sigma_{22}^{-1} (\hat{\bm{\eta}}_{m2} - \hat{\bm{\eta}}_t^*)
+ \hat{\bm{\eta}}_{m1}' A \hat{\bm{\eta}}_{m1}
- 2 (\hat{\bm{\eta}}_{m2} - \bm{\eta}_t^*)' \Sigma_{22}^{-1} \Sigma_{21} A \hat{\bm{\eta}}_{m1}
\nonumber \\
+ (\hat{\bm{\eta}}_{m2} - \hat{\bm{\eta}}_t^*)' \Sigma_{22}^{-1} \Sigma_{21} A \Sigma_{12} \Sigma_{22}^{-1} (\hat{\bm{\eta}}_{m2} - \hat{\bm{\eta}}_t^*)
= 
(\hat{\bm{\eta}}_{m2} - \hat{\bm{\eta}}_t^*)' \Sigma_{22}^{-1} (\hat{\bm{\eta}}_{m2} - \hat{\bm{\eta}}_t^*)
+ (\hat{\bm{\eta}}_{m1} - \bm{\mu})' A (\hat{\bm{\eta}}_{m1} - \bm{\mu}),
\nonumber
\end{align}
where $\bm{\mu}= \Sigma_{12} \Sigma_{22}^{-1} (\hat{\bm{\eta}}_{m2} - \hat{\bm{\eta}}_t^*)$.

Let us characterize the first of these two terms. Since $\sqrt{n} (\hat{\bm{\eta}}_m - \bm{\eta}_m^*) \stackrel{D}{\longrightarrow} N(0,V_m^{-1})$,
with $\bm{\eta}_{m2}^*=\bm{\eta}_t^*$ and $\bm{\eta}_{m1}^*=0$, it follows that
\begin{align}
\Sigma_{22}^{-1/2} (\hat{\bm{\eta}}_{m2} - \hat{\bm{\eta}}_t^*) \stackrel{D}{\longrightarrow} N(0,\Sigma_{22}^{-1/2} V_m^{-1} \Sigma_{22}^{-1/2})
\Rightarrow (\hat{\bm{\eta}}_{m2} - \hat{\bm{\eta}}_t^*)' \Sigma_{22}^{-1} (\hat{\bm{\eta}}_{m2} - \hat{\bm{\eta}}_t^*) = O_p(1/n).
\nonumber
\end{align}
since $V_m$ and $\Sigma_{22}$ are finite, positive-definite matrices.
Note that the same argument shows that $\Sigma_{22}^{-1/2} (\hat{\bm{\eta}}_{m2} - \hat{\bm{\eta}}_t^*)= O_p(1/n)$.
Regarding the second term,
$\sqrt{n} (\hat{\bm{\eta}}_{m1} - \bm{\mu})$ is asymptotically normally-distributed with mean 
$E_{f^*}(\hat{\bm{\eta}}_{m1}) - \Sigma_{12} \Sigma_{22}^{-1} E( \hat{\bm{\eta}}_{m2} - \hat{\bm{\eta}}_t^*)= 0$ and finite, positive-definite covariance
\begin{align}
 \mbox{Cov}_{f^*}(\hat{\bm{\eta}}_{m1} - \bm{\mu})
= \begin{pmatrix} I & - \Sigma_{12} \Sigma_{22}^{-1} \\ 0 & I \end{pmatrix} V_m^{-1}  
\begin{pmatrix} I & 0 \\ - \Sigma_{22}^{-1} \Sigma_{21}  & I \end{pmatrix}.
\nonumber
\end{align}
Therefore, $(\hat{\bm{\eta}}_{m1} - \bm{\mu})' A (\hat{\bm{\eta}}_{m1} - \bm{\mu})= O_p(1/n)$.

In conclusion, this shows that \eqref{eq:quadformdif_logconcave_nested} is
\begin{align}
L_{mt}= n [O_p(1/n) + O_p(1/n) + O_p(1/n)] [1 + o_p(1)]= O_p(1),
\nonumber
\end{align}
as we wanted to prove.

\section{Proof of Proposition \ref{PROP:ZELLNER_PHIKNOWN}}
\label{proof:zeller_phiknown}

\subsection{Zellner's prior and known variance $\phi^*$. Spurious models}
\label{proof:outline_zellner_known_spurious}

Let $m \in S$ be a spurious model. Under $f^*(\by)= N(\by; X_t \btheta_t^*; \phi^* I)$ it is well-known that
$W_{mt}/\phi^* \sim \chi_{p_m-p_t}^2$ (see Lemma \ref{LEM:NESTEDLM_CHITEST} for a formal proof),
and by Lemma \ref{LEM:TAILPROB_TO_MEAN}
\begin{align}
E_{f^*}\left(p(M_m\mid\by)\right) < \int_0^1
P_{f^*} \left( \frac{W_{mt}}{\phi^*} >
\frac{1+\tau n}{\tau n} 2 \log \left[ \frac{(1+\tau n)^{\frac{p_m-p_t}{2}} p(M_t)}{p(M_m) (1/u-1)} \right]
 \right) du.
\label{eq:proof_zellnerknown_expr1}
\end{align}
To bound this integral we use Lemma \ref{LEM:CHISQTAIL_INTBOUND}.
For simplicity we drop the term $(1+\tau n)/(\tau n)$, since it is asymptotically equal to 1 (Condition (B1) implies $\tau n \gg 1$).
Specifically, in Lemma \ref{LEM:CHISQTAIL_INTBOUND} set $\nu=p_m-p_t$, $g= (1+\tau n)^{\frac{p_m-p_t}{2}} p(M_t)/p(M_m)$
and note that $\log(g) \gg  \nu$ as $n\rightarrow \infty$ under Conditions (B1) (which guarantees $\tau n \gg 1$) and (C1).
Then, by Lemma \ref{LEM:CHISQTAIL_INTBOUND},
$$E_{f^*}(p(M_m\mid\by)) \preceq \frac{[\log(g)]^{(p_m-p_t)/2}}{g}  \ll \frac{[p(M_m)/p(M_t)]^\alpha}{(1+\tau n)^{\alpha (p_m-p_t)/2} }$$
 for any fixed $\alpha \in (0,1)$.

As a final remark, inspection of the proof of Lemma \ref{LEM:CHISQTAIL_INTBOUND} shows that the result applies for all $n \geq n_0$, for some fixed $n_0$ that is strictly increasing in $p_m - p_t$. Hence, one may obtain a common $n_0$ that applies to all spurious models by taking that associated to the largest spurious model (i.e. with size $\bar{p}-p_t$).

\subsection{Zellner's prior and known variance $\phi^*$. Non-spurious models}
\label{proof:outline_zellner_known_nonspurious}

Let $m \in S^c$ be a non-spurious model, and define $M_q=M_t \cup M_m$ and $\lambda_{tm}$ as in \eqref{eq:ncp_regression}.
Clearly 
$$W_{mt}=\hat{\btheta}_m'X_m'X_m\hat{\btheta}_m - \hat{\btheta}_q'X_q'X_q\hat{\btheta}_q
+ \hat{\btheta}_q'X_q'X_q\hat{\btheta}_q - \hat{\btheta}_t'X_t'X_t\hat{\btheta}_t =W_{qt} - W_{qm},$$
where, by Lemma \ref{LEM:NESTEDLM_CHITEST},
$W_{qt}/\phi^* \sim \chi_{p_q-p_t}^2$,
$W_{qm}/\phi^* \sim \chi_{p_q-p_m}^2(\lambda_{tm})$.
To ease notation let $b_n(u)= 2 \log ( (1+\tau n)^{\frac{p_m-p_t}{2}} (p(M_t)/p(M_m))/(1/u-1) )$.
Then the integrand in \eqref{eq:proof_zellnerknown_expr1} is equal to
\begin{align}
P_{f^*} \left( \frac{W_{mt}}{\phi^*} > b_n(u)  \right)
= P_{f^*} \left( \frac{W_{qt} - W_{qm}}{\phi^*} > \frac{b_n(u)}{2} + c_n(u) - (c_n(u) - \frac{b_n(u)}{2} )  \right)
\nonumber \\
\leq
P_{f^*} \left( \frac{W_{qm}}{\phi^*} < c_n(u) - \frac{b_n(u)}{2}  \right)
+ P_{f^*} \left( \frac{W_{qt}}{\phi^*} > c_n(u) + \frac{b_n(u)}{2} \right)
\label{eq:tailprob_nonspur_zellner_known}
\end{align}
for any $c_n(u)>0$, the right-hand side following from the union bound. 
The idea is to set $c_n(u)$ to a convenient expression that gives a small value for \eqref{eq:tailprob_nonspur_zellner_known}. 

Specifically we shall take $c_n(u)= 0.5 b_n(u) + \lambda_{tm}/\log(\lambda_{tm})$. Then
\begin{align}
P_{f^*} \left( \frac{W_{qm}}{\phi^*} < c_n(u) - \frac{b_n(u)}{2}  \right)
= P_{f^*} \left( \frac{W_{qm}}{\phi^*} < \frac{\lambda_{tm}}{\log(\lambda_{tm})} \right).
\nonumber
\end{align}
Since $W_{qm}/\phi^* \sim \chi_{p_q-p_m}^2(\lambda_{tm})$ and  $\lambda_{tm} \gg 1$ by assumption, we have that $\lambda_{tm}/\log(\lambda_{tm}) \ll \lambda_{tm}$ and we may apply the non-central chi-square left tail bound in Lemma \ref{LEM:INEQ_NCCHISQ_LEFT}, obtaining
\begin{align}
P_{f^*} \left( \frac{W_{qm}}{\phi^*} < c_n(u) - \frac{b_n(u)}{2}  \right) \leq
\exp \left\{  -\frac{\lambda_{tm}}{2} \left( 1 - \frac{1}{\log(\lambda_{tm})}  \right)^2 \right\}
\ll  e^{- \gamma \lambda_{tm}/2},
\nonumber
\end{align}
for any fixed $\gamma \in (0,1)$.

Regarding the second term in \eqref{eq:tailprob_nonspur_zellner_known},
\begin{align}
P_{f^*} \left( \frac{W_{qt}}{\phi^*} > c_n(u) + \frac{b_n(u)}{2} \right)=
P_{f^*} \left( \frac{W_{qt}}{\phi^*} > 2 \log \left( \frac{(1+\tau n)^{\frac{p_m-p_t}{2}}p(M_t) e^{\frac{\lambda_{tm}}{2\log(\lambda_{tm})}}}{p(M_m)(1/u-1)} \right) \right).
\nonumber
\end{align}
Since $W_{qt}/\phi^* \sim \chi^2_{p_q-p_t}$, this term can be bounded similarly to \eqref{eq:proof_zellnerknown_expr1}.
Specifically, we use the central chi-square provided by Lemma \ref{LEM:CHISQTAIL_INTBOUND}.
Note that to apply Lemma \ref{LEM:CHISQTAIL_INTBOUND} we must have that the threshold is larger than the degrees of freedom $p_q - p_t$, specifically
\begin{align}
 \frac{(1+\tau n)^{\frac{p_m-p_t}{2}}p(M_t) e^{\frac{\lambda_{tm}}{2\log(\lambda_{tm})}}}{p(M_m)} \gg p_q - p_t, 
\nonumber
\end{align}
which holds from Condition (C2), since $p_q - p_t \leq p_m$.
Then, setting $\nu= p_q - p_t$ and 
$g= (1+\tau n)^{\frac{p_m-p_t}{2}}p(M_t) e^{\frac{\lambda_{tm}}{2\log(\lambda_{tm})}} / p(M_m)$
into Lemma \ref{LEM:CHISQTAIL_INTBOUND} gives that
\begin{align}
P_{f^*} \left( \frac{W_{qt}}{\phi^*} > c_n(u) + \frac{b_n(u)}{2} \right)
 \preceq \frac{[\log(g)]^{(p_m-p_t)/2}}{g}  
\ll \left[ \frac{p(M_m) e^{-\frac{\lambda_{tm}}{2\log(\lambda_{tm})}}}{p(M_t) (1+\tau n)^{(p_m-p_t)/2} } \right]^\alpha,
\nonumber
\end{align}
for any fixed $\alpha <1$.
Therefore,
\begin{align}
 E_{f^*}(p(M_m \mid \by)) \ll
e^{- \gamma \lambda_{tm}/2} + e^{-\frac{\lambda_{tm}^\gamma}{2}}
\left[ \frac{p(M_m) }{p(M_t) (1+\tau n)^{(p_m-p_t)/2} } \right]^\gamma,
\nonumber
\end{align}
for any fixed $\gamma <1$. 
The proof is completed by showing that the second term in this expression is of a larger order, since
\begin{align}
\frac{e^{- \gamma \lambda_{tm}/2}}{ e^{-\frac{\lambda_{tm}^\gamma}{2}} \left[ \frac{p(M_m) }{p(M_t) (1+\tau n)^{(p_m-p_t)/2} } \right]^\gamma}
=
\frac{e^{- \frac{\gamma \lambda_{tm}}{2} \left(1 - \frac{1}{\gamma \lambda_{tm}^{1-\gamma}} \right)}}{ \left[ \frac{p(M_m) }{p(M_t) (1+\tau n)^{(p_m-p_t)/2} } \right]^\gamma}
\ll 1,
\label{eq:proof_zellnerknown_largerterm}
\end{align}
the right-hand side following from $\lim_{n \rightarrow \infty} 1/\lambda^{1 - \gamma}=0$ along with Condition (C2).

In conclusion, there exists a finite $n_0$ such that
\begin{align}
 E_{f^*}(p(M_m \mid \by)) \leq
 e^{-\frac{\lambda_{tm}^\gamma}{2}} \left[ \frac{p(M_m) }{p(M_t) (1+\tau n)^{(p_m-p_t)/2} } \right]^\gamma,
\nonumber
\end{align}
for all $n \geq n_0$, as we wished to prove. 
As a final remark, inspection of the proof of Lemma \ref{LEM:INEQ_NCCHISQ_LEFT} and Lemma \ref{LEM:CHISQTAIL_INTBOUND} shows that 
$n_0$ is strictly decreasing in $\lambda_{tm}$, hence one may obtain a common $n_0$ across non-spurious models $M_m$ by taking a minimum of $\lambda_{tm}$ over such $m$.

\section{Tightness of Proposition \ref{PROP:ZELLNER_PHIKNOWN}}
\label{tightness:zellner_phiknown}

To illustrate that the upper-bounds provided by Proposition \ref{PROP:ZELLNER_PHIKNOWN} (i)-(ii) are fairly tight, we shall obtain corresponding lower-bounds.
We assume that $p(M_t \mid \by) \stackrel{P}{\longrightarrow} 1$, since otherwise one does not even attain posterior consistency, and the tightness of the rate becomes less relevant.

The strategy is simple. First, note that the $L_1$ convergence rate of $p(M_m \mid \by)$ to 0 cannot be faster than the rate at which $p(M_m \mid \by)$ converges to 0 in probability (since lack of convergence in probability implies lack of $L_1$ convergence).
Second, $p(M_m \mid \by)$ can be shown to converge to 0 in probability at the same rate as $(1 + B_{tm} p(M_t)/p(M_m))^{-1}$, hence by studying the convergence in probability of $(1 + B_{tm} p(M_t)/p(M_m))^{-1}$ one obtains a lower bound for the $L_1$ convergence rate of $p(M_k \mid \by)$.

To see that $p(M_m \mid \by)$ converges to 0 in probability at the same rate as $(1 + B_{tm} p(M_t)/p(M_m))^{-1}$, note that
\begin{align}
p(M_m \mid \by)= \left( 1 + B_{tm} \frac{p(M_t)}{p(M_m)} + \sum_{k \neq m,t} B_{km} \frac{p(M_k)}{p(M_m)} \right)^{-1}
= \left( 1 + B_{tm} \frac{p(M_t)}{p(M_m)} [1 + o_p(1)] \right)^{-1}
\label{eq:ppsinglemodel_op}
\end{align}
since
\begin{align}
 \frac{\sum_{k \neq m,t} B_{km} p(M_k)/p(M_m)}{B_{tm} (p(M_t)/p(M_m))}=
 \frac{\sum_{k \neq m,t} p(\by \mid M_k) p(M_k)}{p(\by \mid M_t) p(M_t)}
\leq \sum_{k \neq t} B_{kt} \frac{p(M_k)}{p(M_t)} \stackrel{P}{\longrightarrow} 0,
\nonumber
\end{align}
the right-hand side following from the assumption $p(M_t \mid \by)= (1 + \sum_{k \neq t} B_{kt} p(M_k)/p(M_t))^{-1} \stackrel{P}{\longrightarrow} 1$. Therefore, using \eqref{eq:ppsinglemodel_op},
\begin{align}
\frac{(1 + B_{tm} p(M_t)/p(M_m))^{-1}} {p(M_m \mid \by)}
= \frac{\left( 1 + B_{tm} p(M_t)/p(M_m) [1 + o_p(1)] \right)}{1 + B_{tm} p(M_t)/p(M_m)}
= 1 + \frac{\left( B_{tm} p(M_t)/p(M_m) o_p(1) \right)}{1 + B_{tm} p(M_t)/p(M_m)}
\stackrel{P}{\longrightarrow} 1,
\nonumber
\end{align}
since $B_{tm} p(M_t)/p(M_m) \stackrel{P}{\longrightarrow} \infty$.
That is, $p(M_m \mid \by)$ converges to 0 in probability at the same rate as $(1 + B_{tm} p(M_t)/p(M_m))^{-1}$.
All that is left is determining the latter.
Recall that 
\begin{align}
 B_{tm}= \exp \left\{ -\frac{\tau n}{2\phi^*(1+\tau n)} W_{mt} \right\}
(1+\tau n)^{\frac{p_m-p_t}{2}}
\nonumber
\end{align}
where $W_{mt}=\hat{\btheta}_m'X_m'X_m\hat{\btheta}_m - \hat{\btheta}_t'X_t'X_t\hat{\btheta}_t$
is the difference between residual sums of squares under $M_t$ and $M_m$ and $\hat{\btheta}_m= (X_m'X_m)^{-1} X_m' \by$ the least-squares estimate.

Consider first a spurious model $m \in S$, and recall that then $W_{mt}/\phi^* \sim \chi_{p_m-p_t}^2= O_p(1)$ (Lemma \ref{LEM:NESTEDLM_CHITEST}).
Therefore $B_{tm}= (\tau n)^{(p_m-p_t)/2} O_p(1)$ and
\begin{align}
 \left(1 + B_{tm} \frac{p(M_t)}{p(M_m)} \right)^{-1}= (1 + (1+\tau n)^{(p_m-p_t)/2} \frac{p(M_t)}{p(M_m)} O_p(1))^{-1}= 
(\tau n)^{-\frac{p_m-p_t}{2}} \frac{p(M_m)}{p(M_t)} O_p(1),
\nonumber
\end{align}
showing that $E_{f^*}(p(M_m \mid \by)) \succeq (\tau n)^{-(p_m-p_t)/2} p(M_m)/p(M_t)$, as we wished to prove.

Consider now a non-spurious model $m \in S^c$. For simplicity suppose that $M_m$ is a submodel of $M_t$ (the case when $M_m \not\subset M_t$ follows similarly, see Section \ref{proof:zeller_phiknown}), and that $p_t \ll \lambda_{tm}$. Then by Lemma \ref{LEM:NESTEDLM_CHITEST} we have that $-W_{tm}/\phi^* \sim \chi_{p_t-p_m}^2(\lambda_{tm})$.
Using the definition of a non-central chi-square and that $p_t - p_m \ll \lambda_{tm}$ by assumption, it is easy to show that
$-W_{tm}/\phi^*= \lambda_{tm}[1 + o_p(1)]$.
This is because $-W_{tm}/\phi^*$ is equal in distribution to
\begin{align}
(Z_1 + \sqrt{\lambda}_{tm})^2 + \sum_{j=2}^{p_t-p_m} Z_j^2= 
\lambda_{tm} + \sqrt{\lambda}_{tm} Z_1 + \sum_{j=1}^{p_t-p_m} Z_j^2= 
\lambda_{tm} \left( 1 + \frac{Z_1}{\sqrt{\lambda}_{tm}} + \frac{\sum_{j=1}^{p_t-p_m} Z_j^2}{\lambda_{tm}} \right),
\nonumber
\end{align}
where $Z_j \sim N(0,1)$ independently across $j=1,\ldots,p_t-p_m$. Since $\lambda_{tm} \gg p_t-p_m$ by assumption, it follows that
$Z_1/\sqrt{\lambda}_{tm}=o_p(1)$ and $\sum_{j=1}^{p_t-p_m} Z_j^2/\lambda_{tm}= \sum_{j=1}^{p_t-p_m} Z_j^2/(p_t-p_m) o(1)= o_p(1)$.

Therefore
$$
B_{tm}= (1+\tau n)^{(p_m-p_t)/2} e^{\lambda_{tm}[1 + o_p(1)]/2}
$$
and hence $  (1 + B_{tm} p(M_t)/p(M_m))^{-1}=$
\begin{align}
\left(1 + p(M_t)/p(M_m) (1+\tau n)^{(p_m-p_t)/2} e^{\lambda_{tm}[1 + o_p(1)]/2} \right)^{-1}= 
\frac{p(M_m)e^{-\lambda_{tm}/2 [1 + o_p(1)]}}{p(M_t) (\tau n)^{(p_m-p_t)/2}} O_p(1),
\nonumber
\end{align}
which under Condition (C2) is $e^{-0.5 \lambda_{tm} [1+o_p(1)]} O_p(1)$.
Therefore $p(M_m \mid \by)$ converges in probability to 0 at a rate determined by $e^{-0.5 \lambda_{tm}}$, up to lower-order terms.

\section{Proof of Proposition \ref{PROP:ZELLNER_PHIUNKNOWN}}
\label{proof:zellner_phiunknown}

The proof runs analogously to that for the known variance $\phi^*$ case.
Here $\tilde{F}_{mt}$ plays the role of $W_{mt}$ in \eqref{eq:bf_zellner_known} when $\phi^*$ was assumed known.
For precision the right-hand side in \eqref{eq:bf_zellner} does not hold if $p_m=p_t$,
as then $\tilde{F}_{mt}=\infty$,
but our upcoming argument still applies
by lower bounding $\tilde{s}_{tm} \geq \tilde{s}_{tm'}$  in \eqref{eq:bf_zellner},
where $M_{m'}$ adds any single variable to $M_m$.
To apply Lemma \ref{LEM:TAILPROB_TO_MEAN} note that
$P_{f^*}(B_{mt}> (p(M_t)/p(M_m))/(1/u-1))=$
\begin{align}
P_{f^*} \left( \tilde{F}_{mt} > \frac{n-p_m}{p_m-p_t} \left[
\frac{(1+\tau n)^{\frac{p_m-p_t}{n+a_\phi}}}{((\frac{1}{u}-1) [p(M_m)/p(M_t)])^{\frac{2}{n+a_\phi}}} -1
 \right] \right)
< P_{f^*} \left( (p_m-p_t)F_{mt} > b_n(u) \right),
\label{eq:tailprob_zellner_unknown}
\end{align}
where $b_n(u)= 2 [(n-p_m)/(n+a_\phi)] \log ( (1+\tau n)^{\frac{p_m-p_t}{2}} (p(M_t)/p(M_m))/(1/u-1) )$.
The right hand side follows from $\log(z) \leq z-1$ and is written for analogy with
\eqref{eq:tailprob_zellner_known}.

\subsection{Zellner's prior and unknown $\phi^*$. Spurious models}

Let $m \in S$ be a spurious model, then $F_{mt} \sim \mathcal{F}(p_m-p_t,n-p_m)$ where $\mathcal{F}$ is an F-distribution.
In Lemma \ref{LEM:FTAIL_INTBOUND} set
$\nu_1=p_m-p_t$, $\nu_2=n-p_m$, $g=(1+\tau n)^{(p_m-p_t)/2} p(M_t)/p(M_m)$ and $d=2(n-p_m)/(n+a_\phi)$.
Note that Lemma \ref{LEM:FTAIL_INTBOUND} requires that $\log g \gg p_m - p_t$, which holds under under Conditions (B1) and (C1).
Therefore, if $\log(g) \ll n-p_m$, then Lemma \ref{LEM:FTAIL_INTBOUND} gives that $E_{f^*}(p(M_m\mid\by))<$
\begin{align}
 \int_0^1 P_{f^*}((p_m-p_t) F_{mt} > b_n(u)) du
\preceq \left(\frac{1}{g} \right)^{1 - 4 \sqrt{\frac{\log(g)}{n-p_m}}}
\ll \frac{[p(M_m)/p(M_t)]^\alpha}{(\tau n)^{\alpha \frac{p_m-p_t}{2}}}
\label{eq:meanbound_spur_zellunknown}
\end{align}
for any fixed $\alpha <1$.
Also by Lemma \ref{LEM:FTAIL_INTBOUND}, if $\log(g) \gg n-p_m$ then
$$
E_{f^*}(p(M_m\mid\by)) \ll \exp\left\{- \frac{(n-p_t-5)}{2} \log \left( \frac{\log^\gamma (g)}{n-p_t-6} \right) \right\}
\ll e^{-\kappa n}
$$
for any fixed $\gamma <1$, $\kappa >1$, since $p_t \ll n$ by Condition (B2).

\subsection{Zellner's prior and unknown $\phi^*$. Non-spurious models}

Consider now $m \in S^c$ and let $M_q= M_m \cup M_t$. Since $s_q \leq s_m$, then
\begin{align}
  (p_m-p_t) F_{mt} \leq
\frac{n-p_m}{n-p_q}
\left( (p_q-p_t) F_{qt} - (p_q-p_m) F_{qm} \right),
\label{eq:fstat_bound_nonspur}
\end{align}
where marginally $F_{qt} \sim \mathcal{F}(p_q-p_t,n-p_q)$
and $F_{qm} \sim \mathcal{F}(\lambda_{tm},p_q-p_m,n-p_q)$.
Combining \eqref{eq:fstat_bound_nonspur} and \eqref{eq:tailprob_zellner_unknown},
and proceeding analogously to \eqref{eq:tailprob_nonspur_zellner_known}, we obtain that
$P_{f^*} \left( (p_m-p_t) F_{mt} > b_n(u) \right) <$
\begin{align}
P_{f^*} \left( (p_q-p_m) F_{qm} < \left(c_n(u) - \frac{b_n(u)}{2} \right) \frac{n-p_q}{n-p_m}   \right)
+ P_{f^*} \left( (p_q-p_t) F_{qt} > c_n(u) + \frac{b_n(u)}{2} \right),
\label{eq:ftail_bound_nonspur}
\end{align}
for any $c_n(u)>0$. Specifically we set $c_n(u)= 0.5 b_n(u) + \lambda_{tm}/\log(\lambda_{tm})$, therefore

\begin{align}
P_{f^*} \left( (p_q-p_m) F_{qm} < \left(c_n(u) - \frac{b_n(u)}{2} \right) \frac{n-p_q}{n-p_m}   \right)
= P_{f^*} \left( (p_q-p_m) F_{qm} < \frac{\lambda_{tm}}{\log(\lambda_{tm})} \frac{n-p_q}{n-p_m} \right).
\nonumber
\end{align}
Since $F_{qm} \sim \mathcal{F}(\lambda_{tm},p_q-p_m,n-p_q)$ and  $\lambda_{tm} \gg 1$ by assumption, we may apply the left tail bound in Lemma \ref{LEM:INEQ_NCF_LEFT}.

Specifically, in Lemma \ref{LEM:INEQ_NCF_LEFT} set $\nu_1= p_q-p_m$, $\nu_2= n-p_q$,
$w= [\lambda_{tm} / \log(\lambda_{tm}) ] (n-p_q)/(n-p_m)$ and take arbitrarily large but fixed $s$.
Then $sw\ll \lambda_{tm}$ and
\begin{align}
& P_{f^*} \left( (p_q-p_m) F_{qm} < \frac{\lambda_{tm}}{\log(\lambda_{tm})} \frac{n-p_q}{n-p_m} \right)
\nonumber \\
& \leq e^{-\frac{\lambda_{tm}}{2} \left( 1 - \sqrt{\frac{s (n- p_q)}{\log(\lambda_{tm}) (n - p_m)}} \right)^2} + e^{-\frac{n-p_q}{2}(s-1-\log s)}
\ll e^{-\gamma \frac{\lambda_{tm}}{2}} + e^{- \kappa n}
\label{eq:proof_zellunknown_nonspur_term1}
\end{align}
for any fixed $\gamma \in <1$ and $\kappa >0$, since $p_q \ll n$ by Condition (B1).

Regarding the second term in \eqref{eq:ftail_bound_nonspur}, it is equal to
\begin{align}
P_{f^*} \left( (p_q-p_t) F_{qt} > 2 \frac{n-p_m}{n+a_\phi} \log \left( \frac{(1+\tau n)^{\frac{p_m-p_t}{2}} p(M_t) e^{0.5 \lambda_{tm}/\log(\lambda_{tm})}}{ p(M_m) (1/u-1)} \right) \right).
\nonumber
\end{align}
To bound this term we use Lemma \ref{LEM:FTAIL_INTBOUND}.
For simplicity we may drop $(n-p_m)/(n+a_\phi)$, since it converges to 1 as $n \rightarrow \infty$ under Condition (B1).
Specifically, in Lemma \ref{LEM:FTAIL_INTBOUND} set $\nu_1=p_q-p_t$, $\nu_2=n-p_q$, 
and $g= (1 + \tau n)^{\frac{p_m-p_t}{2}} e^{0.5 \lambda_{tm}/\log(\lambda_{tm})} p(M_t)/p(M_m)$.
Then if $p_q - p_t  \log(g) \ll n-p_q$, Lemma \ref{LEM:FTAIL_INTBOUND} Part (i) gives that
\begin{align}
\int_0^1 P_{f^*} \left( (p_q-p_t) F_{qt} > 2 \log \left( \frac{g}{1/u-1} \right)  \right) du 
\ll \frac{1}{g^\alpha}
\ll e^{-\frac{\lambda_{tm}^\gamma}{2}}
\left[ \frac{p(M_m) }{p(M_t) (1+\tau n)^{(p_m-p_t)/2} } \right]^\gamma
\label{eq:proof_zellunknown_nonspur_term2}
\end{align}
for any fixed $\gamma < \alpha < 1$.
In contrast, if $\log(g) \gg n-p_q$ then by Lemma \ref{LEM:FTAIL_INTBOUND} Part (ii) the integral is $\ll e^{-\kappa n}$ for any $\kappa >0$, since $p_q \ll n$.

Combining \eqref{eq:proof_zellunknown_nonspur_term1} and \eqref{eq:proof_zellunknown_nonspur_term2} gives
\begin{align}
 E_{f^*}(p(M_m \mid \by)) \ll
\max \left\{  
e^{- \gamma \lambda_{tm}/2},  e^{-\frac{\lambda_{tm}^\gamma}{2}}
\left[ \frac{p(M_m) }{p(M_t) (1+\tau n)^{(p_m-p_t)/2} } \right]^\gamma, e^{- \kappa n}
\right\},
\nonumber
\end{align}
for any fixed $\gamma <1$, $\kappa >0$.
Then, arguing as in \eqref{eq:proof_zellnerknown_largerterm} to show that the first term is of a smaller order than the second under Condition (C2), gives
\begin{align}
 E_{f^*}(p(M_m \mid \by)) \ll
\max \left\{  e^{-\frac{\lambda_{tm}^\gamma}{2}}
\left[ \frac{p(M_m) }{p(M_t) (1+\tau n)^{(p_m-p_t)/2} } \right]^\gamma, e^{- \kappa n}
\right\},
\nonumber
\end{align}
as we wished to prove.

\section{Proof of Proposition \ref{PROP:NORMALPRIOR}}
\label{proof:normalprior}

The proof proceeds analogously to that for Zellner's prior in Section \ref{proof:zellner_phiunknown}, with suitable adjustments.
First, proceeding as in \eqref{eq:tailprob_zellner_unknown} gives that 
\begin{align}
P_{f^*} \left(B_{mt}> \frac{p(M_t)}{p(M_m) (1/u-1)} \right) < P_{f^*} \left( (p_m-p_t)\tilde{F}_{mt} > b_n(u) \right)
\nonumber
\end{align}
where $b_n(u)= 2 \log ( (p(M_t)/p(M_m)) (c_{mt} \tau n)^{(p_m-p_t)/2}/(1/u-1) )$.

Second, elementary arguments detailed in Lemma \ref{LEM:FSTAT_NORMALPRIOR_BOUND} show that
$(p_m-p_t) \tilde{F}_{mt} \leq (p_m-p_t) F_{mt} + p_m F_{m0}/(1+\tau n\rho_{tp_t})$, hence
\begin{align}
P_{f^*} \left( (p_m-p_t)\tilde{F}_{mt} > b_n(u) \right) <
P_{f^*} \left( (p_m-p_t) F_{mt} > b_n(u) -1 \right)
+ P_{f^*} \left( p_m F_{m0} > 1+ \tau n \rho_{t p_t} \right).
\label{eq:tailprob_normalprior}
\end{align}
The proof is based on noting that the first term in \eqref{eq:tailprob_normalprior} is analogous to \eqref{eq:tailprob_zellner_unknown},
and that under Condition (D2) the second term is bounded by the F right tails in Corollary \ref{COR:INEQ_F_RIGHT}.
Specifically, since $F_{m0} \sim \mathcal{F}_{p_m,n-p_m}(\lambda_{m0})$ where
$\lambda_{m0} \leq \lambda_{t0} \leq \tau n \rho_{tp_t}$ under (D2),
Corollary \ref{COR:INEQ_F_RIGHT}(ii) implies that $P_{f^*}(p_m F_{m0}> 1+ \tau n \rho_{tp_t}) \ll e^{-\tau n \rho_{tp_t}/2}$.

To complete the proof we characterize the first term in \eqref{eq:tailprob_normalprior}.
Briefly, for any spurious model $m \in S$ proceeding as in \eqref{eq:meanbound_spur_zellunknown} gives that
\begin{align}
\int_0^1 P_{f^*} \left( (p_m-p_t) F_{mt} > b_n(u) -1 \right) du \ll
\max \left\{ [p(M_m)/p(M_t)]^\alpha (\tau n)^{\alpha (p_t-p_m)/2}, e^{-\kappa n} \right\}
\label{eq:tail_normalprior_term1}
\end{align}
for any fixed $\alpha <1$, $\kappa >0$.
This gives that
$$E_{f^*}(p(M_m \mid \by)) \ll \max \left\{ [p(M_m)/p(M_t)]^\alpha (\tau n)^{\alpha (p_t-p_m)/2}, e^{-\kappa n}, e^{- \tau n \rho_{tp_t}/2} \right\}$$
for $m \in S$, as we wished to prove.

Consider now a non-spurious model $m \in S^c$. Then, proceeding as in \eqref{eq:ftail_bound_nonspur} shows that
$$
\int_0^1 P_{f^*} \left( (p_m-p_t) F_{mt} > b_n(u)-1 \right) du \ll
\max \left\{  
e^{- \frac{\gamma \lambda_{tm}}{2}},  e^{-\frac{\lambda_{tm}^\gamma}{2}}
\left[ \frac{p(M_m) }{p(M_t) (\tau n)^{(p_m-p_t)/2} } \right]^\gamma, e^{- \kappa n}
\right\},
$$
for any fixed $\gamma < 1$, $\kappa >0$.
Under Condition (C2) the term $e^{- \frac{\gamma \lambda_{tm}}{2}}$ is asymptotically smaller, which gives that
\begin{align}
 E_{f^*}(p(M_m \mid \by)) \ll
\max \left\{  
e^{-\frac{\lambda_{tm}^\gamma}{2}}
\left[ \frac{p(M_m) }{p(M_t) (\tau n)^{(p_m-p_t)/2} } \right]^\gamma, e^{- \kappa n}, e^{-\tau n \rho_{tp_t}/2}
\right\},
\nonumber
\end{align}
as we wished to prove.

\section{Proof of Proposition \ref{PROP:MOMPRIOR_SIMPLE}}
\label{proof:momprior_simple}

The proof rests on two auxiliary results, Proposition \ref{PROP:MOMPRIOR} and Proposition \ref{PROP:PENALTY_NLP}, stated below.

Proposition \ref{PROP:MOMPRIOR} shows that $E_{f^*} p(M_m \mid \by)$ can be bounded by the sum of two terms.
The first term is analogous to the rates obtained under a Normal prior, except that it features an acceleration factor $h_n^{p_m-p_t}$, where $h_n \geq 0$ can in principle be chosen arbitrarily.
Proposition \ref{PROP:PENALTY_NLP} states that, for suitably-chosen $h_n$, the second term is asymptotically negligible.
Proposition \ref{PROP:PENALTY_NLP} requires Conditions (E1)-(E5), stated below, which involve $h_n$ and also parameters measuring the strength of the signal-to-noise.
Recall that we assume Condition (D1) and that $X_k$ has zero column means and unit variances,
so that $V_k= n^{-1} I$, $\rho_{kj}$ are the eigenvalues of $X_k'X_k/n$ and $\tilde{\rho}_{kj}=\rho_{kj}+1/(\tau n)$ those of $\tilde{V}_k^{-1}/n$.

Proposition \ref{PROP:MOMPRIOR_SIMPLE} then follows immediately by setting $h_n= (\tau/n) (\tilde{\rho}_{m p_m}/\tilde{\rho}_{m1})^{1/2}$ in Proposition \ref{PROP:MOMPRIOR}, using that under Condition (E5) the term $P_{f^*}(D_{mt} > 1/h_n^{p_m-p_t})$ is asymptotically negligible, and noting that since $\tau n \gg 1$, 
$$
\frac{(\tilde{\rho}_{m p_m} )^{1/2}}{(\tilde{\rho}_{m 1})^{1/2}}=
\frac{(\rho_{m p_m} + 1/(\tau n))^{1/2}}{(\rho_{m 1} + 1/(\tau n))^{1/2}}
\asymp \frac{\rho_{m p_m}^{1/2}}{\rho_{m 1}^{1/2}}.
$$

We first state Proposition \ref{PROP:MOMPRIOR}, subsequently discuss Conditions (E1)-(E5), and finally state Proposition \ref{PROP:PENALTY_NLP}.

\begin{prop}
Assume that $f^*(\by)= N(\by; X_t \btheta_t^*; \phi^* I)$.
Consider $m \neq t$ and that Conditions (B1), (C1) and (D1) hold.
Let $m \in S$ and $h_n \geq 0$ be an arbitrary sequence. Then
$$E_{f^*}(p(M_m \mid \by)) \ll \max \left\{ [p(M_m)/p(M_t)]^\alpha (h_n \sqrt{\tau n})^{\alpha (p_t-p_m)}, e^{-\kappa n}, e^{-\tau n\rho_{tp_t}/2} \right\}
+ P_{f^*}(D_{mt}>1/h_n^{p_m-p_t})$$
for any fixed $\kappa >0$ and $\alpha \in (0,1)$.
\label{PROP:MOMPRIOR}
\end{prop}


\begin{enumerate}[leftmargin=*,label=(E\arabic*)]
\item As $n \rightarrow \infty$, $\tau n\tilde{\rho}_{mp_m} \gg 1$ and $\tau n\tilde{\rho}_{t1} \gg \sqrt{(\btheta_t^*)'\btheta_t^*}/\min_{j\in M_t} |\theta_{tj}^*|$.

\item As $n \rightarrow \infty$, $h_n^{(1-\delta)(p_m-p_t)/p_t} \gg \phi^*/\min_{j \in M_t} (\theta_{tj}^*)^2$ for fixed $\delta \in (0,1)$.

\item As $n \rightarrow \infty$,
  $h_n^{\epsilon (p_m-p_t)/p_t} \gg  [\max_{j \in M_t} (\theta_{tj}^*)^2 /\phi^*] [\tilde{\rho}_{m1}\tilde{\rho}_{t1}/(\tilde{\rho}_{mp_m}\tilde{\rho}_{tp_t})]^{1/2}$
  for fixed $\epsilon \in (0, 1-\delta)$.

\item As $n \rightarrow \infty$, $n \tilde{\rho}_{mp_m} \gg [\tilde{\rho}_{m1}\tilde{\rho}_{t1}/(\tilde{\rho}_{mp_m}\tilde{\rho}_{tp_t})]^{1/2}$

\item As $n \rightarrow \infty$
\begin{align}
 p_t e^{- \frac{1}{2}\min\{ \gamma \underline{\lambda}_m, (n-p_m)^\gamma  \}} \ll [p(M_m)/p(M_t)]^\alpha (h_n \sqrt{\tau n})^{\alpha (p_t-p_m)}
\nonumber \\
(p_m-p_t) e^{-\frac{1}{2} \min\{h_n^{1-\delta-\epsilon}, (n-p_m)^\gamma\} } \ll [p(M_m)/p(M_t)]^\alpha (h_n \sqrt{\tau n})^{\alpha (p_t-p_m)}.
\end{align}
\end{enumerate}

Conditions (E1)-(E4) are mild and essentially require that $\tau n$ grows quickly enough relative to $\tilde{\rho}_{mp_m}$.
For instance if $\tau n\geq n$, $\btheta_t^*$ is constant in $n$ and $\tilde{\rho}_{m1}/\tilde{\rho}_{mp_m}<c$ for some constant $c$,
then (E1) and (E4) are satisfied when $n \tilde{\rho}_{m,p_m} \gg 1$, (E2) and (E3) hold for any $h_n \gg 1$ and $(\delta,\epsilon)$ can be taken arbitrarily close to 0.
Note that $n \tilde{\rho}_{m,p_m} \gg 1$ is a minimal requirement that posterior variances in $V_m$ converge to 0, which is necessary for the least-squares estimator to be consistent.

Condition (E5) is also minimal, and simply states that the upper-bound on $P_{f^*}(D_{mt} > 1/h_n^{p_m-p_t})$ is asymptotically negligible relative to the first term in Proposition \ref{PROP:MOMPRIOR}. That is, (E5) allows to state a simpler rate in Proposition \ref{PROP:MOMPRIOR_SIMPLE} by simply reporting the first term in Proposition \ref{PROP:MOMPRIOR} as the asymptotically larger of the two.

Finally, we state Proposition \ref{PROP:PENALTY_NLP}.

\begin{prop}
  Let $m \in S$. Let $h_n$, $\tau n$, $\delta$ and $\epsilon$ satisfying (E1)-(E4) and
  $$
  h_n \leq  \frac{1}{2 \sigma} \tau (\tilde{\rho}_{mp_m}/\tilde{\rho}_{m1})^{1/2}
  $$
  where $\sigma$ is the largest diagonal element in $\tilde{V}_m$. Then
\begin{align}
  P_{f^*}(D_{mt}>1/h_n^{p_m-p_t}) \ll
  p_t e^{- \frac{1}{2}\min\{ \gamma \underline{\lambda}_m, (n-p_m)^\gamma  \}}
+ (p_m-p_t) e^{-\frac{1}{2} \min\{h_n^{1-\delta-\epsilon}, (n-p_m)^\gamma\} }
  \nonumber
\end{align}
for any $\gamma <1$,
where $\underline{\lambda}_m=\min_{M_k \subset M_m,p_k=p_m-1} \lambda_{mk}$.

\label{PROP:PENALTY_NLP}
\end{prop}


\section{Proof of Proposition \ref{PROP:MOMPRIOR}}

\subsection{Proof outline}
\label{proof:outline_momprior}

To apply Lemma \ref{LEM:TAILPROB_TO_MEAN} we seek to bound
\begin{align}
P_{f^*}\left(B_{mt} > \frac{(p(M_t)/p(M_m))}{1/u-1} \right)
\leq P_{f^*} \left(
D_{mt} \left( 1 + \frac{p_m-p_t}{n-p_m} \tilde{F}_{mt} \right)^{\frac{a_\phi+n}{2}}
> \frac{p(M_t) (c_{mt} (\tau n))^{\frac{p_m-p_t}{2}} }{p(M_m) (1/u-1)} \right).
\label{eq:tailprob_pmom}
\end{align}

Let $m \in S$,
$h_n \gg 1$ and $b_n(u)=2 (p(M_t)/p(M_m)) (\sqrt{c_{mt}(\tau n)}h_n)^{p_m-p_t}/(1/u-1)$.
From \eqref{eq:tailprob_zellner_unknown},
we have that \eqref{eq:tailprob_pmom} is
\begin{align}
< P_{f^*} \left( (p_m-p_t) \tilde{F}_{mt} > \frac{n-p_m}{n+a_\phi} b_n(u) \right)
+ P_{f^*} \left( D_{mt} > \frac{1}{h_n^{p_m-p_t}} \right).
\label{eq:tailprob_mom_spur}
\end{align}

The first term in \eqref{eq:tailprob_mom_spur} can be bounded as in \eqref{eq:tailprob_normalprior}
since $\lim_{n\rightarrow \infty} (n-p_m)/(n+a_\phi)=1$ under (B1), obtaining
$$
\int_0^1 P_{f^*} \left( (p_m-p_t) \tilde{F}_{mt} > b_n(u) \right) du \ll
\max \left\{ [p(M_m)/p(M_t)]^\alpha (h_n \sqrt{(\tau n)})^{\alpha (p_t-p_m)}, e^{-\kappa n}, e^{-(\tau n)\rho_{tp_t}/2} \right\},
$$
for any fixed $\kappa >0$, as we wished to prove.

\subsection{Derivations for the pMOM prior}
\label{app:derivation_pmom}

The Normal densities in $D_{mt}$ can be lower- and upper-bounded using
the eigenvalues of $\tilde{V}_m$ and $\tilde{V_t}$, giving $D_{mt}\leq$
\begin{align}
&\frac{\tilde{\rho}_{t1}^{p_t/2} \tilde{\rho}_{m1}^{p_m/2}}{\tilde{\rho}_{tp_t}^{p_t/2} \tilde{\rho}_{mp_m}^{p_m/2} \tau^{p_m-p_t}}
\frac{
\int \prod_{j \in M_m}^{} \left( \tilde{\theta}_{mj}^2/\phi + \frac{1}{n\tilde{\rho}_{m,p_m}} \right)
\mbox{IG} \left( \phi; \frac{a_\phi+n}{2}, \frac{\tilde{s}_m}{2} \right) d\phi
}{
\int \prod_{j \in M_t}^{} \left( \tilde{\theta}_{tj}^2/\phi + \frac{1}{n\tilde{\rho}_{t,p_t}} \right)
\mbox{IG} \left( \phi; \frac{a_\phi+n}{2}, \frac{\tilde{s}_t}{2} \right) d\phi
} \leq
\nonumber \\
&
\frac{
\int  \left[ \max_{j \in M_t}^{} \frac{\tilde{\theta}_{mj}^2}{\phi} + \frac{1}{n\tilde{\rho}_{m,p_m}} \right]^{p_t}
\left[ \max_{j \in M_m \setminus M_t}^{} \frac{\tilde{\theta}_{mj}^2}{\phi} + \frac{1}{n\tilde{\rho}_{m,p_m}} \right]^{p_m-p_t}
\mbox{IG} \left( \phi; \frac{a_\phi+n}{2}, \frac{\tilde{s}_m}{2} \right) d\phi
}{
  \tilde{\rho}_{t1}^{-p_t/2} \tilde{\rho}_{m1}^{-p_m/2}
  (\sqrt{2\pi}/e) \tilde{\rho}_{tp_t}^{p_t/2} \tilde{\rho}_{mp_m}^{p_m/2} \tau^{p_m-p_t}
\left[ \min_{j \in M_t} \tilde{\theta}_{tj}^2 (a_\phi +n+2p_t-2)/\tilde{s}_t  \right]^{p_t}
}.
\nonumber
\end{align}
The right hand side follows from Stirling's bound and trivial algebra.
Intuitively the minimum and maximum over $j \in M_t$ both converge to non-zero constants,
whereas that over $j \in M_m\setminus M_t$ vanishes at a rate given by $\tilde{V}_j$.

\section{Proof of Proposition \ref{PROP:PENALTY_NLP}}
\label{proof:penalty_nlp}

The goal is to bound $P(D_{mt}>1/h_n^{p_m-p_t})$, where $D_{mt} \leq$
\begin{align}
\frac{
\tilde{\rho}_{t1}^{p_t/2} \tilde{\rho}_{m1}^{p_m/2}
\int  \left[ \max_{j \in M_t}^{} \frac{\tilde{\theta}_{mj}^2}{\phi} + \frac{1}{n\tilde{\rho}_{m,p_m}} \right]^{p_t}
\left[ \max_{j \in M_m \setminus M_t}^{} \frac{\tilde{\theta}_{mj}^2}{\phi} + \frac{1}{n\tilde{\rho}_{m,p_m}} \right]^{p_m-p_t}
\mbox{IG} \left( \phi; \frac{a_\phi+n}{2}, \frac{\tilde{s}_m}{2} \right) d\phi
}{
  (\sqrt{2\pi}/e) \tilde{\rho}_{tp_t}^{p_t/2} \tilde{\rho}_{mp_m}^{p_m/2} \tau^{p_m-p_t}
\left[ \min_{j \in M_t} \tilde{\theta}_{tj}^2 (a_\phi +n+2p_t-2)/\tilde{s}_t  \right]^{p_t}
  }
  \nonumber
\end{align}

The strategy is to tackle the numerator and denominator separately
using that for any $W_1,W_2$ and any $a>0$, $\delta>0$ it holds that
$P(W_1/W_2>1/a) \leq P(W_1>1/a^{1-\delta}) + P(W_2<a^{\delta})$. Hence
\begin{align}
&P(D_{mt}>1/h_n^{p_m-p_t})  \leq P \left(
  \frac{\sqrt{2\pi}}{4e}
  \left[ \min_{j \in M_t} \frac{\tilde{\theta}_{tj}^2}{\tilde{s}_t} (n-p_t)  \right]^{p_t}
  \leq h_n^{(1-\delta)(p_m-p_t)}
  \right) +
\label{eq:proof_mompen_split} \\
&  P \left(
\frac{
\int  \left[ \max_{j \in M_t}^{} \frac{\tilde{\theta}_{mj}^2}{\phi} + \frac{1}{n\tilde{\rho}_{m,p_m}} \right]^{p_t}
\left[ \max_{j \in M_m \setminus M_t}^{} \frac{\tilde{\theta}_{mj}^2}{\phi} + \frac{1}{n\tilde{\rho}_{m,p_m}} \right]^{p_m-p_t}
\mbox{IG} \left( \phi; \frac{a_\phi+n}{2}, \frac{\tilde{s}_m}{2} \right) d\phi
}{  (1/4) (\tilde{\rho}_{tp_t}/\tilde{\rho}_{t1})^{p_t/2} (\tilde{\rho}_{mp_m}/\tilde{\rho}_{m1})^{p_m/2} \tau^{p_m-p_t} } \right.
                                   \nonumber \\
  &\left. > \frac{1}{h_n^{\delta (p_m-p_t)}}  \right),
  \nonumber
\end{align}
since $a_\phi+n+2p_t-2>n-p_t$.
When $p_t=0$ the first term in \eqref{eq:proof_mompen_split} is trivially 0
and likewise for the second term the upcoming arguments can be trivially modified,
hence we focus on the $p_t\geq 1$ case.
The first term requires the left tail of $\tilde{\theta}_{tj}^2 (n-p_t)/\tilde{s}_m$,
which by Lemma \ref{LEM:NCCHISQ_BAYES_UNIV}
is the ratio of a non-central $\chi^2_1$ and a $\chi^2_{n-p_t}/(n-p_t)$ random variables.
In Lemma \ref{LEM:BAYES_UNIVFTEST_INEQ} we show that when Condition (E1) holds such left tail
is $\ll e^{-\gamma \underline{\lambda}_t/2} + e^{-(n-p_t)^\gamma/2}$,
where $\gamma \in (0,1)$ is fixed
and $\underline{\lambda}_t=\min_{M_k \subset M_t,p_k=p_t-1} \lambda_{tk}$
is the smallest non-centrality parameter due to removing a single variable from $M_t$.
Hence the first term in \eqref{eq:proof_mompen_split} is
\begin{align}
\leq \sum_{j \in M_t}^{} P \left(
  \frac{\tilde{\theta}_{tj}^2}{\tilde{s}_t} (n-p_t)
  \leq \frac{4e}{\sqrt{2\pi}}  h_n^{(1-\delta)\frac{(p_m-p_t)}{p_t}}
  \right)
\ll
p_t (e^{-\gamma \underline{\lambda}_t/2} + e^{-(n-p_t)^\gamma/2}),
\label{eq:proof_mompen_ratedenom}
\end{align}
where the right-hand side follows from Corollary \ref{COR:SUMO}
and noting that the bound $e^{-\gamma \underline{\lambda}_t/2} + e^{-(n-p_t)^\gamma/2}$
holds uniformly across $j \in M_t$.

Consider the second term in \eqref{eq:proof_mompen_split}.
The strategy is to split the integral into four terms, then tackle term separately
using inequalities for ratios of chi-square variables.
Specifically
let $A= \left\{\phi: \max_{j \in M_t} \tilde{\theta}_{mj}^2/\phi > 1/(n\tilde{\rho}_{mp_m})  \right\}$,
$B= \left\{\phi: \max_{j \in M_m \setminus M_t} \tilde{\theta}_{mj}^2/\phi > 1/(n\tilde{\rho}_{mp_m})  \right\}$,
and $A^c$ and $B^c$ be their respective complementary sets.
Then splitting the integral into $A\cap B$, $A^c \cap B$, $A \cap B^c,$ and $A^c \cap B^c$
gives that the second term in \eqref{eq:proof_mompen_split} is
\begin{align}
  \begin{matrix}
< P \Big(
  \left[ \max_{j \in M_t} \frac{\tilde{\theta}_{mj}^2(n-p_m)}{\tilde{s}_m}  \right]^{p_t}
  \left[ \max_{j \in M_m \setminus M_t} \frac{\tilde{\theta}_{mj}^2(n-p_m)}{\tilde{s}_m}  \right]^{p_m-p_t}
  + \left[ \max_{j \in M_t} \frac{\tilde{\theta}_{mj}^2 (n-p_m)}{\tilde{s}_m}  \right]^{p_t}
  \left( \frac{1/n}{\tilde{\rho}_{mp_m}} \right)^{p_m-p_t}
\\
 + \left[ \max_{j \in M_m \setminus M_t} \frac{\tilde{\theta}_{mj}^2 (n-p_m)}{\tilde{s}_m}  \right]^{p_m-p_t}
  \left( \frac{1/n}{\tilde{\rho}_{mp_m}} \right)^{p_t}
  +  \left( \frac{1/n}{\tilde{\rho}_{mp_m}} \right)^{p_m}
> \frac{1/4}{2^{p_m}} \left[\frac{\tilde{\rho}_{tp_t}}{\tilde{\rho}_{t1}}\right]^{p_t/2} \left[\frac{\tilde{\rho}_{mp_m}}{\tilde{\rho}_{m1}}\right]^{p_m/2}
\left[\frac{\tau}{h_n^\delta}\right]^{p_m-p_t}
  \Big)
\end{matrix}
\label{eq:proof_mompen_ratenum}
\end{align}
where we used that, for any $l\geq 1$,
$$
\int \frac{1}{\phi^l} \mbox{IG} \left( \frac{a_\phi +n}{2}, \frac{\tilde{s}_m}{2} \right) d\phi
= \frac{\Gamma \left( \frac{a_\phi +n}{2} +p_t \right)}{\Gamma \left( \frac{a_\phi +n}{2} \right) (\tilde{s}_m/2)^l}
\leq \left( \frac{n-p_m}{\tilde{s}_m} \right)^l.
$$

To bound \eqref{eq:proof_mompen_ratenum} note that
$P(W_1+W_2+W_3+W_4> a/4) \leq \sum_{i=1}^{4} P(W_j>a)$ where $W_1,W_2,W_3,W_4$ are arbitrary random variables,
hence we may split the target probability in \eqref{eq:proof_mompen_ratenum} into four terms.
Let $\epsilon \in (0,1-\delta)$ be a fixed constant, the first term is
\begin{align}
  &P \left(
  \left[ \max_{j \in M_t} \frac{\tilde{\theta}_{mj}^2(n-p_m)}{\tilde{s}_m}  \right]^{p_t}
    \left[ \max_{j \in M_m \setminus M_t} \frac{\tilde{\theta}_{mj}^2(n-p_m)}{\tilde{s}_m}  \right]^{p_m-p_t}
    \right.
    \nonumber \\
&\left.  > \frac{1}{2^{p_m}}\left[\frac{\tilde{\rho}_{tp_t}}{\tilde{\rho}_{t1}}\right]^{p_t/2} \left[\frac{\tilde{\rho}_{mp_m}}{\tilde{\rho}_{m1}}\right]^{p_m/2}
\left[\frac{\tau}{h_n^\delta}\right]^{p_m-p_t}
  \right)
  \nonumber \\
&\leq  P \left(
  \max_{j \in M_t} \frac{\tilde{\theta}_{mj}^2}{\tilde{s}_m/(n-p_m)}
  > \frac{1}{2} \left[ \frac{\tilde{\rho}_{mp_m} \tilde{\rho}_{tp_t}}{\tilde{\rho}_{m1} \tilde{\rho}_{t1}} \right]^{\frac{1}{2}}
  h_n^{\frac{\epsilon (p_m-p_t)}{p_t}}
  \right)
  \nonumber \\
  &+ P \left(
  \max_{j \in M_m \setminus M_t} \frac{\tilde{\theta}_{mj}^2}{\tilde{s}_m/(n-p_m)} >
  \frac{(\tau) \tilde{\rho}_{mp_m}^{1/2}}{2 h_n^{\delta + \epsilon} \tilde{\rho}_{m_1}^{1/2} }
  \right)
    \nonumber \\
&\leq \sum_{j \in M_t} P \left(
   \frac{\tilde{\theta}_{mj}^2}{\tilde{s}_m/(n-p_m)}
  > \frac{1}{2} \left[ \frac{\tilde{\rho}_{mp_m} \tilde{\rho}_{tp_t}}{\tilde{\rho}_{m1} \tilde{\rho}_{t1}} \right]^{\frac{1}{2}}
  h_n^{\frac{\epsilon (p_m-p_t)}{p_t}}
                   \right)
\nonumber \\
 & + \sum_{j \in M_m \setminus M_t} P \left(
   \frac{\tilde{\theta}_{mj}^2}{\tilde{s}_m/(n-p_m)} >
  \frac{\tau \tilde{\rho}_{mp_m}^{1/2}}{2 h_n^{\delta + \epsilon} \tilde{\rho}_{m_1}^{1/2} }
  \right)
  \nonumber \\
&  \ll  p_t \left(e^{-\gamma \underline{\lambda}_m/2} + e^{-(n-p_m)^\gamma/2} \right)
  + (p_m - p_t) \left(e^{-(\gamma/2) \frac{\tau \tilde{\rho}_{mp_m}^{1/2}}{2 \sigma_{jj} h_n^{\delta+\epsilon} \tilde{\rho}_{m1}^{1/2}} } + e^{-(n-p_m)^\gamma/2} \right)
\label{eq:proof_mompen_term1}
\end{align}
where from the definition $h_n \leq  1/(2\sigma) \tau (\tilde{\rho}_{mp_m}/\tilde{\rho}_{m1})^{1/2}$ we have
$$
  \frac{\tau \tilde{\rho}_{mp_m}^{1/2}}{2\sigma_{jj} h_n^{\delta+\epsilon} \tilde{\rho}_{m1}^{1/2}}
\geq h_n^{1-\delta -\epsilon}.
$$

The right hand side in \eqref{eq:proof_mompen_term1} follows from Lemma \ref{LEM:BAYES_UNIVFTEST_INEQ},
Assumption (E3) and Corollary \ref{COR:SUMO} (since the bounds hold uniformly across $j$).
Note that for $j \in M_t$ Lemma \ref{LEM:BAYES_UNIVFTEST_INEQ}(ii) requires that 
$ h_n^{\frac{\epsilon (p_m-p_t)}{p_t}}
(\tilde{\rho}_{mp_m} \tilde{\rho}_{tp_t}/\tilde{\rho}_{m1} \tilde{\rho}_{t1})^{\frac{1}{2}} \gg (\theta_{mj}^*)^2/\phi^*$,
which holds by Assumption (E3).
For $j \in M_m\setminus M_t$ Lemma \ref{LEM:BAYES_UNIVFTEST_INEQ}(iii)
requires that the critical point $\tau \tilde{\rho}_{mp_m}^{1/2}/(h_n^{\delta + \epsilon} \tilde{\rho}_{m_1}^{1/2}) \gg 1$.
Since $\delta + \epsilon <1$ it suffices that
$h_n \leq \tau \tilde{\rho}_{mp_m}^{1/2}/(\tilde{\rho}_{m_1}^{1/2} \sigma_{jj})$,
which is satisfied by definition.

Regarding the second term in \eqref{eq:proof_mompen_ratenum},
\begin{align}
P \left(  \left[ \max_{j \in M_t} \frac{\tilde{\theta}_{mj}^2}{\tilde{s}_m/(n-p_m)}  \right]^{p_t}
  \left( \frac{1}{n\tilde{\rho}_{mp_m}} \right)^{p_m-p_t}
> \frac{1}{2^{p_m}} \left[\frac{\tilde{\rho}_{tp_t}}{\tilde{\rho}_{t1}}\right]^{p_t/2} \left[\frac{\tilde{\rho}_{mp_m}}{\tilde{\rho}_{m1}}\right]^{p_m/2}
\left[\frac{\tau}{h_n^\delta}\right]^{p_m-p_t}
  \right)
\leq
\nonumber \\
P \left( \max_{j \in M_t} \frac{\tilde{\theta}_{mj}^2}{\tilde{s}_m/(n-p_m)} >
\frac{1}{2} \left[ \frac{\tilde{\rho}_{mp_m} \tilde{\rho}_{tp_t}}{\tilde{\rho}_{m1} \tilde{\rho}_{t1}} \right]^{\frac{1}{2}}
  h_n^{\frac{\epsilon (p_m-p_t)}{p_t}}
\right)
+ P \left(  \frac{1}{n\tilde{\rho}_{mp_m}} > \frac{\tau(\tilde{\rho}_{mp_m})^{1/2}}{2 h_n^{\delta+\epsilon} \tilde{\rho}_{m1}^{1/2}}  \right).
\label{eq:proof_mompen_term2}
\end{align}
The first summand in \eqref{eq:proof_mompen_term2} is identical to that in \eqref{eq:proof_mompen_term1},
hence $\ll p_t \left(e^{-\gamma \underline{\lambda}_m/2} + e^{-(n-p_m)^\gamma/2} \right)$.
The second summand is zero since
$h_n \leq n\tilde{\rho}_{mp_m} \frac{\tau(\tilde{\rho}_{mp_m})^{1/2}}{2 \tilde{\rho}_{m1}^{1/2}}$
and $\delta + \epsilon<1$.
To see this note that
$\sigma_{jj}$ is the $j^{th}$ diagonal element in $V_k$,
hence $1/\sigma_{jj} \leq n\tilde{\rho}_{mp_m}$ and by definition
$$h_n \leq  [\tau/(2\sigma)]  (\tilde{\rho}_{mp_m}/\tilde{\rho}_{m1})^{1/2} \leq (\tau n) \tilde{\rho}_{mp_m}^{3/2} /\tilde{\rho}_{m1}^{1/2}.$$

Now, the third term arising from \eqref{eq:proof_mompen_ratenum} is
\begin{align}
P \left( \left[ \max_{j \in M_m \setminus M_t} \frac{\tilde{\theta}_{mj}^2}{\tilde{s}_m/(n-p_m)}  \right]^{p_m-p_t}
  \left( \frac{1}{n\tilde{\rho}_{mp_m}} \right)^{p_t}
> \frac{1}{2^{p_m}} \left[\frac{\tilde{\rho}_{tp_t}}{\tilde{\rho}_{t1}}\right]^{p_t/2} \left[\frac{\tilde{\rho}_{mp_m}}{\tilde{\rho}_{m1}}\right]^{p_m/2}
\left[\frac{\tau}{h_n^\delta}\right]^{p_m-p_t}
\right)
\nonumber \\
=P \left( \max_{j \in M_m \setminus M_t} \frac{\tilde{\theta}_{mj}^2}{\tilde{s}_m/(n-p_m)}
> \frac{(\tau) \tilde{\rho}_{mp_m}^{1/2}}{2 h_n^\delta \tilde{\rho}_{m1}^{1/2}}
\times \left[ \frac{n \tilde{\rho}_{mp_m}^{3/2} \tilde{\rho}_{tp_t}^{1/2}}{\tilde{\rho}_{m1}^{1/2}\tilde{\rho}_{t1}^{1/2}} \right]^{p_t/(p_m-p_t)}
\right)
\nonumber \\
\leq \sum_{j \in M_m \setminus M_t} P \left(  \frac{\tilde{\theta}_{mj}^2}{\tilde{s}_m/(n-p_m)}
> \frac{(\tau) \tilde{\rho}_{mp_m}^{1/2}}{2 h_n^\delta \tilde{\rho}_{m1}^{1/2}}
\times \left[ \frac{n \tilde{\rho}_{mp_m}^{3/2} \tilde{\rho}_{tp_t}^{1/2}}{\tilde{\rho}_{m1}^{1/2}\tilde{\rho}_{t1}^{1/2}} \right]^{p_t/(p_m-p_t)}
\right)
\nonumber \\
\ll
(p_m-p_t) (e^{-\gamma g_n/2} + e^{-(n-p_m)^\gamma/2}),
\label{eq:proof_mompen_term3}
\end{align}
where
$$
g_n= \frac{1}{\sigma_{jj}} \frac{\tau \tilde{\rho}_{mp_m}^{1/2}}{2 h_n^\delta \tilde{\rho}_{m1}^{1/2}}
\times \left[ \frac{n \tilde{\rho}_{mp_m}^{3/2} \tilde{\rho}_{tp_t}^{1/2}}{\tilde{\rho}_{m1}^{1/2}\tilde{\rho}_{t1}^{1/2}} \right]^{p_t/(p_m-p_t)}
.
$$
The right hand side in \eqref{eq:proof_mompen_term3} follows from Lemma \ref{LEM:BAYES_UNIVFTEST_INEQ}(iii)
and noting that
$$
\frac{\tau \tilde{\rho}_{mp_m}^{1/2}}{2 h_n^\delta \tilde{\rho}_{m1}^{1/2}}
\times \left[ \frac{n \tilde{\rho}_{mp_m}^{3/2} \tilde{\rho}_{tp_t}^{1/2}}{\tilde{\rho}_{m1}^{1/2}\tilde{\rho}_{t1}^{1/2}} \right]^{p_t/(p_m-p_t)}
\gg \sigma_{jj},
$$
since by definition $h_n \leq [\tau/(2\sigma)] \tilde{\rho}_{mp_m}^{1/2}/\tilde{\rho}_{m1}^{1/2}$
and by Assumption (E4) $n \tilde{\rho}_{mp_m}^{3/2} \tilde{\rho}_{tp_t}^{1/2} / (\tilde{\rho}_{m1}^{1/2}\tilde{\rho}_{t1}^{1/2}) \gg 1$.
This also implies that $g_n \succeq h_n^{1-\delta}$.

Finally, the fourth term arising from \eqref{eq:proof_mompen_ratenum} is
\begin{align}
&P \left(
  \left( \frac{1}{n\tilde{\rho}_{mp_m}} \right)^{p_m}
> \frac{1}{2^{p_m}} \left[\frac{\tilde{\rho}_{tp_t}}{\tilde{\rho}_{t1}}\right]^{p_t/2} \left[\frac{\tilde{\rho}_{mp_m}}{\tilde{\rho}_{m1}}\right]^{p_m/2}
\left[\frac{\tau}{h_n^\delta}\right]^{p_m-p_t}
\right)
\nonumber \\
&=P \left( h_n^\delta > \frac{\tau n \tilde{\rho}_{mp_m}^{3/2} }{2 \tilde{\rho}_{m1}^{1/2}}
\left[\frac{n\tilde{\rho}_{mp_m}^{3/2}\tilde{\rho}_{tp_t}^{1/2}}{2\tilde{\rho}_{t1}^{1/2}\tilde{\rho}_{m1}^{1/2}}\right]^{p_t/(p_m-p_t)}
 \right)=0,
\label{eq:proof_mompen_term4}
\end{align}
since $1/\sigma_{ii} \leq n \tilde{\rho}_{mp_m}$ and hence
$h_n \leq [\tau/(2\sigma)] \tilde{\rho}_{mp_m}^{1/2}/\tilde{\rho}_{m1}^{1/2} \leq [\tau n/2] \tilde{\rho}_{mp_m}^{3/2}/\tilde{\rho}_{m1}^{1/2}$,
and again recalling that by (E4) $n \tilde{\rho}_{mp_m}^{3/2} \tilde{\rho}_{tp_t}^{1/2} / (\tilde{\rho}_{m1}^{1/2}\tilde{\rho}_{t1}^{1/2}) \gg 1$.

Combining \eqref{eq:proof_mompen_ratedenom}, \eqref{eq:proof_mompen_term1}, \eqref{eq:proof_mompen_term2}, \eqref{eq:proof_mompen_term3} and \eqref{eq:proof_mompen_term4}
 we obtain that $P(D_{mt}>1/h_n^{p_m-p_t}) \ll $
\begin{align}
p_t (e^{-\gamma \underline{\lambda}_t/2} + e^{-(n-p_t)^\gamma/2})
+ p_t \left(e^{-\gamma \underline{\lambda}_m/2} + e^{-(n-p_m)^\gamma/2} \right)
  + (p_m - p_t) \left(e^{-(\gamma/2) h_n^{1-\delta-\epsilon}} + e^{-(n-p_m)^\gamma/2} \right)
\nonumber \\
  + p_t \left(e^{-\gamma \underline{\lambda}_m/2} + e^{-(n-p_m)^\gamma/2} \right)
  + (p_m-p_t) (e^{-\gamma h_n^{1-\delta}/2} + e^{-(n-p_m)^\gamma/2})
\nonumber \\
\ll
 p_t \left(e^{-\gamma \underline{\lambda}_m/2} + e^{-(n-p_m)^\gamma/2} \right)
  + (p_m - p_t) \left(e^{-(\gamma/2) h_n^{1-\delta-\epsilon}} + e^{-(n-p_m)^\gamma/2} \right),
  \nonumber
\end{align}
since $\underline{\lambda}_t \geq \underline{\lambda}_m$ for any $M_m \subset M_t$,
which proves the desired result.

\section{Proof of Proposition \ref{PROP:BF_NONPARAM}}

To proof runs analogously to that for Zellner's prior in the well-specified case shown in Section \ref{proof:zellner_phiunknown}.
The main adjustment is finding a new bound for Bayes factor tail probabilities
under $f^*(\by)=N(\by; W\bbeta^*, \xi^* I)$.

For any $M_k$ denote the KL-optimal variance by $\phi_k^*=\xi^* + (W\bbeta^*)' (I - H_k) W\bbeta^*/n$.
Recall from \eqref{eq:tailprob_zellner_unknown} that
$$P_{f^*}\left( B_{mt}> \frac{p(M_t)/p(M_m)}{(1/u-1)} \right) < P_{f^*}((p_m-p_t)F_{mt}>b_n(u)),$$
where $F_{mt}$ is the F-test statistic defined in \eqref{eq:tailprob_zellner_unknown}
and $b_n(u)= 2 [(n-p_m)/(n+a_\phi)] \log ( (1+ n\tau)^{\frac{p_m-p_t}{2}} (p(M_t)/p(M_m))/(1/u-1) )$.
Further, note that we can re-write
\begin{align}
(p_m-p_t) F_{mt}=
\frac{W_{mt}/\xi^*}{s_m/[(n-p_m)\phi_m^*]} \frac{\xi^*}{\phi_m^*},
\label{eq:fstat_nonparam}
\end{align}
where $W_{mt}=\hat{\btheta}_m' X_m'X_m \hat{\btheta}_m - \hat{\btheta}_t' X_t'X_t \hat{\btheta}_t$,
and $\hat{\btheta}_m$, $\hat{\btheta}_t$ are the least-squares estimates under $M_m$ and $M_t$ respectively.
To bound tail probabilities for \eqref{eq:fstat_nonparam}
note that
\begin{align}
  P_{f^*} \left( (p_m-p_t) F_{mt} > b_n(u) \right)=
  P_{f^*} \left( \frac{U_1}{U_2/(n-p_m)} > \frac{\phi_m^*}{\xi^*} b_n(u) \right),
\nonumber
\end{align}
where $U_1=W_{mt}/\xi^*$ and $U_2=s_m/[(n-p_m)\phi_m^*]$.
In Lemma \ref{LEM:NESTEDLM_NONLINEAR} we prove that
$W_{ls}/\xi^* \sim \chi^2_{p_l-p_s}(\lambda_{sl})$ for any $l$ such that $M_s \subset M_l$,
where $\lambda_{ls}=(W\bbeta^*)'H_l (I-H_s)H_l W\bbeta^*/\xi^*$,
$H_l=X_l (X_l'X_l)^{-1} X_l'$.
Further, since $f^*(\by)$ has Gaussian tails it follows that
$U_2=s_m/[(n-p_m)\phi_m^*]$ has exponential tails with expectation equal to 1 under $f^*$.
This implies that the tail inequalities from Lemma \ref{LEM:INEQ_NCF_RIGHT},
Corollary \ref{COR:INEQ_F_RIGHT}, Lemma \ref{LEM:INEQ_NCF_LEFT} and Lemma \ref{LEM:INEQ_F} still apply up to a constant factor,
since the proofs just require having the ratio of two (possibly dependent) random variables:
$U_1$ being a non-central chi-square and $U_2$ a random variable with exponential tails.
We outline the proof separately for the cases $m \in S$ and $m \in S^c$.

\subsection{Misspecified mean structure. Spurious models}

Let $m \in S$ be a spurious model.
Then $\phi_m^*=\phi_t^*$ and $\lambda_{mt}=(W\bbeta^*)'H_m (I-H_t)H_m W\bbeta^*/\xi^*=0$ by definition.
The result is based on showing that the rate from Lemma \ref{LEM:FTAIL_INTBOUND} still applies.
To see this in \eqref{eq:proof_ftail_bound_nonspur2} define $\nu_1=p_m-p_t$, $\nu_2=n-p_m$,
$d=2 [(n-p_m)/(n+a_\phi)]$ and $g=(1+n \tau)^{\frac{p_m-p_t}{2}} (p(M_t)/p(M_m)) e^{\phi_t^*/\xi^*}$.
If $\log(g) \ll n-p_m$,
since Corollary \ref{COR:INEQ_F_RIGHT} applies up to a constant it follows that
the bound in \eqref{eq:proof_ftail_bound_nonspur1} holds up to a constant,
hence the argument leading to \eqref{eq:proof_ftail_bound_nonspur3} remains valid
and 
$$
E_{f^*}(p(M_m \mid \by)) \ll 1/g^\alpha=
\left[ \frac{p(M_m)}{(n \tau)^{\frac{p_m-p_t}{2}} p(M_t) e^{\phi_t^*/\xi^*}} \right]^\alpha
$$ 
for any $\alpha < d-1$, as we wished to prove.
Conversely if $\log^\gamma(g) \gg \nu_2$ for all fixed $\gamma<1$,
then the bound in Lemma \ref{LEM:INEQ_F} and \eqref{eq:proof_ftail_bound_nonspur4} hold up to a constant
and the argument leading to \eqref{eq:proof_ftail_bound_nonspur5} remains valid,
giving that
$E_{f^*}(p(M_m \mid \by)) \ll  e^{-\kappa n}$ for any fixed $\kappa >0$, as we wished to prove.

\subsection{Misspecified mean structure. Non-spurious models}

Let $m \not\in S$ be a non-spurious model.
We define $M_q=M_m \cup M_t$ and use the bound in \eqref{eq:ftail_bound_nonspur}. 
Specifically, following \eqref{eq:ftail_bound_nonspur}, $E_{f^*}(p(M_m \mid \by)) < \int_0^1 P_{f^*} \left( (p_m-p_t) F_{mt} > b_n(u) \right) du<$
\begin{align}
\int_0^1 P_{f^*} \left( (p_q-p_m) F_{qm} < \left( c_n(u) - \frac{b_n(u)}{2} \right) \frac{n-p_q}{n-p_m}   \right) du
  \nonumber \\
 + \int_0^1 P_{f^*} \left( (p_q-p_t) F_{qt} > c_n(u) +  \frac{b_n(u)}{2} \right) du,
\label{eq:ftail_bound_nonspur_nonpar}
\end{align}
for any $c_n(u)>0$. Specifically we set $c_n(u)= 0.5 b_n(u) + \lambda_{tm}/\log(\lambda_{tm})$.
Note that $\lambda_{qm}=\lambda_{tm}$,
to see this let $M_v= M_m\setminus M_t$ 
and define $\tilde{X}_v= (I-H_t)X_v$, where $H_t= X_t (X_t'X_t)^{-1} X_t$ is the projection matrix onto the column span of $X_t$.
Basic properties of orthogonal projections give $H_q=H_t + \tilde{X}_v(\tilde{X}_v'\tilde{X}_v)^{-1}\tilde{X}_v$, hence
$$\lambda_{qm}=\frac{(W\bbeta^*)'H_t (I-H_m)H_t W\bbeta^*}{\xi^*}=\lambda_{tm},$$
since $\tilde{X}_v(\tilde{X}_v'\tilde{X}_v)^{-1}\tilde{X}_v W\bbeta^*={\bf 0}$ by definition of $M_t$.

Since $\phi_q^*=\phi_t^*$, the first term in \eqref{eq:ftail_bound_nonspur_nonpar} is equal to
$$
\int_0^1 P_{f^*} \left( \frac{U_1}{U_2/(n-p_m)} < \frac{\phi_t^*}{\xi^*} \left[ c_n(u) - \frac{b_n(u)}{2} \right] \frac{n-p_q}{n-p_m}   \right) du.
$$

To bound this integral we argue as in Section \ref{proof:zellner_phiunknown}.
Briefly, we use that $\lambda_{tm} \gg [\lambda_{tm}/\log(\lambda_{tm})] \phi_t^*/\xi^*$ since $\log(\lambda_{tm}) \gg \phi_t^*/\xi^*$ by assumption,
and that $U_1 \sim \chi^2_{p_q-p_m}(\lambda_{tm})$ and $U_2=s_q/[(n-p_q)\phi_t^*]$ so that Lemma \ref{LEM:INEQ_NCF_LEFT} applies up to a constant.
Then the argument in Section \ref{proof:zellner_phiunknown} shows that the first term in \eqref{eq:ftail_bound_nonspur_nonpar} is
$\ll e^{-\gamma \frac{\lambda_{tm}}{2}} + e^{- \kappa n}$.

The second term in \eqref{eq:ftail_bound_nonspur_nonpar} is proven to be
\begin{align}
\ll
\max \left\{  e^{-\frac{\lambda_{tm}^\gamma}{2}}
\left[ \frac{p(M_m) }{p(M_t) (1+\tau n)^{(p_m-p_t)/2} } \right]^\gamma, e^{- \kappa n}
\right\},
\nonumber
\end{align}
for any fixed $\gamma <1$, $\kappa >0$, again following Section \ref{proof:zellner_phiunknown} and using that Lemma \ref{LEM:FTAIL_INTBOUND} applies up to a constant.

Combining the two terms in \eqref{eq:ftail_bound_nonspur_nonpar} gives
\begin{align}
E_{f^*}(p(M_m \mid \by)) \ll
 \max \left\{  
e^{- \gamma \lambda_{tm}/2},  e^{-\frac{\lambda_{tm}^\gamma}{2}}
\left[ \frac{p(M_m) }{p(M_t) (1+\tau n)^{(p_m-p_t)/2} } \right]^\gamma, e^{- \kappa n}
\right\},
\nonumber
\end{align}
for any fixed $\gamma <1$, $\kappa >0$.
Arguing as in \eqref{eq:proof_zellnerknown_largerterm} shows that, under Condition (C2), the first term is of a smaller order than the second term, hence
\begin{align}
E_{f^*}(p(M_m \mid \by)) \ll
 \max \left\{  
 e^{-\frac{\lambda_{tm}^\gamma}{2}} \left[ \frac{p(M_m) }{p(M_t) (1+\tau n)^{(p_m-p_t)/2} } \right]^\gamma, e^{- \kappa n}
\right\},
\nonumber
\end{align}
as we wished to prove.

To complete the proof, we show that $\lambda_{tm} \leq \lambda_m^*$, with equality if and only if $W \bbeta^*= X \btheta_t^*$.
To show that
$$
\lambda_{tm}= \frac{(W\bbeta^*)'H_t (I-H_m)H_t W\bbeta^*}{\xi^*}
\leq \frac{(W\bbeta^*)' (I-H_m) W\bbeta^*}{\xi^*}= \lambda_m^*
$$
it suffices to prove that for any two projection matrices $A,B$ it holds that ${\bf s}'ABA {\bf s} \leq {\bf s}' B {\bf s}$,
i.e. $B-ABA$ is positive semidefinite.
Let $l$ be an eigenvalue of $B-ABA$ and ${\bf v} \in \mathbb{R}^n$ its corresponding eigenvector, that is $(B-ABA){\bf v}= l{\bf v}$.
Since $(B-ABA){\bf v}= B{\bf v} - BA{\bf v}$, it follows that
\begin{align}
(B-ABA)(B-ABA){\bf v}= (B-ABA) (B{\bf v}- BA{\bf v})=
\nonumber \\
B{\bf v} - BA{\bf v} - AB{\bf v} + BA{\bf v}= B{\bf v} - BA{\bf v}= (B-ABA){\bf v}= l{\bf v}.
\nonumber
\end{align}
Hence $l$ is also an eigenvalue of $(B-ABA)(B-ABA)$, so $l \in \{0,1\}$ and $B-ABA$ is positive semidefinite.

\section{Proof of Proposition \ref{PROP:BF_HETEROSKEDASTICITY_KNOWNPHI}}

Using Lemma \ref{LEM:TAILPROB_TO_MEAN} and the expression for the Bayes factor leading to \eqref{eq:tailprob_zellner_known}, the goal is to bound
\begin{align}
E_{f^*}\left(p(M_m\mid\by)\right) < \int_0^1
P_{f^*} \left( \frac{W_{mt}}{\phi^*} >
\frac{1+\tau n}{\tau n} 2 \log \left[ \frac{(1+\tau n)^{\frac{p_m-p_t}{2}} p(M_t)}{(1/u-1) p(M_m)} \right]
 \right) du
\label{eq:proof_goal_hetero_knownphi}
\end{align}
where $W_{mt}=\hat{\btheta}_m'X_m'X_m\hat{\btheta}_m - \hat{\btheta}_t'X_t'X_t\hat{\btheta}_t$
is the difference between residual sums of squares
under $M_t$ and $M_m$ and $\hat{\btheta}_k= (X_k'X_k)^{-1} X_k' \by$ the least-squares estimate for $k \in \{t,m\}$.
For simplicity, since $\tau n \gg 1$ by Assumption (B1), the factor $(1+\tau n)/(\tau n)$ converges to 1 as $n$ grows and may be dropped from the expression.

The proof strategy is to bound $W_{mt}$ by chi-square random variables times certain eigenvalues, which then allows bounding the tail probability in \eqref{eq:proof_goal_hetero_knownphi} and its integral with respect to $u$.
We consider separately the cases where $M_m$ is a spurious and a non-spurious model.

\subsection{Spurious models}

Let $m \in S$ be a spurious model. 
To bound \eqref{eq:proof_goal_hetero_knownphi} we use Lemma \ref{LEM:NESTEDLM_CHITEST_HETERO}, Part (i).
Let $X_{m \setminus t}$ be the columns in $X_m$ that are not included in $X_t$ (i.e. the truly spurious columns in $X_m$).
Then Lemma \ref{LEM:NESTEDLM_CHITEST_HETERO} gives that
\begin{align}
\underline{\omega}_{tm} Z_1 \leq \frac{W_{mt}}{\phi^*} \leq \bar{\omega}_{tm} Z_1
\nonumber
\end{align}
where $Z_1 \sim \chi^2_{p_m-p_t}$, $\underline{\omega}_{tm}$ and $\bar{\omega}_{tm}$ are the smallest and largest eigenvalues of 
$X_{m \setminus t}(I-H_t) \Sigma^* (I-H_t) X_{m \setminus t} (X_{m \setminus t} (I-H_t) X_{m \setminus t})^{-1}$,
and $H_t= X_t(X_t'X_t)^{-1}X_t'$ is the projection matrix onto the column space of $X_t$.
Equivalently, since $(I-H_t)X_m=(I-H_t)(X_t,X_{m \setminus t})=(0,(I-H_t)X_{m \setminus t})$,
$\underline{\omega}_{tm}$ and $\bar{\omega}_{tm}$ are eigenvalues of $X_m(I-H_t) \Sigma^* (I-H_t) X_m (X_m (I-H_t) X_m)^{-1}$.
Hence,
\begin{align}
E_{f^*}\left(p(M_m\mid\by)\right) < \int_0^1
P_{f^*} \left( Z_1 >
\frac{2}{\bar{\omega}_{tm}} \log \left( \frac{(1+\tau n)^{\frac{p_m-p_t}{2}} p(M_t)}{(1/u-1) p(M_m)} \right)
 \right) du
\label{eq:proof_spur_heterosk_mainterm}
\end{align}

To bound this integral we consider separately the cases $\bar{\omega}_{tm} \leq 1$ 
and $\bar{\omega}_{tm} > 1$.
The case $\bar{\omega}_{tm} \leq 1$ is covered by Lemma \ref{LEM:CHISQTAIL_INTBOUND}, Part (i).
Specifically set $\nu=p_m-p_t$, $g= (1+\tau n)^{\frac{p_m-p_t}{2}} p(M_t)/p(M_m)$
and note that $\nu \ll \log(g)$ as $n\rightarrow \infty$ under Conditions (B1) and (C1).
Then by Lemma \ref{LEM:CHISQTAIL_INTBOUND}
$$E_{f^*}(p(M_m\mid\by)) \preceq \frac{[\log(g)]^{(p_m-p_t)/2}}{g}  \ll \frac{[p(M_m)/p(M_t)]^\alpha}{(1+\tau n)^{\alpha (p_m-p_t)/2} }$$
 for any fixed $\alpha \in (0,1)$.

The case $\bar{\omega}_{tm} \in (1,2)$ is covered by Lemma \ref{LEM:CHISQTAIL_INTBOUND}, Part (ii), but it is simpler to 
use Proposition S1 Part (ii) in \cite{rossell:2020}, which covers the wider case $\bar{\omega}_{tm} > 1$.
Setting $c=\phi^*/\bar{\omega}_{tm}$, $\sigma^2= \phi^*$, $d=1$ and $h= (1+\tau n)^{p_m-p_t} [p(M_t)/p(M_m)]^2$ in the statement of Proposition S1 of \cite{rossell:2020}, gives that $c/\sigma^2= 1/\bar{\omega}_{tm}$ and hence
\begin{align}
E_{f^*}(p(M_m\mid\by))
\leq
\frac{2.5  \max \left\{
\left[ \log(h^{1/(\bar{\omega}_{tm}(1+\epsilon))})  \right]^{1/2},
\log \left( h^{\frac{1}{2 \bar{\omega}_{tm} (1+\epsilon) }} \right)
  \right\}}
{h^{\frac{1}{2 \bar{\omega}_{tm} (1 + \epsilon)}}}
\nonumber
\end{align}
where $\epsilon= \sqrt{1/[\bar{\omega}_{tm} \log h]}$.
Since $\bar{\omega}_{tm} > 1$ and $h \gg 1$ by Condition (C1), 
we have that $\bar{\omega}_{tm} \log h \gg 1$ and hence $\lim_{n \rightarrow \infty} \epsilon=0$. Thus, we obtain that
\begin{align}
E_{f^*}(p(M_m\mid\by))
\ll
\left[ \frac{[p(M_m)/p(M_t)]}{(1+\tau n)^{(p_m-p_t)/2} } \right]^{\alpha / \bar{\omega}_{tm}}
\nonumber
\end{align}
for any fixed $\alpha \in (0,1)$.
Combining the cases $\bar{\omega}_{tm} \leq 1$ and $\bar{\omega}_{tm} > 1$ gives
\begin{align}
E_{f^*}(p(M_m\mid\by))
\ll
\left[ \frac{[p(M_m)/p(M_t)]}{(1+\tau n)^{(p_m-p_t)/2} } \right]^{\alpha \min\{1, 1 / \bar{\omega}_{tm} \}},
\nonumber
\end{align}
as we wished to prove.


\subsection{Non-spurious models}

Let $m \in S^c$ be a non-spurious model and $M_q=M_t \cup M_m$ the model with design matrix $X_q$ combining the columns in $X_t$ and $X_m$.
Let $\lambda_{tm}$ be as defined in \eqref{eq:ncp_regression}.

To ease notation in \eqref{eq:proof_goal_hetero_knownphi} denote $b_n(u)= 2 \log ( (1+\tau n)^{\frac{p_m-p_t}{2}} (p(M_t)/p(M_m))/(1/u-1) )$.
Since $W_{mt}=\hat{\btheta}_m'X_m'X_m\hat{\btheta}_m - \hat{\btheta}_q'X_q'X_q\hat{\btheta}_q
+ \hat{\btheta}_q'X_q'X_q\hat{\btheta}_q - \hat{\btheta}_t'X_t'X_t\hat{\btheta}_t =W_{qt} - W_{qm}$, we can bound the tail probability in \eqref{eq:proof_goal_hetero_knownphi} by
\begin{align}
P_{f^*} \left( \frac{W_{mt}}{\phi^*} > b_n(u)  \right)
= P_{f^*} \left( \frac{W_{qt} - W_{qm}}{\phi^*} > \frac{b_n(u)}{2} + c_n(u) - (c_n(u) - \frac{b_n(u)}{2} )  \right)
\nonumber \\
\leq
P_{f^*} \left( \frac{W_{qm}}{\phi^*} < c_n(u) - \frac{b_n(u)}{2}  \right)
+ P_{f^*} \left( \frac{W_{qt}}{\phi^*} > c_n(u) + \frac{b_n(u)}{2} \right)
\label{eq:tailprob_nonspur_zellnerhet_known}
\end{align}
for any $c_n(u)>0$, the right-hand side following from the union bound. 
The idea is to set $c_n(u)$ to a convenient expression that gives a small value for \eqref{eq:tailprob_nonspur_zellnerhet_known}, see below.

We shall bound each term in \eqref{eq:tailprob_nonspur_zellnerhet_known} separately. In both cases we first bound $W_{qm}$ and $W_{qt}$ using Lemma \ref{LEM:NESTEDLM_CHITEST_HETERO} Part (i).
Specifically, for any model $M_k$ nested in $M_q$, denote by $\underline{\omega}_{qk}$ and $\bar{\omega}_{qk}$ the smallest and largest eigenvalues of 
$X_q(I-H_k) \Sigma^* (I-H_k) X_q (X_q (I-H_k) X_q)^{-1}$, where $H_k$ is as usual the projection matrix onto the column space of $X_k$.
Then Lemma \ref{LEM:NESTEDLM_CHITEST_HETERO} Part (i) gives that
\begin{align}
\underline{\omega}_{tq} Z_1 \leq \frac{W_{qt}}{\phi^*} \leq \bar{\omega}_{tq} Z_1
\nonumber \\
\underline{\omega}_{mq} Z_2 \leq \frac{W_{qm}}{\phi^*} \leq \bar{\omega}_{mq} Z_2
\nonumber \\
\frac{\lambda_{qm}}{\bar{\omega}_{mq}} \leq \tilde{\lambda}_{qm} \leq \frac{\lambda_{qm}}{\underline{\omega}_{mq}},
\nonumber
\end{align}
where $Z_1 \sim \chi^2_{p_q-p_t}$, $Z_2 \sim \chi^2_{p_q - p_m}(\tilde{\lambda}_{qm})$,
$\lambda_{qm}= (\tilde{X}_{t \setminus m} \btheta_{t \setminus m}^*)' (I-H_m) \tilde{X}_{t \setminus m} \btheta_{t \setminus m}^*$,
$\tilde{X}_{t \setminus m}= (I - H_m) X_t$, 
$\tilde{\lambda}_{qm}= (\tilde{X}_{t \setminus m} \btheta_{t \setminus m}^*)' W^{-1} \tilde{X}_{t \setminus m} \btheta_{t \setminus m}^*$,
and $W= (\tilde{X}_{t \setminus m}'\tilde{X}_{t \setminus m})^{-1} \tilde{X}_{t \setminus m}' \Sigma^* \tilde{X}_{t \setminus m} (\tilde{X}_{t \setminus m}'\tilde{X}_{t \setminus m})^{-1}$.

We shall take $c_n(u)= 0.5 b_n(u) + \bar{\omega}_{tq} \tilde{\lambda}_{tm}/\log(\tilde{\lambda}_{tm})$, so that
$c_n(u) - b_n(u)/2= \bar{\omega}_{tq} \tilde{\lambda}_{tm}/\log(\tilde{\lambda}_{tm})$ and the first term in \eqref{eq:tailprob_nonspur_zellnerhet_known} leads to
\begin{align}
\int_0^1 P_{f^*} \left( \frac{W_{qm}}{\phi^*} < c_n(u) - \frac{b_n(u)}{2}  \right) du=
P_{f^*} \left( \frac{W_{qm}}{\phi^*} < \frac{\bar{\omega}_{tq} \tilde{\lambda}_{tm}}{\log(\tilde{\lambda}_{tm})}  \right)
\leq
P_{f^*} \left( Z_2 < \frac{\bar{\omega}_{tq} \tilde{\lambda}_{tm}}{\underline{\omega}_{mq} \log(\tilde{\lambda}_{tm})}  \right).
\nonumber
\end{align}
Since $Z_2$ follows a chi-square distribution with non-centrality parameter $\tilde{\lambda}_{tm}$, we shall apply Lemma \ref{LEM:INEQ_NCCHISQ_LEFT} to bound its left tail.
To apply Lemma \ref{LEM:INEQ_NCCHISQ_LEFT} we need that the cutoff is less than the non-centrality parameter,
that is $\bar{\omega}_{tq} \tilde{\lambda}_{tm} / \underline{\omega}_{mq} \log(\tilde{\lambda}_{tm}) \ll \tilde{\lambda}_{tm}$,
which holds under Condition (C2'), which states that $[\underline{\omega}_{mq}/\bar{\omega}_{tq}] \log(\tilde{\lambda}_{tm}) \gg 1$.
We hence set $w= \tilde{\lambda}_{tm}/[\underline{\omega}_{mq} \log (\tilde{\lambda}_{tm})]$ into the statement of Lemma \ref{LEM:INEQ_NCCHISQ_LEFT}, obtaining
\begin{align}
P_{f^*} \left( \underline{\omega}_{mq} Z_2 < \frac{\bar{\omega}_{tq} \tilde{\lambda}_{tm}}{\log(\tilde{\lambda}_{tm})}  \right)
\leq \exp \left\{ -\frac{1}{2} \tilde{\lambda}_{tm} \left(1 - \sqrt{\frac{\bar{\omega}_{tq}}{\underline{\omega}_{mq} \log (\tilde{\lambda}_{tm})}}  \right)^2\right\} \ll e^{-\frac{\gamma \tilde{\lambda}_{tm}}{2}},
\label{eq:term1_nonspur_zellnerhet_known}
\end{align}
for any fixed $\gamma <1$.

Consider now the second term in \eqref{eq:tailprob_nonspur_zellnerhet_known}, since 
$c_n(u) + b_n(u)/2= b_n(u) + \bar{\omega}_{tq} \tilde{\lambda}_{tm}/\log(\tilde{\lambda}_{tm})$, the goal is to bound
\begin{align}
 \int P_{f^*} \left( \frac{W_{qt}}{\phi^*} > b_n(u) + \frac{\bar{\omega}_{tq} \tilde{\lambda}_{tm}}{\log(\tilde{\lambda}_{tm})} \right) du
\leq 
 \int P_{f^*} \left( \bar{\omega}_{tq} Z_1 > b_n(u) + \frac{\bar{\omega}_{tq} \tilde{\lambda}_{tm}}{\log(\tilde{\lambda}_{tm})} \right) du=
\nonumber \\
 \int P_{f^*} \left( Z_1 >
2 \log \left( \left[ \frac{(1+\tau n)^{\frac{p_m-p_t}{2}} p(M_t)}{(1/u-1) p(M_m)} \right]^{1/\bar{\omega}_{tq}}
 e^{\frac{\tilde{\lambda}_{tm}}{2 \log(\tilde{\lambda}_{tm})}} \right) \right) du
\nonumber
\end{align}
where recall that $Z_1 \sim \chi^2_{p_q-p_t}$.

This expression is analogous to \eqref{eq:proof_spur_heterosk_mainterm} obtained in the spurious models case, and can be bounded via the right tails of a central chi-square distribution, provided the cutoff is larger than its degrees of freedom $p_q - p_t$, that is
\begin{align}
 \left[ \frac{(\tau n)^{\frac{p_m-p_t}{2}} p(M_t)}{p(M_m)} \right]^{1/\bar{\omega}_{tq}}
 e^{\frac{\tilde{\lambda}_{tm}}{2 \log(\tilde{\lambda}_{tm})}} \gg p_q - p_t 
\nonumber
\end{align}
which holds under Condition (C2'), since $p_q - p_t \leq p_m$.

Arguing as in \eqref{eq:proof_spur_heterosk_mainterm} shows that the second term in \eqref{eq:tailprob_nonspur_zellnerhet_known} is
\begin{align}
E_{f^*}(p(M_m\mid\by))
\ll
\left[ \frac{[p(M_m)/p(M_t)]}{(1+\tau n)^{(p_m-p_t)/2} } \right]^{\alpha \min\{1, 1 / \bar{\omega}_{tq} \}}
e^{-\frac{\tilde{\lambda}_{tm}}{2 \log(\tilde{\lambda}_{tm})} \min\{\alpha, \alpha \bar{\omega}_{tq} \}  },
\label{eq:term2_nonspur_zellnerhet_known}
\end{align}
for any fixed $\alpha < 1$.

To complete the proof we combine \eqref{eq:term1_nonspur_zellnerhet_known} and \eqref{eq:term2_nonspur_zellnerhet_known}, which gives
\begin{align}
 E_{f^*}(p(M_k \mid \by)) \ll
\max \left\{  e^{ -\gamma \tilde{\lambda}_{tm}/2},
\left[ \frac{[p(M_m)/p(M_t)]}{(\tau n)^{(p_m-p_t)/2} } \right]^{\gamma \min\{1, 1 / \bar{\omega}_{tq} \}}
e^{-\frac{\gamma \min\{1, \bar{\omega}_{tq} \} \tilde{\lambda}_{tm}}{2}  }
\right\}.
\nonumber
\end{align}
for any fixed $\gamma < 1$, as we wished to prove.

\section{Proof of Lemma \ref{LEM:CONSISTENCYCOND_MODELPRIORS}}

Lemma \ref{LEM:CONSISTENCYCOND_MODELPRIORS} gives two sets of results for the cases $p_m \geq p_t$ and $p_m < p_t$.
The former includes both spurious models ($m \in S$) and non-spurious models ($m \not\in S$). We prove each of these sub-cases separately.

\subsection{Spurious models}

Let $m \in S$ be a spurious model.
We wish to prove (C1), that is as $n\rightarrow\infty$, 
$$ 
(n\tau)^{(p_m-p_t)/2} \frac{p(M_t)}{p(M_m)} \gg 1
$$
Note that for $m \in S$ we have $p_m-p_t \geq 1$.

Under a uniform model prior $p(M_m)/p(M_t)=1$, thus (C1) holds if and only if $1 \ll n\tau$.

Under a Beta-Binomial(1,1) prior $p(M_m)/p(M_t)= $
\begin{align}
\frac{{p \choose p_t}}{{p \choose p_m}}= \frac{p_m!(p-p_m)!}{p_t!(p-p_t)!}
\leq \frac{p_m^{p_m-p_t} (p-p_m)!}{(p-p_t)!}
\asymp \frac{(ep_m)^{p_m-p_t} (p-p_m)^{p-p_m+1/2}}{(p-p_t)^{p-p_t+1/2}}
  < \frac{(ep_m)^{p_m-p_t}}{(p-p_t)^{p_m-p_t}}
\label{eq:betabin_rbound}
\end{align}
where the $\asymp$ statement follows from Stirling's bound.
Hence, a sufficient condition for $(n\tau)^{(p_m-p_t)/2} p(M_t)/p(M_m) \gg 1$ to hold is
$$
\frac{\sqrt{\tau n} (ep_m)}{(p-p_t)} \gg 1
$$
Since $p_m \leq \bar{p}$, for (C2) to hold it suffices that $n\tau \gg [e\bar{p}/(p-p_t)]^2$.

Proceeding analogously to \eqref{eq:betabin_rbound}, under the complexity prior
\begin{align}
\frac{p(M_m)}{p(M_t)} \asymp \frac{{p \choose p_t}}{p^{c(p_m-p_t)}  {p \choose p_m}} \preceq
  \left( \frac{ep_m}{p^c(p-p_t)}\right)^{p_m-p_t},
\label{eq:complexprior_rbound}
\end{align}
where note that $ep_m/p^c(p-p_t) \asymp p_m/p^{c+1} \ll 1$, implying that $p(M_m)/p(M_t) \ll 1$.
Hence $n\tau \gg 1$ suffices to guarantee $(n\tau)^{(p_m-p_t)/2} p(M_t)/p(M_m) \gg 1$.

\subsection{Non-spurious models}

Let $m \in S^c$. We wish to prove that (C2) holds, i.e. that as $n\rightarrow\infty$,
$$
\frac{\lambda_{tm}}{2 \log(\lambda_{tm})} + \frac{p_m-p_t}{2} \log(n\tau) - \log \left( \frac{p(M_m)}{p(M_t)} \right) \gg  1.
$$ 
We treat separately the cases $p_m \geq p_t$ and $p_m < p_t$.

Consider first $p_m \geq p_t$.
Under a uniform prior $p(M_m)/p(M_t)=1$.
Since $\lambda_{tm} \gg 1$, $p_m-p_t \geq 0$ and
and $n \tau \gg 1$ from Condition (B1), we obtain
$$
\frac{\lambda_{tm}}{2 \log(\lambda_{tm})} + \frac{p_m-p_t}{2} \log(n\tau) \gg 1
$$
as desired.
Under a Beta-Binomial(1,1) prior, in \eqref{eq:betabin_rbound} we showed that
$p(M_m)/p(M_t) \preceq (ep_m/(p-p_t))^{p_m-p_t}$
and hence
$$
\log\left( \frac{p(M_m)}{p(M_t)} \right) \preceq (p_m-p_t) \log(ep_m/(p-p_t)) \ll \frac{\lambda_{tm}}{2 \log(\lambda_{tm})} + (p_m-p_t)\log(n\tau),
$$
since $\lambda_{tm} \gg 1$ and  $n\tau \gg \bar{p}/(p-p_t) > p_m/(p-p_t)$ by Condition (B1).
Similarly, for the complexity prior \eqref{eq:complexprior_rbound} implies
$\log(p(M_m)/p(M_t)) \preceq (p_m-p_t) \log \left( \frac{ep_m}{p^c(p-p_t)}\right)$.
Since $\lambda_{tm} \gg 1$ and $n\tau \gg 1$ by assumption,
$$
(p_m-p_t) \log \left( \frac{ep_m}{p^c(p-p_t)}\right) \ll \frac{\lambda_{tm}}{2 \log(\lambda_{tm})} + \frac{p_m-p_t}{2} \log(n\tau)
$$
as we wished to prove.

Consider now $p_m < p_t$.
Under a uniform prior $\log(p(M_m)/p(M_t))=0$, hence (C2) holds if and only if
$\lambda_{tm} - (p_t-p_m) \log(n\tau) \gg 1$.
Under a Beta-Binomial(1,1) prior
$\log(p(M_m)/p(M_t)) \preceq (p_m-p_t) \log(ep_m/(p-p_t))$,
thus a sufficient condition for (C2) to hold is that
$$
(p_m-p_t) \log\left(\frac{ep_m}{p-p_t}\right) \ll \frac{\lambda_{tm}}{2 \log(\lambda_{tm})} + \frac{p_m-p_t}{2} \log(n\tau)
\Leftrightarrow
(p_t-p_m) \log\left(\frac{\sqrt{n\tau} (p-p_t)}{ep_m}\right) \ll \frac{\lambda_{tm}}{2 \log(\lambda_{tm})},
$$
as we wished to prove.
Finally, under the complexity prior
we saw that $\log(p(M_m)/p(M_t)) \preceq (p_m-p_t) \log \left( \frac{ep_m}{p^c(p-p_t)}\right)$,
hence a sufficient condition for (C2) is
$$
(p_t-p_m) \log\left(\frac{\sqrt{n\tau} p^c (p-p_t)}{ep_m}\right) \ll \frac{\lambda_{tm}}{2 \log(\lambda_{tm})}
$$
concluding the proof.

\section{Proof of Lemma \ref{LEM:LAMBDABOUND_SMALLMODELS}}

We first prove Part (i). 
Recall that $\lambda_{tm}=(\btheta_t^*)' X_t' (I-H_m) X_t \btheta_t^*/\phi^*$,
where $I-H_m=I-X_m(X_m'X_m)^{-1}X_m'$ is a projection matrix with rank $\geq n-p_m$.
Let $q$ be the rank of $X_t' (I - H_m) X_t$ and note that,
since $X_t$ is assumed full-rank and $I-H_m$ projects onto the orthogonal space to that spanned by the columns of $X_m$,
we have that $q \in [p_t-p_m,p_t]$.
Consider the eigendecomposition $X_t' (I-H_m) X_t/n=E W E'$,
where $E$ is a $p_t \times q$ matrix containing the eigenvectors
and $W=\mbox{diag}(w_1,\ldots,w_q)$ contains its eigenvalues.
Since $w_j \geq v_{tm}$, it follows that
\begin{align}
\phi^* \lambda_{tm}= n (\btheta_t^*)' EWE' \btheta_t^*=
n \sum_{j=1}^{q} w_j (\btheta_t^*)' {\bf e}_j {\bf e}_j' \btheta_t^*
\geq n v_{tm} \sum_{j=1}^{q}  (\btheta_t^*)' {\bf e}_j {\bf e}_j' \btheta_t^*=
n v_{tm} (\btheta_t^*)' EE' \btheta_t^*.
\nonumber
\end{align}

We consider separately the cases where $X_t'(I-H_m)X_t$ has full rank $q=p_t$ and reduced rank $q < p_t$.
If $X_t'(I-H_m)X_t$ has full-rank ($q=p_t$), then $E$ is a square orthonormal matrix satisfying $EE'=I$,
hence 
$$\lambda_{tm} \geq \frac{n v_{tm} (\btheta_t^*)' \btheta_t^*}{\phi^*} \geq n v_{tm} p_t \min_{j} \frac{(\theta_j^*)^2}{\phi} 
\geq n v_{tm} (p_t - p_m) \min_{j} \frac{(\theta_j^*)^2}{\phi},
$$
as we wished to prove.
Consider now the reduced rank $q<p_t$ case, then
\begin{align}
n v_{tm} \sum_{j=1}^{q}  (\btheta_t^*)' {\bf e}_j {\bf e}_j' \btheta_t^*
> n v_{tm} \left( \min_j (\theta_{j}^*)^2 \right) {\bf 1}' EE' {\bf 1},
\nonumber
\end{align}
where ${\bf 1}=(1,\ldots,1)'$ is the $q\times 1$ vector.
Since the columns in $E$ contain eigenvectors and these have unit length,
we obtain ${\bf 1}' E= {\bf 1}'$ and hence ${\bf 1}' EE' {\bf 1}=q \geq p_t-p_m$. Therefore
\begin{align}
\lambda_{tm} \geq \frac{n v_{tm} (\btheta_t^*)' EE' \btheta_t^*}{\phi^*} \geq v_{tm} q \frac{\min_j (\theta_{j}^*)^2}{\phi^*}
\geq n v_{tm} (p_t - p_m) \frac{\min_j (\theta_{j}^*)^2}{\phi^*},
\nonumber
\end{align}
as we wished to prove.

The proof of Part (ii) follows easily by noting that, by definition models $m \in S_{p_m,j}^c$ miss exactly $p_m-j$ columns from $X_t$.
Since the design matrix $X_q$ has full-rank, this implies that the columns in $X_m$ are linearly independent from those in $X_t$, and therefore that $X_t'(I-H_m)X_t$ has rank $q=p_t - j$.
The rest of the proof proceeds as for Part (i), giving
\begin{align}
\lambda_{tm} \geq \frac{n v_{tm} (\btheta_t^*)' EE' \btheta_t^*}{\phi^*} \geq v_{tm} q \frac{\min_j (\theta_{j}^*)^2}{\phi^*}
= n v_{tm} (p_t - j) \frac{\min_j (\theta_{j}^*)^2}{\phi^*},
\nonumber
\end{align}
as we wished to prove.

\section{Proof of Proposition \ref{PROP:VARSEL_NONSPUR}}

Denote by $r_{p_k,p_t}=p(M_k)/p(M_t)$ the ratio of prior probabilities between models of size $p_k$ and $p_t$.
Recall that the uniform prior on the models $p(M_k)$ corresponds to $r_{p_k,p_t}=1$, the Complexity prior with parameter $c > 0$ to
\begin{align}
 r_{p_k,p_t}= \frac{{p \choose p_t}}{{p \choose p_k}} \frac{1}{p^{c(p_k-p_t)}},
\nonumber
\end{align}
and the Beta-Binomial to the latter expression with $c=0$.

\subsection*{\underline{Part (i)}}

Consider non-spurious models of size $l<p_t$ and note that the number of such models is $|S_l^c|={p \choose l}$, i.e. it is equal to the total number of size $l$ models, since models of size $<p_t$ cannot contain $M_t$.
Using that 
$\underline{\lambda}= \min_{p_k < p_t} \lambda_{tk}^\alpha / (p_t-p_k)$ by definition, we obtain
\begin{align}
&E_{f^*}\left( \sum_{p_m=0}^{p_t-1} P(S_l^c \mid \by) \right)
=\sum_{l=0}^{p_t-1} \sum_{k \in S_l^c} E_{f^*} ( p(M_k \mid \by) )
\leq \sum_{l=0}^{p_t-1} \sum_{k \in S_l^c} (n\tau)^{\alpha(p_t-l)/2}  \left(\frac{p(M_k)}{p(M_t)}\right)^\alpha e^{- \frac{\lambda_{tk}^\alpha}{2}}
\nonumber \\
&\leq \sum_{l=0}^{p_t-1} {p \choose l} (n\tau)^{\alpha(p_t-l)/2} r_{l,p_t}^\alpha e^{- \frac{\underline{\lambda} (p_t-l)}{2}} 
\leq 
e^{-\frac{p_t \underline{\lambda}}{2}} (n\tau)^{\frac{\alpha p_t}{2}} 
\sum_{l=0}^{p_t-1} \left(\frac{p e^{\underline{\lambda}/2}}{(n\tau)^{\alpha/2}}\right)^l r_{l,p_t}^\alpha.
\label{eq:bound_nonspur_small}
\end{align}
Evaluating \eqref{eq:bound_nonspur_small} for specific $r_{l,p_t}$ gives rates for any desired model space prior.
For the uniform prior $r_{l,p_t}=1$, using the geometric series gives
\begin{align}
E_{f^*}\left( \sum_{p_m=0}^{p_t-1} P(S_l^c \mid \by) \right) \leq
e^{-\frac{\underline{\lambda} p_t}{2}} (n\tau)^{\frac{\alpha p_t}{2}} \left[ \frac{1-\left(p e^{\underline{\lambda}/2} (n\tau)^{-\alpha/2}  \right)^{p_t}}{1-p e^{\underline{\lambda}/2} (n\tau)^{-\alpha/2} } \right].
\nonumber
\end{align}
Since $\lim_{n \rightarrow \infty} \underline{\lambda}/2 + \log p - 0.5 \log(\tau n)= \infty$ by assumption, we have that
$\lim_{n \rightarrow \infty} p e^{\underline{\lambda}/2} (n\tau)^{-\alpha/2}= \infty$ and
the latter expression converges to
\begin{align}
\asymp 
e^{-\frac{\underline{\lambda} p_t}{2}} (n\tau)^{\frac{\alpha p_t}{2}} \frac{\left( p e^{\underline{\lambda}/2} (n\tau)^{-\alpha/2} \right)^{p_t} }{p e^{\underline{\lambda}/2} (n\tau)^{-\alpha/2}}
= 
\frac{p^{p_t-1} (n\tau)^{\alpha/2}}{e^{\underline{\lambda}/2}}
=
e^{-\underline{\lambda}/2 + (p_t-1) \log p  + \frac{\alpha}{2} \log(n \tau)},
\nonumber
\end{align}
as desired.

Regarding the Complexity and Beta-Binomial priors, plugging in the expression for $r_{l,p_t}$ into \eqref{eq:bound_nonspur_small} gives
\begin{align}
 e^{-\frac{p_t \underline{\lambda}}{2}} (n\tau)^{\frac{\alpha p_t}{2}} 
\sum_{l=0}^{p_t-1} {p \choose l} \left(\frac{e^{\underline{\lambda}/2}}{(n\tau)^{\alpha/2}}\right)^l
\left[ \frac{{p \choose p_t}}{{p \choose l}} p^{c(p_t-l)} \right]^\alpha.
\nonumber
\end{align}
Using that ${p \choose l} < p^l$ and that ${p \choose p_t} < p^{p_t}$, this expression is
\begin{align}
< p^{\alpha c p_t} e^{-\frac{p_t \underline{\lambda}}{2}} (n\tau)^{\frac{\alpha p_t}{2}} 
\sum_{l=0}^{p_t-1} p^{l (1 - \alpha)} \left(\frac{e^{\underline{\lambda}/2}}{(n\tau)^{\alpha/2}}\right)^l \frac{p^{p_t \alpha}}{p^{c \alpha l}}
=
p^{\alpha p_t (c+1)} e^{-\frac{p_t \underline{\lambda}}{2}} (n\tau)^{\frac{\alpha p_t}{2}} 
\sum_{l=0}^{p_t-1}  \left(\frac{p^{1 - \alpha (1+c)} e^{\underline{\lambda}/2}}{(n\tau)^{\alpha/2}}\right)^l 
\nonumber \\
=
p^{\alpha p_t (c+1)} e^{-\frac{p_t \underline{\lambda}}{2}} (n\tau)^{\frac{\alpha p_t}{2}} 
\frac{1 - \left( \frac{p^{1 - \alpha (1+c)} e^{\underline{\lambda}/2}}{(n\tau)^{\alpha/2}} \right)^{p_t}}{1 - \frac{p^{1 - \alpha (1+c)} e^{\underline{\lambda}/2}}{(n\tau)^{\alpha/2}}}.
\label{eq:bound_nonspur_small_complexity}
\end{align}
the right-hand side following from the geometric series.
Since $\lim_{n \rightarrow \infty} \underline{\lambda}/2 - (\alpha (1+c) - 1) \log p - 0.5 \log(n\tau) = \infty$ by assumption, it follows that 
\begin{align}
\lim_{n \rightarrow \infty}
\frac{p^{1 - \alpha (1+c)} e^{\underline{\lambda}/2}}{(n\tau)^{\alpha/2}} = \infty
\nonumber
\end{align}
and hence \eqref{eq:bound_nonspur_small_complexity} converges to
\begin{align}
\asymp
 p^{\alpha p_t (c+1)} e^{-\frac{p_t \underline{\lambda}}{2}} (n\tau)^{\frac{\alpha p_t}{2}} 
\frac{\left( \frac{p^{1 - \alpha (1+c)} e^{\underline{\lambda}/2}}{(n\tau)^{\alpha/2}} \right)^{p_t}}{\frac{p^{1 - \alpha (1+c)} e^{\underline{\lambda}/2}}{(n\tau)^{\alpha/2}}}
=
\frac{p^{p_t - 1 + \alpha (1+c)} (n\tau)^{\alpha/2}}{e^{\underline{\lambda}/2}}=
e^{-\frac{\underline{\lambda}}{2} + [p_t -1 + \alpha (1+c)] \log p + \frac{\alpha}{2} \log(n\tau) }
,
\nonumber
\end{align}
as we wished to prove. Note that the result for the uniform prior corresponds to $c= -1$, and that for the Beta-Binomial to $c=0$.


\subsection*{\underline{Part (ii)}}

We obtain separate bounds for $p_m=p_t$ and $p_m>p_t$, then add them up.
First consider $p_m=p_t$. Since ${p \choose p_t} \leq p^{p_t}$, we obtain that for all $n \geq n_0$
\begin{align}
E_{f^*}\left( \sum_{l=p_t} P(S_l^c \mid \by) \right) \leq
{p \choose p_t} e^{- \min_{p_m=p_t} \frac{\lambda_{tm}^\alpha}{2}}
\leq \left(\frac{pe}{p_t} \right)^{p_t} e^{- \min_{p_m=p_t} \frac{\lambda_{tm}^\alpha}{2}}
\leq e^{- \frac{\bar{\lambda}}{2} + p_t \log(pe)}.
\label{eq:nonspur_varsel_equal}
\end{align}

Now consider models size $p_m>p_t$.
Let $M_m$ be such a model, and denote by $j$ the number of truly active variables selected by $M_m$ out of the $p_t$, so that $p_m-j$ is the number of truly inactive variables in $M_m$. 
Note that there are ${p_t \choose j} {p-p_t \choose p_m-j}$ models selecting the same number of active and inactive variables than $M_m$, 
and that all these models miss $p_t - j$ truly active variables,
so that $\bar{\lambda} \leq \lambda_{tm}^\alpha/(p_t-j)$ for all these models. Hence we obtain that
for any $n \geq n_0$,
\begin{align}
&E_{f^*} \left( \sum_{l=p_t+1}^{\bar{p}} P(S_l^c \mid \by) \right) \leq
\sum_{l=p_t+1}^{\bar{p}} \sum_{j=0}^{p_t-1}  {p_t \choose j} {p-p_t \choose l-j} 
\frac{r_{l,p_t}^\alpha  e^{-(p_t-j) \bar{\lambda}/2}}{(n\tau)^{\alpha (l-p_t)/2}} 
  \nonumber \\
&= (n\tau)^{\alpha p_t/2} e^{-p_t \bar{\lambda}/2} \sum_{j=0}^{p_t-1}  \left(\frac{p_t}{p-p_t}\right)^j e^{j \bar{\lambda}/2}
\sum_{l=p_t+1}^{\bar{p}} \frac{(p-p_t)^l r_{l,p_t}^\alpha}{(n\tau)^{\alpha l/2}}.
\label{eq:nonspur_varsel_large1}
\end{align}

To complete the proof we evaluate \eqref{eq:nonspur_varsel_large1} after setting $r_{p_k,p_t}=p(M_k)/p(M_t)$ to its expression for the Complexity prior with parameter $c$. Using Stirling's bound for the factorial function it is easy to show that
\begin{align}
  r_{l,p_t}= \frac{{p \choose p_t}}{{p \choose l}} \frac{1}{p^{c(l-p_t)}} \asymp \frac{1}{p^{(c+1)(l-p_t)}}.
\nonumber
\end{align}
Recall that the Beta-Binomial corresponds to $c=0$, and note that when setting $c=-1$ above one recovers $r_{l,p_t}=1$, i.e. the uniform prior on the models.
For simplicity we replace $r_{l,p_t}$ by this asymptotic expression, but multiplying it by a constant gives an exact finite $n$ expression.

Using the geometric series, the inner summation in \eqref{eq:nonspur_varsel_large1} becomes
\begin{align}
p^{p_t \alpha (c+1)} \sum_{l=p_t+1}^{\bar{p}} \left( \frac{p-p_t}{(n\tau)^{\alpha/2} p^{\alpha(c+1)}} \right)^l=
\nonumber
& p^{p_t \alpha (c+1)}
\frac{z^{\bar{p}+1} - z^{p_t+1}}{1- z}
\nonumber
\end{align}
where
\begin{align}
 z= \frac{p-p_t}{(n\tau)^{\alpha/2} p^{\alpha(c+1)}}.
\nonumber
\end{align}
It is simple to show that $(z^{\bar{p}+1} - z^{p_t+1})/(1-z) \leq z^{\bar{p}}$ for any $z>0$, hence \eqref{eq:nonspur_varsel_large1} is
\begin{align}
 \leq (n\tau)^{\alpha p_t/2} e^{-p_t \bar{\lambda}/2} \sum_{j=0}^{p_t-1}  \left(\frac{p_t e^{\bar{\lambda}/2}}{p-p_t}\right)^j 
p^{p_t \alpha (c+1)} \left(  \frac{p-p_t}{(n\tau)^{\alpha/2} p^{\alpha(c+1)}} \right)^{\bar{p}}
\nonumber \\
=(n\tau)^{\alpha p_t/2} e^{-p_t \bar{\lambda}/2} p^{p_t \alpha (c+1)} \left(  \frac{p-p_t}{(n\tau)^{\alpha/2} p^{\alpha(c+1)}} \right)^{\bar{p}}
\frac{1 - \left( \frac{p_t e^{\bar{\lambda}/2}}{p-p_t} \right)^{p_t}}{1 - \frac{p_t e^{\bar{\lambda}/2}}{p-p_t}}
\label{eq:nonspur_varsel_large2}
\end{align}

To complete the proof we consider separately the cases where $\bar{\lambda}/2 + \log p_t - \log(p - p_t)$ converges to $\infty$ and $-\infty$.
If $\lim_{n \rightarrow \infty} \bar{\lambda}/2 + \log p_t - \log(p - p_t)= \infty$, then \eqref{eq:nonspur_varsel_large2} converges to
\begin{align}
 =(n\tau)^{\alpha p_t/2} e^{-p_t \bar{\lambda}/2} p^{p_t \alpha (c+1)} \left(  \frac{p-p_t}{(n\tau)^{\alpha/2} p^{\alpha(c+1)}} \right)^{\bar{p}}
\left( \frac{p_t e^{\bar{\lambda}/2}}{p-p_t} \right)^{p_t-1}
\nonumber \\
=
 e^{-\bar{\lambda}/2}
 \frac{(p-p_t)^{\bar{p} - p_t + 1}}{[(n\tau)^{\alpha/2} p^{\alpha(c+1)}]^{\bar{p} - p_t}}
p_t^{p_t-1}
<
 \frac{ e^{-\bar{\lambda}/2 + p_t \log p_t + \log p}}{[(n\tau)^{\alpha/2} p^{\alpha(c+1)-1}]^{\bar{p} - p_t}},
\nonumber
\end{align}
as we wished to prove.

In contrast, if $\lim_{n \rightarrow \infty} \bar{\lambda}/2 + \log p_t - \log(p - p_t)= -\infty$, then the last term in \eqref{eq:nonspur_varsel_large2} converges to 1, hence the limit of \eqref{eq:nonspur_varsel_large2} is
\begin{align}
\leq e^{-p_t\bar{\lambda}/2}
 \left( \frac{1}{n \tau} \right)^{\frac{\alpha (\bar{p} - p_t)}{2}} \frac{1}{p^{\alpha (c+1) (\bar{p}-p_t) -1}},
\nonumber
\end{align}
as we wished to prove.

\section{Proof of Lemma \ref{LEM:INEQ_CHISQ}}

To bound the right tail consider $s \in (0,1/2)$, then by Markov's inequality
$$
P(W>w)= P(e^{sW} > e^{sw}) \leq e^{-ws} E(e^{sW}),
$$
where $E(e^{sW})= (1-2s)^{-\frac{\nu}{2}}$. The log of the right hand side is minimized
for $s^*=\frac{1}{2} - \frac{\nu}{2w}$, where note that $w>\nu$ guarantees that $s^*>0$.
Plugging $s=s^*$ gives the bound stated above.

To bound the left tail consider $s<0$, then by Markov's inequality
$$
P(W<w)= P(e^{sW} > e^{ws}) \leq e^{-ws} E(e^{sW}),
$$
where the right hand side is again minimized by
$s^*=\frac{1}{2} - \frac{\nu}{2w}$ and note that $w<\nu$ guarantees that $s^*<0$.
Plugging $s=s^*$ gives the desired bound.

\section{Proof of Lemma \ref{LEM:INEQ_NCCHISQ_LEFT}}

The result is analogous to Lemma \ref{LEM:INEQ_CHISQ}.
For any $s<0$ Markov's inequality gives $P(W<w)= P(e^{sW} > e^{ws}) \leq e^{-ws} E(e^{sW})$,
where $E(e^{sW})=  \frac{\exp\{\frac{\lambda s}{1-2s}\}}{(1-2s)^{\nu/2}}$
is the non-central chi-square moment generating function.
Straightforward algebra shows that the upper bound is minimized for
$s= \frac{1}{2} - \frac{\nu}{4w} - \frac{1}{2} \sqrt{\frac{\nu^2}{4w^2} + \frac{\lambda}{w}}$.

\section{Proof of Lemma \ref{LEM:INEQ_NCCHISQ_RIGHT}}

The proof is identical to that of Lemma \ref{LEM:INEQ_NCCHISQ_LEFT} except that here we restrict $s \in (0,1/2)$.
The restriction $s<1/2$ arises from the fact that the chi-square moment generating function is undefined for $s \geq 1/2$.

\section{Proof of Lemma \ref{LEM:INEQ_NCF_RIGHT}}

By definition $\nu_1W=U_1\nu_2/U_2$, where $U_1 \sim \chi_{\nu_1}^2(\lambda)$ and $U_2 \sim \chi_{\nu_2}^2$.
Let $s>0$. If $U_1<ws$ and $U_2/\nu_2>s$ then it follows that $U_1 \nu_2/U_2 < w$, hence
$$
P \left( \frac{U_1 \nu_2}{U_2} > w \right) \leq P(U_1 > ws) + P(U_2 < s \nu_2).
$$

The result for the case $\lambda =0$ is obtained by using
the Chernoff bounds in Lemma \ref{LEM:INEQ_CHISQ} to bound $P(U_1 > ws)$ and $P(U_2 < s \nu_2)$.
When $\lambda >0$ use Lemma \ref{LEM:INEQ_NCCHISQ_RIGHT} to bound $P(U_1 > ws)$
and Lemma \ref{LEM:INEQ_CHISQ} to bound $P(U_2 < s \nu_2)$.
Note that Lemma \ref{LEM:INEQ_NCCHISQ_RIGHT} requires $ws > \lambda + \nu_1$
and Lemma \ref{LEM:INEQ_CHISQ} requires $s \nu_2 < \nu_2$, i.e. $s$ must satisfy the constraints $(\lambda + \nu_1)/w < s < 1$.

\section{Proof of Lemma \ref{LEM:INEQ_NCF_LEFT}}

By definition $\nu_1W=U_1\nu_2/U_2$, where $U_1 \sim \chi_{\nu_1}^2(\lambda)$ and $U_2 \sim \chi_{\nu_2}^2$.
Then for any $s \geq 1$
\begin{align}
P(\nu_1W < w) \leq
P \left( U_1 < ws \right) + P(U_2 > \nu_2 s) \leq
\frac{\exp\{\frac{\lambda t}{1-2t} -tws\}}{(1-2t)^{\nu_1/2}}
+ e^{-\frac{\nu_2}{2} (s - 1 -\log(s))}
\nonumber
\end{align}
where the first term is given by Lemma \ref{LEM:INEQ_NCCHISQ_LEFT} for any $t<0$
and the second term by Lemma \ref{LEM:INEQ_CHISQ}.
By Lemma \ref{LEM:INEQ_NCCHISQ_LEFT} the first term is minimized for
$t= \frac{1}{2} - \frac{\nu_1}{4ws} - \frac{1}{2} \sqrt{\frac{\nu_1^2}{4ws^2} + \frac{\lambda}{ws}}$,
and if $ws < \lambda$ then we may set $t=\frac{1}{2} -\frac{1}{2} \sqrt{\frac{\lambda}{ws}}$,
obtaining
$$
\frac{e^{-\frac{1}{2} (\sqrt{\lambda} - \sqrt{ws})^2}}{(\lambda/(ws))^{\frac{\nu_1}{4}}}
+ e^{-\frac{\nu_2}{2} (s - 1 -\log(s))}.
$$

\section{Proof of Lemma \ref{LEM:INEQ_F}}

Let $s>0$, by Markov's inequality $P(W>w) \leq E(W^s)/w^s$.
For $s \in [1,\nu_2/2-2]$ we have that
\begin{align}
E(W^s)= \left( \frac{\nu_2}{\nu_1} \right)^s \frac{\Gamma \left( \frac{\nu_1}{2}+s \right)}{\Gamma \left( \frac{\nu_1}{2} \right)}
\frac{\Gamma \left( \frac{\nu_2}{2} -s \right)}{\Gamma \left( \frac{\nu_2}{2} \right)}.
\label{eq:f_rawmoments}
\end{align}
Expression \eqref{eq:ftail_momentbound} is obtained
by upper-bounding the Gamma functions in the numerator in \eqref{eq:f_rawmoments}
and lower-bounding the denominator using Stirling's formula bounds
$$
\sqrt{2\pi} z^{z+\frac{1}{2}} e^{-z} \leq \Gamma \left( z + 1 \right) \leq e z^{z+\frac{1}{2}} e^{-z},
$$
which hold for any $z>0$, i.e. $\nu_1 > 2$, $\nu_2 > 2$.
The bounds for $\nu_1=1$ and $\nu_1=2$ are obtained similarly by noting that then $\Gamma(\nu_1/2)=\sqrt{\pi}$ and $\Gamma(\nu_1/2)=1$ respectively.

To motivate the choice $s=\min\{ (w-1)\nu_1/2+1, \nu_2/2-2 \}$, we minimize an asymptotic
version of the bound as $\nu_2 \rightarrow \infty$.
Let $c= \lim_{\nu_2 \rightarrow \infty} 2s/\nu_2 \in [0,1)$, then in \eqref{eq:ftail_momentbound}
$$
\lim_{\nu_2 \rightarrow \infty} \frac{\nu_2}{(\nu_2/2-s-1)}= \frac{1}{1/2 - c/2}= 2/(1-c)
$$
and
$$
\left( 1 - \frac{s}{\nu_2/2-1} \right)^{\frac{\nu_2-1}{2}}=
\left( 1 - \frac{s}{\nu_2/2-1} \right)^{\frac{\nu_2}{2}-1} \left( 1 - \frac{s}{\nu_2/2-1} \right)^{\frac{1}{2}}
 \leq e^{-s},
$$
since $(1-s/z)^z \leq e^{-s}$ for all $z \geq 1$.
Hence as $\nu_2 \rightarrow \infty$
\begin{align}
P(W>w) \preceq a \left( \frac{2/(1-c)}{\nu_1w} \right)^s \frac{(s+\nu_1/2-1)^{s+\frac{\nu_1-1}{2}}}{(\nu_1/2-1)^{\frac{\nu_1-1}{2}\mbox{I}(\nu_1>2)}} e^{-s}
\nonumber
\end{align}
Setting the derivative of the log-bound equal to zero gives
$$
\log \left( \frac{2/(1-c)}{\nu_1w} \right) + \log(s+\frac{\nu_1}{2}-1) + \frac{1}{2(s+\frac{\nu_1}{2}-1)}=0.
$$
Note that as $w \rightarrow \infty$ the solution to the equation above must satisfy $s \rightarrow \infty$, hence its third term converges to 0
and an approximate solution is given by
$$
\frac{\nu_1w}{2/(1-c)} = s+\frac{\nu_1}{2}-1 \Rightarrow
s= \frac{\nu_1}{2} \left( \frac{w}{1-c} - 1 \right) + 1.
$$
Recall that $s<\nu_2/2$, so that $c$ is 0 when $s \ll \nu_2/2$ and $c<1$ when $s \asymp \nu_2/2$,
hence the solution must satisfy $s \geq \frac{\nu_1}{2} \left( w - 1 \right) + 1$.
We thus take $s= \min\{ \frac{\nu_1}{2} \left( w - 1 \right) + 1, \nu_2/2-2\}$.

Finally, to obtain the desired bounds note that if $w \leq (\nu_1+\nu_2-6)/\nu_1$
then $s=\frac{\nu_1}{2} \left( w - 1 \right) + 1$,
whereas if $w > (\nu_1+\nu_2-6)/\nu_1$ then $s=\nu_2/2-2$.
Plugging in $s=\frac{\nu_1}{2} \left( w - 1 \right) + 1$ and rearranging terms
gives the expressions in Lemma \ref{LEM:INEQ_F} for the $w \leq (\nu_1+\nu_2-6)/\nu_1$ case.

When $w > (\nu_1+\nu_2-6)/\nu_1$, plug in $s=\nu_2/2-2$ to obtain
\begin{align}
a \left( \frac{\nu_2}{\nu_2-2} \right)^{\frac{\nu_2}{2}-2}
\left( \frac{\nu_1+\nu_2-6}{\nu_1w} \right)^{\frac{\nu_2}{2}-2}
\frac{(\nu_1/2+\nu_2/2-3)^{\frac{\nu_1-1}{2}}}{(\nu_1/2-1)^{\frac{\nu_1-1}{2}\mbox{I}(\nu_1>2)}}  \left( \frac{1}{\nu_2/2-1} \right)^{\frac{3}{2}}.
\label{eq:ftailbound2}
\end{align}
and note that
$$
\left( \frac{\nu_2}{\nu_2-2} \right)^{\frac{\nu_2}{2}-2}=
\left(1+ \frac{2}{\nu_2-2} \right)^{\frac{\nu_2}{2}-2}
< e^{\frac{\nu_2-4}{\nu_2-2}} < e.
$$

\section{Proof of Lemma \ref{LEM:NESTEDLM_CHITEST}}

Let $\tilde{X}_s= X_s - X_m(X_m'X_m)^{-1}X_m'X_s$ be orthogonal to the projection of $X_s$ on $X_m$.
Then clearly $X_m'\tilde{X}_s=0$ and
\begin{align}\hat{\btheta}_q' X_q'X_q \hat{\btheta}_q - \hat{\btheta}_m' X_m'X_m \hat{\btheta}_m=
(\hat{\btheta}_m' X_m'X_m \hat{\btheta}_m
+ \tilde{\btheta}_s' \tilde{X}_s'\tilde{X}_s \tilde{\btheta}_s) - \hat{\btheta}_m' X_m'X_m \hat{\btheta}_m
=\tilde{\btheta}_s' \tilde{X}_s'\tilde{X}_s \tilde{\btheta}_s,
\nonumber
\end{align}
where $\tilde{\btheta}_s= (\tilde{X}_s'\tilde{X}_s)^{-1} \tilde{X}_s' \by$.
Since $E_{f^*}(\by)=X_m\btheta_m^* + X_s\btheta_s^*$, $\tilde{X}_s'X_m={\bf 0}$
and $\tilde{X}_s'X_s= \tilde{X_s}'(\tilde{X}_s +X_m(X_m'X_m)^{-1}X_m'X_s)=\tilde{X}_s'\tilde{X}_s$
we have
that $\tilde{\btheta}_s$ is normally distributed with mean
$$
E_{f^*}(\tilde{\btheta}_s)=(\tilde{X}_s'\tilde{X}_s)^{-1} \tilde{X}_s'(X_m \btheta_m^* + X_s \btheta_s^*)=
(\tilde{X}_s'\tilde{X}_s)^{-1} \tilde{X}_s'X_s \btheta_s^*=
(\tilde{X}_s'\tilde{X}_s)^{-1} \tilde{X}_s'\tilde{X}_s \btheta_s^*= \btheta_s^*
$$
and covariance $\phi^* (\tilde{X}_s'\tilde{X}_s)^{-1}=\phi^* [X_s' (I-X_m(X_m'X_m)^{-1}X_m') X_s]^{-1}$.
That is,
$$
\tilde{\btheta}_s \sim N(\btheta_s^*, \phi^* (\tilde{X}_s'\tilde{X}_s)^{-1}),
$$
implying that
$\tilde{\btheta}_s' \tilde{X}_s'\tilde{X}_s \tilde{\btheta}_s / \phi^*$ follows a $\chi^2_{p_q-p_m}$
with non-centrality $\lambda_{qm}= (\btheta_s^*)' \tilde{X}_s'\tilde{X}_s \btheta_s^*/\phi^*$.

Finally, note that for any $k$ with full-rank $X_k$ we have $\btheta_k^*= (X_k'X_k)^{-1}X_k'X_t\btheta_t^*$ and that
$X_s\btheta_s^*= X_q\btheta_q^* - X_m\btheta_m^*= (H_q-H_m) X_t\btheta_t^*$, hence
$\lambda_{qm}= (\btheta_s^*) X_s' (I-H_m) X_s \btheta_s^*/\phi^*=$
$$
\frac{(\btheta_t^*) X_t' (H_q-H_m) (I-H_m) (H_q-H_m)X_t\btheta_t^*}{\phi^*}=
\frac{(\btheta_t^*) X_t' H_q (I-H_m) H_q X_t\btheta_t^*}{\phi^*}.
$$

\section{Proof of Lemma \ref{LEM:NESTEDLM_NONLINEAR}}

The proof is analogous to that for Lemma \ref{LEM:NESTEDLM_CHITEST}.
Briefly for any full-rank $X_k$ we obtain the KL-optimal $\btheta_k^*= (X_k'X_k)^{-1}X_k'W\bbeta^*$,
hence the data-generating truth can be written as
$$
f^*(\by)= N(X_m\btheta_m^* + X_s\btheta_s^* + W\bbeta^* - X_q\btheta_q^*,\xi^*I)=
N(X_m\btheta_m^* + X_s\btheta_s^* + (I-H_q) W\bbeta^*,\xi^*I).
$$
Following the proof of Lemma \ref{LEM:NESTEDLM_CHITEST},
$$
W_{qm}= \hat{\btheta}_q' X_q'X_q \hat{\btheta}_q - \hat{\btheta}_m' X_m'X_m \hat{\btheta}_m=
\tilde{\btheta}_s' \tilde{X}_s'\tilde{X}_s \tilde{\btheta}_s,
$$
where $\tilde{X}_s=(I-H_m)X_s$ and $\tilde{\btheta}_s=(\tilde{X}_s'\tilde{X}_s)^{-1} \tilde{X}_s'\by$.
Since $\tilde{X}_s'X_m= {\bf 0}$ and $\tilde{X}_s'(I-H_q)= {\bf 0}$, it follows that
\begin{align}
E_{f^*}(\tilde{\btheta}_s)= (\tilde{X}_s'\tilde{X}_s)^{-1} \tilde{X}_s' (X_m\btheta_m^* + X_s\btheta_s^* + (I-H_q) W\bbeta^*)=\btheta_s^*
\nonumber
\end{align}
and $\mbox{Cov}_{f^*}(\tilde{\btheta}_s)= \xi^* (\tilde{X}_s'\tilde{X}_s)^{-1}$.
Therefore
$$
\frac{W_{qm}}{\xi^*} \sim \chi^2_{p_q-p_m}(\lambda_{tm})
$$
where $\lambda_{tm}=(\btheta_s^*)' \tilde{X}_s'\tilde{X}_s \btheta_s^*/\xi^*$.
Since $(\btheta_s^*)' \tilde{X}_s'\tilde{X}_s \btheta_s^*= (\btheta_s^*)' X_s' (I-H_m) X_s \btheta_s^*$
and $X_s \btheta_s^*= X_q\btheta_q^*-X_m\btheta_m^*=(H_q-H_m)W\bbeta^*$,
$$
\lambda_{tm}= \frac{(W\bbeta^*)'(H_q-H_m)(I-H_m)(H_q-H_m)W\bbeta^*}{\xi^*}=
\frac{(W\bbeta^*)'H_q(I-H_m)H_qW\bbeta^*}{\xi^*}.
$$

\section{Proof of Lemma \ref{LEM:NESTEDLM_CHITEST_HETERO}}

\subsection*{\underline{Proof of Part (i)}}

In the proof of Lemma \ref{LEM:NESTEDLM_CHITEST} we showed that
$\hat{\btheta}_q' X_q'X_q \hat{\btheta}_q - \hat{\btheta}_m' X_m'X_m \hat{\btheta}_m=
\tilde{\btheta}_s' \tilde{X}_s'\tilde{X}_s \tilde{\btheta}_s$,
where $\tilde{\btheta}_s= (\tilde{X}_s'\tilde{X}_s)^{-1} \tilde{X}_s' \by$,
and $\tilde{X}_s= X_s - X_m(X_m'X_m)^{-1}X_m'X_s$.
By assumption $f^*(\by)= N(\by;X\btheta^*,\phi^* \Sigma^*)$,
hence $f^*(\tilde{\btheta}_s)= N(\tilde{\btheta}_s;\btheta_s^*,\phi^* W)$ where
$$
W=(\tilde{X}_s'\tilde{X}_s)^{-1} \tilde{X}_s' \Sigma^* \tilde{X}_s (\tilde{X}_s'\tilde{X}_s)^{-1}=
(X_s'H_mX_s)^{-1} X_s' H_m \Sigma^* H_m X_s (X_s'H_mX_s)^{-1},
$$
and $H_m=I-X_m(X_m'X_m)^{-1}X_m'$.
Noting that
\begin{align}
\tilde{\btheta}_s' \tilde{X}_s'\tilde{X}_s \tilde{\btheta}_s=
\tilde{\btheta}_s' \tilde{X}_s'\tilde{X}_s W W^{-1} \tilde{\btheta}_s=
\tilde{\btheta}_s' \tilde{X}_s' \Sigma^* \tilde{X}_s (\tilde{X}_s'\tilde{X}_s)^{-1} W^{-1} \tilde{\btheta}_s.
\nonumber
\end{align}
and recalling that $\underline{\omega}_{mq}$ and $\bar{\omega}_{mq}$ are the smallest and largest eigenvalues of
$\tilde{X}_s' \Sigma^* \tilde{X}_s (\tilde{X}_s'\tilde{X}_s)^{-1}$, we obtain
$$
\underline{\omega}_{mq} Z_1
\leq \frac{\tilde{\btheta}_s' \tilde{X}_s'\tilde{X}_s \tilde{\btheta}_s}{\phi^*} \leq
\bar{\omega}_{mq} Z_1,
$$
where $Z_1=\tilde{\btheta}_s' W^{-1} \tilde{\btheta}_s/\phi^* \sim \chi_{p_q-p_m}^2(\tilde{\lambda}_{qm} )$
and $\tilde{\lambda}_{qm}= (\btheta_s^*)' W^{-1} \btheta_s^*/\phi^*$.
Since $W^{-1}= (\tilde{X}_s'\tilde{X}_s) (\tilde{X}_s' \Sigma^* \tilde{X}_s)^{-1} (\tilde{X}_s'\tilde{X}_s)$,
we also have
$$
\frac{(\btheta_s^*)' (\tilde{X}_s'\tilde{X}_s) \btheta_s^*}{\bar{\omega}_{mq} \phi^*}
\leq \tilde{\lambda}_{qm} \leq
\frac{(\btheta_s^*)' (\tilde{X}_s'\tilde{X}_s) \btheta_s^*}{\underline{\omega}_{mq} \phi^*},
$$
where
$(\btheta_s^*)' (\tilde{X}_s'\tilde{X}_s) \btheta_s^*/\phi^*=(\btheta_s^*)'X_s'(I-X_m(X_m'X_m)^{-1}X_m)X_s \btheta_s^*/\phi^*=\lambda_{qm}$,
completing the proof.

\subsection*{\underline{Proof of Part (ii)}}

For any $T \in L_q$ the matrix $X_F=(X_q,T)$ is full-rank, hence
$\by= X_F \hat{\btheta}_F$ and
$
\by'\by - \hat{\btheta}_q' X_q'X_q \hat{\btheta}_q=
\hat{\btheta}_F' X_F'X_F \hat{\btheta}_F - \hat{\btheta}_q' X_q'X_q \hat{\btheta}_q.
$
By Part (i) this implies that
$$
\underline{\omega}_{qF} Z_2
\leq \frac{\by'\by - \hat{\btheta}_q' X_q'X_q \hat{\btheta}_q}{\phi^*} \leq
\bar{\omega}_{qF} Z_2
$$
where $Z_2 \sim \chi_{n-p_q}^2$,
$\underline{\omega}_{qF}=\underline{\varrho}(\tilde{T}' \Sigma \tilde{T} (\tilde{T}'\tilde{T})^{-1})$,
$\bar{\omega}_{qF}= \bar{\varrho}(\tilde{T}' \Sigma \tilde{T} (\tilde{T}'\tilde{T})^{-1})$
and $\tilde{T}= (I - X_q(X_q'X_q)^{-1}X_q')T$.
The bound applies to any $T \in L_q$ and the tightest result is
$$
\underline{\omega}_q Z_2
\leq \by'\by - \hat{\btheta}_q' X_q'X_q \hat{\btheta}_q \leq
\bar{\omega}_q Z_2,
$$
where $\underline{\omega}_q= \max_{L_q} \underline{\omega}_{qF}$, $\bar{\omega}_q= \min_{L_q} \bar{\omega}_{qF}$.

Combining this with Part (i) gives that
$$
\frac{\underline{\omega}_{mq} Z_1}{\bar{\omega}_q Z_2}
\leq \frac{\hat{\btheta}_q' X_q'X_q \hat{\btheta}_q - \hat{\btheta}_m' X_m'X_m \hat{\btheta}_m}{\by'\by - \hat{\btheta}_q' X_q'X_q \hat{\btheta}_q} \leq
\frac{\bar{\omega}_{mq} Z_1}{\underline{\omega}_q Z_2},
$$
as we wished to prove.

\section{Proof of Lemma \ref{LEM:NCCHISQ_BAYES_UNIV}}

From standard Normal theory $\by \sim N(X_q \btheta_q^*, \phi^* I)$ implies that
$\hat{\btheta}_q \sim N(\bmu, \phi^* \Sigma)$.
Hence marginally $\hat{\theta}_{qi} \sim N(\mu_i, \phi^* \sigma_{ii})$
and thus $\hat{\theta}_{qi}^2/(\phi^* \sigma_{ii}) \sim \chi_1^2(\mu_i^2/(\phi^* \sigma_{ii}))$.
All that remains is to bound $\mu_i^2/\sigma_{ii}$.

We first bound $\sigma_{ii}$.
Let ${\bf e}_i$ be the $i^{th}$ canonical eigenvector with $e_{ii}=1$ and $e_{ij}=0$ for $i \leq j$. Then $\sigma_{ii}=$
$${\bf e}_i' \Sigma {\bf e}_i=
{\bf e}_i' (X_q'X_q)^{-1/2} (X_q'X_q)^{1/2} \Sigma (X_q'X_q)^{1/2} (X_q'X_q)^{-1/2} {\bf e}_i=
{\bf e}_i' (X_q'X_q)^{-1/2} A^2 (X_q'X_q)^{-1/2} {\bf e}_i,
$$
where
$A=(X_q'X_q)^{1/2} (X_q'X_q+ V_q^{-1}/(n\tau))^{-1} (X_q'X_q)^{1/2}$
The eigenvalues of $A$ are the same as the eigenvalues of
$(X_q'X_q+ V_q^{-1}/(n\tau))^{-1} X_q'X_q$ (see proof of Lemma \ref{LEM:FSTAT_NORMALPRIOR_BOUND}),
and these are the inverse of the eigenvalues of $I + (X_q'X_q)^{-1} V_q^{-1}/(n\tau)$,
hence the eigenvalues of $A$ are
$n\tau\rho_{q1}/(n\tau\rho_{q1} +1) \geq \ldots \geq n\tau\rho_{qp_q}/(n\tau\rho_{qp_q} +1)>0$.
Therefore
$$
\frac{(n\tau)^2\rho_{qp_q}^2}{(n\tau)^2\rho^2_{qp_q} +1} \tilde{\sigma}_{ii}
\leq \sigma_{ii} \leq
\frac{(n\tau)^2\rho_{q1}^2}{(n\tau)^2\rho_{q1}^2 +1}  \tilde{\sigma}_{ii}
$$
where
$\tilde{\sigma}_{ii}={\bf e}_i' (X_q'X_q)^{-1} {\bf e}_i$
is the $i^{th}$ diagonal element in $(X_q'X_q)^{-1}$.
Using the blockwise matrix inversion formula,
$\tilde{\sigma}_{ii}= ({\bf x}_{qj}' (I - X_m (X_m'X_m)^{-1}X_m') {\bf x}_{qj})^{-1}$.

We now bound $\mu_i^2$.
Note that $\mu_i^2= (\theta_{qi}^*)^2 + (\mu_i-\theta_{qi}^*)^2 + 2\theta_{qi}^*(\mu_i - \theta_{qi}^*)$,
hence if we can find an upper bound $u$ such that $|\mu_i-\theta_{qi}^*|<u$ then it follows that
\begin{align}
\mu_i^2 \geq (\theta_{qi}^*)^2 - 2u|\theta_{qi}^*|= (\theta_{qi}^*)^2 \left( 1 - \frac{2u}{|\theta_{qi}^*|} \right)
\nonumber \\
\mu_i^2 \leq (\theta_{qi}^*)^2 + u^2 + 2u|\theta_{qi}^*|= (\theta_{qi}^*)^2 \left( 1 + \left[\frac{u}{\theta_{qi}^*}\right]^2 + \frac{2u}{|\theta_{qi}^*|} \right)
\label{eq:proof_ncchisq_univ1}
\end{align}
Let $D=(X_q'X_q+ V_q^{-1}/(n\tau))^{-1} X_q'X_q - I$ and note that $\bmu - \btheta_q^*=  D\btheta_q^*$.
We already saw that the largest eigenvalue of $(X_q'X_q+ V_q^{-1}/(n\tau))^{-1} X_q'X_q$ is
$n\tau\rho_{q1}/(n\tau\rho_{q1} +1)$, hence the largest eigenvalue of $D$ is $1/(n\tau\rho_{q1} +1)$.
Therefore
$$
(\mu_i - \theta_{qi}^*)^2 \leq (\bmu - \btheta_q^*)'(\bmu-\btheta_q^*)= (\btheta_q^*)' D^2 \btheta_q^* \leq \frac{(\btheta_q^*)' \btheta_q^*}{(n\tau\rho_{q1} +1)^2}.
$$
In \eqref{eq:proof_ncchisq_univ1} we may set $u=\sqrt{(\btheta_q^*)' \btheta_q^*}/(n\tau\rho_{q1} +1)$.
Hence
$$
\frac{\mu_i^2}{\sigma_{ii}\phi^*} \in
\left[ \frac{(\theta_{qi}^*)^2}{\tilde{\sigma}_{ii}\phi^*} \left(1 - \frac{2u}{|\theta_{qi}^*|}\right) \left(1 + \frac{1}{(n\tau)^2\rho_{q1}^2} \right),
\frac{(\theta_{qi}^*)^2}{\tilde{\sigma}_{ii} \phi^*} \left( 1 + \frac{2u}{|\theta_{qi}^*|} + \frac{u^2}{(\theta_{qi}^*)^2} \right) \left(1+ \frac{1}{(n\tau)^2\rho_{qp_q}^2} \right)  \right].
$$
Recall that $\tilde{\sigma}_{ii}= ({\bf x}_{qj}' (1 - X_m (X_m'X_m)^{-1}X_m') {\bf x}_{qj})^{-1}$, hence
$(\tilde{\theta_{qi}}^*)^2/(\tilde{\sigma}_{ii}\phi^*)= \lambda_{qi}$ as desired.

\section{Proof of Lemma \ref{LEM:BAYES_UNIVFTEST_INEQ}}

Let $\mu_{qi}$ and $\sigma_{ii}$ be as defined in Lemma \ref{LEM:NCCHISQ_BAYES_UNIV}
and $s_q= \by'\by - \by'X_q(X_q'X_q)^{-1} X_q'\by$ the residual sum of squares under the least squares estimate of $\btheta_q$.
The target probability is
\begin{align}
P \left( \frac{\hat{\theta}_{qi}^2/\sigma_{ii}}{s_q/(n-p_q)} \frac{s_q}{\tilde{s}_q} < \frac{1}{\sigma_{ii} h_n} \right) \leq
P \left( \frac{\hat{\theta}_{qi}^2/\sigma_{ii}}{s_q/(n-p_q)} < \frac{2}{\sigma_{ii} h_n} \right)
+ P \left( \frac{s_q}{\tilde{s}_q} < \frac{1}{2} \right).
\label{eq:proof_bayes_univftest1}
\end{align}
From Lemma \ref{LEM:FSTAT_NORMALPRIOR_BOUND} we have that
$s_q/\tilde{s}_q \geq \left(1+ (s_0-s_q)/(s_q (1+n\tau\rho_{q,p_q}))  \right)^{-1}$,
hence the second term in \eqref{eq:proof_bayes_univftest1} is
\begin{align}
P \left( 1+ \frac{s_0-s_q}{s_q (1+n\tau\rho_{q,p_q})} > 2 \right)=
P \left( p_q F_{q0}  > (n-p_q) (1+n\tau\rho_{q,p_q}) \right)
\label{eq:proof_bayes_univftest2}
\end{align}
where $F_{q0}= [(s_0-s_q)/p_q]/[s_q/(n-p_q)]$ is the F-test statistic to compare $M_q$ with the model that includes no covariates.
To bound this tail probability, in Corollary \ref{COR:INEQ_F_RIGHT}(ii)
set $w= (n-p_q) (1+n\tau\rho_{q,p_q})/p_q$, $\nu_1=p_q$, $\nu_2=n$
and note that $w \gg \nu_2 \gg  (n-p_q)^\gamma$ for any $\gamma <1$, hence \eqref{eq:proof_bayes_univftest2}
is $\ll e^{-(n-p_t)^\gamma}$.

To bound the first term in \eqref{eq:proof_bayes_univftest1},
note that $[\hat{\theta}_{qi}^2/\sigma_{ii}]/[s_q/(n-p_q)]$
is the ratio of a non-central chi-square random variable and central chi-square,
each divided by its respective degrees of freedom (1 and $n-p_q$ respectively).
Hence we may apply the tail inequality in Lemma \ref{LEM:INEQ_NCF_LEFT}, setting  $w= 2/(\sigma_{ii}h_n)$.
Lemma \ref{LEM:NCCHISQ_BAYES_UNIV} guarantees that
$$\lim_{n \rightarrow \infty} \frac{\mu_{qi}^2/(\sigma_{ii}\phi^*)}{\lambda_{qi}}
=\lim_{n \rightarrow \infty} \frac{\mu_{qi}^2/(\sigma_{ii}\phi^*)}{(\theta_{qi})^2/(\tilde{\sigma}_{ii}\phi^*)}
=1,
$$
hence the assumption that $h_n \gg \phi^*/(\theta_{qi}^*)^2$ implies that $w \ll \lambda_{qi}$,
{\it i.e.} that the critical point is below the non-centrality parameter $\lambda_{qi}$.
Then by Lemma \ref{LEM:INEQ_NCF_LEFT} the first term in \eqref{eq:proof_bayes_univftest1} is
$$
\ll e^{-\gamma \lambda_{qi}/2} + e^{-(n-p_q)(s-1-\log(s))/2}
\leq e^{-\gamma \lambda_{qi}/2} + e^{-(n-p_q)c}
$$
for any fixed $s < \lambda_{qi}/w$, $\gamma \in (0,1)$ and $c=s-1-\log(s)$ is an arbitrarily large constant.
Summarizing, \eqref{eq:proof_bayes_univftest1} is $e^{-\gamma \lambda_{qi}/2} + e^{-(n-p_q)^\gamma/2}$,
proving the first part of the lemma.
As a direct implication
\begin{align}
&P \left( \prod_{i=1}^{p_q} \frac{\hat{\theta}_{qi}^2 (n-p_q)}{\tilde{s}_q} < \frac{1}{h_n^{p_q}} \right) \leq
  P \left( \min_{i=1,\ldots,p_q} \frac{\hat{\theta}_{qi}^2 (n-p_q)}{\tilde{s}_q} < \frac{1}{h_n} \right)
  \nonumber \\
&  \leq \sum_{i=1,\ldots,p_q} P \left(  \frac{\hat{\theta}_{qi}^2 (n-p_q)}{\tilde{s}_q} < \frac{1}{h_n} \right),
\nonumber
\end{align}
which by Lemma \ref{LEM:SUMO} is $\ll p_q (e^{-\gamma \lambda_{qi}/2} + e^{-(n-p_q)^\gamma})$,
concluding the proof of Part (i).

Regarding Parts (ii) and (iii), by Corollary \ref{COR:INEQ_F_RIGHT}
\begin{align}
  P \left( \frac{\hat{\theta}_{qi}^2/\sigma_{ii}}{\tilde{s}_q/(n-p_q)} > \frac{\tilde{h}_n}{\sigma_{ii}} \right)
  \ll  e^{-\gamma \tilde{h}_n/(2\sigma_{ii})} + e^{-(n-p_q)^\gamma/2}
  \ll  e^{-\gamma \lambda_{qi}/2} + e^{-(n-p_q)^\gamma/2},
\nonumber
\end{align}
as long as $\tilde{h}_n/\sigma_{ii}$ is asymptotically larger than the non-centrality parameter $\lambda_{qi}$.
For the case $\theta_{qi}^*=0$ this simply requires that $\tilde{h}_n \gg \sigma_{ii}$, which holds by assumption.
For the case $\theta_{qi}^*\neq 0$ this requires that $\tilde{h}_n/\sigma_{ii} \gg (\theta_{qi}^*)^2/(\phi^* \sigma_{ii})$,
that is $\tilde{h}_n \gg (\theta_{qi}^*)^2/\phi^*$, which holds by assumption.
Note that we can equivalently use $\sigma_{ii}$ or $\tilde{\sigma}_{ii}$,
since $\lim_{n\rightarrow\infty} \sigma_{ii}/\tilde{\sigma}_{ii}=1$ by Lemma \ref{LEM:NCCHISQ_BAYES_UNIV}.

\section{Proof of Lemma \ref{LEM:FSTAT_NORMALPRIOR_BOUND}}

Note that $s_m$ is the sum of squared residuals under the least-squares estimate,
hence $\tilde{s}_m > s_m$ and
\begin{align}
\frac{p_m-p_t}{n-p_m} \tilde{F}_{mt}=
\frac{\tilde{s}_t-\tilde{s}_m}{\tilde{s}_m} <
\frac{\tilde{s}_t-s_m}{s_m} =
\frac{s_t-s_m}{s_m} + \frac{\tilde{s}_t - s_t}{s_m}
= \frac{p_m-p_t}{n-p_m} F_{mt} + \frac{\tilde{s}_t - s_t}{s_m}.
\label{eq:bound_f_normalprior}
\end{align}

To prove Lemma \ref{LEM:FSTAT_NORMALPRIOR_BOUND} it suffices to show that we can upper-bound $\tilde{s}_t - s_t=$
\begin{align}
\by'X_t(X_t'X_t)^{-1} X_t'\by - \by'X_t (X_t'X_t+ (n\tau)^{-1} V_t^{-1})^{-1} X_t'\by
\leq \frac{\by'X_t(X_t'X_t)^{-1} X_t'\by}{1+n\tau \rho_{tp_t}}=
\frac{s_0 - s_t}{1+n\tau \rho_{tp_t}},
\label{eq:bound_f_normalprior_step2}
\end{align}
where $s_0= \by'\by$ is the sum of squared residuals under the model with no covariates.
Then $\tilde{s}_t \leq s_t + (s_0-s_t)/(1+n\tau \rho_{tp_t})$, proving the first part of the lemma,
and further \eqref{eq:bound_f_normalprior} is
$$
\leq \frac{p_m-p_t}{n-p_m} F_{mt}
+ \frac{1}{1+n\tau \rho_{tp_t}} \left( \frac{s_0-s_m}{s_m} - \frac{s_t-s_m}{s_m} \right)
=
\frac{n\tau \rho_{tp_t}}{1+n\tau\rho_{tp_t}} \frac{p_m-p_t}{n-p_m} F_{mt}
+ \frac{1}{1+n\tau \rho_{tp_t}} \frac{F_{m0} p_m}{n-p_m},
$$
proving the second part of the lemma.

To prove \eqref{eq:bound_f_normalprior_step2}, the left-hand side is equal to
\begin{align}
&\by'X_t (X_t'X_t)^{-\frac{1}{2}}
\left[ I - (X_t'X_t)^{\frac{1}{2}} (X_t'X_t+ (n\tau)^{-1} V_t^{-1})^{-1} (X_t'X_t)^{\frac{1}{2}} \right]
  (X_t'X_t)^{-\frac{1}{2}} X_t'\by
  \nonumber \\
&< (1-l) \by'X_t (X_t'X_t)^{-1} X_t' \by,
\nonumber
\end{align}
where $l$ is the smallest non-zero eigenvalue of
$(X_t'X_t)^{\frac{1}{2}} (X_t'X_t+ (n\tau)^{-1}V_t^{-1})^{-1} (X_t'X_t)^{\frac{1}{2}}$.
An elementary argument below shows that $l= n\tau \rho_{tp_t}/(1+n\tau\rho_{tp_t})$,
hence $1-l= 1/(1+n\tau \rho_{tp_t})$ as desired.

Consider arbitrary invertible square matrices $A,B$ of equal dimension and let $\rho$ be an eigenvalue of $B^{-\frac{1}{2}} A B^{-\frac{1}{2}}$.
Then $B^{-\frac{1}{2}} A B^{-\frac{1}{2}} {\bf v}= \rho {\bf v}$ where ${\bf v}$ is the corresponding eigenvector,
thus $A B^{-1} B^{\frac{1}{2}} {\bf v}= \rho B^{\frac{1}{2}} {\bf v}$ and $\rho$ is an eigenvalue of $A B^{-1}$ with eigenvector $B^{\frac{1}{2}} {\bf v}$.
Consequently $(X_t'X_t)^{-\frac{1}{2}} (X_t'X_t+ (n\tau)^{-1}V_t^{-1}) (X_t'X_t)^{-\frac{1}{2}}$
has the same eigenvalues as $I+ (n\tau)^{-1}V_t^{-1} (X_t'X_t)^{-1}$,
and trivially the eigenvalues of the latter are $1+\frac{1}{n\tau \rho_{rj}}$ for $j=1,\ldots,p_t$.
Therefore the eigenvalues of $(X_t'X_t)^{\frac{1}{2}} (X_t'X_t+ (n\tau)^{-1}V_t^{-1})^{-1} (X_t'X_t)^{\frac{1}{2}}$
are $n\tau \rho_{rj}/(1+n\tau\rho_{rj})$.

\section{Proof of Lemma \ref{LEM:SUMO}}

We need to show that
\begin{align}
\lim_{n \to \infty} \frac{\sum_{l=l_n^{(0)}}^{l_n^{(1)}} \sigma_{l,n}}{\sum_{l=l_n^{(0)}}^{l_n^{(1)}} g_n^{(l)}}=
\lim_{n \to \infty}
\sum_{l=l_n^{(0)}}^{l_n^{(1)}} \frac{\sigma_{l,n}}{ \bar{g}_n (l_n^{(1)}-l_n^{(0)}+1)} < \infty.
\label{eq:sumo_proofgoal}
\end{align}
To see that this series converges we shall show that it increases at a rate slower than $1/(l_n^{(1)}-l_n^{(0)}+1)$,
i.e. the inverse of its number of terms.
Specifically, define $m=l_n^{(1)}-l_n^{(0)}+1$ and $s_{l,m}= \sigma_{l,n} / (m\bar{g}_n)$,
the increase in the series between $m$ and $m+1$ is
\begin{align}
\sum_{l=1}^{m+1} s_{l,m+1} - \sum_{l=1}^{m} s_{l,m}=
s_{m+1,m+1} + \sum_{l=1}^{m} (s_{l,m+1} - s_{l,m}).
\label{eq:sumo_seriesincrease}
\end{align}
If we can show that this increase is $\ll 1/m$ then by the ratio test
it follows that \eqref{eq:sumo_proofgoal} is a convergent series.
Regarding the first term on the right hand side of \eqref{eq:sumo_seriesincrease},
\begin{align}
s_{m+1,m+1} \ll \frac{1}{m} \Longleftrightarrow
\frac{\sigma_{l_n^{(1)},n}}{ (l_n^{(1)}-l_n^{(0)}+1) \bar{g}_n} \ll \frac{1}{l_n^{(1)}-l_n^{(0)}+1} \Longleftrightarrow
\sigma_{l_n^{(1)},n} \ll \bar{g}_n,
\nonumber
\end{align}
which holds by Assumption (i).
Regarding the second term in \eqref{eq:sumo_seriesincrease},
\begin{align}
\sum_{l=1}^{m} s_{l,m+1} - s_{l,m} \ll \frac{1}{m} \Longleftrightarrow
\sum_{l=l_n^{(0)}}^{l_n^{(1)}} \frac{\sigma_{l,n+1}}{(l_{n+1}^{(1)}-l_{n+1}^{(0)}+1) \bar{g}_{n+1}}
- \frac{\sigma_{l,n}}{(l_n^{(1)}-l_n^{(0)}+1) \bar{g}_n} \ll \frac{1}{(l_n^{(1)}-l_n^{(0)}+1)}.
\nonumber
\end{align}
The right-hand side holds since $l_n^{(1)}-l_n^{(0)}$ is non-decreasing by assumption, hence
\begin{align}
\sum_{l=l_n^{(0)}}^{l_n^{(1)}} \frac{\sigma_{l,n+1}}{(l_{n+1}^{(1)}-l_{n+1}^{(0)}+1) \bar{g}_{n+1}}
- \frac{\sigma_{l,n}}{(l_n^{(1)}-l_n^{(0)}+1) \bar{g}_n}
\leq \frac{1}{l_n^{(1)}-l_n^{(0)}+1} \sum_{l=l_n^{(0)}}^{l_n^{(1)}} \frac{\sigma_{l,n+1}}{\bar{g}_{n+1}} - \frac{\sigma_{l,n}}{\bar{g}_n}
\nonumber
\end{align}
and by Assumption (ii)
\begin{align}
\lim_{n \rightarrow \infty} \sum_{l=l_n^{(0)}}^{l_n^{(1)}} \frac{\sigma_{l,n+1}}{\bar{g}_{n+1}} - \frac{\sigma_{l,n}}{\bar{g}_n} \leq 0.
\nonumber
\end{align}

\section{Proof of Corollary \ref{COR:SUMO}}

The result is an immediately application of Lemma \ref{LEM:SUMO} to the case where $g_n^{(l)}$ does not depend on $l$.
Specifically, in Lemma \ref{LEM:SUMO} set $\bar{g}_n=b_n$
and note that Assumption (i) in Lemma \ref{LEM:SUMO} is satisfied since $\mu_{|A_n|,n} \ll b_n$.
Also Assumption (ii) in Lemma \ref{LEM:SUMO} is equivalent to the assumption made in the statement of Corollary \ref{COR:SUMO} that
$$
\lim_{n \rightarrow \infty} \sum_{k \in A_n} \frac{\mu_{k,n+1}}{b_{n+1}} - \frac{\mu_{k,n}}{b_n} \leq 0.
$$

\section{Proof of Lemma \ref{LEM:SUMO_SIZE}}

The result is a particular case of Lemma \ref{LEM:SUMO}.
Specifically, to prove Lemma \ref{LEM:SUMO_SIZE}(i)
in Lemma \ref{LEM:SUMO} set $A_{l,n}=S_l$, $l_n^{(0)}=p_t$, $l_n^{(1)}=\bar{p}$, $\bar{g}_n=b_n$ and $g_n^{(l)}= a_n^{(l)} |S_l|$,
then the desired result follows.

Analogously, to prove Lemma \ref{LEM:SUMO_SIZE}(ii)
in Lemma \ref{LEM:SUMO} set $A_{l,n}=S_l^c$, $l_n^{(0)}=0$, $l_n^{(1)}=\bar{p}$, $\bar{g}_n=\tilde{b}_n$ and $g_n^{(l)}= \tilde{a}_n^{(l)} |S_l^c|$.

\section{Proof of Lemma \ref{LEM:SUMSPUR_ZELLNERKNOWN}}

We seek to bound
\begin{align}
  \sum_{l=p_t+1}^{\bar{p}} {l \choose p_t} \frac{\left[ (l-p_t) \log((n\tau)^{1/2} (p-p_t)) \right]^{\frac{l-p_t}{2}+1}}{(n\tau)^{\frac{l-p_t}{2}}}
  \nonumber \\
  =  \log((n\tau)^{1/2} (p-p_t))  \sum_{l=p_t+1}^{\bar{p}} {l \choose p_t} \left[ \frac{ (l-p_t)^{1+\frac{2}{l-p_t}} \log((n\tau)^{1/2} (p-p_t))}{n\tau}\right]^{\frac{l-p_t}{2}}.
  \label{eq:proofsumspur_zellknown1}
\end{align}

Let $l_0>p_t$ be an arbitrary fixed integer. The sum in \eqref{eq:proofsumspur_zellknown1} can be split into two sums over $l \leq l_0$ and $l > l_0$.
Let $a>1$ be a fixed constant, then there exists fixed $l_0$ such that
for $l>l_0$ it holds that $(l-p_t)^{1+2/(l-p_t)}< (l-p_t)^a \leq (\bar{p}-p_t)^a$.
Also note that, as $n \rightarrow \infty$, for any fixed $l\leq l_0$ we have
$\max_{l\leq l_0}(l-p_t)^{1+2/(l-p_t)} /(n\tau) \asymp 1/(n\tau) \preceq (\bar{p}-p_t)^a/(n\tau)$.
Hence we have that \eqref{eq:proofsumspur_zellknown1} is
\begin{align}
  \preceq \log((n\tau)^{\frac{1}{2}} (p-p_t)) \sum_{l=p_t+1}^{\bar{p}}
  {l \choose p_t} \left\{  \frac{(\bar{p}-p_t)^a \log((n\tau)^{\frac{1}{2}} (p-p_t))}{(n\tau)}  \right\}^{\frac{l-p_t}{2}}
  \nonumber \\
  < \log((n\tau)^{\frac{1}{2}} (p-p_t))
  \left( \left[ 1 - \frac{(\bar{p}-p_t)^{a/2} \log^{1/2}((n\tau)^{1/2} (p-p_t)) }{(n\tau)^{1/2}}  \right]^{-(p_t+1)} -1 \right),
  \label{eq:proofsumspur_zellknown2}
\end{align}
the right-hand side following from the Binomial coefficient's ordinary generating function.
If $(n\tau)^{1/2} \gg (p_t+1) (\bar{p}-p_t)^{a/2} \log^{3/2}((n\tau)^{1/2} (p-p_t))$,
which is guaranteed by assumption,
from the definition of the exponential function \eqref{eq:proofsumspur_zellknown2} is $\asymp$
\begin{align}
& \log((n\tau)^{\frac{1}{2}} (p-p_t))
  \left( e^{(p_t+1)  (\bar{p}-p_t)^{a/2} \log^{1/2}((n\tau)^{1/2} (p-p_t))/(n\tau)^{1/2}}  -1 \right)
  \nonumber \\
&  \asymp \frac{(p_t+1) (\bar{p}-p_t)^{a/2} \log^{3/2}((n\tau)^{1/2} (p-p_t))}{(n\tau)^{1/2}},
\nonumber
\end{align}
the right-hand side following from $\lim_{z \rightarrow 0} (e^z-1)/z=1$.
This proves the desired result.

\section{Proof of Lemma \ref{LEM:SUMSPUR_ZELLNERUNKNOWN}}

Let $g=(n\tau)^{(p_k-p_t)/2} p(M_t)/p(M_k)$
where $p(M_t)/p(M_k)= p_t!(p-p_t)!/(p_k!(p-p_k)!) < (p-p_t)^{p_k-p_t}$.
Then $\log(g) < (p_k-p_t) \log(n\tau (p-p_t)) \ll n-\bar{p}$ by assumption.
In this case for $m \in S$ we saw in Section \ref{ssec:zellner_unknown_spur} that
$$
E_{f^*}(Z_k) < \int_0^1 P((p_m-p_t) F_{mt} > b_n(u)) du
$$
where $b_n(u)= 2 [(n-p_k)/(n+a_\phi)] \log ( (1+n\tau)^{\frac{p_k-p_t}{2}} [p(M_t)/p(M_k)]/(1/u-1) )$.
Then Lemma \ref{LEM:FTAIL_INTBOUND}(i) gives $E_{f^*}(Z_k) \preceq g^{2 \sqrt{\log(g)/(n-p_k)}}/g$.
Applying Lemma \ref{LEM:SUMO_SIZE}
\begin{align}
  E_{f^*}(P(S \mid \by)) \preceq \sum_{l=p_t+1}^{\bar{p}} {p - p_t \choose l-p_t}
  \frac{[(n\tau)^{\frac{l-p_t}{2}} (p-p_t)^{l-p_t}]^{2 \sqrt{(l-p_t)[\log((n\tau)^{1/2} (p-p_t))]/(n-l)}} }
  {(n\tau)^{\frac{l-p_t}{2}}}
  \nonumber \\
  < \sum_{l=p_t+1}^{\bar{p}} {l \choose p_t}
  \left( \frac{[(n\tau)^{1/2} (p-p_t)]^{2c_n}}{(n\tau)^{1/2}} \right)^{l-p_t}
  < \left[ 1 - \frac{[(n\tau)^{1/2} (p-p_t)]^{2c_n}}{(n\tau)^{1/2}} \right]^{-(p_t+1)}  -1
  \label{eq:proofsumspur_zellunknown1}
\end{align}
where $c_n= \sqrt{ (\bar{p}-p_t) [\log((n\tau)^{1/2} (p-p_t))]/(n-\bar{p})}$
and the right-hand side follows from the Binomial coefficient's ordinary generating function.

If $[(n\tau)^{1/2} (p-p_t)]^{2c_n} \ll (n\tau)^{1/2}$, which holds by assumption, then
\eqref{eq:proofsumspur_zellunknown1} is
\begin{align}
&  \asymp \exp \left\{  \frac{(p_t+1) [(n\tau)^{1/2} (p-p_t)]^{2c_n}}{(n\tau)^{1/2}} \right\} -1
  \asymp \frac{(p_t+1) [(n\tau)^{1/2} (p-p_t)]^{2c_n}}{(n\tau)^{1/2}}
  \nonumber \\
  &=   \frac{(p_t+1)}{(n\tau)^{1/2}} e^{2 [\log^{3/2} ((n\tau)^{1/2} (p-p_t))] \sqrt{\frac{p-p_t}{n-\bar{p}}  }},
\nonumber
\end{align}
since $\lim_{z \rightarrow 0} (e^z-1)/z=1$, as we wished to prove.

\section{Proof of Corollary \ref{COR:INEQ_F_RIGHT}}

The result is obtained by setting
$s=1 + \frac{w}{\nu_2} (1-\sqrt{1+2\nu_2/w})$ in Lemma \ref{LEM:INEQ_NCF_RIGHT} and applying basic inequalities.
We first check that when $w \in ((\nu_1+\lambda)/(2-\sqrt{3}),\nu_2)$ then $s \in ((\nu_1+\lambda)/w,1)$ as required by Lemma \ref{LEM:INEQ_NCF_RIGHT}.
Clearly, $s<1$.
Further note that $s$ is decreasing in $w/\nu_2 \leq 1$, thus plugging $w/\nu_2=1$ into the expression of $s$ gives
that $s \geq 2-\sqrt{3}$.
Hence $2-\sqrt{3} > (\nu_1+\lambda)/w $ implies that $s >\nu_1/w$, as desired.

To motivate $s$ asymptotically as $w \rightarrow \infty$
we set $s$ such that the two terms for $P(\nu_1W>w)$ in Lemma \ref{LEM:INEQ_NCF_RIGHT}
under $\lambda=0$ are approximately equal.
The leading factor in the first term is $e^{-ws/2}$
and the second term is $e^{-\nu_2(s-1-\log(s))/2}$,
hence we seek $s$ such that $ws= \nu_2(s-1-\log(s))$, i.e. $(1+\log(s))/s= 1-w/\nu_2$.
To find such $s$ suppose that $\nu_2$ is a function of $w$
such that $\lim_{w \rightarrow \infty} w/\nu_2 =0$,
then $\lim_{w \rightarrow \infty} (1+\log(s))/s=1$ and hence $\lim_{w \rightarrow \infty} s=1$.
Plugging in the second order Taylor expansion $\log(s) \approx s-1 - (s-1)^2/2$ around $s=1$
gives $\nu_2(s-1-\log(s)) \approx \nu_2(s-1)^2/2$.
After simple algebra solving the second order equation
$ws= \nu_2(s-1)^2/2$ under the restriction that $s<1$ gives
$s= 1+ \frac{w}{\nu_2} \left( 1 - \sqrt{1 + 2\nu_2/w}  \right)$, as desired.

To complete the proof we plug our choice of $s$ into Lemma \ref{LEM:INEQ_NCF_RIGHT}.
Consider first Lemma \ref{LEM:INEQ_NCF_RIGHT}(i).
Since $s<1$ the first term is
$$
\leq \left( \frac{ew}{\nu_1} \right)^{\frac{\nu_1}{2}} e^{-ws/2}
= \left( \frac{ew}{\nu_1} \right)^{\frac{\nu_1}{2}} e^{-\frac{w}{2}\left(1+ \frac{w}{\nu_2} \left( 1 - \sqrt{1 + 2\nu_2/w}  \right)\right)},
$$
where considering that $w<\nu_2$ and that $z(1-\sqrt{1+2/z}) \geq -\sqrt{2z}$ for any $z \in (0,1)$,
$$
e^{-\frac{w}{2} \left( 1 + \frac{w}{\nu_2} \left( 1-\sqrt{1+2\nu_2/w} \right) \right)}
\leq e^{-\frac{w}{2} (1 - \sqrt{2w/\nu_2})}.
$$
Regarding the second term, since $s-1-\log(s) > (s-1)^2/2$ for all $s \in (0,1)$,
$$
e^{-\frac{\nu_2}{2}(s-1-\log(s))} <
e^{-\frac{\nu_2}{2} \frac{(s-1)^2}{2}}
= e^{-\frac{\nu_2}{4}\frac{w^2}{\nu_2^2} \left( 1 - \sqrt{1 + 2\nu_2/w}  \right)^2}
< e^{-w/2},
$$
where the last inequality follows from
$z^2(1-\sqrt{1+2/z})^2 > 2z$ for all $z \in (0,1)$ and plugging in $z=w/\nu_2$.

Consider now Part (ii). Then
$$
P(\nu_1W>w) \leq e^{-\frac{ws}{2} \left( 1- \sqrt{\frac{\lambda}{ws}} \right)^2 } \left( \frac{ws}{\lambda} \right)^{\frac{\nu_1}{4}}
+ (es)^{\nu_2/2} e^{- \frac{s\nu_2}{2}}.
$$
The second term is identical to the $\lambda=0$ case, hence $<e^{-w/2}$.
Regarding the first term we just showed that
$s=1 + \frac{w}{\nu_2} (1-\sqrt{1+2\nu_2/w})> 1 - \sqrt{2w/\nu_2}$, hence the first term is
$$
< \left( \frac{w}{\lambda} \right)^{\frac{\nu_1}{4}}
e^{-\frac{w}{2} (1 - \sqrt{2w/\nu_2}) \left( 1 - \sqrt{\frac{\lambda}{w (1 - \sqrt{2w/\nu_2})}} \right)^2},
$$
as we wished to prove.

\section{Proof of Lemma \ref{LEM:INTBOUND_EXPTAILS}}

Since $g/(1/u-1)$ is increasing in $u$,
$$
  \int_{\underline{u}}^{\bar{u}} P\left(W > d \log \left( \frac{g}{1/u-1} \right)   \right) du <
 \left[d \log \left(\frac{g}{1/\bar{u}-1}\right) \right]^c \frac{b}{g^{ld}}
\int_{\underline{u}}^{\bar{u}} (1/u -1)^{ld} du.
$$
The case $ld=1$ follows trivially from
$$
\int_{\underline{u}}^{\bar{u}} (1/u -1)^{ld} du=
\log(\bar{u}) - \log(\underline{u}) - (\bar{u} - \underline{u})
< \log(1/\underline{u}),
$$
since $0 < \underline{u} < \bar{u} < 1$ by assumption.

For the case $ld \neq 1$, applying the change of variables $v=1/u-1$ gives
$$
\int_{\underline{u}}^{\bar{u}} (1/u -1)^{ld} du=
\int_{1/\bar{u}-1}^{1/\underline{u}-1} \frac{v^{ld}}{(v+1)^2} dv <
\int_{1/\bar{u}-1}^{1/\underline{u}-1} v^{ld-2} dv=
\frac{1}{1-ld} \left[
\left( \frac{1}{\bar{u}} -1 \right)^{ld-1}
- \left( \frac{1}{\underline{u}} -1 \right)^{ld-1}
 \right].
$$
If $ld < 1$ then $1-ld>0$ and the right hand side is
$< (\bar{u}/(1-\bar{u}))^{1-ld} / (1-ld)$, as desired.
If $ld > 1$ then $1-ld<0$ and the right hand side is
$< (1/\underline{u} -1)^{ld-1}/ (ld-1)$.

\section{Proof of Lemma \ref{LEM:INTBOUND_POLYTAILS}}

Clearly
\begin{align}
  \int_{\underline{u}}^{\bar{u}} P\left(W > d \log \left( \frac{g}{1/u-1} \right)   \right) du <
  \int_{\underline{u}}^{\bar{u}} \frac{b}{d^c \left[\log \left( \frac{g}{1/u-1} \right) \right]^c} du.
\label{eq:intpolytail1}
\end{align}

A first trivial bound is found by noting that the integrand is decreasing in $u$,
hence \eqref{eq:intpolytail1} is
\begin{align}
< \frac{b}{d^c \left[\log \left( \frac{g}{1/\underline{u}-1} \right) \right]^c}.
\label{eq:polytail_bound1}
\end{align}
Applying the change of variables $v= \log(g) - \log(1/u-1)$ gives that the right hand side
of \eqref{eq:intpolytail1} is

\begin{align}
=  \int_{\underline{v}}^{\bar{v}} \frac{b}{d^c v^c}
  \frac{ge^{-v}}{(1+ge^{-v})^2} dv.
\label{eq:intpolytail2}
\end{align}
where $\underline{v}=\log\left(\frac{g}{1/\underline{u}-1}\right)$
and $\bar{v}=\log\left(\frac{g}{1/\bar{u}-1}\right)$.
Notice that $ge^{-v}/(1+ge^{-v})$ is decreasing in $v$ and
$1/(1+ge^{-v})$ is increasing in $v$,
hence $ge^{-v}/(1+ge^{-v})^2 \leq ge^{-\underline{v}} /[(1+ ge^{-\underline{v}}) (1+ge^{-\bar{v}})]$.
Thus, \eqref{eq:intpolytail2} is
\begin{align}
< \frac{b}{d^c} \frac{ge^{-\underline{v}}}{(1+ ge^{-\underline{v}}) (1+ge^{-\bar{v}})}
\int_{\underline{v}}^{\bar{v}} \frac{1}{v^c} dv=
\frac{b}{d^c} \frac{\bar{u} (1 - \underline{u})}{c-1} \left(
\frac{1}{\log^{c-1}\left(\frac{g}{1/\underline{u}-1}\right)}
- \frac{1}{\log^{c-1}\left(\frac{g}{1/\bar{u}-1}\right)}
 \right)
\nonumber \\
<
\frac{b}{d^c (c-1)}
\frac{1}{\log^{c-1}\left(\frac{g}{1/\underline{u}-1}\right)}.
\label{eq:polytail_bound2}
\end{align}
Both \eqref{eq:polytail_bound1} and \eqref{eq:polytail_bound2} give valid bounds for the target
integral. Clearly, \eqref{eq:polytail_bound2} is smaller than \eqref{eq:polytail_bound1} if and only if
$$
c-1 > \log \left( \frac{g}{1/\underline{u}-1} \right).
$$

\section{Proof of Lemma \ref{LEM:CHISQTAIL_INTBOUND}}

Chernoff's bound for chi-square tails (Lemma \ref{LEM:INEQ_CHISQ}) gives that
\begin{align}
P \left( W > d \log \left( \frac{g}{1/u -1}   \right) \right) \leq
\left( \frac{1/u-1}{g} \right)^{d/2} \left[ \frac{de}{\nu} \log \left( \frac{g}{1/u-1} \right) \right]^{\frac{\nu}{2}}
\nonumber
\end{align}
for any $u$ such that $d \log (g/(1/u -1) > \nu$, that is for $u>\underline{u}$ where we define $\underline{u}=(1+ge^{-\nu/d})^{-1}$.
Further define $\bar{u}=1-\underline{u}$.

These exponential tails are of the form required by Lemma \ref{LEM:INTBOUND_EXPTAILS}.
To prove Part (i), set $l=1/2$ and $c=\nu/2$ in Lemma \ref{LEM:INTBOUND_EXPTAILS}, when $d=2$ then $ld=1$ and hence
\begin{align}
\int_0^1 P \left( W > 2 \log \left( \frac{g}{1/u -1}   \right) \right) du <
2\underline{u} + \frac{1}{g} \log(1/\underline{u}) \left[ \frac{2e}{\nu} \log \left( \frac{g}{1/\bar{u}-1} \right) \right]^{\frac{\nu}{2}}
\nonumber \\
=\frac{2}{1+ge^{-\nu/2}} + \frac{1}{g} \log(1+ge^{-\nu/2}) \left[ \frac{4e}{\nu} \log \left( \frac{g}{e^{\nu/4}} \right) \right]^{\frac{\nu}{2}}
\label{eq:proof_chisqtail_intbound1}
\end{align}
Given that $\log(g) \gg  \nu$ by assumption, the first term in \eqref{eq:proof_chisqtail_intbound1} is of a smaller order than the second term.
The second term in \eqref{eq:proof_chisqtail_intbound1} is $\asymp \frac{1}{g} (4e/\nu)^{\nu/2} \log^{\nu/2+1} (g/e^{\nu/4})$, giving the desired result.
Further, since $\log(g) \gg \frac{\nu}{2} \log\log(1+g)$ by assumption we have that $\frac{1}{g} (4e/\nu)^{\nu/2} \log^{\nu/2+1} (g/e^{\nu/4}) \ll 1/g^{\alpha}$
as $g \rightarrow \infty$ for any fixed $\alpha \in (0,1)$, as we wished to prove.
The case $d >2$ follows trivially since $P( W > d \log ( g/(1/u-1) ) )$ decreases in $d$.
For precision's sake, when $d>2$ Lemma \ref{LEM:INTBOUND_EXPTAILS} and trivial algebra give the finite-sample bound
$$
\int_0^1 P \left( W > 2 \log \left( \frac{g}{1/u -1}   \right) \right) du <
\frac{2}{1+ge^{-\nu/d}} + \frac{1}{g} \left[ \frac{2de^{2/d}}{\nu} \log \left( \frac{g}{e^{\nu/(2d)}} \right) \right]^{\frac{\nu}{2}} \frac{1}{d/2-1}
\ll  1/g^\alpha.
$$

To prove Part (ii), set $l=1/2$ and $c=\nu/2$ in Lemma \ref{LEM:INTBOUND_EXPTAILS}. Since $d<2$ then $ld<1$ and hence
\begin{align}
\int_0^1 P \left( W > 2 \log \left( \frac{g}{1/u -1}   \right) \right) du <
2\underline{u} + \left( \frac{1}{g} \right)^{\frac{d}{2}} \left[ \frac{de}{\nu} \log \left( \frac{g}{1/\bar{u}-1} \right) \right]^{\frac{\nu}{2}}
\frac{1}{1-\frac{d}{2}} \left( \frac{\bar{u}}{1-\bar{u}} \right)^{1-\frac{d}{2}}.
\nonumber
\end{align}
Since $\bar{u}= 1 - \underline{u}$, it holds that $1/\bar{u}-1= \underline{u}/(1-\underline{u})>\underline{u}$ and that $\bar{u}/(1-\bar{u}) < 1/\underline{u}$.
Hence the right-hand side above is
\begin{align}
 < 2\underline{u} + \frac{1}{g^{\frac{d}{2}}} \left( \frac{1}{\underline{u}} \right)^{1-\frac{d}{2}}
\left[ \frac{de}{\nu} \log \left( \frac{g}{\underline{u}} \right) \right]^{\frac{\nu}{2}}
\frac{1}{1-\frac{d}{2}}.
\nonumber
\end{align}
To conclude, take $\underline{u}$ such that $\underline{u}= (1/\underline{u})^{1-\frac{d}{2}} g^{-d/2}$, that is $\underline{u}= g^{-(2-d/2)d/2}$, to obtain the expression 
\begin{align}
 \underline{u} \left( 2 + \left[ \frac{de}{\nu} \log( g )^{(2-d/2)d/2+1} \right]^{\frac{\nu}{2}}
\frac{1}{1-\frac{d}{2}} \right),
\nonumber
\end{align}
as we wished to prove.

\section{Proof of Lemma \ref{LEM:FTAIL_INTBOUND}}

We seek to bound
\begin{align}
\int_{0}^{1} P \left( \nu_1 F > d \log \left( \frac{g}{1/u-1} \right) \right) du \leq
\underline{u} + (1-\bar{u}) + \int_{\underline{u}}^{\bar{u}} P \left( \nu_1 F > d \log \left( \frac{g}{1/u-1} \right) \right) du,
\label{eq:proof_ftail_bound_nonspur2}
\end{align}
for suitably defined $\underline{u}<\bar{u}$.
The strategy is that when $\log(g) \ll \nu_2$
one may use the exponential bound in Corollary \ref{COR:INEQ_F_RIGHT} for the integrand in \eqref{eq:proof_ftail_bound_nonspur2},
then Lemma \ref{LEM:INTBOUND_EXPTAILS} to bound its integral.
For the case $\log^\gamma(g) \gg \nu_2$ for all fixed $\gamma < 1$,
one may use the polynomial bound in Lemma \ref{LEM:INEQ_F} for the integrand,
then Lemma \ref{LEM:INTBOUND_POLYTAILS} to bound \eqref{eq:proof_ftail_bound_nonspur2}.

Consider first the case where $\log(g) \ll \nu_2$.
Corollary \ref{COR:INEQ_F_RIGHT} gives
\begin{align}
  P\left(\nu_1 F > d \log \left( \frac{g}{1/u-1} \right)  \right) <
2 \exp\left\{- \frac{d}{2} \log\left(\frac{g}{1/u-1} \right) \left(1 - \sqrt{\frac{2 d \log(g/(1/u-1))}{\nu_2}} \right)\right\}
  \label{eq:proof_ftail_bound_nonspur1}
\end{align}
for any $d \log (g/(1/u-1)) \in (\nu_1/(2-\sqrt{3}), \nu_2)$.
Equivalently, in terms of $u$, \eqref{eq:proof_ftail_bound_nonspur1} holds for
$$
u \in \left( \frac{1}{1+g e^{-\nu_1 /(d(2-\sqrt{3}))} }, \frac{1}{1 + g e^{-\nu_2/d} } \right).
$$
Denote the left endpoint of this interval by $\underline{u}$, define $\bar{u}= 1 - \underline{u}$,
and for future reference note that $1/\bar{u}-1= e^{\nu_1/(d(2-\sqrt{3}))}/g$.
Clearly, for any $\omega \in \left(0,\sqrt{\frac{2d}{\nu_2} \log (g/(1-\bar{u}-1))}\right)$
Expression \eqref{eq:proof_ftail_bound_nonspur1} is $< 2 e^{-d \log(g/(1/u-1)) (1-\omega)}$,
which is of the form required by Lemma \ref{LEM:INTBOUND_EXPTAILS}.
Before applying the lemma we must check that \eqref{eq:proof_ftail_bound_nonspur1} holds for $u \in (\underline{u},\bar{u})$,
i.e. that $u < \bar{u}$ implies $d \log(g/(1/u-1)) < \nu_2$.
Since $\log(g/(1/u-1))$ is increasing in $u$, $u< \bar{u}$ implies that
$$
d \log \left( \frac{g}{1/u-1} \right) <
d \log \left( \frac{g}{1/\bar{u}-1} \right)=
2d \log \left( \frac{g}{e^{\nu_1/(2d(2-\sqrt{3}))}}   \right),
$$
hence $d \log(g/(1/u-1))<\nu_2$ for large enough $\nu_2$ and thus the bound in \eqref{eq:proof_ftail_bound_nonspur1} is valid for $u \in (\underline{u},\bar{u})$.

In Lemma \ref{LEM:INTBOUND_EXPTAILS} set $b=2$, $c=0$ and $l= (1-\omega)/2$, where $\omega$ is as defined above.
We focus on the case where $ld<1$, the $ld\geq 1$ case follows then trivially by observing the integrand decreases in $d$.
By Lemma \ref{LEM:INTBOUND_EXPTAILS} the integral in \eqref{eq:proof_ftail_bound_nonspur2} is
\begin{align}
  &\leq 2\underline{u}+ \frac{2 \left[ \bar{u}/(1-\bar{u}) \right]^{1- d(1-\omega)/2}}{g^{d(1-\omega)/2} (1- d(1-\omega)/2))}
    \nonumber \\
&  =\frac{2}{1+g e^{-\nu_1 /(d(2-\sqrt{3}))} } + \frac{2 e^{-\nu_1 (1-d(1-\omega)/2)/(d(2-\sqrt{3}))}}{\left[ 1- d(1-\omega)/2 \right] } \left( \frac{1}{g} \right)^{d(1-\omega) -1}
\nonumber \\
&\preceq \frac{1}{g} e^{\frac{\nu_1}{d(2-\sqrt{3})}}
+ e^{-\nu_1} \left( \frac{1}{g} \right)^{d-1 - d\sqrt{ \frac{4}{\nu_2} \log \left(\frac{g}{e^{\nu_1/2}} \right)  }}
\label{eq:proof_ftail_bound_nonspur3}
\end{align}
since $\omega < \sqrt{\frac{2d}{\nu_2} \log (g/(1-\bar{u}-1))}$,
$\nu_1/(d(2-\sqrt{3})) \ll \nu_1$ and $d>1$ by assumption.
In particular since $\log(g) \ll \nu_2$ by assumption
one can take any fixed $\omega \in (0,1)$ that is arbitrarily close to 0,
then \eqref{eq:proof_ftail_bound_nonspur3} is $\ll 1/g^\alpha$  as $\nu_2 \rightarrow \infty$
for any fixed $\alpha < d-1$, as we wished to prove.

Consider now the case $\log^\gamma(g) \gg \nu_2$ for all fixed $\gamma < 1$.
Lemma \ref{LEM:INEQ_F} gives that, for all $d \log(g/(1/u-1)) > \nu_1+\nu_2-6$,
\begin{align}
P \left( \nu_1 F > d \log \left( \frac{g}{1/u-1} \right) \right) <
\left( \frac{(\nu_1+\nu_2-6)^{\frac{\nu_1+\nu_2-5}{\nu_2-4}}}{d \log(g/(1/u-1))} \right)^{\frac{\nu_2}{2}-2}
\frac{1}{2^{\frac{\nu_1-1}{2}}} \frac{ae}{(\nu_2/2 -1)^{\frac{3}{2}}}.
\label{eq:proof_ftail_bound_nonspur4}
\end{align}
It is easy to check that $d \log(g/(1/u-1)) > \nu_1+\nu_2-6$ if and only if
$u > 1/ (1+ g/e^{\nu_1+\nu_2-6})$.
Define $\underline{u}= 1/ (1+ g^\omega/e^{\omega (\nu_1+\nu_2-6)})$
for any fixed $\omega <1$, and $\bar{u}=1-\underline{u}$.
Then clearly
$u>\underline{u}$ implies that $u > 1/ (1+ g/e^{\nu_1+\nu_2-6})$
and hence the bound in \eqref{eq:proof_ftail_bound_nonspur4} applies to all $u>\underline{u}$.

The polynomial tails in \eqref{eq:proof_ftail_bound_nonspur4}
are of the form required by Lemma \ref{LEM:INTBOUND_POLYTAILS} for $c=\nu_2/2-2$.
We apply Lemma \ref{LEM:INTBOUND_POLYTAILS}
(note that $g \geq 1/\underline{u}-1= g^\omega/e^{\omega (\nu_1+\nu_2-6)}$ as required by the lemma),
which gives that \eqref{eq:proof_ftail_bound_nonspur2} is
\begin{align}
< 2\underline{u} +
\left( \frac{(\nu_1+\nu_2-6)^{\frac{\nu_1+\nu_2-5}{\nu_2-4}}}
{d \log \left( \frac{g}{1/\underline{u} -1} \right)} \right)^{\frac{\nu_2}{2}-2}
\frac{1}{2^{\frac{\nu_1-1}{2}}} \frac{ae}{(\nu_2/2 -1)^{\frac{3}{2}}}
\nonumber \\
= \frac{2}{1+ g/e^{\nu_1+\nu_2-6}} +
ae \left( \frac{(\nu_1+\nu_2-6)^{\frac{\nu_1+\nu_2-5}{\nu_2-4}}}
{ d[(1-\omega) \log(g) + \omega (\nu_1+\nu_2-6)]} \right)^{\frac{\nu_2-4}{2}}
\nonumber \\
= \frac{2}{1+ g/e^{\nu_1+\nu_2-6}} + ae
\left( \frac{\nu_1+\nu_2-6}{ \left( d[(1-\omega) \log(g) + \omega (\nu_1+\nu_2-6)] \right)^{\frac{\nu_2-4}{\nu_1+\nu_2-5}} }  \right)^{\frac{\nu_1+\nu_2-5}{2}}
\label{eq:proof_ftail_bound_nonspur5}
\end{align}
where $(\nu_2-4)/(\nu_1+\nu_2-5)$ is arbitrarily close to 1 for large enough $\nu_2$
from the assumption that $\nu_1 \ll \nu_2$.
The first term in \eqref{eq:proof_ftail_bound_nonspur5} is $<2e^{\nu_1+\nu_2-6 - \log(g)} \ll 2e^{-\log^\alpha(g)}$
for any $\alpha <1$, whereas the second in term in \eqref{eq:proof_ftail_bound_nonspur5} is
$$
\ll \exp\left\{-\frac{\nu_1+\nu_2-5}{2} \log \left( \frac{\log^\alpha(g)}{\nu_1+\nu_2-6} \right) \right\},
$$
for any $\alpha< 1$ as $\nu_2 \rightarrow \infty$, as we wished to prove.

\bibliographystyle{plainnat}
\bibliography{references}

\end{document}